\documentclass{article}
\usepackage{amsmath}
\usepackage{array}
\usepackage{amsmath,amssymb,amsfonts,dsfont, latexsym,cancel, wasysym}
\usepackage{pb-diagram}
\usepackage{calrsfs}
\usepackage{amssymb}
\usepackage{amsthm}
\usepackage{float}
\usepackage{amsfonts}
\usepackage{latexsym}
\usepackage[dvips]{graphicx}
\usepackage{graphics}
\usepackage{mathrsfs}
\usepackage{enumerate}
\usepackage{multirow}
\usepackage{layout}
\usepackage{nccpic}
\usepackage[all]{xy}
\usepackage{fancyhdr}
\usepackage{color}
\usepackage{pstricks}
\usepackage{tikz}
\usepackage{wasysym}
\usepackage{subfigure}
\usepackage{enumitem}
\usepackage{subfig}
\usepackage{vmargin}

\definecolor{verdeobscuro}{rgb}{0.05, 0.5, 0.06}
\setlistdepth{9}

\newlist{myEnumerate}{enumerate}{9}
\setlist[myEnumerate,1]{label=(\alph*)}
\setlist[myEnumerate,2]{label=($\alph* _\arabic*$)}
\setlist[myEnumerate,3]{label=(${\alph*_{\arabic* \arabic*}}$)}
\setlist[myEnumerate,4]{label=(${\alph*_{\arabic* \arabic* \arabic*}}$)}
\setlist[myEnumerate,5]{label=(${\alph*_{\arabic* \arabic* \arabic*        \arabic*}}$)}
\setlist[myEnumerate,6]{label=(${\alph*_{\arabic* \arabic* \arabic*        \arabic* \arabic*}}$)}
\setlist[myEnumerate,7]{label=(${\alph*_{\arabic* \arabic* \arabic*        \arabic* \arabic*        \arabic*}}$)}
\setlist[myEnumerate,8]{label=(${\alph*_{\arabic* \arabic* \arabic*   \arabic* \arabic* \arabic* \arabic*}}$)}
\setlist[myEnumerate,9]{label=(${\alph*_{\arabic* \arabic* \arabic*    \arabic*  \arabic*    \arabic*  \arabic*    \arabic*}}$)}

\def\keywords#1{\par\medskip
\noindent\textbf{Keywords.} #1}
\def\subjclass#1{
\noindent{\emph{Mathematics Subject Classification:} #1}}

\theoremstyle{definition}
\newtheorem{definition}{Definition}
\newtheorem{remark}{Remark}
\newtheorem{theorem}[definition]{Theorem}
\newtheorem{proposition}[definition]{Proposition}
\newtheorem{lemma}[definition]{Lemma}
\newtheorem{conjecture}[definition]{Conjecture}

\begin{document}

\title{Arc Representations}
\author{Salom\'on Dom\'inguez,\\
Facultad de Econom\'ia y Negocios,\\
 Universidad An\'ahuac Campus Sur,\\
           {\tt salomon.dominguez3@gmail.com}}

\date{}
\maketitle

\thispagestyle{empty}

\begin{abstract}
This paper was inspired by four articles: surface cluster algebras studied by Fomin-Shapiro-Thurston \cite{fst}, the mutation theory of quivers with potentials initiated by Derksen-Weyman-Zelevinsky \cite{dwz},  string modules associated to arcs on unpunctured surfaces by Assem-Br$\ddot{u}$stle-Charbonneau-Plamondon \cite{acbp} and Quivers with potentials associated to triangulated surfaces, part II: Arc representations by Labardini-Fragoso. \cite{lf2}. For a surface with marked points ($\Sigma,M$) Labardini-Fragoso associated a quiver with potential $(Q(\tau),S(\tau))$ then for an ideal triangulation of ($\Sigma,M$) and an ideal arc Labardini-Fragoso defined an arc representation of $(Q(\tau),S(\tau))$. This paper focuses on extent the definition of arc representation to a more general context by considering a tagged triangulation and a tagged arc. We associate in an explicit way a representation of the quiver with potential constructed Labardini-Fragoso and prove that the Jacobian relations are met.

\\


\end{abstract}

\tableofcontents
\setcounter{section}{-1}
\setcounter{page}{-1}
\section{Introduction.}
\setcounter{page}{1}

Fomin and Zelevinsky introduced cluster algebras \cite{fz} around 2001 in order to deal with questions related with dual canonical bases and total positivity in algebraic Lie theory. Soon afterward it was discovered that this kind of structure appears in many areas of mathematics like algebraic geometry, representation theory, and Teichm$\ddot{u}$ller theory.

By definition, the cluster algebras ${\mathcal A}_{B}$ associated to a skewsymmetrizable matrix with integer coefficients $B$ is a subalgebra of a field of rational functions. Such subalgebra is recursively generated by some sets called clusters. Begining with the matrix $B$ and an initial set called {\it seed} we obtain an other seed through the mutations of the initial seed. Through this paper, the matrix $B$ will be a skewsymmetric matrix.

The set of matrices with integer coefficients correspond bijectively with the set of 2-acyclic quivers. If the matrix $B$ corresponds to the quiver $Q$ let's denote by ${\mathcal A}_{Q}:={\mathcal A}_{B}$ to the cluster algebra associated to the matrix $B$.

This paper is about quivers with potential $(Q,S)$ and the representations of the Jacobian algebra ${\mathcal P}(Q,S)$. Where $S$ is a non degenerated potential of a 2-acyclic quiver $Q$ and ${\mathcal P}(Q,S)$ is the quotient of the complete path algebra of $Q$ module the closure of the ideal generated by the cyclic derivatives of the potential $S$.

We mean by {\it potential} ( \cite[Definition 3.1]{dwz}) in $Q$ an element in $R \langle \langle Q \rangle \rangle$ (see \cite[Definition 2.2]{dwz}) in which all the summands are cycles of positive length. We denote by $R \langle \langle Q \rangle \rangle_{cyc}$ the set of all the potentials in $Q$, which is a closed subspace of $R \langle \langle Q \rangle \rangle$. 

For each arrow $a$ in $Q_{1}$ and every cycle $a_{1}...a_{d}$ in $Q$ we define the {\it cyclic derivative} as follows:

$$\partial_{a}(a_{1}...a_{d})=\displaystyle{ \sum_{i=1}^{d}}\delta_{a,a_{i}}a_{i+1}...a_{d}a_{1}...a_{i-1},$$ 

(where $\delta_{a,a_{i}}$ is the Kronecker delta) It extends linearly and continuously to obtain a morphism $\partial_{a}:R \langle \langle Q \rangle \rangle_{cyc} \xrightarrow{} R \langle \langle Q \rangle \rangle$(\cite[Definition 3.1]{dwz}). 

It was defined in \cite{fst} the cluster algebra ${\mathcal A}_{(\Sigma,M)}$ which comes from the tagged triangulations of a Riemann surface with marked points $(\Sigma,M)$. In addition, they proved that there is a bijection between the clusters of ${\mathcal A}_{(\Sigma,M)}$ and the tagged triangulations of Riemann surfaces $(\Sigma,M)$.

Let's consider an ideal triangulation $\tau$ without self-folded triangles of a surfaces with marked points ($\Sigma,M$). Combining the articles \cite{dwz} and \cite{fst}, Labardini-Fragoso defined a non degenerate potential $S(\tau)$ for $Q(\tau)$.

For a surface $\Sigma$ without marked points in the interior and an ideal triangulation $\tau$ of $\Sigma$. It is defined the arc representation of an arc $i$ (which does not belong to the triangulation $\tau$) in a natural way following the ideas in \cite{acbp} and \cite{ccs}. 

Now let ($\Sigma,M)$ be a surface with marked points in the interior (punctures), $\tau$ an ideal triangulation (triangulation without tags) of $(\Sigma,M)$ without self-folded triangles, $i$ an arc (without tags) which does not belong to $\tau$ and $(Q(\tau),S(\tau))$ the quiver with potential associated to the triangulation $\tau$. An easy calculus shows that the punctures cause a problem if we define in a natural way the representation $m(\tau,i)$. That is, the representation $m(\tau,i)$ does not met the Jacobian relations of $(Q(\tau),S(\tau))$.

Since the representation $m(\tau,i)$ does not met the Jacobian relations. It is defined in \cite{lf2} the detour curves (denoted by $d^{\Delta}_{(a_{1},a_{2})}$) and the detour matrices (denoted by $D^{\Delta^{\alpha}}_{i,j}$). Where $\alpha$ is an arrow in $Q(\tau)$ such that $t(\alpha)=j$.

Using the detour matrices it is defined in \cite{lf2} the arc representation $M(\tau,i)$ and it is to proved that the arc representation $M(\tau,i)$ meets the Jacobian relations.

The motivation of this paper is to extend the arc representations defined by Labardini-Fragoso to tagged triangulations context. The concepts that we need will be defined on a triangulation $\tau$ with no negative signature (see \cite[Definition 3.9]{lf3}) of a surface with marked points $(\Sigma,M)$. Later we will generalize to any tagged triangulation. It is worth mentioning that the definition of arc representation is based on the definition of Labardini-Fragoso and it was deduced by doing several mutations of representations.

Section 1 is divided into 5 subsections. Throughout the five subsections, we remain some central definitions from \cite{fst} as a quiver, a path in a quiver, the quiver mutation in a vertex, the path algebra of a quiver, the complete path algebra and the Jacobian algebra. In addition, we remind a result from \cite{mosher} that relate any pair of ideal triangulations of a surface with marked points through the flip. We also enunciate some concepts and results from \cite{dwz} about quivers with potential for ideal and tagged triangulations defined in \cite{lf1} and \cite{lf3}.

In section 2, we recall some concepts and results from \cite{dwz} about representations of quivers with potential.

Section 3 begins with the definition of the triangle of type 1, 2 and 3, then it is shown how the quiver $\widehat{Q}(\tau^{\circ})$ is drawn on the surface. For each triangle of type 2 and 3, we assign a label to each side. 

Using the concepts mentioned in the previous paragraph it is defined when a curve surround a puncture. This concept will be useful to define the detour curves and the auxiliary curves. With those curves it is defined the string representation $m(\tau,i)$.

On the other hand, for a tagged arc $i$ with ends $p$ and $q$ which does not belong to the triangulations $\tau^{\circ}$, it is drawn a tagged curve $i'$. The tagged curve $i'$ is useful to define the string representation $m(\tau,i)$ in a graph way. 

Now, given an arc $j$ in $\tau^{\circ}$ it is defined the concept of the crossing point of $i'$ with the arc $j$. The section 3 ends defining the string representation $m(\tau,i)$ which is the first approximation to the arc representation.

Following Labardini-Fragoso's ideas in section 4, we define the sets ${\mathcal B}^{\Delta,1}_{i,j}$, with these sets it is drawn the $1$-detour and the $1$-auxiliary curves. We will use those curves to define (in a recursive way) the sets ${\mathcal B}^{\Delta,n}_{i,j}$ and then the $n$-detour and $n$-auxiliary curves. The section 4 ends defining the detour matrices the auxiliary matrices and arc representation $M(\tau,i)$. 

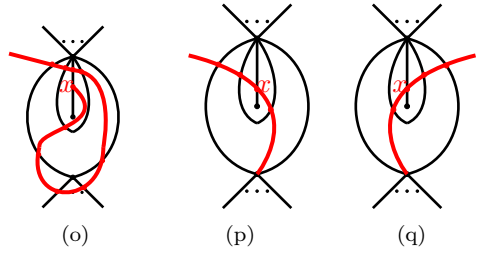
Section 5 begins with a proposition that says that if $x$ is a crossing which it is not an initial point of any detour, then $x\in W_{l}:=\{z\in M(\tau,i) | Jz=0\}$, where $J$ is the Jacobian ideal. Then we state and prove the principal theorem, such theorem state that the arc representation $M(\tau,i)$ satisfies the Jacobian relations and the module is nilpotent.

We finish the section 5 with the conjecture that the mutation of arc representation is compatible with the mutation of  quivers with potential. More precisely if $\tau$ and $ \sigma$ are tagged triangulations of a surface with marked points $(\Sigma,M)$ such that $\sigma$ is obtained from $\tau$ doing a flip in the arc $j$. Then the representations $\mu_{j}(M(\tau,i))$ and $M(\sigma,i)$ are right equivalent.

To finish this paper it is important to mention that the arc representations are interesting because:
\begin{itemize}
\item{} Once established the compatibility with the mutation and using the results of \cite[Section 13]{msw}. It can be calculated in a combinatoric way, the Euler characteristic of the quiver Grassmannians of $M(\tau,i)$, although the arc representation can be so difficult to calculate.

\item{} Even though the Jacobian algebras which are obtained from surfaces are gentle, it is very complicated to describe this class of representation with $E$- invariant equal to $0$ and $-1$. 

\end{itemize}
\section{Background on triangulations of surfaces and quivers with potential.}

\subsection{Quiver.}

A {\it quiver} $Q$ is the directed graph, that is, $Q=(Q_{0},Q_{1},s,t)$, where $Q_{0}$ is the vertex set of $Q$, $Q_{1}$ the arrow set and $s,t: Q_{1}\longrightarrow {}{Q_{0}}$ are functions defined as follows:\\
Given an arrow $\alpha:i \longrightarrow {}{j}$ in $Q$ let's define $s(\alpha):=i$ and $t(\alpha):=j$. If there are no arrows 
$\alpha:i\longrightarrow{}{j}$ in $Q_{1}$ such that $i=j$ then we eill say that $Q$ is a {\it loop-free quiver}. We will always deal only with loop-free quivers.

A {\it path of length $d\textgreater 0$} is a sequence $a_{d}a_{d-1}\cdots a_{2}a_{1}$ of arrows with $t(a_{j})=s(a_{j+1})$ for $j=1,...,d-1$, a path $a_{d}a_{d-1}\cdots a_{2}a_{1}$ of length $d\textgreater 0$ is a {\it d-cycle} if $t(a_{d})=s(a_{1})$. We say that a quiver $Q$ is {\it 2-acyclic} if $Q$ has no 2-cycles.

Given two paths $a=a_{d'}a_{d'-1}\cdots a_{2}a_{1}$ and $b_{d}b_{d-1}\cdots b_{2}b_{1}$, if $t(b_{d})=s(a_{1})$ then the composition $ab$ is given by $ab=a_{d'}a_{d'-1}\cdots a_{2}a_{1}b_{d}b_{d-1}\cdots b_{2}b_{1}$.

$$\bullet \xrightarrow {b_{1}}{\bullet}\xrightarrow {b_{2}}{\bullet}\cdots \bullet \xrightarrow {b_{d-1}}{\bullet}\xrightarrow {b_{d}}{\bullet} \xrightarrow {a_{1}}{\bullet}\xrightarrow {a_{2}}{\bullet}\cdots \bullet\xrightarrow {a_{d'-1}}{\bullet}\xrightarrow {a_{d'}}{\bullet}$$

For $i\in Q$, an {\it i-hook} is a path $ab$ of length $2$, such that $t(b)=i=s(a)$.

\begin{figure}[H]
\centering 
\begin{tikzpicture}[scale=0.55]

\filldraw [black] 
(-2, -3) circle (2pt)
(2,-3) circle (2pt);

\draw (0,-1)node{$i$};

\draw[->][line width=1pt] (-1.9,-2.9) -- (-0.15, -1.1)
node[pos=0.4,above] {$b$}; 
\draw[->][line width=1pt] (0.15,-1.1) -- (1.9, -2.9)
node[pos=0.6,above] {$a$}; 

\end{tikzpicture}
\caption{$i$-hook} 
\end{figure} 

\begin{definition}\label{carcaj}
Given a quiver $Q$ and a vertex $i$ in $Q$ such that $Q$ does not have incident 2-cycles at the vertex $i$ . We define the {\it mutation of $Q$ in direction $i$} as the quiver $\mu_{i}(Q)$ with vertex set $Q_{0}$ that results after applying the following three step procedure.

\begin{itemize}
\item[Step 1.] For each $i$-hook $ab$ we intruduce a new arrow $[ab]:s(b)\longrightarrow {}{t(a)}$.
\item[Step 2.]Replace each arrow $a: i \longrightarrow{} {t(a)}$ of $Q$ by an arrow $a^{*}: t(a) \longrightarrow {}{i}$ and each arrow $b: s(b) \longrightarrow{} {i} $ by $b^{*}: i \longrightarrow {}{s(b)}$. 
\item[Step 3.] Choose a maximal collection of disjoint $2$-cycles and remove them.
\end{itemize}

The quiver obtained after the steps 1 and 2 is called {\it Premutation} of $Q$ in direction $i$, and we will denote it by $\tilde{\mu}_{i}(Q)$.

\end{definition}

\begin{remark}\

\begin{itemize}
\item[1)] The mutation is defined for non necessarily $2$-acyclic quiver but in order to be able to perform mutation at any vertex of $Q$ we need it to be $2$-acyclic.
\item[2)] Choosing of the maximal collection in the third step of the Definition \ref{carcaj} is not given following a canonical procedure. However, up to this choose the mutation is an involution on the class of $2$-acyclic quivers, that is, $\mu^{2}_{i}(Q)\simeq Q$ for every $2$-acyclic quiver.
\end{itemize}

\end{remark}

\subsection{Ideal triangulations of surfaces and their flips.}\label{carcaj de triangulacion}
We briefly review the material on triangulation of surfaces and their signed adjacency matrices and flips. The reader could find more information in \cite{fst}.

\begin{definition}\label{no sigma}
{\it A border surface with marked points} is a pair $(\Sigma,M)$, where $\Sigma$ is a compact connected oriented Riemann surface with no empty boundary, and $M$ is a finite set on $\Sigma$, call {\it marked points}. Such that $M$ is a non empty set and has at least one point from each connected component of the boundary of $\Sigma$. The marked points that lie in the interior of $\Sigma$ will be called punctures, and the set of punctures of $(\Sigma,M)$ will be denoted by $P$. We will always assume that $(\Sigma,M)$ is none of the following:

\begin{enumerate}
\item{}A sphere with less than 5 punctures,
\item{}A unpunctured monogon, digon or triangle,
\item{}A once-punctured monogon.

\end{enumerate}
\end{definition}

Here by a monogon (resp. digon, triangle), we mean a disk with exactly one (resp. two, three) marked point(s) on the boundary.

\begin{definition}
We define a {\it curve} $\gamma$ in $(\Sigma,M)$ as a continuous function $\gamma: [0,1]\longrightarrow {}{\Sigma}$. Usually we will denote by $I$ to the interval $[0,1]$.

\end{definition}

Let $\gamma_{1}$ and $\gamma_{2}$ be two curves in $(\Sigma ,M)$, we say that $\gamma_{1}$ is isotopic to $\gamma_{2}$ if there is a continuos function $H:I\times \Sigma \xrightarrow {}{\Sigma}$ such that the following properties are met:
\begin{itemize}
\item[1)]For all $t\in I$ the function $H(0,t):\Sigma \longrightarrow {}{\Sigma}$ is a homeomorphism,
\item[2)]$H(0,x)=x$ for all $x\in \Sigma$,
\item[3)]$H(1,\gamma_{1})=\gamma_{2}$,
\item[4)]$H(t,m)=m$ for all $m\in M$ and $t\in I$.
\end{itemize}
We call the function $H$ isotopy.

\begin{definition}\label{arco}
Let $(\Sigma,M)$ be a border surface with marked points. A curve $\gamma$ in $\Sigma$ (up to isotopy) is an {\it arc} if the following four conditions are met;
\begin{itemize}
\item[1)] The ends of $\gamma$ are points in $M$,
\item[2)]The curve $\gamma$ does not intersect itself, except possibly at the ends,
\item[3)]The relative interior of $\gamma$ is a disjoint of $M$ and $\partial \Sigma$, where $\partial \Sigma$ is the boundary of $\Sigma$,
\item[4)]$\gamma$ does not cut to a monogon without punctures or a digon without punctures.

\end{itemize}
\end{definition}
We will call {\it loop} to the arc that the ends coincide, and we say that two arcs are compatible if there are arcs 
in their respective isotopy classes whose relative interior does not intersect \cite[Definition 2.4]{fst}.

\begin{proposition}
Given any collection of pairwise compatible arcs, it is always possible to find representatives in their isotopy classes whose relative interiors do not intersect each other.
\end{proposition}
For more details \cite[Proposition 2.5]{fst}.

\begin{definition}
{\it An ideal tringulation $\tau$} of $(\Sigma,M)$ is any maximal collection of pairwise campatible arcs such that the relative interior of those arcs do not intersect each other. 
\end{definition}

If $\tau$ is an ideal triangulation of $(\Sigma,M)$ and we take a connected component of the complement in $\Sigma$ of the union of the arcs in $\tau$, the closure $\Delta$ of this component will be called an ideal triangle of $\tau$. An ideal triangle is called {\it interior} if its intersection with the boundary of $\Sigma$ consisting only of (possibly none) marked points. Otherwise, it will be called {\it non interior}. An interior ideal triangle $\Delta$ is {\it self-folded} if it contains exactly two arcs of $\tau$ (note that every interior ideal triangle contains at least two and at most three arcs of $\tau$. While non interior ideal triangle contains at least one and at most two arcs).

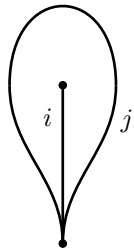
\begin{figure}[H]
\centering 
\begin{tikzpicture}[scale=0.7]

\draw[line width=1pt] (0, 0) .. controls(0, 2) and (2, 2) .. (2, 0);
\draw[line width=1pt] (0,0) to [out=-90, in=90] (1,-3);
\draw[line width=1pt] (1,-3) to [out=90, in=-90] node[pos=0.8, right] {$j$} (2,0);
\filldraw [black] (1,-3) circle (2pt)
(1,0) circle (2pt);
\draw[line width=1pt] (1,-3) -- (1,0) node[pos=0.8, left] {$i$};

\end{tikzpicture}
\caption{Self-folded triangle.\label{td}}
\end{figure} 
We will call the arc $i$ the folded side of $\triangle$.

The number $n$ of arcs in a triangulation $\tau$ is determined by the genus $g$ of $\Sigma$, the number $b$ of component of the boundary of $\Sigma$, the number $p$ of punctures and the number $c$ of points in the boundary of $\Sigma$, according to the formula $n=6g+3b+3p+c-6$. This formula could be proved using the definition and basic properties of the Euler characteristic. Hence the number $n$ is an invariant of $(\Sigma,M)$, called the rank of $(\Sigma,M)$ because $n$ coincides with the rank of the cluster algebra associated to $(\Sigma,M)$.

\begin{remark}
The reason for excluding the surfaces in the Definition \ref{no sigma} is the fact that their triangulations are empty or there is only one triangulation. The reason for excluding the spheres with less than five punctures is a bit more subtle (self-folded triangles present some unpleasant properties in these surface). 
\end{remark}

The following was first remarked in \cite{mosher}. Let $\tau$ be an ideal triangulation of $(\Sigma,M)$, we consider an arc $i\in \tau$ such that $i$ is not the folded side of any self-folded triangle. Then there is an unique arc $j$ such that $\sigma= \{ \tau - \{ i \} \} \cup \{ j \}$ is an ideal triangulation. To the procedure of change the arc $i$ by the arc $j$ is called the {\t flip in the arc $i$}. An important remark is that it is not possible to do the flip in the arc $i$ of the Figure \ref{td}. The case when the arc $i$ is the folded side of a self-folded triangle will be considered in the subsection $5$ of this section.

\begin{proposition}
Any pair of ideal triangulations are related by a sequence of flips.
\end{proposition}

Let $\tau$ be an ideal triangulation of $(\Sigma,M)$, now we are going to remind how it is defined the quiver $Q(\tau)$ in \cite[Definition 4.1]{fst}.

Given an ideal triangulation $\tau$ we associate a skewsymmetric matrix $B(\tau)$ of $n\times n$ whose arrows and columns correspond to the arcs of the triangulation $\tau$. Let $\pi_{\tau}:\tau \xrightarrow {}{\tau}$ be a function which is the identity function in the arcs $k$ such that $k$ is not a side of any self-folded triangle. Otherwise if $i$ and $j$ are the unique arcs of a self-folded triangle $\Delta$, where $i$ is the folded side of $\Delta$, then $\pi_{\tau}(i):=j$ and $\pi_{\tau}(j):=j$. For each non self-folded triangle $\Delta$ of $\tau$ let $B^{\triangle}=(b^{\triangle}_{i,j})$ be the matrix of $n\times n$ defined as follow:

$$ b_{i,j}^{\triangle} = \left\{ 
\begin{tabular}{cc}
1 & If $\pi_{\tau}(i)$, $\pi_{\tau}(j)$ are sides of $\triangle$ and $\pi_{\tau}(j)$ follows to $\pi_{\tau}(i)$ in clockwise order.\\
-1 & If the condition above is met but, in counterclockwise order. \\
0 & Any other case.
\end{tabular}
\right. $$

The {\it adjacent matrix} $B(\tau)$ is defined by:

$$B(\tau)=\displaystyle{\sum_{\triangle}}B^{\triangle},$$

where the sum runs over all the self-folded triangles $\tau$.

The matrix $B(\tau)$ defines an {\it adjacency quiver} $Q(\tau)$, whose vertex are the arcs of $\tau$, with $b_{i,j}$ arrows from $i$ to $j$ whenever $b_{i,j} \textgreater 0$. Since $B(\tau)$ is skewsymmetric, $Q(\tau)$ is a $2$-acyclic quiver. 

The next proposition is a direct consequence of the previous concepts and was remarked in \cite[Proposition 4.8]{fst}.

\begin{proposition}
Let $\tau$ and $\sigma$ be two ideal triangulations of $(\Sigma,M)$. If $\sigma$ is obtained from $\tau$ doing a flip in an arc $i$, then $Q(\sigma)=\mu_{i}(Q(\tau))$. That is, if $\tau$ and $\sigma$ are two ideal triangulations such that they are related by the flip in the arc $i$, then the associated quivers to $\tau$ and $\sigma$ are related by the mutation in the vertex $i$.
\end{proposition}

\begin{definition}
If a puncture is incident to exactly two arcs $i_{1}$ and $i_{2}$ (as is shown in Figure \ref{dc}) of an ideal triangulation $\tau$, then the quiver $Q(\tau)$ has neither arrows from $i_{1}$ to $i_{2}$ nor from $i_{2}$ to $i_{1}$. For $i_{1}$ and $i_{2}$ we introduce a new arrow from $i_{1}$ to $i_{2}$ and an other arrow from $i_{2}$ to $i_{1}$. The quiver which is obtained by gluing both arrows  will be called {\it non reduced adjacency quiver} and it will be denote it by $\widehat{Q}(\tau)$.
\end{definition}
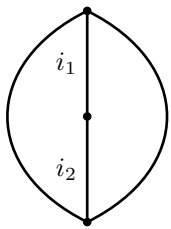
\begin{figure}[H]
\centering 
\begin{tikzpicture}[scale=0.7]
\filldraw [black] (0,0) circle (2pt)
(0,-2) circle (2pt)
(0,-4) circle (2pt);

\draw[line width=1pt] (0, 0) .. controls(-2, -1) and (-2, -3) .. (0, -4);
\draw[line width=1pt] (0, 0) .. controls(2, -1) and (2, -3) .. (0, -4);
\draw[line width=1pt] (0, 0) -- (0, -2)
node[pos=0.5,left] {$i_{1}$}; 
\draw[line width=1pt] (0, -2) -- (0, -4)
node[pos=0.5,left] {$i_{2}$};

\end{tikzpicture}
\caption{Puncture incident only at the arcs $i_{1}$ and $i_{2}$.\label{dc}}
\end{figure}
Clearly to obtain $Q(\tau)$ from $\widehat{Q}(\tau)$ is enough to remove the 2-cycles.

\subsection{Quivers with potentials.}

This section is about some concepts and results about of quivers with potential, for more details see \cite{dwz}.

Given a quiver $Q$ we denote by $R$ to the $K$-vector space with basis $\{e_{i} | i\in Q_{0}\}$, now $R$ with the operation $e_{i}e_{j}=\delta_{ij}e_{i}$ has structure of commutative $K$-algebra and we will call $R$ {\it generated space by the vertex} of $Q$. We will think to each vertex $e_{i}$ as a path of zero length in $i$. Analogously we define the {\it generated space by the arrows} of $Q$ as the $K$-vector space $A$ with basis the set of arrow of $Q$. We note that $A$ has structure of $R$-bimodule if we define $e_{i}a=\delta_{i,t(a)}a$ and $ae_{j}=a\delta_{s(a),j}$ for $i\in Q_{0}$ and $a\in Q_{1}$. For $d\geq 0$ we denote  by $A^{d}$ to the $K$-vector space with basis the set of all the paths of length $d$ in $Q$, this vector space has structure of $R$ bimodule too. For $d=0$ and $d=1$, the $R$-bimodule $A^{d}$ coincide with $R$ and $A$ respectively.

Now we will remember two definitions ({\it path and complete path algebra}) which are very important for this work.

\begin{definition}\label{algebra completa}
The {\it path algebra} of a quiver $Q$ is defined as:
$$R \langle Q \rangle =\displaystyle{\bigoplus_{d=0}^{\infty}}A^{d},$$

where $A^{d}$ is the generated space by the paths of length $d$.

\end{definition}

For the spaces $A^{n}$ and $A^{m}$ is easy to verify that $A^{n}\cdot A^{m}\subseteq A^{n+m}$, with the multiplication operation induced by the concatenation of paths (\cite[Definition 2.2]{dwz}). This is because if $c$ is a path of length $n$
and $c'$ is a path of length $m$ then $c\cdot c'$ is $0$ or is a path of length $m+n$.

\begin{definition}\label{algebra completa}
The {\it complete path algebra} of a quiver $Q$ is defined as:

$$R \langle \langle Q \rangle \rangle = \displaystyle{\prod_{d=0}^{\infty}}A^{d}.$$

This means that the elements of $R \langle \langle Q \rangle \rangle$ are $K$-linear combinations (possibly infinite) of elements in $R \langle Q \rangle $. The multiplication in $R \langle \langle Q \rangle \rangle$ is defined extending in a natural way the multiplication (concatenation) of $R \langle Q \rangle $.

\end{definition}

If in $R \langle \langle Q \rangle \rangle = \displaystyle{\prod_{d=0}^{\infty}}A^{d}$ we consider the topology $m$-adica which fundamental system of open neighborhoods around the zero is given by the powers of $\mathfrak{m}=\mathfrak{m}(A)=\displaystyle{\prod_{d\geq 1}}A^{d}$ (the ideal of $R \langle \langle Q \rangle \rangle $ generated by the arrows) then $R \langle Q \rangle $ is dense in $R \langle \langle Q \rangle \rangle$.

\begin{remark}
A sequence $(x_{n})_{n\in \mathds{N}}$ of elements in $R \langle \langle Q \rangle \rangle$ converge if only if for all $d\geq 0$, the sequence stabilizes when $n \xrightarrow {}{\infty}$, in which case $\displaystyle{\lim_{n \to \infty}}x_{n}=\displaystyle{\sum_{d\geq 0}\lim_{n \to \infty}}x_{n}^{(d)}$ where $x_{n}^{(d)}$ denote the component of degree $d$ of $x_{n}$.
\end{remark}

Although the action of $R$ in $R \langle \langle Q \rangle \rangle$ ($R \langle Q \rangle$) is not central (in the sense that if $a$ and $b$ are two paths in $Q$, then $e_{t(a)}ab=ae_{s(a)}b=abe_{s(b)}$), the action is compatible with the multiplication of $R \langle \langle Q \rangle \rangle$. Therefore we say that $R \langle \langle Q \rangle \rangle$ ($R \langle Q \rangle$) is a $R$-algebras. In consequence if $\varphi$ is a homomorphism of $K$-algebras between path algebras (or complete path algebras), we say that $\varphi$ is a homomorphism of $R$-algebras if the subyacent quivers have the same set of vertex and $\varphi(r)=r$ for all $r\in R$. It is easy to verify that each homomorphism of $R$-algebras between complete paths algebras is continuos.

The following is a useful way to define homomorphism between complete path algebras.

Each pair $(\varphi^{1},\varphi^{2})$ of homomorphism of $R$-bimodules $\varphi^{1}:A\xrightarrow {}A'$ and $\varphi^{2}:A \xrightarrow {}\mathfrak{m}(A')^{2}$ can be extended uniquely to a continuous homomorphism of $R$-algebras $\varphi:R \langle \langle Q \rangle \rangle \xrightarrow{} R \langle \langle Q' \rangle \rangle$ such that $\varphi |_{A}=(\varphi^{1},\varphi^{2})$. In addition $\varphi$ is a isomorphism of $R$ algebras if only if $\varphi^{1}$ is a isomorphism of $R$-bimodules.

A {\it potential} in $A$ (or $Q$) is an element of $R \langle \langle Q \rangle \rangle$ such that all the terms are cycles of positive length ( \cite[Definition 3.1]{dwz}). We denote by $R \langle \langle Q \rangle \rangle_{cyc}$ the set of all the potentials in $A$ and the set $R \langle \langle Q \rangle \rangle_{cyc}$ is a closed vector subspace of $R \langle \langle Q \rangle \rangle$. Two potentials $S$ and $S'$ in $R \langle \langle Q \rangle \rangle_{cyc}$ are cyclically equivalent if $S-S'$ belong to the closure of the vector subspace of $R \langle \langle Q \rangle \rangle$ generated by all the elements $a_{1}...a_{d} - a_{2}...a_{d}a_{1}$. Where $a_{1}...a_{d}$ is a cylce of positive length (\cite[Definition 3.2]{dwz}).

A quiver with potential is a pair $(A,S)$ (or $(Q,S)$), where $S$ is a potential in $A$ such that given two different terms of $S$, those terms are not cyclically equivalent. (\cite[Definition 3.2]{dwz}).
We will use the notation $QP$ instead of quiver with potential. The {\it direct sum} of two QPs $(A,S)$ and $(A',S')$ in the same vertex set is the QP $(A,S)\oplus(A',S')=(A\oplus A',S+ S' )$.

If $(A,S)$ and $(A',S')$ are two $QPs$ in the same vertex set, we say that $(A,S)$ is right equivalent to $(A',S')$ if there is a right equivalence between $(A,S)$ and $(A',S')$. This means that there is an isomorphism of $R$-algebras $\varphi:R \langle \langle Q \rangle \rangle \xrightarrow{} R \langle \langle Q' \rangle \rangle$ and $\varphi(S)$ is cyclically equivalent to $S'$(\cite[Definition 4.2]{dwz}).

For each arrow $a$ in $Q_{1}$ and each cycle in $Q$ we define the cyclic derivative as follows:

$$\partial_{a}(a_{1}...a_{d})=\displaystyle{ \sum_{i=1}^{d}}\delta_{a,a_{i}}a_{i+1}...a_{d}a_{1}...a_{i-1},$$

(where $\delta_{a,a_{i}}$ is the Kronecker delta) and extend lineally and continuously to obtain a morphism $\partial_{a}:R \langle \langle Q \rangle \rangle_{cyc} \xrightarrow{} R \langle \langle Q \rangle \rangle$(\cite[Definition 3.1]{dwz}). If the potentials $S$ and $S'$ are cyclically equivalent, then $\partial_{a}(S)= \partial_{a}(S')$.

The {\it Jacobian ideal } $J(S)$ is the closure of the bilateral ideal of $R \langle \langle Q \rangle \rangle$ generated by $\{\partial_{a}(S) | a\in Q_{1}\}$ and the {\it Jacobian algebra} ${\mathcal P}(Q,S)$ is the quotient algebra $R \langle \langle Q \rangle \rangle / J(S)$(\cite[Definition 3.1]{dwz}). The Jacobian ideal and the Jacobian algebra are invariant under right equivalences, in the sense that if $\varphi:R \langle \langle Q \rangle \rangle \xrightarrow{} R \langle \langle Q' \rangle \rangle$ is a right equivalence between $(A,S)$ and $(A',S')$, then $\varphi$ send to $J(S)$ in $J(S')$ and therefore induces an isomorphism between ${\mathcal P}(Q,S)$ and ${\mathcal P}(Q',S')$ (\cite[Proposition 3.7]{dwz}).

A $QP$ is {\it trivial} if $S\in A^{2}$ and $A$ is generated by $\{\partial_{a}(S) | a\in Q_{1}\}$(\cite[Definition 4.4 and see Proposition 4.4]{dwz}). Otherwise, we say that a QP is reduced if the component of degree $2$ of $S$ is zero, this means that in $S$ do not appear $2$-cycles. We noticed that the subjacent quiver of a QP could have 2-cycles and we say that a QP is 2-acyclic if the subjacent quiver does not have 2-cycles we say that is 2-acyclic.

\begin{theorem}\label{split}
(Splitting theorem, \cite[theorem 4.6]{dwz})For every $QP$ $(A,S)$ there is a trivial $(A_{triv},S_{triv})$ an a reduced $QP$ $(A_{red},S_{red})$ such that $(A,S)$ is right equivalent to the direct sum $(A_{triv},S_{triv})\oplus (A_{red},S_{red})$. Furthermore the right equivalence class of the $QPs$ $(A_{triv},S_{triv})$ and $(A_{red},S_{red})$ is determined by the right equivalence class of $(A,S)$.
\end{theorem}

In the situation of the Theorem \ref{split} the $QP$ $(A_{triv},S_{triv})$ is called trivial part of $(A,S)$, respectively $(A_{red},S_{red})$ is called reduce part of $(A,S)$(\cite[Definition 4.13]{dwz}).

\subsection{The potential of a triangulation}

Given a triangulation $\tau$ of a bordered surface with marked points, we associate a quiver $Q(\tau)$ in the section \ref{carcaj de triangulacion}, now we will define a potential. This potential was defined in \cite[Definition 3.1]{lf1}.

Let $(\Sigma,M)$ be a bordered surface with marked points and $P\subseteq M$ the puncture set of $(\Sigma,M)$. For each $p\in P$ we choose a scalar non zero $x_{p}$ in $K$. This choice will be fixed for each triangulation of $(\Sigma,M)$.

\begin{definition}\label{potencialD}
Let $\tau$ be an ideal triangulation of $(\Sigma,M)$. Based on the choice of the scalars $(x_{p})_{p\in P}$ we associate to $\tau$ a potential $S(\tau)\in R \langle \langle Q(\tau) \rangle \rangle$ as follows. Let $\widehat{A}(\tau)$ be the generated space by the arrows of $\widehat{Q}(\tau)$.

\begin{enumerate}
\item{}For each ideal triangle which is not a self-folded triangle $\Delta$ in $\tau$, we denote by $\widehat{S}^{\triangle}$ the triangle oriented in clockwise in $\widehat{Q}(\tau)$ (up to cyclical equivalence).

\item{}If $\Delta$ is a triangle as in item $1.$ and the sides of $\Delta$ are $j,k,l$ such that $\Delta$ is adjacent to two self-folded triangles as is shown in the next Figure.

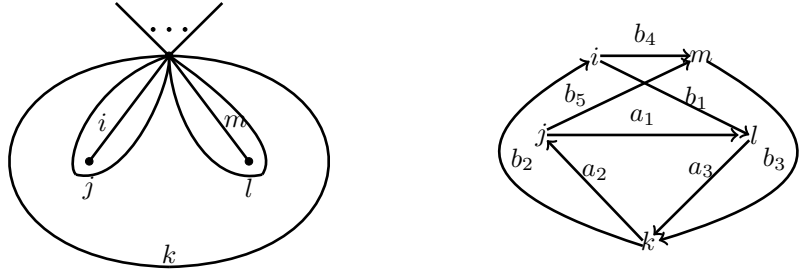
\begin{figure}[H]
\centering 
\begin{tikzpicture}[scale=0.7]

\filldraw [black] (-7,2.5) circle (2pt)
(-8.5,0.5) circle (2pt)
(-5.5,0.5) circle (2pt)
(-7,3) circle (1pt)
(-6.7,3) circle (1pt)
(-7.3,3) circle (1pt);

\draw[line width=1pt] (-7,2.5) -- (-8, 3.5);
\draw[line width=1pt] (-7,2.5) -- (-6, 3.5);
\draw[line width=1pt] (-7,2.5) -- (-8.5,0.5);
\draw[line width=1pt] (-7,2.5) -- (-5.5,0.5);
\draw[line width=1pt] (-7, 2.5) .. controls(-11,2.4) and (-11,-1.4) .. (-7, -1.5);
\draw[line width=1pt] (-7, 2.5) .. controls(-3,2.4) and (-3,-1.4) .. (-7, -1.5);
\draw[line width=1pt] (-7, 2.5) .. controls(-8.5,2) and (-9,0.5) .. (-8.75, 0.25);
\draw[line width=1pt] (-7, 2.5) .. controls(-6.95,2) and (-7.75,0) .. (-8.75, 0.25);
\draw[line width=1pt] (-7, 2.5) .. controls(-7,1.5) and (-6.3,0) .. (-5.25, 0.25);
\draw[line width=1pt] (-7, 2.5) .. controls(-7,2.5) and (-4.75,1) .. (-5.25, 0.25);

\draw (-8.5,0) node{$j$};
\draw (-8.25,1.25)node{$i$};
\draw (-5.75,1.25)node{$m$};
\draw (-5.5,0)node{$l$};
\draw (-7,-1.25)node{$k$};

\draw (0,1)node{$j$};
\draw (1,2.5)node{$i$};
\draw (3,2.5)node{$m$};
\draw (4,1)node{$l$};
\draw (2,-1)node{$k$};

\draw[->][line width=1pt] (1.1,2.4) -- node[pos=0.5,right]{$b_{1}$}(3.9,1);
\draw[->][line width=1pt] (1.1,2.5) -- node[pos=0.5,above]{$b_{4}$}(2.8,2.5);
\draw[->][line width=1pt] (0.1,1) -- node[pos=0.5,above]{$a_{1}$}(3.7,1);
\draw[->][line width=1pt] (0.1,1.1) -- node[pos=0.2,above]{$b_{5}$}(2.8,2.4);
\draw[->][line width=1pt] (3.9,.9) -- node[pos=0.5,above]{$a_{3}$}(2.1,-.9);
\draw[->][line width=1pt] (1.9,-1) -- node[pos=0.5,above]{$a_{2}$}(0.1,.9);
\draw[->][line width=1pt] (1.9,-1.1) .. controls (-1.5, 0) and (-1.5, 1) .. node[pos=0.5,right]{$b_{2}$}(.9,2.4);
\draw[->][line width=1pt] (3.1,2.4) .. controls (5.5, 1) and (5.5, 0) .. node[pos=0.5,left]{$b_{3}$}(2.2,-1);

\end{tikzpicture}

\end{figure}

Then it is defined $\widehat{T}^{\triangle}=\frac{b_{2}b_{3}b_{4}}{x_{p}x_{q}}$ (up to cyclical equivalence), where $p$ and $q$ are the punctures which are enclosed by the non folded sides of the self-folded triangles that are adjacent to $\Delta$. Otherwise, if $\Delta$ is adjacent to less than two self-folded triangles then it is defined $\widehat{T}^{\triangle}=0$.

\item{} If the puncture $p$ is adjacent to exactly one arc $i$ of $\tau$ then $i$ is the folded side of some self-folded triangle in $\tau$ and we have the following configuration around the arc $i$.

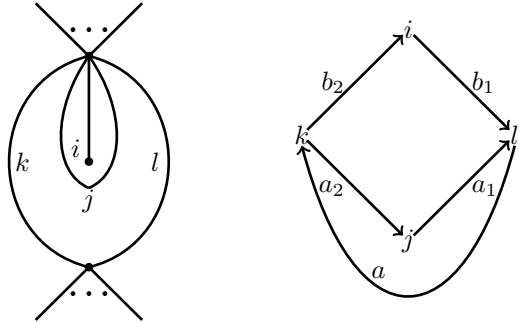
$$\begin{tikzpicture}[scale=0.7]

\filldraw [black] (0,2.5) circle (2pt)
(0,-1.5) circle (2pt)
(0,0.5) circle (2pt)
(0,3) circle (1pt)
(0.3,3) circle (1pt)
(-0.3,3) circle (1pt)
(0,-2) circle (1pt)
(0.3,-2) circle (1pt)
(-0.3,-2) circle (1pt);

\draw[line width=1pt] (0,2.5) -- (-1, 3.5);
\draw[line width=1pt] (0,2.5) -- (1, 3.5);
\draw[line width=1pt] (0,-1.5) -- (1, -2.5);
\draw[line width=1pt] (0,-1.5) -- (-1, -2.5);
\draw[line width=1pt] (0,2.5) -- (0, 0.5);
\draw[line width=1pt] (0,2.5) .. controls(-0.7,1.5) and (-0.7,0.3) .. (0, 0);
\draw[line width=1pt] (0,2.5) .. controls(0.7,1.5) and (0.7,0.3) .. (0, 0);
\draw[line width=1pt] (0,2.5) .. controls(-2, 2) and (-2,-1) .. (0, -1.5);
\draw[line width=1pt] (0,2.5) .. controls(2, 2) and (2,-1) .. (0, -1.5);

\draw (-1.25,0.5)node{$k$};
\draw (1.25,0.5)node{$l$};
\draw (-0.25,0.75)node{$i$};
\draw (0,-0.25)node{$j$};

\draw[<-][line width=1pt](4,0.8) .. controls (5, -3) and (7, -3) .. node[pos=0.3, right] {$a$}(8,0.8);
\draw (4,1)node{$k$};
\draw (8,1)node{$l$};
\draw (6,3)node{$i$};
\draw (6,-1)node{$j$};
\draw[->][line width=1pt] (4.1,1.1) -- node[pos=0.5,left]{$b_{2}$}(5.9,2.9);
\draw[->][line width=1pt] (6.1,2.9) -- node[pos=0.5,right]{$b_{1}$}(7.9,1.1);
\draw[->][line width=1pt] (4.1,0.9) -- node[pos=0.5,left]{$a_{2}$}(5.9,-0.9);
\draw[->][line width=1pt] (6.1,-0.9) -- node[pos=0.5,right]{$a_{1}$}(7.9,0.9);

\end{tikzpicture}$$

The case when $k$ and $l$ are arcs of $\tau$ that do not belong to the boundary of $\Sigma$, then we define $\widehat{S}^{p}=-\frac{ab_{1}b_{2}}{x_{p}}$ up to cyclical equivalence.

\item{}Let $p$ be a puncture which is adjacent to more than one arc, we delete all incident loops at the puncture $p$. The arrows between the remaining adjacents arcs to $p$ form a cycle $a_{1}^{p},...,a_{d}^{p}$ (up to cyclical equivalence) oriented counter clockwise around $p$. Under this situation we define $\widehat{S}^{p}=x_{p}a_{1}^{p},...,a_{d}^{p}$. 

\end{enumerate}

The {\it no reduced potential} $\widehat{S}(\tau)\in R\langle \langle \widehat{Q }\rangle \rangle$ of $\tau$ is defined as:

$$\widehat{S}(\tau)=\displaystyle{\sum_{\triangle}}(\widehat{S}^{\triangle} + \widehat{T}^{\triangle})+ \displaystyle{\sum_{p \in P}}\widehat{S}^{p}$$

where the first sum run over all the interior non self-folded triangles.

Finally, we define $(Q(\tau),S(\tau))$ as the reduced part of $(\widehat{Q}(\tau),\widehat{S}(\tau))$(up to cyclical equivalence).

\end{definition}

Now by \cite[Proposition 10.2]{lf2} it is possible to choose the coefficients $x_{p}$ of the potential $S(\tau)$  as $-1$.

\subsection{Tagged triangulations.}

For an arc $i$ in $\tau$ which is not the folded side of a self-folded triangle, it is defined the flip in the arc $i$. Now to extend the concept to any arc we remind some definitions and results of \cite[Section 9.3]{fst}.

We will think that each arc $i$  is divided into three equal parts, which are the central part that consists only of interior of the arc $i$ and the other two parts which contain one end in each part. In the parts which contain one end of the arc $i$ we will set one label "plain" or "notched".

\begin{definition}
Let ($\Sigma$,M) be a bordered surface with marked points and $\tau$ an ideal triangulation of $(\Sigma,)$. A {\it tagged arc} is an arc which has a label at each end, such label is in the set \{"plain","notched"\} and the following conditions are met.

\end{definition}

\begin{itemize}
\item[i)]The arc is not a loop which encloses only one puncture,
\item[ii)]If the end of one arc is in the boundary of $\Sigma$ such end has been labeled "plain", 
\item[iii)] If the arc is a loop then each end of the arc has the same label.
\end{itemize}

The label "plain" will be omitted in the graphic representation  of the arc $i$ while the label "notched" will be represented as $\Bowtie$ . Now, each  ideal triangulation $\tau$ of $(\Sigma,M)$  could be thought as a tagged triangulation as follows;

If $\gamma$ is not a loop which encloses only one puncture then $\gamma$ thought as a tagged arc  will be the same arc with labels "plain" at both ends. Otherwise, if $\gamma$ is a loop based in the puncture $a$ and $\gamma$ enclosed only to the puncture $b$, as is shown in the following figure.

$$
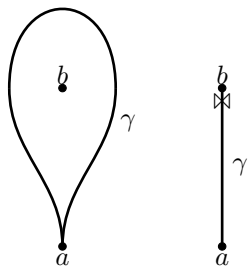
\begin{tikzpicture}[scale=0.7]
\draw[line width=1pt] (0, 0) .. controls(0, 2) and (2, 2) .. (2, 0);
\draw[line width=1pt] (0,0) to [out=-90, in=90] (1,-3);
\draw[line width=1pt] (1,-3) to [out=90, in=-90] node[pos=0.8, right] {$\gamma$} (2,0);
\filldraw [black] (1,-3) circle (2pt)
(1,0) circle (2pt);
\draw (1,-3.25)node{$a$};
\draw (1,0.25)node{$b$};
\filldraw [black] (4,-3) circle (2pt)
(4,0) circle (2pt);
\draw (4,-3.25)node{$a$};
\draw (4,0.25)node{$b$};
\draw[line width=1pt] (4,-3) -- (4,0)
node[pos=0.5, right] {$\gamma$}
node[pos=0.8, above] {$\Join$};

\end{tikzpicture}$$

Then we replace the arc $\gamma$ by the tagged arc which ends $a$ and $b$, the labels will be "plain" at the ends $a$ and $\Bowtie$ at the end $b$. Without losing generality we will denote this arc by $\gamma$ and $\bold{A^{\Join}}(\Sigma,M)$ will be the set of all the tagged arcs.

The following definitions were introduced in \cite[Secci\'on 2]{lf2}.

\begin{definition}\label{signatura}
Let $\varepsilon:P\longrightarrow \{1,-1\}$ be a function, where $P$ is the set of punctures, we define the {\it function} $t_{\varepsilon}:\bold{A}^{o}(\Sigma,M) \longrightarrow \bold{A^{\Join}}(\Sigma,M)$ as follows;

\begin{itemize}
\item[i)] If the arc $i$ is not a loop which encloses only one puncture, then the version without labels of $t_{\varepsilon}(i)$ is $i$. An end $p$ of $t_{\varepsilon}(i)$ will be tagged "notched" if only if $\varepsilon(p)$ is $-1$.
\item[ii)] If $i$ is a loop based in a marked point $q$ which encloses only one puncture $p$ then $t_{\varepsilon}(i)$ is the arc which connects the marked points $p$ and $q$, such that the label at the end $q$ is "notched" if only if $\varepsilon(q)=-1$. On the other hand the end $p$ is labeled "notched" if only if $\varepsilon(p)=1$.
\end{itemize}
\end{definition}

Let $\tau$ be a tagged triangulation of $(\Sigma,M)$, we define the {\it signature} of $\tau$ as the function $\delta_{\tau}:P\longrightarrow \{1,-1,0\}$ given by:

\[ \delta_{\tau}(p) = \left\{ \begin{array}{ll}
1 & \mbox{If the end of each arc incident at the puncture $p$ is labeled $"\ "$.}\\
-1 & \mbox{If the end of each arc incident at the puncture $p$ is labeled $\Join$.}\\
0 & \mbox{else.}\end{array} \right. \]

\begin{definition}
Let $\tau$ be a tagged triangulation of $(\Sigma,M)$, we define the {\it weak signature} of $\tau$ as the function:

\[
\varepsilon_{\tau}(p) = \left \{ \begin{array}{lcccl}
1 & \mbox{If $\delta_{\tau}(p)\in \{ 0,1\}$}\\
-1 & \mbox{else.} \\
\end{array}
\right.
\]
\end{definition}

Now given a tagged triangulation $\tau$, we can obtain an ideal triangulation denoted by $\tau^{\circ}$ as follows.

\begin{itemize}
\item[i)] Remove all the tags at the puncture $p$ with signature non zero.
\item[ii)] For each puncture $p$ with signature  equal to zero, we replace the tagged arc $i\in \tau$ which is labeled $\Bowtie$ at $p$ by the loop which enclosed only the puncture $p$ and  the arc $i$.
\end{itemize}

\begin{definition}
Let $\tau$ be a {\it tagged triangulation} of ($\Sigma$, M).
\begin{itemize}
\item[i)] We define the no reduced potential $\widehat{S}(\tau)\in R \langle \langle \widehat{Q}(\tau) \rangle \rangle$ associated to $\tau$ as;
$$\widehat{S}(\tau)=t_{\varepsilon_{\tau}}(\displaystyle{\sum_{\triangle}}(\widehat{S}^{\triangle}(\tau^{\circ}) + \widehat{T}^{\triangle}(\tau^{\circ}))+ \displaystyle{\sum_{p \in P}}(\varepsilon_{\tau}(p)\widehat{S}^{p}(\tau^{\circ})))$$

where the first sum runs over all the interior non self-folded triangles of $\tau^{\circ}$ and $\varepsilon_{\tau}$ is the weak signature of $\tau$.

\item[ii)] We define $(Q(\tau),S(\tau))$ as the reduced part of $(\widehat{Q}(\tau),\widehat{S}(\tau))$.
\end{itemize}
\end{definition}

\section{Background on representations of quivers with potential.}

Remembering that $R$ denotes the $K$-vector space with basis the set $\{e_{i} | i\in Q_{0}\}$, this vector space is in fact a commutative ring with the multiplication operation $e_{i}e_{j}=\delta_{ij}e_{i}$.

\begin{definition}
(\cite[Definition 10.1]{dwz}) Let $(Q,S)$ be a QP, a {\it representation of quiver with potential} is a third ${\mathcal M}=(Q,S,M)$, where $M$ consists of the following two families. 
\end{definition}

\begin{itemize}
\item[1)] A family $(M_{i})_{i\in Q_{0}}$ of $K$-vector spaces of finite dimension,
\item[2)] A family $(a_{M}:M_{s(a)}\longrightarrow {}{M_{t(a)}})_{a\in Q_{0}}$ of $K$-linear transformations such that $\partial_{a_{M}}(S)=0$ for all $a_{M} $.
\end{itemize}

A representation ${\mathcal M}=(Q,S,M)$ of a QP ($Q,S$) is {\it nilpotent} if there is an integer $r\geqslant 1$ such that any path of length at least $r$ is zero.

\begin{remark}

It is not true that any representation annihilated by the cyclic derivatives of $S(\tau)$ is nilpotent. That is, it is possible to construct $R\langle Q(\tau)\rangle$-modules that satisfy the cyclic derivatives of $S(\tau)$ but can not be given the structure of ${\mathcal P}(Q(\tau),S(\tau))$-module. One example of this is given by the representation shown in the Figure \ref{rep no nil}.

\begin{figure}[htp]
\begin{center}
\subfigure[Quiver]{
\hspace*{-1em}{\begin{tikzpicture}[scale=1]

\draw (0,0)node{$3$};
\draw (2,-2)node{$2$};
\draw (-2,-2)node{$1$};

\draw[->][line width=1.5pt] (-1.9,-1.75) -- (-0.2, 0.05)
node[pos=0.5,above] {$a_{1}$};
\draw[->][line width=1.5pt] (-1.75,-1.9) -- (-0.1, -0.15)
node[pos=0.65,below] {$a_{2}$}; 

\draw[<-][line width=1.5pt] (-1.6,-1.9) -- (1.75, -1.9)
node[pos=0.5,above] {$b_{2}$};
\draw[<-][line width=1.5pt] (-1.6,-2.1) -- (1.75, -2.1)
node[pos=0.5,below] {$b_{1}$};
\draw[<-][line width=1.5pt] (2.1,-1.75) -- (0.2, 0.05)
node[pos=0.5,above] {$c_{1}$};
\draw[<-][line width=1.5pt] (1.9,-1.85) -- (0.1, -0.15)
node[pos=0.65,below] {$c_{2}$};

\end{tikzpicture} }
}
\subfigure[Representation]{
\hspace*{0.6em}{\begin{tikzpicture}[scale=1]

\draw (0,0.15)node{$K$};
\draw (2,-2)node{$K$};
\draw (-2,-2)node{$K$};

\draw[->][line width=1.5pt] (-1.9,-1.75) -- (-0.2, 0.05)
node[pos=0.5,above] {$1$};
\draw[->][line width=1.5pt] (-1.75,-1.9) -- (-0.1, -0.15)
node[pos=0.65,below] {$1$}; 

\draw[<-][line width=1.5pt] (-1.6,-1.9) -- (1.75, -1.9)
node[pos=0.5,above] {$1$};
\draw[<-][line width=1.5pt] (-1.6,-2.1) -- (1.75, -2.1)
node[pos=0.5,below] {$1$};
\draw[<-][line width=1.5pt] (2.1,-1.75) -- (0.2, 0.05)
node[pos=0.5,above] {$1$};
\draw[<-][line width=1.5pt] (1.9,-1.85) -- (0.1, -0.15)
node[pos=0.65,below] {$1$};

\end{tikzpicture}}

} 
\end{center}
\caption{Representation which is not nilpotent \label{rep no nil}}

\end{figure}
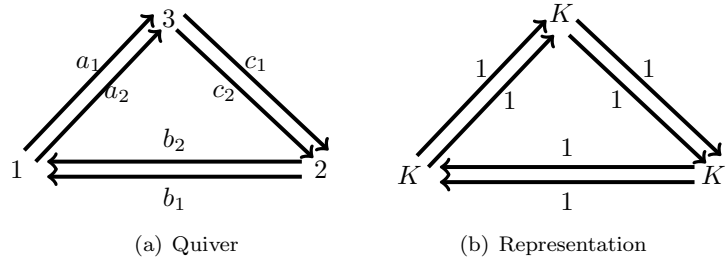

An easy calculus shows that the representation of Figure \ref{rep no nil} satisfy the cyclic derivatives of the potential $S=a_{1}b_{1}c_{1} + a_{2}b_{2}c_{2}-a_{1}b_{2}c_{1}a_{2}b_{1}c_{2}$, but is not nilpotent.

\end{remark}

\begin{definition}\label{equivalencia derecha QP}
(\cite[Definition 10.2]{dwz}) Let $(Q,S)$ and $(Q',S')$ be two quivers with potential in the same vertex set and let $\mathcal M= (Q,S,M)$ , ${\mathcal M'}= (Q',S',M')$ be the representations of $(Q,S)$ and $(Q',S')$ respectively. We define a {\it right equivalence} between ${\mathcal M}$ y ${\mathcal M'}$ as a pair $(\varphi,\psi)$ that satisfy the following two properties:

\begin{itemize}
\item{} $\varphi: R \langle \langle Q \rangle \rangle \longrightarrow {}{R \langle \langle Q' \rangle \rangle}$ is a right equivalence of quivers with potential between $(Q,S)$ and $(Q',S')$,
\item{} $\psi:M_{i}\longrightarrow {}{M'}_{i}$ is an isomorphism of vector spaces such that $\psi \circ u_{M}=\varphi(u)_{M'}\circ \psi$ for all elements $u\in R \langle \langle Q \rangle \rangle$.

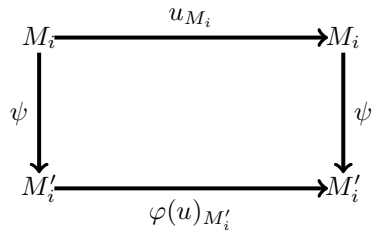
\begin{figure}[htp]
\begin{center}
\begin{tikzpicture}[scale=1]

\draw (2,0)node{$M_{i}$}; 
\draw (-2,0)node{$M_{i}$};
\draw (2,-2)node{$M'_{i}$};
\draw (-2,-2)node{$M'_{i}$};

\draw[->][line width=1.5pt] (-1.8,0) -- (1.8, 0)
node[pos=0.5,above] {$u_{M_{i}}$};
\draw[->][line width=1.5pt] (-1.8,-2) -- (1.8, -2)
node[pos=0.5,below] {$\varphi(u)_{M'_{i}}$}; 

\draw[->][line width=1.5pt] (-2,-0.2) -- (-2, -1.8)
node[pos=0.5,left] {$\psi$};
\draw[->][line width=1.5pt] (2,-0.2) -- (2, -1.8)
node[pos=0.5,right] {$\psi$};

\end{tikzpicture} 
\end{center}

\end{figure}

\end{itemize}

\end{definition}

Let's see an example of the Definition \ref{equivalencia derecha QP}. Let's consider the QP $(Q,0)$, where $Q$ is the following quiver.

\begin{figure}[H]
\begin{center}
\begin{tikzpicture}[scale=1]

\draw (0,0)node{$3$};
\draw (2,-2)node{$2$};
\draw (-2,-2)node{$1$};

\draw[->][line width=1.5pt] (-1.75,-1.9) -- (-0.1, -0.15)
node[pos=0.65,above] {$c$}; 

\draw[->][line width=1.5pt] (-1.6,-2) -- (1.75, -2)
node[pos=0.5,below] {$a$};

\draw[<-][line width=1.5pt] (1.9,-1.85) -- (0.1, -0.15)
node[pos=0.65,above] {$b$};

\end{tikzpicture} 
\end{center}
\caption{Quiver \label{ed}}

\end{figure}

For any $\lambda$ in $K$, the QP-representations of Figure \ref{ed} are right equivalents
\begin{figure}[htp]
\begin{center}
\subfigure[]{
\hspace*{-1em}{\begin{tikzpicture}[scale=1]

\draw (0,0)node{$K$};
\draw (2,-2)node{$K$};
\draw (-2,-2)node{$K$};

\draw[->][line width=1.5pt] (-1.75,-1.9) -- (-0.1, -0.15)
node[pos=0.65,above] {$1_{K}$}; 

\draw[->][line width=1.5pt] (-1.6,-2) -- (1.75, -2)
node[pos=0.5,below] {$\lambda 1_{K}$};

\draw[<-][line width=1.5pt] (1.9,-1.85) -- (0.1, -0.15)
node[pos=0.65,above] {$1_{K}$};

\end{tikzpicture} }
}
\subfigure[]{
\hspace*{0.6em}{\begin{tikzpicture}[scale=1]

\draw (0,0)node{$K$};
\draw (2,-2)node{$K$};
\draw (-2,-2)node{$K$};

\draw[->][line width=1.5pt] (-1.75,-1.9) -- (-0.1, -0.15)
node[pos=0.65,above] {$1_{K}$}; 

\draw[->][line width=1.5pt] (-1.6,-2) -- (1.75, -2)
node[pos=0.5,below] {$0$};

\draw[<-][line width=1.5pt] (1.9,-1.85) -- (0.1, -0.15)
node[pos=0.65,above] {$1_{K}$};

\end{tikzpicture}}

} 
\end{center}
\caption{Right equivalence representations \label{ed}}

\end{figure}
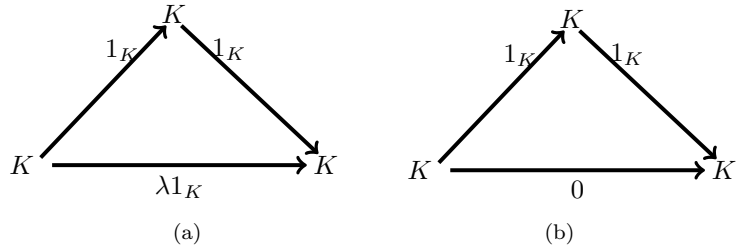

with the pair $(\varphi,\psi)$, where $\varphi:R \langle \langle Q \rangle \rangle \longrightarrow {}{R \langle \langle Q \rangle \rangle}$ is the isomorphism of $R$-algebras with the action in the arrows  given by $a\longmapsto \lambda bc$, $b\longmapsto b$, $c\longmapsto c$ and $\psi$ is the identity in each copy of $K$. This example shows that there are representations which are right equivalents but are not isomorphic. 

Remembering that each QP is right equivalent to the direct sum of its reduced part and its trivial part, and it is well defined up to right equivalence (Theorem \ref{split}). Let $(Q,S)$ be a QP and $\varphi:R \langle \langle (Q_{red\oplus C}) \rangle \rangle \longrightarrow {}{R \langle \langle Q \rangle \rangle}$ a right equivalence between $(Q_{red},S_{red}) \oplus (C,T)$ and $(Q,S)$, where $(Q_{red},S_{red})$ is a reduced QP and $(C,T)$ is a trivial QP. Let ${\mathcal M}=(Q,S,M)$ be a representation and set $M^{\varphi}=M$ as $K$-vector space. It is define the action of $R \langle \langle Q_{red} \rangle \rangle $ on $M^{\varphi}$ setting $u_{M^{\varphi}}=\varphi(u)_{M}$ for all elements $u\in R \langle \langle Q_{red} \rangle \rangle $.

\begin{proposition}
(\cite[Proposition 4.5]{dwz}) With the action of $R \langle \langle Q_{red} \rangle \rangle $ on $M^{\varphi}$ just defined, the third $(Q_{red},S_{red},M^{\varphi})$ becomes a QP representation. Moreover, the right-equivalence class of $(Q_{red},S_{red},M^{\varphi})$ is determined by the right equivalence class of ${\mathcal M}$.

\end{proposition}

\begin{definition}
(\cite[Definition 10.4]{dwz})The (right equivalence class of the) QP representation ${\mathcal M}_{red}=(Q_{red},S_{red},M^{\varphi})$ is the reduced part of ${\mathcal M}$.
\end{definition}

\begin{remark}\label{nota 5.0.15}
The construction of a right equivalence between a QP ($Q,S$) and the direct sum of a reduced QP with a trivial is not given by a canonical procedure in any obvious way. That is, there is no canonical way to construct $(Q_{red},S_{red})$ nor a right-equivalence $(Q,S) \longrightarrow {}{(Q_{red},S_{red})\oplus(Q_{triv},S_{triv})}$ (even when the QP is $(Q_{red},S_{red})$ is already known). A specific right equivalence $\varphi$ is defined in \cite{dwz}, satisfying the following condition. $\varphi$  acts as the identity on all the arrows of $Q$ that do not appear in the component of degree 2 (denoted $S^{(2)}$) of $S$ as long as no arrow appearing in $S^{(2)}$ appears in different summands of $S^{(2)}$. Using this property of $\varphi$ , it is easy to see that given a representation ${\mathcal M}=(Q,S,M)$, the action of $R \langle \langle Q_{red} \rangle \rangle $ on $M^{\varphi}$ coincides with its action on $M$ induced by the inclusion of quivers $Q_{red} \hookrightarrow Q$. That is, restricting the action of $R \langle \langle Q \rangle \rangle $ on $M$ to its subalgebra $R \langle \langle Q_{red} \rangle \rangle $ gives us the reduced part ${\mathcal M}_{red}$.

\end{remark}

\section{String representation $m(\tau,i)$.}

\setpapersize{A4} 
\setmargins{2.5cm} 
{1.5cm} 
{0.5cm} 
{23.42cm} 
{10pt} 
{1cm} 
{0pt} 
{2cm} 

The aim of this article is to extend the arc representation defined in $\cite{lf2}$ to a context of tagged arcs and tagged triangulations. The Sections 3 and 4 describe the concepts which we need defined to enunciate and prove the main theorem of this paper, that is:

\begin{theorem}\label{relaciones}
If $\tau$ is a tagged triangulation of a bordered surface with marked points $(\Sigma, M)$ and $i$ is a tagged arc which does not belong to $\tau$, then the arc representation $M(\tau,i)$ satisfies the Jacobian relatios and the module $M(\tau,i)$is nilpotent.
\end{theorem} 

Then we will start with the definitions and constructions that we need.

\begin{definition}
Let $\tau$ be a tagged triangulation of a bordered surface with marked points $(\Sigma,M)$ with non negative signature and $\Delta$ a non self-folded triangle of $\tau^{\circ}$. We say  that  $\Delta$ is a triangle of;
\begin{itemize}
\item{} {\it Type 1} if any side of $\Delta$ is side of a self-folded triangle.
\item{} {\it Type 2} if $\Delta$ shares exactly one side with a self-folded triangle $\Delta'$.
\item{} {\it Type 3} if $\Delta$ shares exactly two sides with two different self-folded triangles $\Delta'$ and $\Delta''$. 
\end{itemize}
\end{definition}

\begin{figure}[H]
\centering 
\begin{tikzpicture}[scale=0.7]

\filldraw [black] (-2,0) circle (2pt)
(-4, -2) circle (2pt)
(-2,-4) circle (2pt);

\draw (-3,-2)node{$\Delta$};

\draw[color=blue][line width=1pt] (-2,0) -- (-4, -2);
\draw[color=blue][line width=1pt] (-2,-4) -- (-4, -2);
\draw[color=blue][line width=1pt] (-2,0) -- (-2, -4);



\filldraw [black] (0,0) circle (2pt)
(0,-4) circle (2pt)
(0,-2) circle (2pt)
(0,0.5) circle (1pt)
(0.3,0.5) circle (1pt)
(-0.3,0.5) circle (1pt)
(0,-4.5) circle (1pt)
(0.3,-4.5) circle (1pt)
(-0.3,-4.5) circle (1pt);
\draw (0,-3)node{$\Delta$};
\draw (-0.25,-1.5)node{$\Delta'$};

\draw[line width=1pt] (0,0) -- (-1, 1);
\draw[line width=1pt] (0,0) -- (1, 1);
\draw[line width=1pt] (0,-4) -- (1, -5);
\draw[line width=1pt] (0,-4) -- (-1, -5);
\draw[line width=1pt] (0,0) -- (0, -2);
\draw[color=blue][line width=1pt] (0,0) .. controls(-0.7,-1) and (-0.7,-2.2) .. (0, -2.5);
\draw[color=blue][line width=1pt] (0,0) .. controls(0.7,-1) and (0.7,-2.2) .. (0, -2.5);
\draw[color=blue][line width=1pt] (0,0) .. controls(-2, -0.5) and (-2,-3.5) .. (0, -4);
\draw[color=blue][line width=1pt] (0,0) .. controls(2, -0.5) and (2,-3.5) .. (0, -4);

\filldraw [black] (5,0) circle (2pt)
(3.5,-2) circle (2pt)
(6.5,-2) circle (2pt)
(5,0.5) circle (1pt)
(5.3,0.5) circle (1pt)
(4.7,0.5) circle (1pt);
\draw (5,-3)node{$\Delta$};
\draw (4.25,-1.5)node{$\Delta'$};
\draw (5.75,-1.5)node{$\Delta''$};
\draw[line width=1pt] (5,0) -- (4, 1);
\draw[line width=1pt] (5,0) -- (6, 1);
\draw[line width=1pt] (5,0) -- (3.5,-2);
\draw[line width=1pt] (5,0) -- (6.5,-2);
\draw[color=blue][line width=1pt] (5, 0) .. controls(1,-0.1) and (1,-3.9) .. (5, -4);
\draw[color=blue][line width=1pt] (5, 0) .. controls(9,-0.1) and (9,-3.9) .. (5, -4);
\draw[color=blue][line width=1pt] (5, 0) .. controls(3.5,-0.5) and (3,-2) .. (3.25, -2.25);
\draw[color=blue][line width=1pt] (5, 0) .. controls(5.05,-0.5) and (4.25,-2.5) .. (3.25, -2.25);
\draw[color=blue][line width=1pt] (5, 0) .. controls(5,-1) and (5.7,-2.5) .. (6.75, -2.25);
\draw[color=blue][line width=1pt] (5, 0) .. controls(5,0) and (7.25,-1.5) .. (6.75, -2.25);

\end{tikzpicture}
\caption{Triangles of type 1,2 and 3.}
\end{figure}
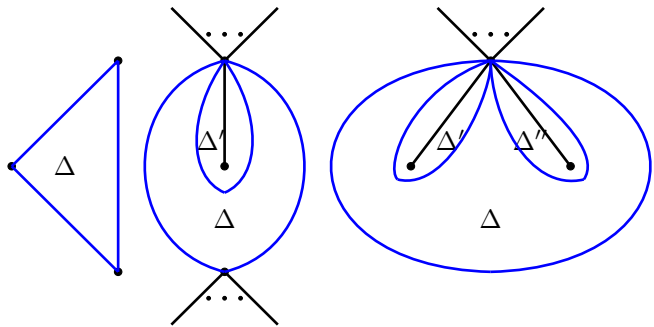

\begin{figure}[H]
\begin{center}
\subfigure[]{
\hspace*{-1em}{\begin{tikzpicture}[scale=0.7]

\filldraw [black] (5,0) circle (2pt)
(3.5,-2) circle (2pt)
(6.5,-2) circle (2pt)
(5,0.5) circle (1pt)
(5.3,0.5) circle (1pt)
(4.7,0.5) circle (1pt);

\draw[line width=1pt] (5,0) -- (4, 1);
\draw[line width=1pt] (5,0) -- (6, 1);
\draw[line width=1pt] (5,0) -- (3.5,-2);
\draw[line width=1pt] (5,0) -- (6.5,-2);
\draw[line width=1pt] (5, 0) .. controls(1,-0.1) and (1,-3.9) .. (5, -4);
\draw[line width=1pt] (5, 0) .. controls(9,-0.1) and (9,-3.9) .. (5, -4);
\draw[line width=1pt] (5, 0) .. controls(3.5,-0.5) and (3,-2) .. (3.25, -2.25);
\draw[line width=1pt] (5, 0) .. controls(5.05,-0.5) and (4.25,-2.5) .. (3.25, -2.25);
\draw[line width=1pt] (5, 0) .. controls(5,-1) and (5.7,-2.5) .. (6.75, -2.25);
\draw[line width=1pt] (5, 0) .. controls(5,0) and (7.25,-1.5) .. (6.75, -2.25);

\draw[->][color=red][line width=1.5pt] (5,-3.9) -- (3.25, -2.35);
\draw[->][color=red][line width=1.5pt] (3.5, -2.35) ..controls (4, -2.75) and (6, -2.75) .. (6.45, -2.35);
\draw[->][color=red][line width=1.5pt] (6.75, -2.35) -- (5.1,-3.9);
\draw[->][color=red][line width=1.5pt] (4.9,-3.9) ..controls(2,-3.9) and (2, -1.25) .. (4, -1.25)
node[pos=0.8,left] {$\alpha_{1}$};
\draw[->][color=red][line width=1.5pt] (4.15, -1.25) -- (5.85, -1.25)
node[pos=0.5,above] {$\beta$};
\draw[->][color=red][line width=1.5pt] (6.1, -1.25) ..controls (8.5, -1.25) and (7.6,-3.9) .. (5.4,-3.8)
node[pos=0.2,right] {$\alpha_{2}$};
\draw[->][color=red][line width=1.5pt] (4.15, -1.35) -- (6.1, -2.35)
node[pos=0.85,left] {$\alpha_{3}$};
\draw[->][color=red][line width=1.5pt] (3.7, -2.35) -- (5.9, -1.45)
node[pos=0.15,right] {$\alpha_{4}$};

\end{tikzpicture} }
}
\subfigure[]{
\hspace*{0.6em}{\begin{tikzpicture}[scale=0.7]

\filldraw [black] (0,0) circle (2pt)
(0,-4) circle (2pt)
(0,-2) circle (2pt)
(0,0.5) circle (1pt)
(0.3,0.5) circle (1pt)
(-0.3,0.5) circle (1pt)
(0,-4.5) circle (1pt)
(0.3,-4.5) circle (1pt)
(-0.3,-4.5) circle (1pt);

\draw[line width=1pt] (0,0) -- (-1, 1);
\draw[line width=1pt] (0,0) -- (1, 1);
\draw[line width=1pt] (0,-4) -- (1, -5);
\draw[line width=1pt] (0,-4) -- (-1, -5);
\draw[line width=1pt] (0,0) -- (0, -2);
\draw[line width=1pt] (0,-4) -- (0, -2);
\draw[line width=1pt] (0,0) .. controls(-2, -0.5) and (-2,-3.5) .. (0, -4);
\draw[line width=1pt] (0,0) .. controls(2, -0.5) and (2,-3.5) .. (0, -4);

\draw[->][color=red][line width=1.5pt] (0.15,-3) ..controls(0.3,-2) and (0.3,-2).. (0.15, -1);
\draw[->][color=red][line width=1.5pt] (0.2, -0.9) -- (1.4,-1.9);
\draw[<-][color=red][line width=1.5pt] (0.2,-3.1) -- (1.4,-2.1) ; 

\draw[<-][color=red][line width=1.5pt] (-0.15,-3) ..controls(-0.3,-2) and (-0.3,-2).. (-0.15, -1);
\draw[<-][color=red][line width=1.5pt] (-0.2, -0.9) -- (-1.4,-1.9);
\draw[->][color=red][line width=1.5pt] (-0.2,-3.1) -- (-1.4,-2.1) ;

\end{tikzpicture}}

} 
\end{center}
\caption{The quiver $Q(\tau^{\circ})$ drawn in $(\Sigma,M).$\label{piezas}}

\end{figure}
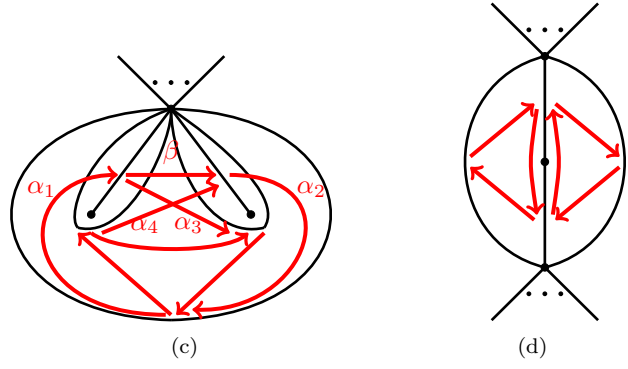

Now, we are going to draw the quiver $\widehat{Q}(\tau^{\circ})$ on $(\Sigma,M)$. According to \cite[Remark 4.2]{fst} every ideal triangulation can be obtained gluing puzzle pieces which appear in Figure \ref{piezas} (see \cite[Remark 4.2]{fst}). So it is enough to draw the quiver in each one of those pieces.

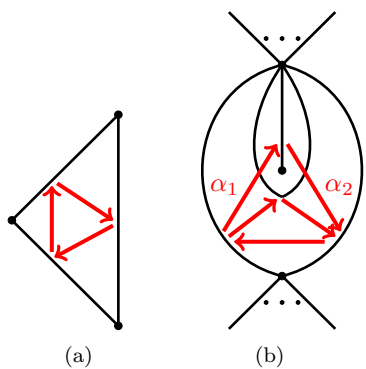
\begin{figure}[H]
\begin{center}
\subfigure[]{
\hspace*{-1em}{\begin{tikzpicture}[scale=0.7]

\filldraw [black] (-2,0) circle (2pt)
(-4, -2) circle (2pt)
(-2,-4) circle (2pt);

\draw[line width=1pt] (-2,0) -- (-4, -2);
\draw[line width=1pt] (-2,-4) -- (-4, -2);
\draw[line width=1pt] (-2,0) -- (-2, -4);


\draw[->][color=red][line width=1.5pt] (-3.25,-2.6) -- (-3.25, -1.35);
\draw[->][color=red][line width=1.5pt] (-3.15, -1.3) -- (-2.1,-2);
\draw[->][color=red][line width=1.5pt] (-2.1,-2.1) -- (-3.2,-2.7) ;

\end{tikzpicture} }
}
\subfigure[]{
\hspace*{0.6em}{\begin{tikzpicture}[scale=0.7] 

\filldraw [black] (0,0) circle (2pt)
(0,-4) circle (2pt)
(0,-2) circle (2pt)
(0,0.5) circle (1pt)
(0.3,0.5) circle (1pt)
(-0.3,0.5) circle (1pt)
(0,-4.5) circle (1pt)
(0.3,-4.5) circle (1pt)
(-0.3,-4.5) circle (1pt);

\draw[line width=1pt] (0,0) -- (-1, 1);
\draw[line width=1pt] (0,0) -- (1, 1);
\draw[line width=1pt] (0,-4) -- (1, -5);
\draw[line width=1pt] (0,-4) -- (-1, -5);
\draw[line width=1pt] (0,0) -- (0, -2);
\draw[line width=1pt] (0,0) .. controls(-0.7,-1) and (-0.7,-2.2) .. (0, -2.5);
\draw[line width=1pt] (0,0) .. controls(0.7,-1) and (0.7,-2.2) .. (0, -2.5);
\draw[line width=1pt] (0,0) .. controls(-2, -0.5) and (-2,-3.5) .. (0, -4);
\draw[line width=1pt] (0,0) .. controls(2, -0.5) and (2,-3.5) .. (0, -4);


\draw[->][color=red][line width=1.5pt] (-1,-3.25) -- (-0.1, -2.55);
\draw[->][color=red][line width=1.5pt] (0, -2.55) -- (1,-3.25);
\draw[->][color=red][line width=1.5pt] (0.8,-3.35 )-- (-0.9,-3.35);
\draw[->][color=red][line width=1.5pt] (-1.1,-3.15) -- (-0.1, -1.5)
node[pos=0.5,left] {$\alpha_{1}$};
\draw[->][color=red][line width=1.5pt] (0.1, -1.5) -- (1.1,-3.15)
node[pos=0.5,right] {$\alpha_{2}$};

\end{tikzpicture}}

} 
\end{center}

\end{figure}

\begin{remark}The next two remarks are about the quiver $\widehat{Q}({\tau^{\circ}})$ drawn on $(\Sigma,M)$.
\begin{myEnumerate}
\item[i)] The arrows $\alpha_{i}$ of $\widehat{Q}(\tau^{\circ})$ in the configuration (b) and (c) of the Figure \ref{piezas} have  non empty intersection with exactly one arc of $\tau^{\circ}$. 
\item[ii)]The arrow $\beta$ of $\widehat{Q}(\tau^{\circ})$ in the configuration (c) of the Figure \ref{piezas} has non empty intersection with exactly two arcs of $\tau^{\circ}$. 
\end{myEnumerate}
\end{remark}

Now, for each triangle $\Delta$ of type 1 or 2 we will assign a label to each side of $\Delta$, this is useful to define when a curve $\gamma$ surrounds a puncture.

\begin{itemize}
\item[a)] Let $\Delta$ be a triangle of type 2, we denote by $m$ the unique side of $\Delta$ which is a side of a self-folded triangle $\Delta'$.
We assign the label $l_{1}$ to $m$ in $\Delta$, and going through the sides of $\Delta$ in  clockwise,  we assign the labels $l_{2}$ and $l_{3}$ in $\Delta$ respectively.

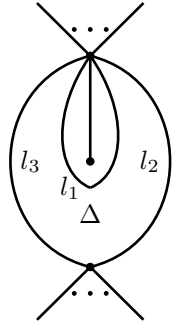
\begin{figure}[H]
\centering 
\begin{tikzpicture}[scale=0.7]


\filldraw [black] (0,0) circle (2pt)
(0,-4) circle (2pt)
(0,-2) circle (2pt)
(0,0.5) circle (1pt)
(0.3,0.5) circle (1pt)
(-0.3,0.5) circle (1pt)
(0,-4.5) circle (1pt)
(0.3,-4.5) circle (1pt)
(-0.3,-4.5) circle (1pt);
\draw (0,-3)node{$\Delta$};

\draw[line width=1pt] (0,0) -- (-1, 1);
\draw[line width=1pt] (0,0) -- (1, 1);
\draw[line width=1pt] (0,-4) -- (1, -5);
\draw[line width=1pt] (0,-4) -- (-1, -5);
\draw[line width=1pt] (0,0) -- (0, -2);
\draw[line width=1pt] (0,0) .. controls(-0.7,-1) and (-0.7,-2.2) .. (0, -2.5)
node[pos=1,left] {$l_{1}$};
\draw[line width=1pt] (0,0) .. controls(0.7,-1) and (0.7,-2.2) .. (0, -2.5);
\draw[line width=1pt] (0,0) .. controls(-2, -0.5) and (-2,-3.5) .. (0, -4)
node[pos=0.5,right] {$l_{3}$};
\draw[line width=1pt] (0,0) .. controls(2, -0.5) and (2,-3.5) .. (0, -4)
node[pos=0.5,left] {$l_{2}$};

\end{tikzpicture}
\caption{Labels of a triangle of type 2.}
\end{figure} 

\item[b)] Let $\Delta$ be a triangle of type 3, we denote by $m$ the unique side of $\Delta$ which is not a side of any self-folded triangle $\Delta'$.
We assign the label $l_{1}$ to $m$ in $\Delta$ and going through the sides of $\Delta$ in  clockwise, we assign the labels $l_{2}$ and $l_{3}$ in $\Delta$ respectively.

\begin{figure}[H]
\centering 
\vspace{0.2em}\hspace*{-3em}{\begin{tikzpicture}[scale=0.7]

\filldraw [black] (0,0) circle (2pt)
(-1.5,-2) circle (2pt)
(1.5,-2) circle (2pt)
(0,0.5) circle (1pt)
(0.3,0.5) circle (1pt)
(-0.3,0.5) circle (1pt);
\draw (0,-3)node{$\Delta$};
\draw[line width=1pt] (0,0) -- (-1, 1);
\draw[line width=1pt] (0,0) -- (1, 1);
\draw[line width=1pt] (0,0) -- (-1.5,-2);
\draw[line width=1pt] (0,0) -- (1.5,-2);
\draw[line width=1pt] (0, 0) .. controls(-4,-0.1) and (-4,-3.9) .. (0, -4)
node[pos=1,above] {$l_{1}$};;
\draw[line width=1pt] (0, 0) .. controls(4,-0.1) and (4,-3.9) .. (0, -4);
\draw[line width=1pt] (0, 0) .. controls(-1.5,-0.5) and (-2,-2) .. (-1.75, -2.25);
\draw[line width=1pt] (0, 0) .. controls(0.05,-0.5) and (-0.75,-2.5) .. (-1.75, -2.25)
node[pos=0.8,right] {$l_{2}$};
\draw[line width=1pt] (0, 0) .. controls(0,-1) and (0.7,-2.5) .. (1.75, -2.25)
node[pos=0.8,left] {$l_{3}$};
\draw[line width=1pt] (0, 0) .. controls(0,0) and (2.25,-1.5) .. (1.75, -2.25);

\end{tikzpicture}} 
\caption{Labels of a triangle of type 3.}

\end{figure}
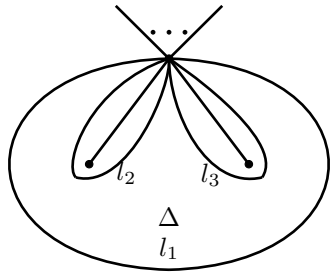 
\end{itemize}

Now, we are going to define the concept of surround a puncture.

\begin{definition}\label{rodear}
Let $\gamma$ be a curve in $(\Sigma,M)$ with the followings properties;

\begin{itemize}
\item[1)] $\gamma$ is a non selfintersecting oriented curve on $\Sigma$,
\item[2)] The intersection of $\gamma$ and $M$ is empty,
\item[3)] $\gamma$ minimizes the crossing points with the arcs of $\tau^{\circ}$ in his isotopy class,
\item[4)] The extremes $q_{0}$ and $q_{1}$ of $\gamma$ lie in the interior of an arc $j$ which belong to $\tau^{\circ}$,
\item[5)]$\gamma$ is oriented from $q_{0}$ to $q_{1}$.
\end{itemize}
In this case we divide $\gamma$ into segments $[r_{0},r_{1}]_{\gamma}$,$[r_{1},r_{2}]_{\gamma}$,...,$[r_{l},r_{l+1}]_{\gamma}$ with $r_{0}=q_{0}$ and $r_{l+1}=q_{1}$ such that $[r_{k},r_{k+1}]_{\gamma} \cap \tau^{\circ}=\{r_{k},r_{k+1}\}$ for $k \in \{0,...,l\}$. Given a puncture $p\in M$ we will say that $\gamma$ {\it surrounds} $p$ {\it with respect} $\tau^{\circ}$ (in counterclockwise), if $\gamma \cup [q_{1},q_{0}]_{j}$ is a closed curve, oriented in counterclockwise sense and it is one of the two conditions described below is satisfied.

\begin{itemize}
\item[a)] $\delta_{\tau}(p)\neq 0$.

Depending on if the arc $j$ is the folded side of a self-folded triangle or not, we are going to consider the following two cases $a_{1})$ and $a_{2})$.

\begin{itemize}
\item[$a_{1})$] $j$ is not a side of any self-folded triangle.

If the arc $j$ is not a side of any self-folded triangle then $j$ could be a side of a triangle of type $1$,$2$ or $3$. We are going to consider the three possibilities mentioned above. 
\begin{itemize}
\item[$a_{11})$]The arc $j$ is a side of a triangle $\Delta$ of type $1$ and the next four conditions $i)$-$iv)$ are met. See Figure \ref{figura a11}.

\begin{myEnumerate}
\item[i)] $p$ is the vertex of $\Delta$ opposite to $j$, 
\item[ii)] For $k=1,..,l-1$ the segments $[r_{k},r_{k+1}]_{\gamma}$ are contractible to the puncture $p$, with the homotopy that avoids $M$ and each of whose intermediate maps are segments with endpoints in the arcs of $\tau^{\circ}$ to which $r_{k}$ and $r_{k+1}$ belong.

\item[iii)] The segments $[q_{0},r_{1}]_{\gamma}$ and $[r_{l},q_{1}]_{\gamma}$ are contained in $\Delta$,
\item[iv)] $\gamma \cup [q_{0},q_{1}]_{j}$ divides $\Sigma$ in two regions and one of those two regions is homeomorphic to a disk that contains one puncture namely $p$.
\begin{figure}[H]
\centering 
\begin{tikzpicture}[scale=0.8]

\filldraw [black] (3,0) circle (2pt)
(3,-4) circle (2pt)
(-0.5,-2) circle (2pt)
(-1.5,-2) circle (1pt)
(-1.5,-2.5) circle (1pt)
(-1.5,-1.5) circle (1pt);
\draw (2,-2)node{$\Delta$};
\draw[line width=1pt] (3, 0) -- (3, -4)
node[pos=0.5,left] {$j$};
\draw[line width=1pt] (3, 0) -- (-0.5, -2);
\draw[line width=1pt] (3, -4) -- (-0.5, -2);
\draw[line width=1pt] (-0.5, -2) -- (-0.5,0);
\draw[line width=1pt] (-0.5, -2) -- (-2.5,-0.5);
\draw[line width=1pt] (-0.5, -2) -- (-2.5,-3.5);
\draw[line width=1pt] (-0.5, -2) -- (-0.5,-4)
node[pos=0.25,left] {$p$};

\filldraw [red] (3,-1) circle (2pt)
(3,-3) circle (2pt)
(1.65,-3.2) circle (2pt)
(-0.5,-3.3) circle (2pt)
(-1.8,-3) circle (2pt)
(-1.8,-1) circle (2pt)
(-0.5,-0.45) circle (2pt)
(1.95,-0.6) circle (2pt);

\draw[color=red][line width=1.5pt] (3,-1) .. controls(0, 0.5) and (-2.5,-1) .. (-2.5,-2)
node[pos=0,right] {$q_{0}=r_{0}$}
node[pos=0.11,above] {$r_{1}$}
node[pos=0.4,above] {$r_{2}$}
node[pos=0.77,above] {$r_{3}$}; 
\draw[color=red][line width=1.5pt] (3,-3) .. controls(0, -3.5) and (-2.5,-3.5) .. (-2.5,-2)
node[pos=0,right] {$q_{1}=r_{l+1}$}
node[pos=0.15,below] {$r_{l}$}
node[pos=0.38,below] {$r_{l-1}$}
node[pos=0.65,below] {$r_{l-2}$}; 

\end{tikzpicture}
\caption{Case $a_{11}$) of the Definition \ref{rodear}.}\label{figura a11}
\end{figure} 

\end{myEnumerate}
\item[$a_{12})$]$j$ is a side of a triangle $\Delta$ of type $2$. We denote by $m$ the side of $\Delta$ which is a non folded side of $\Delta'$, where $\Delta'$ is the unique self-folded triangle which shares the side $m$ with $\Delta$. In this case the following five conditions $i)$- $v)$ are met. See Figure \ref{figura a12}.

\begin{myEnumerate}

\item[i)] $m$ is based on the puncture $p$,
\item[ii)] If $r_{0}$ belongs to the side labeled with $l_{3}$ in $\Delta$ then $r_{1}$ belongs to the side with label $l_{2}$ in $\Delta$, $r_{l-3}$ belongs to the side with label $l_{3}$, $r_{l-1}$ belongs to the folded side of $\Delta'$ ,$r_{l}$ belongs to $m$ and the segment $[r_{l-3},r_{l-1}]_{\gamma}$ is contractible to $p$. On the other hand if $r_{0}$ belongs to the side labeled $l_{2}$ in $\Delta$ then $r_{4}$ belongs to the side with label $l_{2}$ in $\Delta$, $r_{1}$ belongs to $m$, $r_{2}$ belongs to the folded side of $\Delta'$, $r_{l}$ belongs to the side with label $l_{3}$ in $\Delta$ and the segment $[r_{2},r_{4}]_{\gamma}$ is contractible to $p$, 
\item[iii)] The segments $[q_{0},r_{1}]_{\gamma}$ and $[r_{l},q_{1}]_{\gamma}$ are contained in $\Delta$,
\item[iv)] For $k=1,..,l-1$ the segments $[r_{k},r_{k+1}]_{\gamma}$ are contractible to the puncture $p$, with the homotopy that avoids $M$ and each of whose intermediate maps are segments with endpoints in the arcs of $\tau^{\circ}$ to which $r_{k}$ and $r_{k+1}$ belong,

\item[v)]$\gamma \cup [q_{0},q_{1}]_{j}$ divides $\Sigma$ in two regions and one of those two regions is homeomorphic to a disk that contains one puncture namely $p$.
\end{myEnumerate}

\begin{figure}[H]
\centering 
$$ \begin{tikzpicture}[scale=0.8]

\filldraw [black] (0,0) circle (2pt)
(0,-4) circle (2pt)
(0,-2) circle (2pt)
(0,0.5) circle (1pt)
(0.3,0.5) circle (1pt)
(-0.3,0.5) circle (1pt)
(0,-4.5) circle (1pt)
(0.3,-4.5) circle (1pt)
(-0.3,-4.5) circle (1pt);
\draw (0,-3)node{$\Delta$};
\draw (4,-3.25)node{$\Delta$};
\draw[line width=1pt] (0,0) -- (-1, 1)
node[pos=0.3,right] {$p$};
\draw[line width=1pt] (0,0) -- (1, 1);
\draw[line width=1pt] (0,-4) -- (1, -5);
\draw[line width=1pt] (0,-4) -- (-1, -5);
\draw[line width=1pt] (0,0) -- (0, -2);
\draw[line width=1pt] (0,0) .. controls(-0.7,-1) and (-0.7,-2.2) .. (0, -2.5)
node[pos=0.5,left] {$m$};
\draw[line width=1pt] (0,0) .. controls(0.7,-1) and (0.7,-2.2) .. (0, -2.5);
\draw[line width=1pt] (0,0) .. controls(-2, -0.5) and (-2,-3.5) .. (0, -4)
node[pos=0.5,right] {$j$};
\draw[line width=1pt] (0,0) .. controls(2, -0.5) and (2,-3.5) .. (0, -4);

\filldraw [black] 
(4,0) circle (2pt)
(4,-4) circle (2pt)
(4,-2) circle (2pt)
(4,0.5) circle (1pt)
(4.3,0.5) circle (1pt)
(3.7,0.5) circle (1pt)
(4,-4.5) circle (1pt)
(4.3,-4.5) circle (1pt)
(3.7,-4.5) circle (1pt);
\draw[line width=1pt] (4,0) -- (3, 1)
node[pos=0.25,right] {$p$};
\draw[line width=1pt] (4,0) -- (5, 1);
\draw[line width=1pt] (4,-4) -- (5, -5);
\draw[line width=1pt] (4,-4) -- (3, -5);
\draw[line width=1pt] (4,0) -- (4, -2);
\draw[line width=1pt] (4,0) .. controls(2, -0.5) and (2,-3.5) .. (4, -4);
\draw[line width=1pt] (4,0) .. controls(6, -0.5) and (6,-3.5) .. (4, -4)
node[pos=0.5,left] {$j$};
\draw[line width=1pt] (4,0) .. controls(3.3,-1) and (3.3,-2.2) .. (4, -2.5)
node[pos=0.5,left] {$m$};
\draw[line width=1pt] (4,0) .. controls(4.7,-1) and (4.7,-2.2) .. (4, -2.5);



\filldraw [red] (4.55,-3.8) circle (2pt)
(5.25,-3) circle (2pt)
(3.55,-2) circle (2pt)
(4.4,-0.75) circle (2pt)
(4,-1.05) circle (2pt)
(4.9,-0.55) circle (2pt)
(4.8,0.75) circle (2pt)
(3.5,0.55) circle (2pt)
(2.65,-1.2) circle (2pt);
\draw[color=red][line width=1.5pt] (5.25,-3) .. controls(3,-3.5) and (3,-1) .. (5,-0.5)
node[pos=0,right] {$q_{0}=r_{0}$}
node[pos=0.5,right] {$r_{1}$}
node[pos=0.75,right] {$r_{2}$}
node[pos=0.93,below] {$r_{3}$}
node[pos=1.04,below] {$r_{4}$}; 
\draw[color=red][line width=1.5pt] (5,-0.5) .. controls(7,0.5) and (3,2) .. (2.65,-1)
node[pos=0.38,above] {$r_{5}$}
node[pos=0.7,above] {$r_{l-1}$}
node[pos=1,left] {$r_{l}$};
\draw[color=red][line width=1.5pt] (4.55,-3.8) .. controls(3.45,-3.8) and (2.5,-1.7) .. (2.65,-1)
node[pos=0,right] {$q_{1}=r_{l+1}$}; 

\filldraw [red] (-0.55,-3.8) circle (2pt)
(-1.25,-3) circle (2pt)
(1.3,-1.05) circle (2pt)
(0.7,0.7) circle (2pt)
(-0.75,0.75) circle (2pt)
(-0.9,-0.5) circle (2pt)
(-0.4,-0.75) circle (2pt)
(0,-1.05) circle (2pt)
(0,-1.05) circle (2pt)
(0.35,-2.1) circle (2pt);

\draw[color=red][line width=1.5pt] (-0.55,-3.8) .. controls(1,-3.5) and (1.5,-0.75) .. (1.25,-0.75)
node[pos=0,left] {$q_{0}=r_{0}$}
node[pos=0.9,right] {$r_{1}$};
\draw[color=red][line width=1.5pt] (1.3,-0.75) .. controls(1.25,2) and (-2.25,0.75) .. (-1,-0.5)
node[pos=0.25,right] {$r_{2}$}
node[pos=0.6,left] {$r_{l-4}$}
node[pos=1,left] {$r_{l-3}$};
\draw[color=red][line width=1.5pt] (-1.25,-3) .. controls(1,-3) and (0.75,-1) .. (-1,-0.5)
node[pos=0.85,left] {$r_{l-2}$}
node[pos=0.75,right] {$r_{l-1}$}
node[pos=0.3,right] {$r_{l}$}
node[pos=0,left] {$r_{l+1}=q_{1}$};

\end{tikzpicture}$$
\caption{Case $a_{12}$) of the Definition \ref{rodear}.}\label{figura a12}
\end{figure} 
\item[$a_{13})$]$j$ is a side of a triangle $\Delta$ of type $3$. We denote by $m$ and $m'$ the sides of $\Delta $ which are the not folded sides of $\Delta'$ and $\Delta''$. Where $\Delta'$ and $\Delta''$ are the two unique self-folded triangles which share the sides $m$ and $m'$ respectively. See Figure \ref{figura a13}.
\end{itemize}

\begin{myEnumerate} 
\item[i)] $m$ and $m'$ are based on the punctue $p$, 
\item[ii)] $r_{2}$ belongs to the folded side of $\Delta''$ and $r_{l-1}$ belongs to the folded side of $\Delta'$, 
\item[iii)] The segments $[q_{0},r_{1}]_{\gamma}$ and $[r_{l},q_{1}]_{\gamma}$ are contained in $\Delta$. 
\item[iv)] For $k=0,..,l$ the segments $[r_{k},r_{k+1}]_{\gamma}$ are contractible to the puncture $p$, with the homotopy that avoids $M$ and each of whose intermediate maps are segments with endpoints in the arcs of $\tau^{\circ}$ to which $r_{k}$ and $r_{k+1}$ belong, 
\item[v)] $\gamma \cup [q_{0},q_{1}]_{j}$ divides $\Sigma$ in two regions and one of those two regions is homeomorphic to a disk that contains one puncture namely $p$.
\end{myEnumerate}

\begin{figure}[H]
\centering 
$$ \begin{tikzpicture}[scale=0.8]
\filldraw [black] (0,0) circle (2pt)
(-1.5,-2) circle (2pt)
(1.5,-2) circle (2pt)
(0,0.5) circle (1pt)
(0.3,0.5) circle (1pt)
(-0.3,0.5) circle (1pt);
\draw (0,-3)node{$\Delta$};
\draw[line width=1pt] (0,0) -- (-1, 1)
node[pos=0.3,right] {$p$};
\draw[line width=1pt] (0,0) -- (1, 1);
\draw[line width=1pt] (0,0) -- (-1.5,-2);
\draw[line width=1pt] (0,0) -- (1.5,-2);
\draw[line width=1pt] (0, 0) .. controls(-4,-0.1) and (-4,-3.9) .. (0, -4);
\draw[line width=1pt] (0, 0) .. controls(4,-0.1) and (4,-3.9) .. (0, -4)
node[pos=1,below] {$j$};
\draw[line width=1pt] (0, 0) .. controls(-1.5,-0.5) and (-2,-2) .. (-1.75, -2.25);
\draw[line width=1pt] (0, 0) .. controls(0.05,-0.5) and (-0.75,-2.5) .. (-1.75, -2.25)
node[pos=0.8,right] {$m$}
node[pos=0.6,left] {$\Delta'$};
\draw[line width=1pt] (0, 0) .. controls(0,-1) and (0.7,-2.5) .. (1.75, -2.25)
node[pos=0.8,left] {$m'$}
node[pos=0.5,right] {$\Delta''$};
\draw[line width=1pt] (0, 0) .. controls(0,0) and (2.25,-1.5) .. (1.75, -2.25);


\filldraw [red] (-1,-3.9) circle (2pt)
(-1.25,-2.15) circle (2pt)
(-1.3,-1.7) circle (2pt)
(-1.23,-0.8) circle (2pt)
(-1.05,-0.15) circle (2pt)
(-0.65,0.6) circle (2pt)
(0.65,0.6) circle (2pt)
(1.05,-0.15) circle (2pt)
(1.23,-0.95) circle (2pt)
(1.3,-1.7) circle (2pt)
(1.25,-2.2) circle (2pt)
(1,-3.9) circle (2pt);

\draw[color=red][line width=1.5pt] (-1,-3.9) .. controls(-2.5,2.5) and (2.5,2.5) .. (1,-3.9)
node[pos=1,below] {$q_{0}=r_{0}$}
node[pos=0,below] {$q_{1}=r_{l+1}$}
node[pos=0.08,left] {$r_{l}$}
node[pos=0.14,left] {$r_{l-1}$}
node[pos=0.21,left] {$r_{l-2}$}
node[pos=0.29,left] {$r_{l-3}$}
node[pos=0.37,left] {$r_{l-4}$}
node[pos=0.62,right] {$r_{5}$}
node[pos=0.72,right] {$r_{4}$}
node[pos=0.82,right] {$r_{3}$}
node[pos=0.87,right] {$r_{2}$}
node[pos=0.92,right] {$r_{1}$};

\end{tikzpicture}$$
\caption{Case $a_{13}$) of the Definition \ref{rodear}.}\label{figura a13}
\end{figure}

\item[$a_{2})$] $j$ is the folded side of some self-folded triangle $\Delta'$. We denote by $m$ the not folded side of $\Delta'$.

For this case $m$ could be a side of a triangle of type $2$ or type $3$. We are going to consider these two subcases:
\begin{itemize}
\item[$a_{21})$]$m$ is a side of a triangle $\Delta$ of type 2, in this case the followings five conditions $i)$- $v)$ are met. See Figure \ref{figura a21}.

\begin{myEnumerate}
\item[i)] $p$ is vertex of $\Delta$ oppsite to $m$,
\item[ii)] $r_{2}$ belongs to the side of $\Delta$ with label $l_{3}$ in $\Delta$ and $r_{l-1}$ belongs to the side with label $l_{2}$ in $\Delta$,
\item[iii)] The segments $[r_{1},r_{2}]_{\gamma}$ and $[r_{l-1},r_{l}]_{\gamma}$ are contained in $\Delta$,
\item[iv)] For $k=2,..,l-2$ the segments $[r_{k},r_{k+1}]_{\gamma}$ are contractible to the puncture $p$, with the homotopy that avoids $M$ and each of whose intermediate maps are segments with endpoints in the arcs of $\tau^{\circ}$ to which $r_{k}$ and $r_{k+1}$ belong, 
\item[v)] $\gamma \cup [q_{0},q_{1}]_{j}$ divides $\Sigma$ in two regions and one of those two regions is homeomorphic to a disk that contains one puncture namely $p$.

\end{myEnumerate} 

\begin{figure}[H]
\centering 
$$ \begin{tikzpicture}[scale=0.8]

\filldraw [black] (0,0) circle (2pt)
(0,-4) circle (2pt)
(0,-2) circle (2pt)
(0,0.5) circle (1pt)
(0.3,0.5) circle (1pt)
(-0.3,0.5) circle (1pt)
(0,-4.5) circle (1pt)
(0.3,-4.5) circle (1pt)
(-0.3,-4.5) circle (1pt);
\draw (0,-3)node{$\Delta$};
\draw (4,-3.25)node{$\Delta$};
\draw (3.7,-1.45)node{$\Delta'$};
\draw (-0.25,-1.8)node{$\Delta'$};
\draw (4.7,-2)node{$m$};
\draw (4.2,-1)node{$j$};
\draw (-0.7,-1)node{$m$};
\draw (-0.2,-1)node{$j$};
\draw (0,-4.3)node{$p$};
\draw[line width=1pt] (0,0) -- (-1, 1);
\draw[line width=1pt] (0,0) -- (1, 1);
\draw[line width=1pt] (0,-4) -- (1, -5);
\draw[line width=1pt] (0,-4) -- (-1, -5);
\draw[line width=1pt] (0,0) -- (0, -2);
\draw[line width=1pt] (0,0) .. controls(-0.7,-1) and (-0.7,-2.2) .. (0, -2.5);
\draw[line width=1pt] (0,0) .. controls(0.7,-1) and (0.7,-2.2) .. (0, -2.5);
\draw[line width=1pt] (0,0) .. controls(-2, -0.5) and (-2,-3.5) .. (0, -4);
\draw[line width=1pt] (0,0) .. controls(2, -0.5) and (2,-3.5) .. (0, -4);

\filldraw [black] 
(4,0) circle (2pt)
(4,-4) circle (2pt)
(4,-2) circle (2pt)
(4,0.5) circle (1pt)
(4.3,0.5) circle (1pt)
(3.7,0.5) circle (1pt)
(4,-4.5) circle (1pt)
(4.3,-4.5) circle (1pt)
(3.7,-4.5) circle (1pt);
\draw (4,-4.3)node{$p$};
\draw[line width=1pt] (4,0) -- (3, 1);
\draw[line width=1pt] (4,0) -- (5, 1);
\draw[line width=1pt] (4,-4) -- (5, -5);
\draw[line width=1pt] (4,-4) -- (3, -5);
\draw[line width=1pt] (4,0) -- (4, -2);
\draw[line width=1pt] (4,0) .. controls(2, -0.5) and (2,-3.5) .. (4, -4);
\draw[line width=1pt] (4,0) .. controls(6, -0.5) and (6,-3.5) .. (4, -4);
\draw[line width=1pt] (4,0) .. controls(3.3,-1) and (3.3,-2.2) .. (4, -2.5);
\draw[line width=1pt] (4,0) .. controls(4.7,-1) and (4.7,-2.2) .. (4, -2.5);



\filldraw [red] (4,-0.5) circle (2pt)
(3.55,-0.85) circle (2pt)
(2.77,-3) circle (2pt)
(3.45,-4.55) circle (2pt)
(4.75,-4.75) circle (2pt)
(4.9,-3.5) circle (2pt)
(3.55,-1.9) circle (2pt)
(4,-1.5) circle (2pt) ;
\draw[color=red][line width=1.5pt] (4,-0.5) .. controls(3.25,-1) and (2.75,-2) .. (2.75,-2.9)
node[pos=0,right] {$q_{0}=r_{0}$}
node[pos=0.2,left] {$r_{1}$}
node[pos=1,left] {$r_{2}$};
\draw[color=red][line width=1.5pt] (2.75,-2.9) .. controls(3,-5.5) and (5.5,-5.5) .. (4.9,-3.5)
node[pos=0.3,left] {$r_{3}$}
node[pos=0.7,right] {$r_{l-2}$}
node[pos=1,right] {$r_{l-1}$};
\draw[color=red][line width=1.5pt] (4.9,-3.5) .. controls(4.8,-2.8) and (2.5,-2.5) .. (4,-1.5)
node[pos=0.8,left] {$r_{l}$}
node[pos=1,right] {$q_{1}=r_{l+1}$};


\filldraw [red] (0,-0.5) circle (2pt)
(0.4,-0.75) circle (2pt)
(1,-3.45) circle (2pt)
(0.55,-4.55) circle (2pt)
(-0.85,-4.85) circle (2pt)
(-0.95,-3.5) circle (2pt)
(0.45,-2) circle (2pt)
(0,-1.5) circle (2pt) ;

\draw[color=red][line width=1.5pt] (0,-0.5) .. controls(1,-1) and (1.5,-2) .. (1,-3.5)
node[pos=0.01,left] {$r_{l+1}=q_{1}$}
node[pos=0.1,right] {$r_{l}$}
node[pos=1,right] {$r_{l-1}$};
\draw[color=red][line width=1.5pt] (1,-3.5) .. controls(0.5,-6) and (-2,-5) .. (-1,-3.5)
node[pos=0.2,right] {$r_{l-2}$}
node[pos=0.6,left] {$r_{3}$}
node[pos=1,left] {$r_{2}$};
\draw[color=red][line width=1.5pt] (-1,-3.5) .. controls(-0.5,-3) and (1.25,-2.5) .. (0,-1.5)
node[pos=0.8,right] {$r_{1}$}
node[pos=1,left] {$r_{0}=q_{0}$};

\end{tikzpicture}$$
\caption{Case $a_{21}$) of the Definition \ref{rodear}.}\label{figura a21}
\end{figure} 
\item[$a_{22})$]$m$ is a side of a triangle $\Delta$ of type $3$. We denote by $m'$ the side of $\Delta$ which is different from $m$ that is not the folded side of a self-folded triangle $\Delta''$. In this case the following eleven conditions are met.

Under these conditions the arc $m$ has label $l_{2}$ o $m$ has label $l_{3}$ in $\Delta$. First we are going to consider that $m$ has label $l_{2}$ in $\Delta$. See Figure \ref{figura a22}.

\begin{myEnumerate}
\item[i)] $m$ and $m'$ are based on the puncture $p$,
\item[ii)] $r_{2},r_{l-7}$ belong to the arc with label $l_{1}$ in $\Delta$ and $r_{l},r_{l-4},\ m$ , $r_{l-5}\in j$,
\item[iii)] $[p,r_{r-5}]_{j}\cap[q_{0},q_{1}]_{j}=\emptyset$,
\item[iv)] $[r_{l},r_{l-4}]_{\gamma} \cup [r_{l-4},r_{l}]_{m}$ is a closed curve which cuts to the arc $m'$ in exactly two points and to the folded side of $\Delta''$ in exactly one point,
\item[v)] The segment $[r_{l},r_{l-4}]_{\gamma} \cup [r_{l-4},r_{l}]_{m}$ contains in the interior only the puncture which is enclosed by $m'$,
\item[vi)] The segment $[r_{l-2},q_{1}]_{\gamma}$ is not contractible to $p$,
\item[vii)] $\gamma \cup [q_{0},q_{1}]_{j}$ divides $\Sigma$ in two regions. One of those two regions is homeomorphic to a disk that contains one puncture namely $p$, and possibly to the one puncture which is enclosed by $m$,
\item[viii)] The segments $[r_{1},r_{2}]_{\gamma}$ and $[r_{l-1},r_{l}]_{\gamma}$ are contained in $\Delta$,
\item[ix)] For $k=0,..,l$ the segments $[r_{k},r_{k+1}]_{\gamma}$ are contractible to the puncture $p$, with the homotopy that avoids $M$ and each of whose intermediate maps are segments with endpoints in the arcs of $\tau^{\circ}$ to which $r_{k}$ and $r_{k+1}$ belong,
\item[x)]The interior of the segment $[r_{0},r_{2}]_{\gamma}$ only cuts $m$ in the point $r_{1}$, 
\item[xi)] The interior of the segment $[r_{l-7},r_{l-4}]_{\gamma}$ only cuts the arcs $m$ and $j$ in one point respectively. 
\end{myEnumerate}
\end{itemize}
\end{itemize}
\begin{figure}[H]
\centering 
\hspace*{-6em}{\begin{tikzpicture}[scale=0.8]

\filldraw [black] (0,0) circle (2pt)
(-1.5,-2) circle (2pt)
(1.5,-2) circle (2pt)
(0,0.5) circle (1pt)
(0.3,0.5) circle (1pt)
(-0.3,0.5) circle (1pt);
\draw (0,-3)node{$\Delta$};
\draw (-1,-1)node{$\Delta'$};
\draw (-1.5,-1.6)node{$j$};
\draw (0.8,-1.5)node{$\Delta''$};
\draw[line width=1pt] (0,0) -- (-1, 1)
node[pos=0.3,right] {$p$};
\draw[line width=1pt] (0,0) -- (1, 1);
\draw[line width=1pt] (0,0) -- (-1.5,-2);
\draw[line width=1pt] (0,0) -- (1.5,-2);
\draw[line width=1pt] (0, 0) .. controls(-4,-0.1) and (-4,-3.9) .. (0, -4);
\draw[line width=1pt] (0, 0) .. controls(4,-0.1) and (4,-3.9) .. (0, -4);
\draw[line width=1pt] (0, 0) .. controls(-1.5,-0.5) and (-2,-2) .. (-1.75, -2.25);
\draw[line width=1pt] (0, 0) .. controls(0.05,-0.5) and (-0.75,-2.5) .. (-1.75, -2.25)
node[pos=0.8,right] {$m$};
\draw[line width=1pt] (0, 0) .. controls(0,-1) and (0.7,-2.5) .. (1.75, -2.25)
node[pos=0.8,left] {$m'$};
\draw[line width=1pt] (0, 0) .. controls(0,0) and (2.25,-1.5) .. (1.75, -2.25);


\filldraw [black] (8,0) circle (2pt)
(6.5,-2) circle (2pt)
(9.5,-2) circle (2pt)
(8,0.5) circle (1pt)
(8.3,0.5) circle (1pt)
(7.7,0.5) circle (1pt);
\draw (8,-3)node{$\Delta$};
\draw (7.15,-1.6)node{$\Delta'$};
\draw (6.8,-1.3)node{$j$};
\draw (8.8,-1.5)node{$\Delta''$};
\draw[line width=1pt] (8,0) -- (7, 1)
node[pos=0.3,right] {$p$};
\draw[line width=1pt] (8,0) -- (9, 1);
\draw[line width=1pt] (8,0) -- (6.5,-2);
\draw[line width=1pt] (8,0) -- (9.5,-2);
\draw[line width=1pt] (8, 0) .. controls(4,-0.1) and (4,-3.9) .. (8, -4);
\draw[line width=1pt] (8, 0) .. controls(12,-0.1) and (12,-3.9) .. (8, -4);
\draw[line width=1pt] (8, 0) .. controls(6.5,-0.5) and (6,-2) .. (6.25, -2.25);
\draw[line width=1pt] (8, 0) .. controls(8.05,-0.5) and (7.25,-2.5) .. (6.25, -2.25)
node[pos=0.8,right] {$m$};
\draw[line width=1pt] (8, 0) .. controls(8,-1) and (8.7,-2.5) .. (9.75, -2.25)
node[pos=0.8,left] {$m'$};
\draw[line width=1pt] (8, 0) .. controls(8,0) and (10.25,-1.5) .. (9.75, -2.25);


\filldraw [red] (-0.95,-1.4) circle (2pt)
(-1.45,-1.1) circle (2pt)
(-3,-1.7) circle (2pt)
(0.85,0.9) circle (2pt)
(-0.75,0.8) circle (2pt)
(-0.95,-0.1) circle (2pt)
(-0.8,-0.4) circle (2pt)
(-0.5,-0.7) circle (2pt)
(-0.2,-0.8) circle (2pt)
(0.2,-0.9) circle (2pt)
(0.75,-1) circle (2pt)
(1.4,-1.2) circle (2pt)
(-1.8,-1.9) circle (2pt)
(-1.35,-1.7) circle (2pt);
\draw[color=red][line width=1.5pt] (-0.95,-1.4) ..controls(-1.1,-1) and (-2.5,-0.8) .. (-3,-1.7)
node[pos=-0,right] {{\scriptsize $q_{0}=r_{0}$}}
node[pos=0.4,above] {{\scriptsize $r_{1}$}}
node[pos=1,left] {{\scriptsize $r_{2}$}}; 
\draw[color=red][line width=1.5pt] (-3,-1.7) .. controls(-4,-3) and (-2,-5) .. (1.5,-4.15);
\draw[color=red][line width=1.5pt] (1.5,-4.15) .. controls(7,-2.5) and (-1,3.5) .. (-1,0)
node[pos=0.66,right] {$r_{3}$}
node[pos=0.9,left] {{\scriptsize $r_{l-8}$}}
node[pos=1,left] {{\scriptsize$r_{l-7}$}};
\draw[color=red][line width=1.5pt] (-1,0) .. controls(-0.8,-1) and (0.5,-0.8) .. (1.4,-1.2)
node[pos=1.01,right] {{\scriptsize $r_{l-1}$}};
\draw[color=red][line width=1.5pt] (1.4,-1.2) .. controls(5,-4) and (-5.5,-3.5) .. (-1.35,-1.7)
node[pos=0.95,left] {{\scriptsize $r_{l}$}}
node[pos=1,right] {{\scriptsize $q_{1}=r_{l+1}$}};

\filldraw [red] (7.25,-1) circle (2pt)
(9.15,-3.85) circle (2pt)
(8.85,0.9) circle (2pt)
(7.25,0.8) circle (2pt)
(7.05,-0.1) circle (2pt)
(7.2,-0.4) circle (2pt)
(7.5,-0.7) circle (2pt)
(7.8,-0.8) circle (2pt)
(8.2,-0.9) circle (2pt)
(8.75,-1) circle (2pt)
(9.4,-1.2) circle (2pt)
(7.45,-1.55) circle (2pt)
(6.75,-1.6) circle (2pt)
(6.85,-2.2) circle (2pt);
\draw[color=red][line width=1.5pt] (6.75,-1.6) .. controls(7,-4) and (9.3,-3.9) .. (9.5,-3.8)
node[pos=-0,left] {{\scriptsize $q_{0}$}}
node[pos=0.15,left] {{\scriptsize $r_{1}$}}
node[pos=1,below] {{\scriptsize $r_{2}$}};
\draw[color=red][line width=1.5pt] (9.5,-3.8) .. controls(15,-2.5) and (7,3.5) .. (7,0) 
node[pos=0.66,right] {$r_{3}$}
node[pos=0.9,left] {{\scriptsize $r_{l-8}$}}
node[pos=1,left] {{\scriptsize$r_{l-7}$}};
\draw[color=red][line width=1.5pt] (7,0) .. controls(7.2,-1) and (8.5,-0.8) .. (9.4,-1.2);
\draw[color=red][line width=1.5pt] (9.4,-1.2) .. controls(12,-3) and (8,-3.7) .. (7.25,-1) 
node[pos=0.9,below] {{\scriptsize $r_{l}$}}
node[pos=0.99,left] {{\scriptsize $q_{1}$}};

\end{tikzpicture}}

\end{figure} 
\begin{figure}[H]
\centering 
\vspace{-8em}\hspace*{-4em}{\begin{tikzpicture}[scale=0.8]

\filldraw [black] (0,0) circle (2pt)
(-1.5,-2) circle (2pt)
(1.5,-2) circle (2pt)
(0,0.5) circle (1pt)
(0.3,0.5) circle (1pt)
(-0.3,0.5) circle (1pt);
\draw (0,-3)node{$\Delta$};
\draw (-1.75,-1)node{$\Delta'$};
\draw (-1.15,-1.9)node{$j$};
\draw (0.8,-1.5)node{$\Delta''$};
\draw[line width=1pt] (0,0) -- (-1, 1)
node[pos=0.3,right] {$p$};
\draw[line width=1pt] (0,0) -- (1, 1);
\draw[line width=1pt] (0,0) -- (-1.5,-2);
\draw[line width=1pt] (0,0) -- (1.5,-2);
\draw[line width=1pt] (0, 0) .. controls(-4,-0.1) and (-4,-3.9) .. (0, -4);
\draw[line width=1pt] (0, 0) .. controls(4,-0.1) and (4,-3.9) .. (0, -4);
\draw[line width=1pt] (0, 0) .. controls(-1.5,-0.5) and (-2,-2) .. (-1.75, -2.25);
\draw[line width=1pt] (0, 0) .. controls(0.05,-0.5) and (-0.75,-2.5) .. (-1.75, -2.25)
node[pos=0.8,right] {$m$};
\draw[line width=1pt] (0, 0) .. controls(0,-1) and (0.7,-2.5) .. (1.75, -2.25)
node[pos=0.8,left] {$m'$};
\draw[line width=1pt] (0, 0) .. controls(0,0) and (2.25,-1.5) .. (1.75, -2.25);


\filldraw [black] (8,0) circle (2pt)
(6.5,-2) circle (2pt)
(9.5,-2) circle (2pt)
(8,0.5) circle (1pt)
(8.3,0.5) circle (1pt)
(7.7,0.5) circle (1pt);
\draw (8,-3)node{$\Delta$};
\draw (7.15,-1.6)node{$\Delta'$};
\draw (6.8,-1.3)node{$j$};
\draw (8.8,-1.5)node{$\Delta''$};
\draw[line width=1pt] (8,0) -- (7, 1)
node[pos=0.3,right] {$p$};
\draw[line width=1pt] (8,0) -- (9, 1);
\draw[line width=1pt] (8,0) -- (6.5,-2);
\draw[line width=1pt] (8,0) -- (9.5,-2);
\draw[line width=1pt] (8, 0) .. controls(4,-0.1) and (4,-3.9) .. (8, -4);
\draw[line width=1pt] (8, 0) .. controls(12,-0.1) and (12,-3.9) .. (8, -4);
\draw[line width=1pt] (8, 0) .. controls(6.5,-0.5) and (6,-2) .. (6.25, -2.25);
\draw[line width=1pt] (8, 0) .. controls(8.05,-0.5) and (7.25,-2.5) .. (6.25, -2.25)
node[pos=0.8,right] {$m$};
\draw[line width=1pt] (8, 0) .. controls(8,-1) and (8.7,-2.5) .. (9.75, -2.25)
node[pos=0.8,left] {$m'$};
\draw[line width=1pt] (8, 0) .. controls(8,0) and (10.25,-1.5) .. (9.75, -2.25);

\filldraw [red] (-0.75,-1) circle (2pt)
(-1.8 ,-1.55) circle (2pt)
(1.15,-3.85) circle (2pt)
(0.85,0.9) circle (2pt)
(-0.75,0.8) circle (2pt)
(-0.95,-0.1) circle (2pt)
(-0.8,-0.4) circle (2pt)
(-0.5,-0.7) circle (2pt)
(-0.2,-0.8) circle (2pt)
(0.2,-0.9) circle (2pt)
(0.75,-1) circle (2pt)
(1.4,-1.2) circle (2pt)
(-0.55,-1.5) circle (2pt)
(-1.25,-1.6) circle (2pt);
\draw[color=red][line width=1.5pt] (-1.25,-1.6) .. controls(-5,-1) and (-0.7,-3.9) .. (1.5,-3.8)
node[pos=-0,right] {{\scriptsize $q_{0}$}}
node[pos=0.06,above] {{\scriptsize $r_{1}$}}
node[pos=1,below] {{\scriptsize $r_{2}$}};
\draw[color=red][line width=1.5pt] (1.5,-3.8) .. controls(7,-2.5) and (-1,3.5) .. (-1,0) 
node[pos=0.66,right] {$r_{3}$}
node[pos=0.9,left] {{\scriptsize $r_{l-8}$}}
node[pos=1,left] {{\scriptsize$r_{l-7}$}};
\draw[color=red][line width=1.5pt] (-1,0) .. controls(-0.8,-1) and (0.5,-0.8) .. (1.4,-1.2)
node[pos=1.01,right] {{\scriptsize $r_{l-1}$}};
\draw[color=red][line width=1.5pt] (1.4,-1.2) .. controls(4,-3) and (0,-3.7) .. (-0.75,-1)
node[pos=0.9,below] {{\scriptsize $r_{l}$}}
node[pos=0.99,left] {{\scriptsize $q_{1}$}};

\filldraw [red] (7.25,-1) circle (2pt)
(6.55,-1) circle (2pt)
(9.15,-3.85) circle (2pt)
(8.85,0.9) circle (2pt)
(7.25,0.8) circle (2pt)
(7.05,-0.1) circle (2pt)
(7.2,-0.4) circle (2pt)
(7.5,-0.7) circle (2pt)
(7.8,-0.8) circle (2pt)
(8.2,-0.9) circle (2pt)
(8.75,-1) circle (2pt)
(9.4,-1.2) circle (2pt)
(7.05,-2) circle (2pt)
(6.75,-1.6) circle (2pt);
\draw[color=red][line width=1.5pt] (7.25,-1) .. controls(3,-1) and (7.3,-3.9) .. (9.5,-3.8)
node[pos=-0,right] {{\scriptsize $q_{0}$}}
node[pos=0.06,above] {{\scriptsize $r_{1}$}}
node[pos=1,below] {{\scriptsize $r_{2}$}};
\draw[color=red][line width=1.5pt] (9.5,-3.8) .. controls(15,-2.5) and (7,3.5) .. (7,0) 
node[pos=0.66,right] {$r_{3}$}
node[pos=0.9,left] {{\scriptsize $r_{l-8}$}}
node[pos=1,left] {{\scriptsize$r_{l-7}$}};
\draw[color=red][line width=1.5pt] (7,0) .. controls(7.2,-1) and (8.5,-0.8) .. (9.4,-1.2)
node[pos=1.01,right] {{\scriptsize $r_{l-1}$}};
\draw[color=red][line width=1.5pt] (9.4,-1.2) .. controls(12,-3) and (8,-3.7) .. (6.75,-1.6) 
node[pos=0.9,below] {{\scriptsize $r_{l}$}}
node[pos=0.99,left] {{\scriptsize $q_{1}$}};
\end{tikzpicture}} 
\caption{Case $a_{22}$) of the Definition \ref{rodear}.}\label{figura a22}

\end{figure} 
On the other hand, we consider the case when $m $ has label $l_{3}$ in $\Delta$. See Figure \ref{figura a22p}.

\begin{myEnumerate}
\item[i)] $m$ and $m'$ are based on the puncture $p$,
\item[ii)] $r_{8},r_{l-1}$ belongs to the arc with label $l_{1}$ in $\Delta$, $r_{2},r_{4}$ belongs to the arc with label $l_{2}$ in $\Delta$, $r_{1}$, $r_{5}$ belongs to $m$ and $r_{6}\in j$, 
\item[iii)]$[p,r_{6}]_{j}\cap[q_{0},q_{1}]_{j}=\emptyset$, 
\item[iv)] $[r_{1},r_{5}]_{\gamma} \cup [r_{5},r_{1}]_{m}$ is a closed curve that cuts  $m'$ in exactly two points and the folded side of $\Delta''$ in only one point,
\item[v)] $[r_{1},r_{5}]_{\gamma} \cup [r_{5},r_{1}]_{m}$ contains in the interior only the unique puncture which is enclosed by $m$,
\item[vi)] $\gamma \cup [q_{0},q_{1}]_{j}$ divides $\Sigma$ in two regions. One of those two regions is homeomorphic to a disk that contains one puncture namely $p$, and possibly to the one puncture which is enclosed by $m$,
\item[vii)]For $k=0,..,l$ the segments $[r_{k},r_{k+1}]_{\gamma}$ are contractibles to $p$ with the homotopy that avoids $M$ and each of whose intermediate maps are segments with endpoints in the arcs of $\tau^{\circ}$ to which $r_{k}$ and $r_{k+1}$ belong,
\item[viii)] The segments $[r_{1},r_{2}]_{\gamma}$ and $[r_{l-1},r_{l}]_{\gamma}$ are contained to $\Delta$, 
\item[ix)] The segment $[r_{5},r_{8}]_{\gamma}$ is contractible to $p$,
\item[x)] The interior of the segment $[q_{1},r_{l-1}]_{\gamma}$ cuts only $m$ in one point, 
\item[xi)] The segment $[r_{4},r_{6}]_{\gamma}$ is contractible to $p$. 
\end{myEnumerate}

\begin{figure}[H]
\centering 
\vspace{0em}\hspace*{-6em}{\begin{tikzpicture}[scale=0.8]

\filldraw [black] (0,0) circle (2pt)
(-1.5,-2) circle (2pt)
(1.5,-2) circle (2pt)
(0,0.5) circle (1pt)
(0.3,0.5) circle (1pt)
(-0.3,0.5) circle (1pt);
\draw (0,-2.5)node{$\Delta$};
\draw (-1,-0.8)node{$\Delta''$};
\draw (1.4,-1.4)node{$j$};
\draw (0.8,-1.5)node{$\Delta'$};
\draw[line width=1pt] (0,0) -- (-1, 1)
node[pos=0.3,right] {$p$};
\draw[line width=1pt] (0,0) -- (1, 1);
\draw[line width=1pt] (0,0) -- (-1.5,-2);
\draw[line width=1pt] (0,0) -- (1.5,-2);
\draw[line width=1pt] (0, 0) .. controls(-4,-0.1) and (-4,-3.9) .. (0, -4);
\draw[line width=1pt] (0, 0) .. controls(4,-0.1) and (4,-3.9) .. (0, -4);
\draw[line width=1pt] (0, 0) .. controls(-1.5,-0.5) and (-2,-2) .. (-1.75, -2.25);
\draw[line width=1pt] (0, 0) .. controls(0.05,-0.5) and (-0.75,-2.5) .. (-1.75, -2.25)
node[pos=0.8,right] {$m'$};
\draw[line width=1pt] (0, 0) .. controls(0,-1) and (0.7,-2.5) .. (1.75, -2.25)
node[pos=0.8,left] {$m$};
\draw[line width=1pt] (0, 0) .. controls(0,0) and (2.25,-1.5) .. (1.75, -2.25);


\filldraw [black] (8,0) circle (2pt)
(6.5,-2) circle (2pt)
(9.5,-2) circle (2pt)
(8,0.5) circle (1pt)
(8.3,0.5) circle (1pt)
(7.7,0.5) circle (1pt);
\draw (8,-3.5)node{$\Delta$};
\draw (7.15,-1.6)node{$\Delta''$};
\draw (9.4,-1.4)node{$j$};
\draw (8.6,-1.2)node{$\Delta'$};
\draw[line width=1pt] (8,0) -- (7, 1)
node[pos=0.3,right] {$p$};
\draw[line width=1pt] (8,0) -- (9, 1);
\draw[line width=1pt] (8,0) -- (6.5,-2);
\draw[line width=1pt] (8,0) -- (9.5,-2);
\draw[line width=1pt] (8, 0) .. controls(4,-0.1) and (4,-3.9) .. (8, -4);
\draw[line width=1pt] (8, 0) .. controls(12,-0.1) and (12,-3.9) .. (8, -4);
\draw[line width=1pt] (8, 0) .. controls(6.5,-0.5) and (6,-2) .. (6.25, -2.25);
\draw[line width=1pt] (8, 0) .. controls(8.05,-0.5) and (7.25,-2.5) .. (6.25, -2.25)
node[pos=0.8,right] {$m'$};
\draw[line width=1pt] (8, 0) .. controls(8,-1) and (8.7,-2.5) .. (9.75, -2.25);
\draw[line width=1pt] (8, 0) .. controls(8,0) and (10.25,-1.5) .. (9.75, -2.25)
node[pos=0.8,right] {$m$};

\filldraw [red] (1.3,-1.8) circle (2pt)
(1.75,-1.75) circle (2pt)
(-1.6,-1.35) circle (2pt)
(-1,-1.25) circle (2pt)
(-0.3,-1.05) circle (2pt)
(0.15,-0.85) circle (2pt)
(0.45,-0.6) circle (2pt)
(0.6,-0.4) circle (2pt)
(0.7,-0.05) circle (2pt)
(0.6,0.6) circle (2pt)
(-0.8,0.8) circle (2pt)
(1.1,-0.85) circle (2pt)
(2.75,-1.2) circle (2pt)
(0.75,-1) circle (2pt);
\draw[color=red][line width=1.5pt] (0.8,-1) .. controls(1.2,-0.5) and (3.3,-1).. (3.3,-2)
node[pos=0,left] {$q_{1}$}
node[pos=0.2,above] {$r_{l}$}; 
\draw[color=red][line width=1.5pt] (-3.5,-2) ..controls(-3.6,-5) and (3.05,-5) .. (3.3,-2);
\draw[color=red][line width=1.5pt] (-3.5,-2) ..controls(-3.1,1) and (0.2,1.2) .. (0.6,0.6);
\draw[color=red][line width=1.5pt] (-1.75,-1.35) ..controls(-1.15,-1.35) and (1.3,-1) .. (0.6,0.6)
node[pos=0.3,below] {$r_{3}$}; 
\draw[color=red][line width=1.5pt] (-1.65,-1.35) ..controls (-3.5,-2) and (-0.95,-3.5) .. (1.05,-2.75)
node[pos=0,above] {$r_{2}$};
\draw[color=red][line width=1.5pt] (1.05,-2.75).. controls (1.45,-2.7) and (3,-1.5) .. (1.3,-1.8)
node[pos=1,left] {$q_{0}$}
node[pos=0.87,above] {$r_{1}$};

\filldraw [red] (9.1,-1.5) circle (2pt)
(8.65,-1.85) circle (2pt)
(6.35,-1.35) circle (2pt)
(7,-1.3) circle (2pt)
(7.7,-1.1) circle (2pt)
(8.1,-0.9) circle (2pt)
(8.45,-0.6) circle (2pt)
(8.6,-0.45) circle (2pt)
(8.7,-0.05) circle (2pt)
(8.6,0.6) circle (2pt)
(7.2,0.8) circle (2pt)
(5,-2) circle (2pt)
(9.25,-2.2) circle (2pt)
(9.4,-1.9) circle (2pt);

\draw[color=red][line width=1.5pt] (5,-2) .. controls (5.2,-4.25) and (9,-3.5) .. (9.4,-1.9)
node[pos=1,right] {$q_{1}$}
node[pos=0.9,below] {$r_{l}$};
\draw[color=red][line width=1.5pt] (5,-2) ..controls(4.9,1) and (8.2,1.2) .. (8.6,0.6);
\draw[color=red][line width=1.5pt] (6.35,-1.35) ..controls(6.85,-1.35) and (9.3,-1) .. (8.6,0.6);
\draw[color=red][line width=1.5pt] (6.35,-1.35) ..controls (4.5,-2) and (7.05,-3.5) .. (9.05,-1.5)
node[pos=1,above] {$q_{0}$}
node[pos=0.87,above] {$r_{1}$}
node[pos=0.05,below] {$r_{2}$};

\end{tikzpicture}} 

\end{figure} 
\begin{figure}[H]
\centering 
\vspace{-4em}\hspace*{-2em}{\begin{tikzpicture}[scale=0.8]

\filldraw [black] (0,0) circle (2pt)
(-1.5,-2) circle (2pt)
(1.5,-2) circle (2pt)
(0,0.5) circle (1pt)
(0.3,0.5) circle (1pt)
(-0.3,0.5) circle (1pt);
\draw (0,-3)node{$\Delta$};
\draw (-1,-1)node{$\Delta''$};
\draw (-1.5,-1.6)node{$j$};
\draw (0.8,-1.5)node{$\Delta'$};
\draw[line width=1pt] (0,0) -- (-1, 1)
node[pos=0.3,right] {$p$};
\draw[line width=1pt] (0,0) -- (1, 1);
\draw[line width=1pt] (0,0) -- (-1.5,-2);
\draw[line width=1pt] (0,0) -- (1.5,-2);
\draw[line width=1pt] (0, 0) .. controls(-4,-0.1) and (-4,-3.9) .. (0, -4);
\draw[line width=1pt] (0, 0) .. controls(4,-0.1) and (4,-3.9) .. (0, -4);
\draw[line width=1pt] (0, 0) .. controls(-1.5,-0.5) and (-2,-2) .. (-1.75, -2.25);
\draw[line width=1pt] (0, 0) .. controls(0.05,-0.5) and (-0.75,-2.5) .. (-1.75, -2.25)
node[pos=0.8,right] {$m'$};
\draw[line width=1pt] (0, 0) .. controls(0,-1) and (0.7,-2.5) .. (1.75, -2.25)
node[pos=0.8,left] {$m$};
\draw[line width=1pt] (0, 0) .. controls(0,0) and (2.25,-1.5) .. (1.75, -2.25);

\filldraw [black] (8,0) circle (2pt)
(6.5,-2) circle (2pt)
(9.5,-2) circle (2pt)
(8,0.5) circle (1pt)
(8.3,0.5) circle (1pt)
(7.7,0.5) circle (1pt);
\draw (8,-3)node{$\Delta$};
\draw (7,-1)node{$\Delta''$};
\draw (6.5,-1.6)node{$j$};
\draw (8.8,-1.5)node{$\Delta'$};
\draw[line width=1pt] (8,0) -- (7, 1)
node[pos=0.3,right] {$p$};
\draw[line width=1pt] (8,0) -- (9, 1);
\draw[line width=1pt] (8,0) -- (6.5,-2);
\draw[line width=1pt] (8,0) -- (9.5,-2);
\draw[line width=1pt] (8, 0) .. controls(4,-0.1) and (4,-3.9) .. (8, -4);
\draw[line width=1pt] (8, 0) .. controls(12,-0.1) and (12,-3.9) .. (8, -4);
\draw[line width=1pt] (8, 0) .. controls(6.5,-0.5) and (6,-2) .. (6.25, -2.25);
\draw[line width=1pt] (8, 0) .. controls(8.05,-0.5) and (7.25,-2.5) .. (6.25, -2.25)
node[pos=0.8,right] {$m'$};
\draw[line width=1pt] (8, 0) .. controls(8,-1) and (8.7,-2.5) .. (9.75, -2.25)
node[pos=0.8,left] {$m$};
\draw[line width=1pt] (8, 0) .. controls(8,0) and (10.25,-1.5) .. (9.75, -2.25);
\filldraw [red] (1.3,-1.8) circle (2pt)
(0.5,-1.75) circle (2pt)
(-1.6,-1.35) circle (2pt)
(-1,-1.25) circle (2pt)
(-0.3,-1.05) circle (2pt)
(0.15,-0.85) circle (2pt)
(0.45,-0.6) circle (2pt)
(0.6,-0.4) circle (2pt)
(0.7,-0.05) circle (2pt)
(0.6,0.6) circle (2pt)
(-0.8,0.8) circle (2pt)
(1.5,-1.25) circle (2pt)
(2.75,-1.2) circle (2pt)
(0.75,-1) circle (2pt);
\draw[color=red][line width=1.5pt] (1.3,-1.8) .. controls(1.2,-0.5) and (3.3,-1).. (3.3,-2)
node[pos=0,left] {$q_{1}$}
node[pos=0.3,above] {$r_{l}$}; 
\draw[color=red][line width=1.5pt] (-3.5,-2) ..controls(-3.6,-5) and (3.05,-5) .. (3.3,-2);
\draw[color=red][line width=1.5pt] (-3.5,-2) ..controls(-3.1,1) and (0.2,1.2) .. (0.6,0.6);

\draw[color=red][line width=1.5pt] (-1.75,-1.35) ..controls(-1.15,-1.35) and (1.3,-1) .. (0.6,0.6)
node[pos=0.3,below] {$r_{3}$}; 
\draw[color=red][line width=1.5pt] (-1.65,-1.35) ..controls (-3.5,-2) and (-1.2,-2.9) .. (-1,-2.75)
node[pos=0,above] {$r_{2}$};
\draw[color=red][line width=1.5pt] (-1,-2.75).. controls (0.45,-2.7) and (0.5,-1.5) .. (0.8,-1)
node[pos=1,right] {$q_{0}$}
node[pos=0.6,left] {$r_{1}$};
\filldraw [red] (9.3,-1.8) circle (2pt)
(8.7,-1.9) circle (2pt)
(6.4,-1.35) circle (2pt)
(7,-1.25) circle (2pt)
(7.7,-1.05) circle (2pt)
(8.15,-0.85) circle (2pt)
(8.45,-0.6) circle (2pt)
(8.6,-0.4) circle (2pt)
(8.7,-0.05) circle (2pt)
(8.6,0.6) circle (2pt)
(7.2,0.8) circle (2pt)
(9.1,-0.85) circle (2pt)
(10.75,-1.2) circle (2pt)
(8.75,-1) circle (2pt);
\draw[color=red][line width=1.5pt] (8.8,-1) .. controls(9.2,-0.5) and (11.3,-1).. (11.3,-2)
node[pos=0,left] {$q_{1}$}
node[pos=0.3,above] {$r_{l}$}; 
\draw[color=red][line width=1.5pt] (4.5,-2) ..controls(4.4,-5) and (11.05,-5) .. (11.3,-2);
\draw[color=red][line width=1.5pt] (4.5,-2) ..controls(4.9,1) and (8.2,1.2) .. (8.6,0.6);

\draw[color=red][line width=1.5pt] (6.25,-1.35) ..controls(6.85,-1.35) and (9.3,-1) .. (8.6,0.6)
node[pos=0.3,below] {$r_{3}$}; 
\draw[color=red][line width=1.5pt] (6.35,-1.35) ..controls (4.5,-2) and (6.8,-2.9) .. (7,-2.75)
node[pos=0,above] {$r_{2}$};
\draw[color=red][line width=1.5pt] (7,-2.75).. controls (8.45,-2.7) and (8.5,-1.5) .. (9.3,-1.8)
node[pos=1,right] {$q_{0}$}
node[pos=0.6,left] {$r_{1}$};

\end{tikzpicture}} 
\caption{Case $a_{22}$) of the Definition \ref{rodear}.}\label{figura a22p} 

\end{figure}

\item[b)] $\delta_{\tau}(p) = 0$.

The condition $\delta_{\tau}(p) = 0$ implies that there is only one arc $j'$ in $\tau^{\circ}$ which is adjacent to the puncture $p$ and $j'$ is the folded side of a self-folded triangle $\Delta'$. We denote by $m$ the not folded side of $\Delta'$. As in the case $a)$ we have the followings two subcases for the arc $j$.
\begin{itemize}
\item[$b_{1})$]The arc $j$ is not a side of any self-folded triangle.

The hypothesis $\delta_{\tau}(p)=0$ implies that $l=3$ and we have to consider the following two possibilities:
\begin{itemize}
\item[$b_{11})$]The arc $j$ is a side of a triangle $\Delta$ wich is of type 2 and $\Delta$ shares the side $m$ with $\Delta'$. See Figure \ref{figura b11}.

\begin{myEnumerate}
\item[i)] $j'$ is adjacent to $p$,
\item[ii)] $\gamma$ cuts the arc $m$ only in $r_{1}$ and $r_{3}$,
\item[iii)] $\gamma$ cut only in $r_{2}$ to $j'$,
\item[iv)] The curve $\gamma$ does not cut any other arc of $\tau^{\circ}$, except those mentioned in the previous points, 
\item[v)] $\gamma \cup [q_{0},q_{1}]_{j}$ divides $\Sigma$ in two regions and one of those two regions is homeomorphic to a disk that contains one puncture namely $p$. 
\end{myEnumerate}
\begin{figure}[H]
\centering 
$$ \begin{tikzpicture}[scale=0.8]

\filldraw [black] (0,0) circle (2pt)
(0,-4) circle (2pt)
(0,-2) circle (2pt)
(0,0.5) circle (1pt)
(0.3,0.5) circle (1pt)
(-0.3,0.5) circle (1pt)
(0,-4.5) circle (1pt)
(0.3,-4.5) circle (1pt)
(-0.3,-4.5) circle (1pt);
\draw (0,-3)node{$\Delta$};
\draw (4,-3.25)node{$\Delta$};
\draw (3.7,-1.25)node{$\Delta'$};
\draw (-0.25,-1.8)node{$\Delta'$};
\draw (4.7,-2)node{$m$};
\draw (5.8,-2)node{$j$};
\draw (-0.7,-1.5)node{$m$};
\draw (-1.7,-2)node{$j$};
\draw (0,-2.3)node{$p$};
\draw[line width=1pt] (0,0) -- (-1, 1);
\draw[line width=1pt] (0,0) -- (1, 1);
\draw[line width=1pt] (0,-4) -- (1, -5);
\draw[line width=1pt] (0,-4) -- (-1, -5);
\draw[line width=1pt] (0,0) -- (0, -2);
\draw[line width=1pt] (0,0) .. controls(-0.7,-1) and (-0.7,-2.2) .. (0, -2.5);
\draw[line width=1pt] (0,0) .. controls(0.7,-1) and (0.7,-2.2) .. (0, -2.5);
\draw[line width=1pt] (0,0) .. controls(-2, -0.5) and (-2,-3.5) .. (0, -4);
\draw[line width=1pt] (0,0) .. controls(2, -0.5) and (2,-3.5) .. (0, -4);

\filldraw [black] 
(4,0) circle (2pt)
(4,-4) circle (2pt)
(4,-2) circle (2pt)
(4,0.5) circle (1pt)
(4.3,0.5) circle (1pt)
(3.7,0.5) circle (1pt)
(4,-4.5) circle (1pt)
(4.3,-4.5) circle (1pt)
(3.7,-4.5) circle (1pt);
\draw (4,-2.3)node{$p$};
\draw[line width=1pt] (4,0) -- (3, 1);
\draw[line width=1pt] (4,0) -- (5, 1);
\draw[line width=1pt] (4,-4) -- (5, -5);
\draw[line width=1pt] (4,-4) -- (3, -5);
\draw[line width=1pt] (4,0) -- (4, -2);
\draw[line width=1pt] (4,0) .. controls(2, -0.5) and (2,-3.5) .. (4, -4);
\draw[line width=1pt] (4,0) .. controls(6, -0.5) and (6,-3.5) .. (4, -4);
\draw[line width=1pt] (4,0) .. controls(3.3,-1) and (3.3,-2.2) .. (4, -2.5);
\draw[line width=1pt] (4,0) .. controls(4.7,-1) and (4.7,-2.2) .. (4, -2.5);



\filldraw [red] (5.1,-0.75) circle (2pt)
(4.45,-1) circle (2pt)
(4,-1.25) circle (2pt)
(3.5,-1.95) circle (2pt)
(5.1,-3.25) circle (2pt);
\draw[color=red][line width=1.5pt] (5.1,-0.75) ..controls(3,-1.25) and (3,-2.75) .. (5.1,-3.25)
node[pos=0,right] { $q_{0}=r_{0}$}
node[pos=0.05,below] { $r_{1}$}
node[pos=0.28,right] { $r_{2}$}
node[pos=0.45,left] { $r_{3}$}
node[pos=1,right] { $q_{1}=r_{4}$};

\filldraw [red] (-1.1,-0.75) circle (2pt)
(-0.45,-0.95) circle (2pt)
(0,-1.17) circle (2pt)
(0.5,-1.95) circle (2pt)
(-1.1,-3.25) circle (2pt);

\draw[color=red][line width=1.5pt] (-1.1,-0.75) .. controls(1,-1.25) and (1,-2.75) .. (-1.1,-3.25)
node[pos=0,left] { $q_{1}=r_{4}$}
node[pos=0.1,above] { $r_{3}$}
node[pos=0.20,right] { $r_{2}$}
node[pos=0.45,right] { $r_{1}$}
node[pos=1,left] { $q_{0}=r_{0}$};

\end{tikzpicture}$$
\caption{Case $b_{11}$) of the Definition \ref{rodear}.}\label{figura b11}
\end{figure}
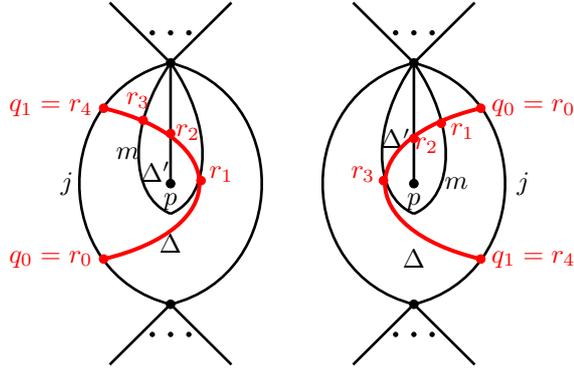
\item[$b_{12})$]The arc $j$ has been labeled $l_{1}$ in $\Delta$ and $\Delta$ is triangle of type $3$. We are going to consider to $\Delta'$ and $\Delta''$ as the unique self-folded triangles which share one side with $\Delta$, where $m$ and $m'$ are those sides respectively. We denote by $j'$ and $j''$ the folded sides of $\Delta'$ and $\Delta''$ respectively. See Figure \ref{figura b12}.

\begin{myEnumerate}
\item[i)] $j'$ is adjacent to $p$,
\item[ii)] $\gamma$ cuts $m$ only in two points $r_{1}$ and $r_{3}$,
\item[iii)] $\gamma$ cuts $j'$ only in $r_{2}$,
\item[iv)] The curve $\gamma$ does not cut any other arc of $\tau^{\circ}$, except those mentioned in the previous points, 
\item[v)]$\gamma \cup [q_{0},q_{1}]_{j}$ divides $\Sigma$ in two regions and one of those two regions is homeomorphic to a disk that contains one puncture namely $p$. 
\end{myEnumerate}
\begin{figure}[H]
\begin{center}
\subfigure[]{
\hspace*{-1em}{\begin{tikzpicture}[scale=0.82]

\filldraw [black] (0,0) circle (2pt)
(-1.5,-2) circle (2pt)
(1.5,-2) circle (2pt)
(0,0.5) circle (1pt)
(0.3,0.5) circle (1pt)
(-0.3,0.5) circle (1pt);
\draw (0,-3)node{$\Delta$};
\draw (-1,-1)node{$\Delta'$};
\draw (-3.2,-2)node{$j$};
\draw (0.8,-1.5)node{$\Delta''$};

\draw (-1.6,-2.15)node{$p$};

\draw[line width=1pt] (0,0) -- (-1, 1);
\draw[line width=1pt] (0,0) -- (1, 1);
\draw[line width=1pt] (0,0) -- (-1.5,-2);
\draw[line width=1pt] (0,0) -- (1.5,-2);
\draw[line width=1pt] (0, 0) .. controls(-4,-0.1) and (-4,-3.9) .. (0, -4);
\draw[line width=1pt] (0, 0) .. controls(4,-0.1) and (4,-3.9) .. (0, -4);
\draw[line width=1pt] (0, 0) .. controls(-1.5,-0.5) and (-2,-2) .. (-1.75, -2.25);
\draw[line width=1pt] (0, 0) .. controls(0.05,-0.5) and (-0.75,-2.5) .. (-1.75, -2.25)
node[pos=0.6,right] {$m$};
\draw[line width=1pt] (0, 0) .. controls(0,-1) and (0.7,-2.5) .. (1.75, -2.25)
node[pos=0.8,left] {$m'$};
\draw[line width=1pt] (0, 0) .. controls(0,0) and (2.25,-1.5) .. (1.75, -2.25);


\filldraw [red] 
(-2.6,-1) circle (2pt)
(-1.6,-1.35) circle (2pt)
(-1.25,-1.555) circle (2pt)
(-1,-2.05) circle (2pt)
(-2.6,-3) circle (2pt);

\draw[color=red][line width=1.5pt] (-2.6,-1) ..controls(-0.5,-1.5) and (-0.5,-2.5) .. (-2.6, -3) 
node[pos=1,left] {$r_{0}=q_{0}$}
node[pos=0.55,right] {$r_{1}$}
node[pos=0.32,right] {$r_{2}$}
node[pos=0.08,right] {$r_{3}$}
node[pos=0,left] {$r_{4}=q_{1}$};

\end{tikzpicture} }
}
\subfigure[]{
\hspace*{0.6em}{\begin{tikzpicture}[scale=0.82] 

\filldraw [black] (0,0) circle (2pt)
(-1.5,-2) circle (2pt)
(1.5,-2) circle (2pt)
(0,0.5) circle (1pt)
(0.3,0.5) circle (1pt)
(-0.3,0.5) circle (1pt);
\draw (0,-3)node{$\Delta$};
\draw (-1,-1)node{$\Delta''$};
\draw (-3.2,-2)node{$j$};
\draw (0.8,-1.5)node{$\Delta'$};

\draw (1.7,-2.15)node{$p$};

\draw[line width=1pt] (0,0) -- (-1, 1);
\draw[line width=1pt] (0,0) -- (1, 1);
\draw[line width=1pt] (0,0) -- (-1.5,-2);
\draw[line width=1pt] (0,0) -- (1.5,-2);
\draw[line width=1pt] (0, 0) .. controls(-4,-0.1) and (-4,-3.9) .. (0, -4);
\draw[line width=1pt] (0, 0) .. controls(4,-0.1) and (4,-3.9) .. (0, -4);
\draw[line width=1pt] (0, 0) .. controls(-1.5,-0.5) and (-2,-2) .. (-1.75, -2.25);
\draw[line width=1pt] (0, 0) .. controls(0.05,-0.5) and (-0.75,-2.5) .. (-1.75, -2.25)
node[pos=0.6,right] {$m'$};
\draw[line width=1pt] (0, 0) .. controls(0,-1) and (0.7,-2.5) .. (1.75, -2.25)
node[pos=0.6,left] {$m$};
\draw[line width=1pt] (0, 0) .. controls(0,0) and (2.25,-1.5) .. (1.75, -2.25);


\filldraw [red] 
(2.6,-1) circle (2pt)
(1.6,-1.35) circle (2pt)
(1.25,-1.555) circle (2pt)
(1,-2.05) circle (2pt)
(2.6,-3) circle (2pt);

\draw[color=red][line width=1.5pt] (2.6,-1) ..controls(0.5,-1.5) and (0.5,-2.5) .. (2.6, -3) 
node[pos=1,right] {$r_{1}=q_{1}$}
node[pos=0.55,left] {$r_{3}$}
node[pos=0.33,above] {$r_{2}$}
node[pos=0.2,above] {$r_{1}$}
node[pos=0,right] {$r_{0}=q_{0}$};
\end{tikzpicture}}

} 
\end{center}
\caption{Case $b_{12}$) of the Definition \ref{rodear}.}\label{figura b12}

\end{figure}
\end{itemize}
\item[$b_{2})$]The arc $j$ is the folded side of some self-folded triangle $\Delta'$. We denote by $m$ the not folded side of $\Delta'$.

According to with conditions over $j$ and $m$, it is easy to see that $m$ must have been a side of a triangle of type $3$. Let $m'$ be the unique side of $\Delta$ which is different to $m$ and in addition $m'$ is not a folded side of a self-folded triangle $\Delta''$, $p'$ is the adjacent puncture to $j$ with signature $0$ and $j'$ is the folded side of $\Delta''$. As the puncture $p$ has signature $0$ then $l=5$. See Figure \ref{figura b2}.

\begin{myEnumerate}
\item[i)] $j'$ is adjacent to $p$,
\item[ii)]If $m$ has a label $l_{2}$ in $\Delta$ then $q_{0}\in [p',q_{1}]_{j}$. On the other hand, if $m$ has a label $l_{3}$ in $\Delta$ then $q_{1}\in[p',q_{0}]_{j}$,

\item[iii)] $\gamma$ cuts $m$ only in $r_{1}$ and $r_{5}$,
\item[iv)] $\gamma$ cuts $m$ only in $r_{2}$ and $r_{4}$,
\item[v)] $\gamma$ cuts $j'$ only in $r_{3}$,
\item[vi)] $\gamma \cup [q_{0},q_{1}]_{j}$ divides $\Sigma$ in two regions and one of those two regions is homeomorphic to a disk that contains one puncture namely $p$ and possibly to the puncture $p'$, 
\item[vii)] If $m$ has a label $l_{2}$ in $\Delta$ then $r_{0}\notin[q,q_{1}]_{j'}$ where $q$ is the opposite puncture to $p$ in $j'$. On the other hand, if $m$ has a label $l_{3}$ in $\Delta$ then $r_{l+1}\notin[q,r_{0}]_{j'}$.

\end{myEnumerate}

\begin{figure}[H]
\centering 

\hspace*{-6em}{\begin{tikzpicture}[scale=0.8]

\filldraw [black] (0,0) circle (2pt)
(-1.5,-2) circle (2pt)
(1.5,-2) circle (2pt)
(0,0.5) circle (1pt)
(0.3,0.5) circle (1pt)
(-0.3,0.5) circle (1pt);
\draw (0,-3)node{$\Delta$};
\draw (-0.8,-1.5)node{$\Delta'$};
\draw (-1.2,-1.3)node{$j$};
\draw (0.9,-1.75)node{$\Delta''$};
\draw (1.3,-2.05)node{$p$};
\draw (-1.5,-2.25)node{$p'$};

\draw[line width=1pt] (0,0) -- (-1, 1);
\draw[line width=1pt] (0,0) -- (1, 1);
\draw[line width=1pt] (0,0) -- (-1.5,-2);
\draw[line width=1pt] (0,0) -- (1.5,-2);
\draw[line width=1pt] (0, 0) .. controls(-4,-0.1) and (-4,-3.9) .. (0, -4);
\draw[line width=1pt] (0, 0) .. controls(4,-0.1) and (4,-3.9) .. (0, -4);
\draw[line width=1pt] (0, 0) .. controls(-1.5,-0.5) and (-2,-2) .. (-1.75, -2.25);
\draw[line width=1pt] (0, 0) .. controls(0.05,-0.5) and (-0.75,-2.5) .. (-1.75, -2.25)
node[pos=0.6,right] {$m$};
\draw[line width=1pt] (0, 0) .. controls(0,-1) and (0.7,-2.5) .. (1.75, -2.25)
node[pos=0.6,left] {$m'$};
\draw[line width=1pt] (0, 0) .. controls(0,0) and (2.25,-1.5) .. (1.75, -2.25);


\filldraw [red] (-1.35,-1.8) circle (2pt)
(-0.9,-1.96) circle (2pt)
(1.5,-2.25) circle (2pt)
(1,-1.3) circle (2pt)
(0.2,-1.1) circle (2pt)
(-0.3,-1) circle (2pt)
(-0.8,-1) circle (2pt);

\draw[color=red][line width=1.5pt] (-1.35,-1.8) ..controls(2.65,-3.5) and (2.65,-1) .. (-0.8,-1)
node[pos=0,left] {$r_{0}=q_{0}$}
node[pos=0.05,below] {$r_{1}$}
node[pos=0.4,below] {$r_{2}$}
node[pos=0.75,above] {$r_{3}$}
node[pos=0.87,above] {$r_{4}$}
node[pos=0.96,above] {$r_{5}$}
node[pos=,left] {$r_{6}=q_{1}$};

\filldraw [black] (8,0) circle (2pt)
(6.5,-2) circle (2pt)
(9.5,-2) circle (2pt)
(8,0.5) circle (1pt)
(8.3,0.5) circle (1pt)
(7.7,0.5) circle (1pt);
\draw (8,-2.5)node{$\Delta$};
\draw (6.7,-1.4)node{$\Delta'$};
\draw (7.2,-1.4)node{$j$};
\draw (9.3,-2.05)node{$p$};
\draw (8.8,-1.5)node{$\Delta''$};

\draw (6.5,-2.25)node{$p'$};
\draw[line width=1pt] (8,0) -- (7, 1);
\draw[line width=1pt] (8,0) -- (9, 1);
\draw[line width=1pt] (8,0) -- (6.5,-2);
\draw[line width=1pt] (8,0) -- (9.5,-2);
\draw[line width=1pt] (8, 0) .. controls(4,-0.1) and (4,-3.9) .. (8, -4);
\draw[line width=1pt] (8, 0) .. controls(12,-0.1) and (12,-3.9) .. (8, -4);
\draw[line width=1pt] (8, 0) .. controls(6.5,-0.5) and (6,-2) .. (6.25, -2.25);
\draw[line width=1pt] (8, 0) .. controls(8.05,-0.5) and (7.25,-2.5) .. (6.25, -2.25)
node[pos=0.9,right] {$m$};
\draw[line width=1pt] (8, 0) .. controls(8,-1) and (8.7,-2.5) .. (9.75, -2.25)
node[pos=0.8,left] {$m'$};
\draw[line width=1pt] (8, 0) .. controls(8,0) and (10.25,-1.5) .. (9.75, -2.25);


\filldraw [red] 
(6.65,-1.8) circle (2pt)
(6.2,-1.95) circle (2pt)
(9.8,-2) circle (2pt)
(8.2,-1.1) circle (2pt)
(8.95,-1.3) circle (2pt)
(7.75,-1) circle (2pt)
(7.2,-1) circle (2pt);

\draw[color=red][line width=1.5pt] (6.65,-1.8) ..controls(4.5,-2.5) and (7.1,-3.3) .. (8,-3.25)
node[pos=0,right] {$q_{0}=r_{0}$}
node[pos=0.1,left] {$r_{1}$};

\draw[color=red][line width=1.5pt] (8,-3.25) ..controls(10.65,-3.5) and (10.65,-1) .. (7.2,-1)
node[pos=0.5,right] {$r_{2}$}
node[pos=0.75,above] {$r_{3}$}
node[pos=0.87,above] {$r_{4}$}
node[pos=0.96,above] {$r_{5}$}
node[pos=,left] {$r_{6}=q_{1}$};

\end{tikzpicture}}
\end{figure}

\begin{figure}[H]
\centering 

\hspace*{-6em}{\begin{tikzpicture}[scale=0.8]

\filldraw [black] (0,0) circle (2pt)
(-1.5,-2) circle (2pt)
(1.5,-2) circle (2pt)
(0,0.5) circle (1pt)
(0.3,0.5) circle (1pt)
(-0.3,0.5) circle (1pt);
\draw (0,-3)node{$\Delta$};
\draw (-0.8,-1.5)node{$\Delta''$};
\draw (1.2,-1.3)node{$j$};
\draw (0.9,-1.75)node{$\Delta'$};
\draw (1.3,-2.15)node{$p'$};
\draw (-1.3,-2.05)node{$p$};

\draw[line width=1pt] (0,0) -- (-1, 1);
\draw[line width=1pt] (0,0) -- (1, 1);
\draw[line width=1pt] (0,0) -- (-1.5,-2);
\draw[line width=1pt] (0,0) -- (1.5,-2);
\draw[line width=1pt] (0, 0) .. controls(-4,-0.1) and (-4,-3.9) .. (0, -4);
\draw[line width=1pt] (0, 0) .. controls(4,-0.1) and (4,-3.9) .. (0, -4);
\draw[line width=1pt] (0, 0) .. controls(-1.5,-0.5) and (-2,-2) .. (-1.75, -2.25);
\draw[line width=1pt] (0, 0) .. controls(0.05,-0.5) and (-0.75,-2.5) .. (-1.75, -2.25)
node[pos=0.6,right] {$m'$};
\draw[line width=1pt] (0, 0) .. controls(0,-1) and (0.7,-2.5) .. (1.75, -2.25)
node[pos=0.6,left] {$m$};
\draw[line width=1pt] (0, 0) .. controls(0,0) and (2.25,-1.5) .. (1.75, -2.25);


\filldraw [red] (1.35,-1.8) circle (2pt)
(0.85,-1.96) circle (2pt)
(-1.5,-2.25) circle (2pt)
(-1,-1.3) circle (2pt)
(-0.2,-1.1) circle (2pt)
(0.3,-1) circle (2pt)
(0.8,-1) circle (2pt);

\draw[color=red][line width=1.5pt] (1.35,-1.8) ..controls(-2.65,-3.5) and (-2.65,-1) .. (0.8,-1)
node[pos=0,right] {$r_{6}=q_{1}$}
node[pos=0.05,below] {$r_{5}$}
node[pos=0.4,below] {$r_{4}$}
node[pos=0.75,above] {$r_{3}$}
node[pos=0.87,above] {$r_{2}$}
node[pos=0.95,above] {$r_{1}$}
node[pos=,right] {$r_{0}=q_{0}$};


\filldraw [black] (8,0) circle (2pt)
(6.5,-2) circle (2pt)
(9.5,-2) circle (2pt)
(8,0.5) circle (1pt)
(8.3,0.5) circle (1pt)
(7.7,0.5) circle (1pt);
\draw (8,-2.5)node{$\Delta$};
\draw (7.2,-1.5)node{$\Delta''$};
\draw (9.3,-1.4)node{$j$};
\draw (8.8,-1.5)node{$\Delta'$};

\draw (6.7,-2.05)node{$p$};
\draw (9.3,-2.15)node{$p'$};
\draw[line width=1pt] (8,0) -- (7, 1);
\draw[line width=1pt] (8,0) -- (9, 1);
\draw[line width=1pt] (8,0) -- (6.5,-2);
\draw[line width=1pt] (8,0) -- (9.5,-2);
\draw[line width=1pt] (8, 0) .. controls(4,-0.1) and (4,-3.9) .. (8, -4);
\draw[line width=1pt] (8, 0) .. controls(12,-0.1) and (12,-3.9) .. (8, -4);
\draw[line width=1pt] (8, 0) .. controls(6.5,-0.5) and (6,-2) .. (6.25, -2.25);
\draw[line width=1pt] (8, 0) .. controls(8.05,-0.5) and (7.25,-2.5) .. (6.25, -2.25)
node[pos=0.9,right] {$m'$};
\draw[line width=1pt] (8, 0) .. controls(8,-1) and (8.7,-2.5) .. (9.75, -2.25)
node[pos=0.8,left] {$m$};
\draw[line width=1pt] (8, 0) .. controls(8,0) and (10.25,-1.5) .. (9.75, -2.25);


\filldraw [red] 
(9.35,-1.8) circle (2pt)
(9.8,-1.95) circle (2pt)
(6.2,-2) circle (2pt)
(7.7,-1.1) circle (2pt)
(7.05,-1.3) circle (2pt)
(8.25,-1) circle (2pt)
(8.8,-1) circle (2pt);

\draw[color=red][line width=1.5pt] (9.35,-1.8) ..controls(11.5,-2.5) and (8.9,-3.3) .. (8,-3.25)
node[pos=0,left] {$q_{1}=r_{6}$}
node[pos=0.1,right] {$r_{5}$};

\draw[color=red][line width=1.5pt] (8,-3.25) ..controls(5.35,-3.5) and (5.35,-1) .. (8.8,-1)
node[pos=0.5,left] {$r_{4}$}
node[pos=0.75,above] {$r_{3}$}
node[pos=0.87,above] {$r_{2}$}
node[pos=0.95,above] {$r_{1}$}
node[pos=,right] {$r_{0}=q_{0}$};

\end{tikzpicture}}
\caption{Case $b_{2}$) of the Definition \ref{rodear}.}\label{figura b2}
\end{figure}
\end{itemize}
\end{itemize}
\end{definition}

\begin{remark}
Analyzing the Definition \ref{rodear}, it could be possible to define the concept \textit{surround a puncture} choosing the non folded side of a self-folded triangle, but to make the writing simpler we choose the folded side.
\end{remark}

\begin{definition}
Let $\tau$ be a tagged triangulation of $(\Sigma,M)$, we say that a loop $k$ \textit{ encloses closely the arc } $i\in \tau^{\circ}$ if the following two conditions are met.

\begin{myEnumerate}
\item[i)] $k$ minimizes the intersection points with the arcs of $\tau^{\circ}$,
\item[ii)] $i$ and $k$ are the two only sides of a self-folded triangle $\Delta'$, where $i$ is the folded side. It is important to observe that $\Delta'$ is not necessarily a triangle in $\tau^{\circ}$.
\end{myEnumerate}
\end{definition}

\begin{definition}\label{curva i'}(Derivative curve $i'=i'(\tau^{\circ},i)$)
Let $\tau$ be a tagged triangulation of $(\Sigma,M)$ with non negative signature and $i$ be a tagged arc which does not belong to $\tau$. Let $\{p,q\}\subseteq M$ be the set of the ends of an arc $j\in \tau^{\circ}$. We are going to define a tagged curve $i'=i'(\tau^{\circ},i)$ ( in two steps) in $\Sigma$ from an arc $i$. First we define a curve without labels $i_{0}$, next we assign a label at each end of $i_{0}$. 

The construction depends on the followings three parameters.

\begin{myEnumerate}
\item[a)] $p=q$ \'o $p\neq q$. 
\item[b)] The signature in the ends of $i$.
\item[c)] The label at the ends of $i$.
\end{myEnumerate}

Let's consider the three parameters above to construct $i'$.

\begin{myEnumerate}

\item[$\bold{a)}$] $p\neq q$.
\begin{myEnumerate}
\item[$a_{1}$)] $\delta_{\tau}(p) \neq 0\neq \delta_{\tau}(q)$.
\begin{myEnumerate}
\item[$a_{11}$)] The arc $i$ has a label $"plain"$ at both ends as shown in the Figure \ref{i' a11)}.
Let's define $i_{0}:=i$ and we assign the label $"plain"$ at both ends of $i_{0}$. Under this situation the derivative curve $i'$ is equal to the tagged arc $i$. 

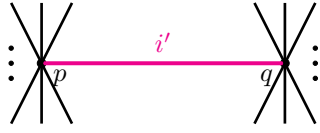
\begin{figure}[H]
\centering 
$$ \begin{tikzpicture}[scale=0.8]
\filldraw [black] (-2,0) circle (2pt)
(2,0) circle (2pt)
(-2.5,0) circle (1pt)
(-2.5,0.25) circle (1pt)
(-2.5,-0.25) circle (1pt)
(2.5,0) circle (1pt)
(2.5,0.25) circle (1pt)
(2.5,-0.25) circle (1pt); 
\draw (-1.7,-0.25)node{$p$};
\draw (1.7,-0.25)node{$q$};

\draw[color=magenta][line width=1.5pt] (-2,0) -- (2, 0)
node[pos=0.5,above] {$i'$};
\draw[line width=1pt] (-2,0) -- (-2, 1);
\draw[line width=1pt] (-2,0) -- (-1.5, 1);
\draw[line width=1pt] (-2,0) -- (-1.5, -1);
\draw[line width=1pt] (-2,0) -- (-2, -1);
\draw[line width=1pt] (-2,0) -- (-2.5, -1);
\draw[line width=1pt] (-2,0) -- (-2.5, 1);

\draw[line width=1pt] (2,0) -- (2, 1);
\draw[line width=1pt] (2,0) -- (2.5, 1);
\draw[line width=1pt] (2,0) -- (2.5, -1);
\draw[line width=1pt] (2,0) -- (2, -1);
\draw[line width=1pt] (2,0) -- (1.5, -1);
\draw[line width=1pt] (2,0) -- (1.5, 1);
\end{tikzpicture}$$
\caption{Case $a_{11}$) of the Definition \ref{curva i'}.}\label{i' a11)}
\end{figure} 
\item[$a_{12}$)]The arc $i'$ has a label $\Bowtie$ only at one end. Without loss of generality we can suppose that the end with label $\Bowtie$ let be $p$.

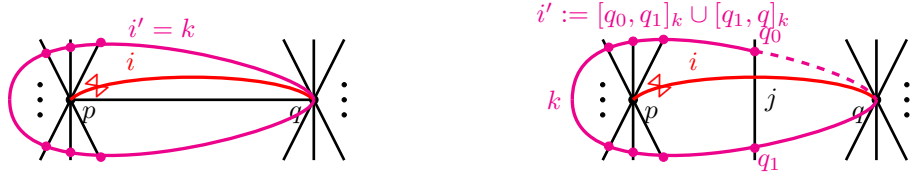
\begin{figure}[htp]
\begin{center}
\subfigure[The version without label of $i$ belong to $\tau^{\circ}$.]{
\hspace*{4em}{\begin{tikzpicture}[scale=0.8]

\filldraw [black] (-2,0) circle (2pt)
(2,0) circle (2pt)
(-2.5,0) circle (1pt)
(-2.5,0.25) circle (1pt)
(-2.5,-0.25) circle (1pt)
(2.5,0) circle (1pt)
(2.5,0.25) circle (1pt)
(2.5,-0.25) circle (1pt);
\filldraw [white] (4.5,0) circle (2pt)
(-4.5,0) circle (2pt);
\draw (-1.7,-0.25)node{$p$};
\draw (1.7,-0.25)node{$q$};
\draw[line width=1pt] (-2,0) -- (2, 0);

\draw[line width=1pt] (-2,0) -- (-2, 1);
\draw[line width=1pt] (-2,0) -- (-1.5, 1);
\draw[line width=1pt] (-2,0) -- (-1.5, -1);
\draw[line width=1pt] (-2,0) -- (-2, -1);
\draw[line width=1pt] (-2,0) -- (-2.5, -1);
\draw[line width=1pt] (-2,0) -- (-2.5, 1);

\draw[line width=1pt] (2,0) -- (2, 1);
\draw[line width=1pt] (2,0) -- (2.5, 1);
\draw[line width=1pt] (2,0) -- (2.5, -1);
\draw[line width=1pt] (2,0) -- (2, -1);
\draw[line width=1pt] (2,0) -- (1.5, -1);
\draw[line width=1pt] (2,0) -- (1.5, 1);
\draw [color=red] (-1.6,0.2)node {$\rotatebox{-50}{\textbf{\Bowtie}}$}; 

\draw[color=red][line width=1.25pt] (-2,0) ..controls(-1.5,0.5) and(1.5,0.5) .. (2, 0)
node[pos=0.3,above] {$i$}; 

\filldraw [magenta] (-1.5,0.95) circle (2pt)
(-2,0.87) circle (2pt)
(-2.4,0.77) circle (2pt)
(-2.4,-0.77) circle (2pt)
(-2,-0.87) circle (2pt)
(-1.5,-0.95) circle (2pt);
\draw[color=magenta][line width=1.25pt] (2,0) ..controls(2,-0.5) and(-3,-1.8) .. (-3, 0);
\draw[color=magenta][line width=1.25pt] (2,0) ..controls(2,0.5) and(-3,1.8) .. (-3, 0)
node[pos=0.5,above] {$i'=k$};

\end{tikzpicture} 

}
}
\subfigure[The version without label of $i$ does not belong to $\tau^{\circ}$. ]{
\hspace*{0.8em}{\begin{tikzpicture}[scale=0.8] 
\filldraw [black] (4,0) circle (2pt)
(8,0) circle (2pt)
(3.5,0) circle (1pt)
(3.5,0.25) circle (1pt)
(3.5,-0.25) circle (1pt)
(8.5,0) circle (1pt)
(8.5,0.25) circle (1pt)
(8.5,-0.25) circle (1pt);

\draw[line width=1pt] (4,0) -- (4, 1);
\draw[line width=1pt] (4,0) -- (4.5, 1);
\draw[line width=1pt] (4,0) -- (4.5, -1);
\draw[line width=1pt] (4,0) -- (4, -1);
\draw[line width=1pt] (4,0) -- (3.5, -1);
\draw[line width=1pt] (4,0) -- (3.5, 1);

\draw[line width=1pt] (8,0) -- (8, 1);
\draw[line width=1pt] (8,0) -- (8.5, 1);
\draw[line width=1pt] (8,0) -- (8.5, -1);
\draw[line width=1pt] (8,0) -- (8, -1);
\draw[line width=1pt] (8,0) -- (7.5, -1);
\draw[line width=1pt] (8,0) -- (7.5, 1);

\draw[line width=1pt] (6,1) -- (6,-1)
node[pos=0.5,right] {$j$}; 
\draw (4.3,-0.25)node{$p$};
\draw (7.7,-0.25)node{$q$};

\draw [color=red] (4.4,0.2)node {$\rotatebox{-50}{\textbf{\Bowtie}}$};
\draw[color=red][line width=1.25pt] (4,0) ..controls(4.5,0.5) and(7.5,0.5) .. (8, 0)
node[pos=0.3,above] {$i$};
\filldraw [magenta] (6,0.8) circle (2pt)
(4.5,1) circle (2pt)
(4,0.97) circle (2pt)
(3.6,0.87) circle (2pt)
(3.6,-0.77) circle (2pt)
(4,-0.87) circle (2pt)
(4.5,-0.95) circle (2pt)
(6,-0.8) circle (2pt);
\draw[color=magenta][line width=1.25pt] (8,0) ..controls(8,-0.5) and(3,-1.8) .. (3, 0)
node[pos=0.4,below] {$q_{1}$};
\draw[color=magenta][line width=1.25pt] (6,0.8) ..controls(6,0.8) and(3,1.6) .. (3, 0)
node[pos=0.5,above] {$i':= [q_{0},q_{1}]_{k} \cup [q_{1},q]_{k}$}
node[pos=1,left] {$k$};
\draw[dashed][color=magenta][line width=1.25pt] (6,0.8) ..controls(6.2,0.8) and (7.8,0.4) .. (8, 0)
node[pos=0.2,above] {$q_{0}$};
\end{tikzpicture}}

} 

\end{center}
\caption{Case $a_{12}$) of the Definition \ref{curva i'}.}\label{i' a12)}

\end{figure}

Let $k$ be the only loop based on the puncture $q$ which encloses closely the arc $i$. If the version without label of the arc $i$ belongs to $\tau^{\circ}$, then we define $i_{0}:=k$ and we assign the label "plain" at both ends of $i_{0}$ as is shown in (a) of Figure \ref{i' a12)}. Otherwise we denote by $\gamma:=[q_{0},q_{1}]_{k}$ the segment of $k$ which surrounds the puncture $p$, where the ends $q_{0}$ and $q_{1}$ of $\gamma$ belong to an arc $j$ of $\tau^{\circ}$ (See Definition \ref{rodear}). It is defined $i_{0}:=\gamma' \cup \gamma$ where $\gamma'$ is the segment of the loop $k$ which connects to $q$ with $q_{1}$ such that $q_{0}$ does not belong to $\gamma'$. We assign the label "plain" at both ends as is shown in (b) of the Figure \ref{i' a12)}.

\item[$a_{13}$)] The arc $i$ has a label $"\Bowtie"$ at both ends.

Let $i_{1}$ be the arc which is obtained from $i$ by changing the label at the end $p$ and $i_{2}$ be the arc obtained from $i$ by changing the label at the end $q$. We denote by $i_{0}^{1}$ and $i_{0}^{2}$ the curves without label associated to $i_{1}$ and $i_{2}$ respectively. If the version without label of the arc $i$ belongs to $\tau^{\circ}$ then we define $i_{0}:=i_{0}^{1}\cup i_{0}^{2}$ (See (a) of Figure \ref{i' a12)}) and we assign the label "plain" at both ends of $i_{0}^{1}$ and $i_{0}^{2}$, as is shown in Figure \ref{i' a13)}. Otherwise let $\gamma':[q'_{0},q'_{1}]_{i_{0}^{1}}$ be the segment of $i_{0}^{1}$ which surrounds the puncture $q$ where the ends $q'_{0}$ and $q'_{1}$ of $\gamma'$ belong to an arc $j_{1}$ of $\tau^{\circ}$. Analogously for $i_{0}^{2}$ we have the segment $\gamma:=[q_{0},q_{1}]_{i_{0}^{2}}$ which surrounds the puncture $p$, where the ends $q_{0}$ and $q_{1}$ of $\gamma$ belong to an arc $j_{2}$ of $\tau^{\circ}$. With this information we observe the following three points.

\begin{myEnumerate}
\item[1)]The arc $i$ intersects the arcs $j_{1}$ and $j_{2}$. Let $c_{1}$ and $c_{2}$ be the intersection points of the arc $i$ with $j_{1}$ and $j_{2}$ respectively.
\item[2)] Let $x$ be the intersection point of $i_{0}^{1}$ with $j_{2}$ and $z$ the intersection point of $i_{2}$ with $j_{1}$. The segments $[c_{2},c_{1}]_{i}$, $[x,q'_{1}]_{i_{0}^{1}}$ and $[q_{1},z]_{i_{0}^{2}}$ are homotopics with the homotopy that avoids $M$ and each of whose intermediate maps are segments with endpoints in the arcs $j_{1}$ and $j_{2}$ of $\tau^{\circ}$.

\item[3)] According to the points 1) and 2) we can identify the points $x$,$c_{2}$ with $q_{1}$ and $z$,$c_{1}$ with $q'_{1}$.

As a consequence of the three observations described above, we define $i_{0}:=\gamma \cup [q_{1},q'_{1}]\cup \gamma'$ and assign the label "plain" to both ends of $i_{0}$. See Figure \ref{i' a13)b}.

\begin{figure}[H]
\centering 
$$ \begin{tikzpicture}[scale=0.8]

\filldraw [black] (-2,0) circle (2pt)
(2,0) circle (2pt)
(-2.5,0) circle (1pt)
(-2.5,0.25) circle (1pt)
(-2.5,-0.25) circle (1pt)
(2.5,0) circle (1pt)
(2.5,0.25) circle (1pt)
(2.5,-0.25) circle (1pt);
\draw (-1.7,-0.25)node{$p$};
\draw (1.7,-0.25)node{$q$};
\draw[line width=1pt] (-2,0) -- (-2, 1);
\draw[line width=1pt] (-2,0) -- (-1.5, 1);
\draw[line width=1pt] (-2,0) -- (-1.5, -1);
\draw[line width=1pt] (-2,0) -- (-2, -1);
\draw[line width=1pt] (-2,0) -- (-2.5, -1);
\draw[line width=1pt] (-2,0) -- (-2.5, 1);

\draw[line width=1pt] (2,0) -- (2, 1);
\draw[line width=1pt] (2,0) -- (2.5, 1);
\draw[line width=1pt] (2,0) -- (2.5, -1);
\draw[line width=1pt] (2,0) -- (2, -1);
\draw[line width=1pt] (2,0) -- (1.5, -1);
\draw[line width=1pt] (2,0) -- (1.5, 1);
\draw[line width=1pt] (-2,0) -- (2, 0);

\draw [color=red] (-1.4,0.4)node {$\rotatebox{77}{\textbf{\Bowtie}}$};
\draw [color=red] (1.5,0.4)node {$\rotatebox{50}{\textbf{\Bowtie}}$};
\draw[color=red][line width=1.25pt] (-2,0) ..controls(-1.7,0.75) and (1.7,0.75) .. (2, 0)
node[pos=0.7,above] {$i$};


\filldraw [magenta] (-1.5,0.9) circle (2pt)
(-2,0.84) circle (2pt)
(-2.35,0.7) circle (2pt)
(-2.35,-0.7) circle (2pt)
(-2,-0.84) circle (2pt)
(-1.5,-0.9) circle (2pt)
(1.55,0.9) circle (2pt)
(2,0.84) circle (2pt)
(2.35,0.7) circle (2pt)
(2.35,-0.7) circle (2pt)
(2,-0.84) circle (2pt)
(1.55,-0.9) circle (2pt);
\draw[color=magenta][line width=1.25pt] (-2,0) ..controls(-2,0) and (2.5,2) .. (2.75, 0)
node[pos=0.9,right] {$i^{1}_{0}$}
node[pos=0.5,above] {$i'=i^{1}_{0} \cup i^{2}_{0}$};
\draw[color=magenta][line width=1.25pt] (-2,0) ..controls(-2,0) and (2.5,-2) .. (2.75, 0);
\draw[color=magenta][line width=1.25pt] (-2.75,0) ..controls(-2.5,-2) and (2,0) .. (2, 0);
\draw[color=magenta][line width=1.25pt] (-2.75,0) ..controls(-2.5,2) and (2,0) .. (2, 0)
node[pos=0,left] {$i^{2}_{0}$};
\end{tikzpicture}$$

\caption{Case $a_{13}$) of the definition \ref{curva i'}.}\label{i' a13)}
\end{figure} 
\begin{figure}[H]

\centering 
$$ \begin{tikzpicture}[scale=0.8]

\filldraw [black] (4,0) circle (2pt)
(8,0) circle (2pt)
(3.5,0) circle (1pt)
(3.5,0.25) circle (1pt)
(3.5,-0.25) circle (1pt)
(8.5,0) circle (1pt)
(8.5,0.25) circle (1pt)
(8.5,-0.25) circle (1pt);
\draw[line width=1pt] (4,0) -- (4, 1);
\draw[line width=1pt] (4,0) -- (4.5, 1);
\draw[line width=1pt] (4,0) -- (4.5, -1);
\draw[line width=1pt] (4,0) -- (4, -1);
\draw[line width=1pt] (4,0) -- (3.5, -1);
\draw[line width=1pt] (4,0) -- (3.5, 1);

\draw[line width=1pt] (8,0) -- (8, 1);
\draw[line width=1pt] (8,0) -- (8.5, 1);
\draw[line width=1pt] (8,0) -- (8.5, -1);
\draw[line width=1pt] (8,0) -- (8, -1);
\draw[line width=1pt] (8,0) -- (7.5, -1);
\draw[line width=1pt] (8,0) -- (7.5, 1);

\draw[line width=1pt] (5.25,1) -- (5.25,-1)
node[pos=0.4,right] {$j_{2}$};
\draw[line width=1pt] (6.75,1) -- (6.75,-1)
node[pos=0.6,left] {$j_{1}$};

\draw (4.3,-0.25)node{$p$};
\draw (7.7,-0.25)node{$q$};

\draw[color=red][line width=1.25pt] (4,0) ..controls(4.5,-0.5) and (7.5,0.5) .. (8, 0)
node[pos=0.7,above] {$i$};
\draw [color=red] (4.6,-0.2)node {$\rotatebox{50}{\textbf{\Bowtie}}$};
\draw [color=red] (7.7,0.1)node {$\rotatebox{50}{\textbf{\Bowtie}}$};
\filldraw [magenta] (5.24,0.8) circle (2pt)
(4.5,0.95) circle (2pt)
(4,0.92) circle (2pt)
(3.65,0.75) circle (2pt)
(3.65,-0.7) circle (2pt)
(4,-0.84) circle (2pt)
(4.5,-0.9) circle (2pt)
(5.24,-0.85) circle (2pt)
(6.75,0.85) circle (2pt)
(7.55,0.9) circle (2pt)
(8,0.84) circle (2pt)
(8.35,0.7) circle (2pt)
(8.35,-0.75) circle (2pt)
(8,-0.9) circle (2pt)
(7.55,-0.95) circle (2pt)
(6.75,-0.85) circle (2pt);
\draw[color=magenta][line width=1.25pt] (4,0) ..controls(4,0) and (6.25,0.8) .. (6.75,0.8);
\draw[color=magenta][line width=1.25pt] (6.75,0.85) ..controls(6.75,0.85) and (8.75,1.25) .. (8.75, 0)
node[pos=1,right] {$i^{1}_{0}$}
node[pos=-0,above] {$q'_{1}$};
\draw[color=magenta][line width=1.25pt] (6.75,-0.8) ..controls(6.75,-0.8) and (8.5,-1.5) .. (8.75, 0)
node[pos=-0.1,left] {$q'_{0}$};

\draw[color=magenta][line width=1.25pt] (5.24,-0.85) ..controls(5.5,-0.85) and (8,0) .. (8,0);
\draw[color=magenta][line width=1.25pt] (3.25,0) ..controls(3.25,-1.25) and (5.24,-0.85) .. (5.24,-0.85)
node[pos=1,below] {$q_{1}$}
node[pos=0.1,left] {$i^{2}_{0}$};
\draw[color=magenta][line width=1.25pt] (3.25,0) ..controls(3.5,1.5) and (5.25,0.8) .. (5.25, 0.8)
node[pos=0.9,right] {$q_{0}$};

\end{tikzpicture}$$

$$ \begin{tikzpicture}[scale=0.8]
\filldraw [black] (4,0) circle (2pt)
(8,0) circle (2pt)
(3.5,0) circle (1pt)
(3.5,0.25) circle (1pt)
(3.5,-0.25) circle (1pt)
(8.5,0) circle (1pt)
(8.5,0.25) circle (1pt)
(8.5,-0.25) circle (1pt);
\draw[line width=1pt] (4,0) -- (4, 1);
\draw[line width=1pt] (4,0) -- (4.5, 1);
\draw[line width=1pt] (4,0) -- (4.5, -1);
\draw[line width=1pt] (4,0) -- (4, -1);
\draw[line width=1pt] (4,0) -- (3.5, -1);
\draw[line width=1pt] (4,0) -- (3.5, 1);

\draw[line width=1pt] (8,0) -- (8, 1);
\draw[line width=1pt] (8,0) -- (8.5, 1);
\draw[line width=1pt] (8,0) -- (8.5, -1);
\draw[line width=1pt] (8,0) -- (8, -1);
\draw[line width=1pt] (8,0) -- (7.5, -1);
\draw[line width=1pt] (8,0) -- (7.5, 1);

\draw[line width=1pt] (5.25,1) -- (5.25,-1)
node[pos=0.4,right] {$j_{2}$};
\draw[line width=1pt] (6.75,1) -- (6.75,-1)
node[pos=0.6,left] {$j_{1}$};

\draw (4.3,-0.25)node{$p$};
\draw (7.7,-0.25)node{$q$};

\draw[color=red][line width=1.5pt] (4,0) ..controls(4.5,-0.5) and (7.5,0.5) .. (8, 0)
node[pos=0.7,above] {$i$};
\draw [color=red] (4.6,-0.2)node {$\rotatebox{50}{\textbf{\Bowtie}}$};
\draw [color=red] (7.7,0.1)node {$\rotatebox{50}{\textbf{\Bowtie}}$};
\filldraw [magenta] (5.24,0.8) circle (2pt)
(4.5,0.95) circle (2pt)
(4,0.92) circle (2pt)
(3.65,0.75) circle (2pt)
(3.65,-0.7) circle (2pt)
(4,-0.84) circle (2pt)
(4.5,-0.9) circle (2pt)
(5.24,-0.85) circle (2pt)
(6.75,0.85) circle (2pt)
(7.55,0.9) circle (2pt)
(8,0.84) circle (2pt)
(8.35,0.7) circle (2pt)
(8.35,-0.75) circle (2pt)
(8,-0.9) circle (2pt)
(7.55,-0.95) circle (2pt)
(6.75,-0.85) circle (2pt);
\draw[color=magenta][line width=1.5pt] (5.25,-0.8) .. controls(5.75,-1) and (6.25,1) .. (6.75,0.8);
\draw[dashed][color=magenta][line width=1.5pt] (4,0) ..controls(4,0) and (6.25,0.8) .. (6.75,0.8);
\draw[color=magenta][line width=1.5pt] (6.75,0.85) ..controls(6.75,0.85) and (8.75,1.25) .. (8.75, 0)
node[pos=0.8,right] {$i^{1}_{0}$}
node[pos=0,above] {$q'_{1}$};
\draw[color=magenta][line width=1.25pt] (6.75,-0.8) ..controls(6.75,-0.8) and (8.5,-1.5) .. (8.75, 0)
node[pos=-0.1,left] {$q'_{0}$};

\draw[dashed][color=magenta][line width=1.25pt] (5.24,-0.85) ..controls(5.5,-0.85) and (8,0) .. (8,0);
\draw[color=magenta][line width=1.25pt] (3.25,0) ..controls(3.25,-1.25) and (5.24,-0.85) .. (5.24,-0.85)
node[pos=1,below] {$q_{1}$};
\draw[color=magenta][line width=1.25pt] (3.25,0) ..controls(3.5,1.5) and (5.25,0.8) .. (5.25, 0.8)
node[pos=0.9,right] {$q_{0}$}
node[pos=0.2,left] {$i^{2}_{0}$};

\draw [color=magenta] (6,-1.5)node {$i':=[q_{0},q_{1}]_{i^{2}_{0}}\cup[q_{1},q'_{1}]\cup[q'_{0},q'_{1}]_{i^{1}_{0}}$};

\filldraw [black] (-2,0) circle (2pt)
(2,0) circle (2pt)
(-2.5,0) circle (1pt)
(-2.5,0.25) circle (1pt)
(-2.5,-0.25) circle (1pt)
(2.5,0) circle (1pt)
(2.5,0.25) circle (1pt)
(2.5,-0.25) circle (1pt);
\draw[line width=1pt] (-2,0) -- (-2, 1);
\draw[line width=1pt] (-2,0) -- (-2.5, 1);
\draw[line width=1pt] (-2,0) -- (-2.5, -1);
\draw[line width=1pt] (-2,0) -- (-2, -1);
\draw[line width=1pt] (-2,0) -- (-1.5, -1);
\draw[line width=1pt] (-2,0) -- (-1.5, 1);

\draw[line width=1pt] (2,0) -- (2, 1);
\draw[line width=1pt] (2,0) -- (2.5, 1);
\draw[line width=1pt] (2,0) -- (2.5, -1);
\draw[line width=1pt] (2,0) -- (2, -1);
\draw[line width=1pt] (2,0) -- (1.5, -1);
\draw[line width=1pt] (2,0) -- (1.5, 1);

\draw[line width=1pt] (-0.75,1) -- (-0.75,-1);
\draw[line width=1pt] (0.75,1) -- (0.75,-1);

\draw (-1.7,-0.25)node{$p$};
\draw (1.7,-0.25)node{$q$};
\draw [magenta] (-1,-0.25)node{$c_{2}$};
\draw [magenta] (1,0.35)node{$c_{1}$};

\draw[color=red][line width=1.25pt] (-2,0) ..controls(-1.5,-0.5) and (0.5,0.5) .. (2, 0)
node[pos=0.7,above] {$i$};
\draw [color=red] (-1.4,-0.2)node {$\rotatebox{50}{\textbf{\Bowtie}}$};
\draw [color=red] (1.7,0.1)node {$\rotatebox{50}{\textbf{\Bowtie}}$};
\filldraw [magenta] (-0.76,0.8) circle (2pt)
(-1.5,0.95) circle (2pt)
(-2,0.92) circle (2pt)
(-2.35,0.75) circle (2pt)
(-2.35,-0.7) circle (2pt)
(-2,-0.84) circle (2pt)
(-1.5,-0.9) circle (2pt)
(-0.76,-0.85) circle (2pt)
(0.75,0.85) circle (2pt)
(1.55,0.9) circle (2pt)
(2,0.84) circle (2pt)
(2.35,0.7) circle (2pt)
(2.35,-0.75) circle (2pt)
(2,-0.9) circle (2pt)
(1.55,-0.95) circle (2pt)
(0.75,-0.85) circle (2pt)
(-0.76,0.45) circle (2pt)
(-0.76,-0.05) circle (2pt)
(0.76,0.1) circle (2pt)
(0.76,-0.4) circle (2pt);
\draw[color=magenta][line width=1.25pt] (-0.75,-0.8) .. controls(-0.75,-1) and (0.25,1) .. (0.75,0.8)
node[pos=0,below] {$q_{1}$};
\draw[dashed][color=magenta][line width=1.25pt] (-2,0) ..controls(-2,0) and (0.25,0.8) .. (0.75,0.8)
node[pos=0.45,above] {$x$};

\draw[color=magenta][line width=1.25pt] (0.75,0.85) ..controls(0.75,0.85) and (2.75,1.25) .. (2.75, 0)
node[pos=0,above] {$q'_{1}$};
\draw[color=magenta][line width=1.25pt] (0.75,-0.8) ..controls(0.75,-0.8) and (2.5,-1.5) .. (2.75, 0);

\draw[dashed][color=magenta][line width=1.25pt] (-0.76,-0.85) ..controls(-0.5,-0.85) and (2,0) .. (2,0)
node[pos=0.55,below] {$z$};
\draw[color=magenta][line width=1.25pt] (-2.75,0) ..controls(-2.75,-1.25) and (-0.76,-0.85) .. (-0.76,-0.85);
\draw[color=magenta][line width=1.25pt] (-2.75,0) ..controls(-2.5,1.5) and (-0.75,0.8) .. (-0.75, 0.8);

\end{tikzpicture}$$

\caption{Case $a_{13}$) of the Definition \ref{curva i'}.}\label{i' a13)b}
\end{figure}

\end{myEnumerate}

\item[$a_{2}$)] $\delta_{\tau}(p) =0$ and $ \delta_{\tau}(q)\neq 0$ or $\delta_{\tau}(p) \neq0$ and $ \delta_{\tau}(q) = 0$. 

Without loss of generality it is enough to consider the case $\delta_{\tau}(p) =0$ and $ \delta_{\tau}(q)\neq 0$ since the case $\delta_{\tau}(p) \neq0$ and $ \delta_{\tau}(q) = 0$ is symmetric. Let's denote by $j'$ the unique arc of $\tau^{\circ}$ which is incident at $i$ and by $m$ to the loop based on the opposite puncture to $p$ in $j'$  which encloses closely to $j'$.

\item[$a_{21}$)] The arc $i$ has a label "plain" at both ends.

Under these conditions we define the curve $i_{0}:= i$ and assign the label "plain" at both ends, as is shown in Figure \ref{i' a21)}.

\begin{figure}[H]
\centering 
$$ \begin{tikzpicture}[scale=0.8]

\filldraw [black] (-2,0) circle (2pt)
(2,0) circle (2pt)
(-2,-2) circle (2pt)
(2.3,0) circle (1pt)
(2.3,-0.2) circle (1pt)
(2.3,0.2) circle (1pt);

\draw[line width=1pt] (2,0) -- (2, 1);
\draw[line width=1pt] (2,0) -- (2.5, 1);
\draw[line width=1pt] (2,0) -- (2.5, -1);
\draw[line width=1pt] (2,0) -- (2, -1);
\draw[line width=1pt] (2,0) -- (1.5, -1);
\draw[line width=1pt] (2,0) -- (1.5, 1);

\draw (-1.7,-0.25)node{$p$};
\draw (1.7,-0.25)node{$q$};

\draw[line width=1pt] (-2,-2) -- (-2, 0);
\draw[line width=1pt] (-2,-2) ..controls(-2,-2) and (-3.5,0.5) .. (-2, 0.5);
\draw[line width=1pt] (-2,-2) ..controls(-2,-2) and (-0.5,0.5) .. (-2, 0.5);
\filldraw [magenta] (-1.45,0.2) circle (2pt);
\draw[color=magenta][line width=1.5pt] (-2,0) ..controls(-1.5,0.5) and (1.5,0.5) .. (2, 0)
node[pos=0.5,above] {$i'$};

\end{tikzpicture}$$
\caption{Case $a_{21}$) of the Definition \ref{curva i'}.}\label{i' a21)}
\end{figure}
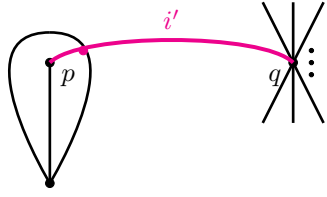

\item[$a_{22}$)] The arc $i$ has a label "plain" at the end $p$ and label $"\Bowtie"$ at the end $q$.

Let $k$ be the loop based on $p$ which encloses closely to $i$, if the version without labels of $i$ belongs to $\tau^{\circ}$ then it defines $i':=k$. Otherwise, let's consider the next two cases:

\begin{itemize}
\item[$a_{221}$)] The arc $i$ itersects more than one arc of $\tau^{\circ}$.
Under this condition $i_{0}$ is defined as  in $a_{12})$ since the signature at the end $p$ does not affect the construction of $i_{0}$. The curve $i_{0}$ is  labeled "plain" at both ends, as is shown in the Figure \ref{i' a22)}.

\begin{figure}[H]
\centering 
$$ \begin{tikzpicture}[scale=0.8]

\filldraw [black] (-2,0) circle (2pt)
(2,0) circle (2pt)
(-2,-2) circle (2pt)
(2.3,0) circle (1pt)
(2.3,-0.2) circle (1pt)
(2.3,0.2) circle (1pt);

\draw[line width=1pt] (2,0) -- (2, 1);
\draw[line width=1pt] (2,0) -- (2.5, 1);
\draw[line width=1pt] (2,0) -- (2.5, -1);
\draw[line width=1pt] (2,0) -- (2, -1);
\draw[line width=1pt] (2,0) -- (1.5, -1);
\draw[line width=1pt] (2,0) -- (1.5, 1);
\draw[line width=1pt] (0.5,-2) -- (0.5,1)
node[pos=0.7,right] {$j$};
\draw (-1.7,-0.25)node{$p$};
\draw (1.7,-0.25)node{$q$};

\draw[line width=1pt] (-2,-2) -- (-2, 0);
\draw[line width=1pt] (-2,-2) ..controls(-2,-2) and (-3.5,0.5) .. (-2, 0.5);
\draw[line width=1pt] (-2,-2) ..controls(-2,-2) and (-0.5,0.5) .. (-2, 0.5);
\filldraw [magenta] (-1.45,0.2) circle (2pt)
(0.5,0.75) circle (2pt)
(1.55,0.9) circle (2pt)
(2,0.85) circle (2pt)
(2.4,0.75) circle (2pt)
(2.4,-0.85) circle (2pt)
(2,-0.95) circle (2pt)
(1.55,-0.9) circle (2pt)
(0.5,-0.5) circle (2pt);

\draw[color=magenta][line width=1.5pt] (-2,0) ..controls(-2,0) and (2.75,2) .. (3, 0)
node[pos=0.35,above] {$i'$}
node[pos=0.49,below] {$q_{1}$};
\draw[color=magenta][line width=1.5pt] (0.5,-0.5) ..controls(0.5,-0.5) and (2.75,-1.75) .. (3, 0)
node[pos=0,left] {$q_{0}$};

\end{tikzpicture}$$
\caption{Case $a_{221}$) of the Definition \ref{curva i'}.}\label{i' a22)}
\end{figure}

\item[$a_{222}$)] The arc $i$ intersects only the arc $m$ of $\tau^{\circ}$ in one point.

Touring $k$ clockwise, we denote by $t_{r}$ the $r^{th}$ intersection point of $k$ with some arc of $\tau^{\circ}$ for $r\in \{1,...,n+1\}$. The intersection points $t_{1}$ and $t_{n+1}$ are the unique two intersection points of $k$ with $m$. It is defined $i_{0}:=[p,t_{1}]_{k}\cup [t_{1},t_{n+1}]_{k}$ with the property that $t_{n+1}\notin [p,t_{1}]_{k}$ and it is  labeled "plain" at both ends. See Figure \ref{i' a222)}.

\begin{figure}[H]
\centering 
$$ \begin{tikzpicture}[scale=0.8]

\filldraw [black] (-2,0) circle (2pt)
(2,0) circle (2pt)
(-2,-2) circle (2pt)
(2.3,0) circle (1pt)
(2.3,-0.2) circle (1pt)
(2.3,0.2) circle (1pt);

\draw[line width=1pt] (2,0) -- (2, 1);
\draw[line width=1pt] (2,0) -- (2.5, 1);
\draw[line width=1pt] (2,0) -- (2.5, -1);
\draw[line width=1pt] (2,0) -- (2, -1);
\draw[line width=1pt] (2,0) -- (1.5, -1);
\draw[line width=1pt] (2,0) -- (1.5, 1);
\draw[line width=1pt] (-2,0) -- (-2, -2); 

\draw (-1.7,-0.25)node{$p$};
\draw (1.7,-0.25)node{$q$};

\draw[line width=1pt] (-2,-2) ..controls(-2,-2) and (-3.5,0.5) .. (-2, 0.5);
\draw[line width=1pt] (-2,-2) ..controls(-2,-2) and (-0.5,0.5) .. (-2, 0.5);
\filldraw [magenta] (-1.45,0.2) circle (2pt)
(1.55,0.9) circle (2pt)
(2,0.85) circle (2pt)
(2.4,0.75) circle (2pt)
(2.4,-0.85) circle (2pt)
(2,-0.95) circle (2pt)
(1.55,-0.9) circle (2pt)
(-1.45,-0.5) circle (2pt);

\draw[color=magenta][line width=1.5pt] (-2,0) ..controls(-2,0) and (2.75,2) .. (3, 0)
node[pos=0.35,above] {$i'$}
node[pos=0.25,above] {$t_{1}$};
\draw[color=magenta][line width=1.5pt] (-1.4,-0.5) ..controls(0.5,-0.5) and (2.75,-1.75) .. (3, 0)
node[pos=0.1,below] {$t_{n+1}$};

\end{tikzpicture}$$
\caption{Case $a_{221}$) of the Definition \ref{curva i'}.}\label{i' a222)}
\end{figure}

\end{itemize} 
\item[$a_{23}$)] The arc $i$ has a label "plain" at the end $q$ and label $"\Bowtie"$ at the end $p$.

Under this condition it is defined $i_{0}$ as the version without labels of the arc $i$. The curve $i_{0}$ is labeled "plain" at the ends $q$ and $\Bowtie$  at the end $p$, as is shown in Figure \ref{i' a23)}. 
\begin{figure}[H]
\centering 
$$ \begin{tikzpicture}[scale=0.8]

\filldraw [black] (-2,0) circle (2pt)
(2,0) circle (2pt)
(-2,-2) circle (2pt);

\draw[line width=1pt] (2,0) -- (2, 1);
\draw[line width=1pt] (2,0) -- (2.5, 1);
\draw[line width=1pt] (2,0) -- (2.5, -1);
\draw[line width=1pt] (2,0) -- (2, -1);
\draw[line width=1pt] (2,0) -- (1.5, -1);
\draw[line width=1pt] (2,0) -- (1.5, 1);

\draw (-1.7,-0.25)node{$p$};
\draw (1.7,-0.25)node{$q$};

\draw[line width=1pt] (-2,-2) -- (-2, 0);
\draw[line width=1pt] (-2,-2) ..controls(-2,-2) and (-3.5,0.5) .. (-2, 0.5);
\draw[line width=1pt] (-2,-2) ..controls(-2,-2) and (-0.5,0.5) .. (-2, 0.5);
\filldraw [magenta] (-1.45,0.2) circle (2pt);
\draw [color=red] (-1.85,0.1)node {$\rotatebox{-50}{\textbf{\Bowtie}}$}; 
\draw[color=magenta][line width=1.5pt] (-2,0) ..controls(-1.5,0.5) and (1.5,0.5) .. (2, 0)
node[pos=0.5,above] {$i'$};

\end{tikzpicture}$$
\caption{Case $a_{23}$) of the Definition \ref{curva i'}.}\label{i' a23)}
\end{figure}
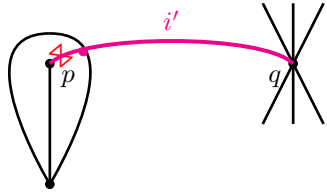

\item[$a_{24}$)] The arc $i$ has a label $"\Bowtie"$ at both ends.

Let $i^{1}$ be the arc which is obtained from $i$ by changing the label at the end $p$. For $i^{1}$ we have defined the curve $i_{0}^{1}$ (see $a_{22}$)). The curve $i_{0}:=i_{0}^{1}$ and it is assigned the label $\Bowtie$ at the end  $p$  and the label "plain" at the opposite end to $p$ in $i_{0}$, as is shown in Figure \ref{i' a24)}.

\begin{figure}[H]
\centering 
$$ \begin{tikzpicture}[scale=0.8]

\filldraw [black] (-2,0) circle (2pt)
(2,0) circle (2pt)
(-2,-2) circle (2pt)
(2.3,0) circle (1pt)
(2.3,-0.2) circle (1pt)
(2.3,0.2) circle (1pt);

\draw[line width=1pt] (2,0) -- (2, 1);
\draw[line width=1pt] (2,0) -- (2.5, 1);
\draw[line width=1pt] (2,0) -- (2.5, -1);
\draw[line width=1pt] (2,0) -- (2, -1);
\draw[line width=1pt] (2,0) -- (1.5, -1);
\draw[line width=1pt] (2,0) -- (1.5, 1);
\draw[line width=1pt] (0.5,-2) -- (0.5,1)
node[pos=0.65,right] {$j$};
\draw (-1.7,-0.25)node{$p$};
\draw (1.7,-0.25)node{$q$};

\draw[line width=1pt] (-2,-2) -- (-2, 0);
\draw[line width=1pt] (-2,-2) ..controls(-2,-2) and (-3.5,0.5) .. (-2, 0.5);
\draw[line width=1pt] (-2,-2) ..controls(-2,-2) and (-0.5,0.5) .. (-2, 0.5);
\filldraw [magenta] (-1.45,0.2) circle (2pt)
(0.5,0.75) circle (2pt)
(1.55,0.9) circle (2pt)
(2,0.85) circle (2pt)
(2.4,0.75) circle (2pt)
(2.4,-0.85) circle (2pt)
(2,-0.95) circle (2pt)
(1.55,-0.9) circle (2pt)
(0.5,-0.5) circle (2pt);

\draw[color=magenta][line width=1.5pt] (-2,0) ..controls(-2,0) and (2.75,2) .. (3, 0)
node[pos=0.35,above] {$i_{0}$}
node[pos=0.49,below] {$q_{1}$};
\draw[color=magenta][line width=1.5pt] (0.5,-0.5) ..controls(0.5,-0.5) and (2.75,-1.75) .. (3, 0)
node[pos=0,left] {$q_{0}$}; 

\filldraw [black] (4,0) circle (2pt)
(8,0) circle (2pt)
(4,-2) circle (2pt)
(8.3,0) circle (1pt)
(8.3,-0.2) circle (1pt)
(8.3,0.2) circle (1pt);

\draw[line width=1pt] (8,0) -- (8, 1);
\draw[line width=1pt] (8,0) -- (8.5, 1);
\draw[line width=1pt] (8,0) -- (8.5, -1);
\draw[line width=1pt] (8,0) -- (8, -1);
\draw[line width=1pt] (8,0) -- (7.5, -1);
\draw[line width=1pt] (8,0) -- (7.5, 1);
\draw[line width=1pt] (6.5,-2) -- (6.5,1);
\draw (4.3,-0.25)node{$p$};
\draw (7.7,-0.25)node{$q$};

\draw[line width=1pt] (4,-2) -- (4, 0);
\draw[line width=1pt] (4,-2) ..controls(4,-2) and (2.5,0.5) .. (4, 0.5);
\draw[line width=1pt] (4,-2) ..controls(4,-2) and (5.5,0.5) .. (4, 0.5);
\filldraw [magenta] (4.55,0.2) circle (2pt)
(6.5,0.75) circle (2pt)
(7.55,0.9) circle (2pt)
(8,0.85) circle (2pt)
(8.4,0.75) circle (2pt)
(8.4,-0.85) circle (2pt)
(8,-0.95) circle (2pt)
(7.55,-0.9) circle (2pt)
(6.5,-0.5) circle (2pt);

\draw [color=red] (4.15,0.1)node {$\rotatebox{-50}{\textbf{\Bowtie}}$}; 
\draw[color=magenta][line width=1.5pt] (4,0) ..controls(4,0) and (8.75,2) .. (9, 0)
node[pos=0.35,above] {$i'$};
\draw[color=magenta][line width=1.5pt] (6.5,-0.5) ..controls(6.5,-0.5) and (8.75,-1.75) .. (9, 0);

\end{tikzpicture}$$
\caption{Case $a_{24}$) of the Definition \ref{curva i'}.}\label{i' a24)}
\end{figure}
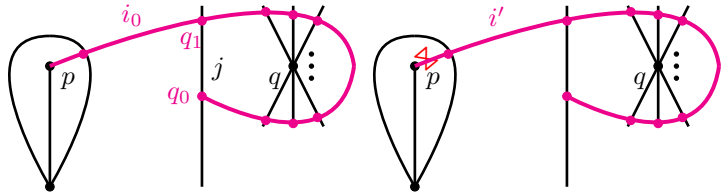 
\end{myEnumerate}

\item[$a_{3}$)] $\delta_{\tau}(p) = 0= \delta_{\tau}(q)$.

For this case,  the curve $i_{0}$ is defined as the version without label of $i$. The labels are assigned  at each end of $i_{0}$ depending on of the label at each end of $i$:

\begin{myEnumerate}
\item[$a_{31}$)] If the arc $i$ has a label "plain" at both ends.

The curve $i_{0}$ is labeled "plain" at both ends, as is shown in Figure \ref{i' a31)}.

\begin{figure}[H]
\centering 
$$ \begin{tikzpicture}[scale=0.8]

\filldraw [black] (-2,0) circle (2pt)
(2,0) circle (2pt)
(-2,-2) circle (2pt)
(2,-2) circle (2pt);

\draw (-1.7,-0.25)node{$p$};
\draw (1.7,-0.25)node{$q$};

\draw[line width=1pt] (-2,-2) -- (-2, 0);
\draw[line width=1pt] (-2,-2) ..controls(-2,-2) and (-3.5,0.5) .. (-2, 0.5);
\draw[line width=1pt] (-2,-2) ..controls(-2,-2) and (-0.5,0.5) .. (-2, 0.5);

\draw[line width=1pt] (2,-2) -- (2, 0);
\draw[line width=1pt] (2,-2) ..controls(2,-2) and (0.5,0.5) .. (2, 0.5);
\draw[line width=1pt] (2,-2) ..controls(2,-2) and (3.5,0.5) .. (2, 0.5);
\filldraw [magenta] (-1.45,0.2) circle (2pt)
(1.45,0.25) circle (2pt);
\draw[color=magenta][line width=1.5pt] (-2,0) ..controls(-1.5,0.5) and (1.5,0.5) .. (2, 0)
node[pos=0.5,above] {$i'$};

\end{tikzpicture}$$
\caption{Case $a_{31}$) of the Definition \ref{curva i'}.}\label{i' a31)}
\end{figure}
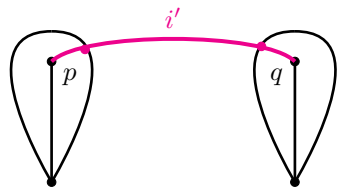 

\item[$a_{32}$)] The arc $i$ has a label $\Bowtie$ at only one end. Without loss of generality, we can assume that the end with label $\Bowtie$ is the puncture $q$. 

The curve $i_{0}$  is labeled  "plain" at the end $p$ and $\Bowtie$ at the $q$, as is shown in Figure \ref{i' a32)}.

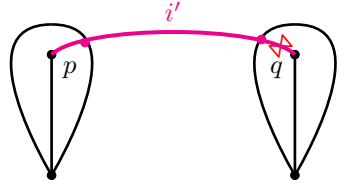
\begin{figure}[H]
\centering 
$$ \begin{tikzpicture}[scale=0.8]

\filldraw [black] (-2,0) circle (2pt)
(2,0) circle (2pt)
(2,-2) circle (2pt)
(-2,-2) circle (2pt);

\draw (-1.7,-0.25)node{$p$};
\draw (1.7,-0.25)node{$q$};

\draw[line width=1pt] (-2,-2) -- (-2, 0);
\draw[line width=1pt] (-2,-2) ..controls(-2,-2) and (-3.5,0.5) .. (-2, 0.5);
\draw[line width=1pt] (-2,-2) ..controls(-2,-2) and (-0.5,0.5) .. (-2, 0.5);

\draw[line width=1pt] (2,-2) -- (2, 0);
\draw[line width=1pt] (2,-2) ..controls(2,-2) and (0.5,0.5) .. (2, 0.5);
\draw[line width=1pt] (2,-2) ..controls(2,-2) and (3.5,0.5) .. (2, 0.5);
\filldraw [magenta] (-1.45,0.2) circle (2pt)
(1.45,0.25) circle (2pt);
\draw [color=red] (1.8,0.1)node {$\rotatebox{40}{\textbf{\Bowtie}}$}; 
\draw[color=magenta][line width=1.5pt] (-2,0) ..controls(-1.5,0.5) and (1.5,0.5) .. (2, 0)
node[pos=0.5,above] {$i'$};

\end{tikzpicture}$$
\caption{Case $a_{32}$) of the Definition \ref{curva i'}.}\label{i' a32)}
\end{figure} 

\item[$a_{33}$)]If the arc $i$ has a label $\Bowtie$ at both ends.

The curve $i_{0}$ is labeled  $\Bowtie$ at both ends, as is shown in Figure \ref{i' a33)}.

\begin{figure}[H]
\centering 
$$ \begin{tikzpicture}[scale=0.8]

\filldraw [black] (-2,0) circle (2pt)
(2,0) circle (2pt)
(2,-2) circle (2pt)
(-2,-2) circle (2pt);
\draw (-1.7,-0.25)node{$p$};
\draw (1.7,-0.25)node{$q$};

\draw[line width=1pt] (-2,-2) -- (-2, 0);
\draw[line width=1pt] (-2,-2) ..controls(-2,-2) and (-3.5,0.5) .. (-2, 0.5);
\draw[line width=1pt] (-2,-2) ..controls(-2,-2) and (-0.5,0.5) .. (-2, 0.5);

\draw[line width=1pt] (2,-2) -- (2, 0);
\draw[line width=1pt] (2,-2) ..controls(2,-2) and (0.5,0.5) .. (2, 0.5);
\draw[line width=1pt] (2,-2) ..controls(2,-2) and (3.5,0.5) .. (2, 0.5);
\filldraw [magenta] (-1.45,0.2) circle (2pt)
(1.45,0.25) circle (2pt);
\draw [color=red] (-1.85,0.1)node {$\rotatebox{-50}{\textbf{\Bowtie}}$}; 
\draw [color=red] (1.8,0.1)node {$\rotatebox{40}{\textbf{\Bowtie}}$}; 
\draw[color=magenta][line width=1.5pt] (-2,0) ..controls(-1.5,0.5) and (1.5,0.5) .. (2, 0)
node[pos=0.5,above] {$i'$};

\end{tikzpicture}$$
\caption{Case $a_{33}$) of the Definition \ref{curva i'}.}\label{i' a33)}
\end{figure}
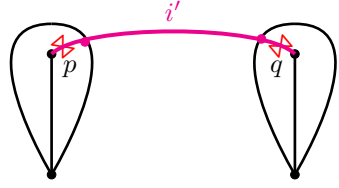 
\end{myEnumerate}
\end{myEnumerate}

\item[b)] $p= q$.

Since $p=q$ and the arc  $i$ belongs to an other triangulation $\tau_{0}$ of $(\Sigma,M)$, then the labels at the ends should be the same.

\begin{myEnumerate}
\item[$b_{1}$)] $\delta_{\tau}(p) \neq 0$.
\begin{myEnumerate}
\item[$b_{11}$)] The arc $i$ has a label "plain" at both ends.

The curve  $i_{0}:=i$  and it is labeled  "plain" at both ends of $i_{0}$. See Figure \ref{i' b11)}

\begin{figure}[H]
\centering 
$$ \begin{tikzpicture}[scale=0.8]
\filldraw [black] (-2,0) circle (2pt)
(-2,1) circle (2pt)
(-2.5,0) circle (1pt)
(-2.5,-0.5) circle (1pt)
(-2.5,-0.25) circle (1pt)
(-1.5,0) circle (1pt)
(-1.5,-0.5) circle (1pt)
(-1.5,-0.25) circle (1pt)
(-2.25,0.5) circle (1pt)
(-2,0.5) circle (1pt)
(-1.75,0.5) circle (1pt)
(-2.25,1.65) circle (1pt)
(-2,1.65) circle (1pt)
(-1.75,1.65) circle (1pt)
(-2.25,2.25) circle (1pt)
(-2,2.25) circle (1pt)
(-1.75,2.25) circle (1pt);
\draw[line width=1pt] (-2,1) -- (-1, 2.5);
\draw[line width=1pt] (-2,0) -- (-1, 0.25);
\draw[line width=1pt] (-2,0) -- (-3, 0.25);
\draw[line width=1pt] (-2,0) -- (-1.5, -1);
\draw[line width=1pt] (-2,0) -- (-2.5, -1);
\draw[line width=1pt] (-2,1) -- (-3, 2.5);

\draw (-2,-0.45)node{$p$};


\draw[color=magenta][line width=1.5pt] (-2,0) ..controls(-3.5,0.5) and (-3.5,2) .. (-2, 2)
node[pos=0.5,left] {$i'$};
\draw[color=magenta][line width=1.5pt] (-2,0) ..controls(-0.5,0.5) and (-0.5,2) .. (-2, 2);

\end{tikzpicture}$$

\caption{Case $b_{11}$) of the Definition \ref{curva i'}.}\label{i' b11)}.
\end{figure}
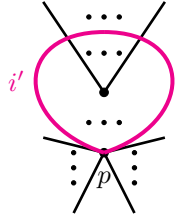 

\item[$b_{12}$)] The arc $i$ has a label $"\Bowtie"$ at both ends. 

\begin{myEnumerate}
\item[$b_{121}$)] The version without label of $i$ does not belong to the triangulation $\tau^{\circ}$. Let $i^{1}$ be the curve in $(\Sigma,M)$ which is obtained from $i$ by changing the label of $i$ and $i^{2}$ the curve in $(\Sigma,M)$ which is obtained from $i$ by changing the label of $i$ at the opposite end to the one used to obtained $i^{1}$. The curve $i'$ is obtained by doing the same construction that we made in $a_{13})$ forgetting the hypothesis that $i_{1}$ and $i_{2}$ are  arcs of $\tau^{\circ}$, as is shown in Figure \ref{i' b12)b}.
\begin{figure}[H]
\centering 
$$ \begin{tikzpicture}[scale=0.7]
\filldraw [black] (2,0) circle (2pt)
(1.5,0) circle (1pt)
(1.5,-0.5) circle (1pt)
(1.5,-0.25) circle (1pt)
(2.5,0) circle (1pt)
(2.5,-0.5) circle (1pt)
(2.5,-0.25) circle (1pt);
\draw[line width=1pt] (2,0) -- (2.5, 1);
\draw[line width=1pt] (2,0) -- (3, 0.25);
\draw[line width=1pt] (2,0) -- (0.25, 0.25);
\draw[line width=1pt] (2,0) -- (2.5, -1);
\draw[line width=1pt] (2,0) -- (1.5, -1);
\draw[line width=1pt] (2,0) -- (1.5, 1);

\draw[line width=1pt] (0.25,1) -- (1.25, 1.5);

\draw[line width=1pt] (3,1.5) -- (4, 1);

\draw (2,-0.45)node{$p$};
\draw (-0.65,1)node{$j_{1}$};
\draw (-3.05,1)node{$j_{2}$};


\draw[dashed][color=magenta][line width=1pt] (2,0) ..controls(2,0) and (3.25,0.25).. (3.5, 1.25)
node[pos=0.97,left] {$z$};
\draw[color=magenta][line width=1pt] (3.5, 1.25) ..controls(3.25,2.5) and (0.75,2.5).. (1,1.35);
\draw[color=magenta][line width=1pt] (1,1.35) ..controls(1.25, 0.25) and (3.5,0.75) ..(2.75,-0.75)
node[pos=1,right] {$i_{0}^{2}$}
node[pos=0,right] {$q_{1}$};
\draw[color=magenta][line width=1pt] (0.35,1.05) ..controls(0.75,-0.5) and (2.25,-1.5)..(2.75,-0.75)
node[pos=0,above] {$q_{0}$};

\filldraw [magenta] (3.5,1.25) circle (2pt)
(1,1.4) circle (2pt)
(1.65,0.7) circle (2pt)
(2.25,0.45) circle (2pt)
(2.7,0.2) circle (2pt)
(2.5,-0.95) circle (2pt)
(1.6,-0.75) circle (2pt)
(0.7,0.2) circle (2pt)
(0.35,1.1) circle (2pt);

\filldraw [black] (-2,0) circle (2pt)
(-2.5,0) circle (1pt)
(-2.5,-0.5) circle (1pt)
(-2.5,-0.25) circle (1pt)
(-1.5,0) circle (1pt)
(-1.5,-0.5) circle (1pt)
(-1.5,-0.25) circle (1pt);
\draw[line width=1pt] (-2,0) -- (-1.5, 1);
\draw[line width=1pt] (-2,0) -- (-0.5, 0.25);
\draw[line width=1pt] (-2,0) -- (-3, 0.25);
\draw[line width=1pt] (-2,0) -- (-1.5, -1);
\draw[line width=1pt] (-2,0) -- (-2.5, -1);
\draw[line width=1pt] (-2,0) -- (-2.5, 1);

\draw[line width=1pt] (-3.75,1) -- (-2.75, 1.5);

\draw[line width=1pt] (-1,1.5) -- (0, 1);

\draw (-2,-0.45)node{$p$};
\draw (3.65,1)node{$j_{1}$};
\draw (0.75,1)node{$j_{2}$};


\filldraw [magenta] (-3.45,1.15) circle (2pt)
(-0.15,1.1) circle (2pt)
(-0.75,0.2) circle (2pt)
(-1.7,-0.6) circle (2pt)
(-2.5,-0.95) circle (2pt)
(-2.8,0.2) circle (2pt)
(-2.3,0.5) circle (2pt)
(-1.6,0.75) circle (2pt)
(-0.9,1.4) circle (2pt);

\draw[dashed][color=magenta][line width=1pt] (-2,0) ..controls(-2,0) and (-3.55,0.35) .. (-3.45, 1.15)
node[pos=0.9,left] {$x$};
\draw[color=magenta][line width=1pt] (-3.45, 1.1)..controls(-3.45,2.5) and (-0.35,3).. (-0.15, 1.15);
\draw[color=magenta][line width=1pt] (-0.15, 1.15) ..controls(-0.15,0.5) and (-2.9,-2) .. (-3,-0.5)
node[pos=1,left] {$i_{0}^{1}$}
node[pos=0.03,below] {$q_{1}^{'}$};
\draw[color=magenta][line width=1pt] (-0.9, 1.45) ..controls(-0.9,0.45) and(-3,1) .. (-3,-0.5)
node[pos=0,above] {$q_{0}^{'}$};

\end{tikzpicture}$$
$$ \begin{tikzpicture}[scale=0.7]

\filldraw [black] (-2,0) circle (2pt)
(-2.5,0) circle (1pt)
(-2.5,-0.5) circle (1pt)
(-2.5,-0.25) circle (1pt)
(-1.5,0) circle (1pt)
(-1.5,-0.5) circle (1pt)
(-1.5,-0.25) circle (1pt);
\draw[line width=1pt] (-2,0) -- (-1.5, 1);
\draw[line width=1pt] (-2,0) -- (-0.5, 0.25);
\draw[line width=1pt] (-2,0) -- (-3.5, 0.25);
\draw[line width=1pt] (-2,0) -- (-1.5, -1);
\draw[line width=1pt] (-2,0) -- (-2.5, -1);
\draw[line width=1pt] (-2,0) -- (-2.5, 1);

\draw[line width=1pt] (-3.75,1) -- (-2.75, 1.5);

\draw[line width=1pt] (-1,1.5) -- (0, 1);

\draw (-2,-0.45)node{$p$};
\draw (-0.65,1)node{$j_{1}$};
\draw (-3.05,1)node{$j_{2}$};
\filldraw [magenta] 
(-3,1.4) circle (2pt)
(-2.35,0.7) circle (2pt)
(-1.75,0.45) circle (2pt)
(-1.25,0.15) circle (2pt)
(-1.5,-0.95) circle (2pt)
(-2.4,-0.75) circle (2pt)
(-3.3,0.2) circle (2pt)
(-3.65,1.1) circle (2pt);
\filldraw [magenta] 
(-0.15,1.1) circle (2pt)
(-0.75,0.2) circle (2pt)
(-1.7,-0.6) circle (2pt)
(-2.5,-0.95) circle (2pt)
(-2.8,0.15) circle (2pt)
(-2.3,0.5) circle (2pt)
(-1.6,0.75) circle (2pt)
(-0.9,1.4) circle (2pt);

\draw[color=magenta][line width=1pt] (-0.15, 1.15) ..controls(-0.15,0.5) and (-2.9,-2) .. (-3,-0.5)
node[pos=0.03,below] {$q_{1}^{'}$};
\draw[color=magenta][line width=1pt] (-0.9, 1.45) ..controls(-0.9,0.45) and(-3,1) .. (-3,-0.5)
node[pos=0,above] {$q_{0}^{'}$};

\draw[color=magenta][line width=1pt] (-3,1.35) ..controls(-2.75, 0.25) and (-0.5,0.75) ..(-1.25,-0.75)
node[pos=0.03,right] {$q_{1}$};
\draw[color=magenta][line width=1pt] (-3.65,1.05) ..controls(-3.25,-0.5) and (-1.75,-1.5)..(-1.25,-0.75)
node[pos=0,left] {$q_{0}$};;

\draw[color=magenta][line width=1pt] (-0.15, 1.15) ..controls(-0.5,2.75) and (-2.75,2.75) .. (-3,1.35);

\draw [color=magenta] (-2,-1.5)node {$i':=[q_{0},q_{1}]_{i^{2}_{0}}\cup[q_{1},q'_{1}]\cup[q'_{0},q'_{1}]_{i^{1}_{0}}$};
\end{tikzpicture}$$

\caption{Case $b_{121}$) de la Definici\'on \ref{curva i'}.}\label{i' b12)b}
\end{figure}
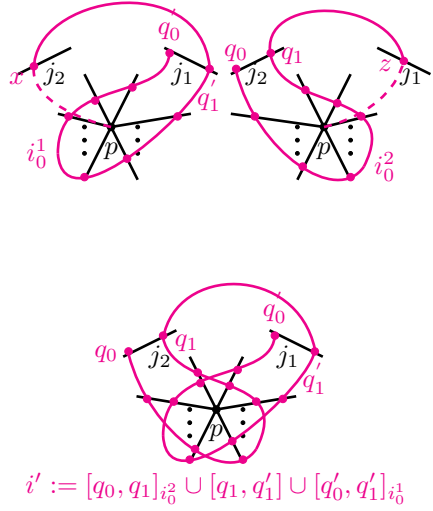

\item[$b_{122}$)]The version without labels of the arc $i$ belong to $\tau^{\circ}$. 

Let $\Delta_{1}$ and $\Delta_{2}$ be the unique non self-folded triangles of $\tau^{\circ}$ which share the version without label of the arc $i$ as one of its sides. Where  $\Delta_{1}$ is the  triangle such that  the version without labels of the arc $i$ encloses the other two sides of $\Delta_{1}$. We denote by $i,j,k$ the sides of the triangle $\Delta_{1}$ touring it clockwise and denoted by $i,j',k'$ to the sides of $\Delta_{2}$ touring it clockwise.

If the arc $j$ is not the non folded side of any self-folded triangle, then we draw a curve  $\gamma_{1}$ which surrounds the puncture $p$ and the ends of $\gamma_{1}$ belong to the arc $j$. Otherwise, if the arc $j$ is the non folded side of a self-folded triangle $\Delta'$ whit folded side $j''$, then we draw a curve  $\gamma_{1}$ which surrounds the puncture $p$ and its ends belong to the arc $j''$. So that the segment $[r_{l+1},r_{l-1}]_{\gamma_{1}}$ is not contractible to the puncture $p$. See Figure \ref{b121a} and Figure \ref{b121b}. 

\begin{figure}[H]
\begin{center}

\begin{tikzpicture}[scale=0.75]
\filldraw [black] (0,0) circle (2pt)
(0,2) circle (2pt)
(0,1) circle (1pt)
(-0.5,1) circle (1pt)
(0.5,1) circle (1pt)
(0,-0.5) circle (1pt)
(-0.5,-0.5) circle (1pt)
(0.5,-0.5) circle (1pt)
(-2.5,-0.25) circle (1pt)
(-2.5,0) circle (1pt)
(-2.5,0.25) circle (1pt)
(2.5,-0.25) circle (1pt)
(2.5,0) circle (1pt)
(2.5,0.25) circle (1pt);

\draw (0,2.5)node{$\Delta^{1}$}; 
\draw (0,4.5)node{$\Delta^{2}$}; 
\draw (0,3.75)node{$i$}; 
\draw (1,1.25)node{$j$}; 
\draw (-1,1.25)node{$k$}; 
\draw[line width=1pt] (0,0) ..controls(-4,0.5) and (-4,3.5) .. (0, 4); 
\draw[line width=1pt] (0,0) ..controls(4,0.5) and (4,3.5) .. (0, 4);

\draw[line width=1pt] (0,0) ..controls(-1.75,0.75) and (-1.75,1.25) .. (0, 2); 
\draw[line width=1pt] (0,0) ..controls(1.75,0.75) and (1.75,1.25) .. (0, 2);

\draw[line width=1pt] (0,0) ..controls (-2,0) and (-4,0.5) .. (-4, 3);
\draw[line width=1pt] (0,0) -- (-4, -1);

\draw[line width=1pt] (0,0) ..controls (2,0) and (4,0.5) .. (4, 3);
\draw[line width=1pt] (0,0) -- (4, -1);

\filldraw [magenta] (1,1.5) circle (2pt)
(-2.32,0.8) circle (2pt)
(-2.03,0.25) circle (2pt)
(-1.35,-0.32) circle (2pt)
(1.3,-0.35) circle (2pt)
(1.55,0.15) circle (2pt)
(1.53,0.45) circle (2pt)
(1.18,0.7) circle (2pt)
(-1.2,1.25) circle (2pt)
(0.35,1.85) circle (2pt);

\draw[color=magenta][line width=1.5pt] (0, 3.5) .. controls (1,3.5) and (1.5,2.75) .. (1,1.5);
\draw[color=magenta][line width=1.5pt] (0,-0.75) ..controls(-3,-0.5) and (-3.5,3) .. (0, 3.5);
\draw[color=magenta][line width=1.5pt] (0,-0.75) ..controls (2,-0.75) and (2.2,1.2) .. (0, 0.75);
\draw[color=magenta][line width=1.5pt] (0,0.75) ..controls(-1.5,0.25) and (-2,3) .. (0, 2.75); 
\draw[color=magenta][line width=1.5pt] (0, 2.75) ..controls (1,2.75) and (0.85,1.75) .. (0.35,1.85);
\end{tikzpicture}

\end{center}
\caption{The arc $j$ is not a side of any self-folded triangle. \label{b121a}}

\end{figure}
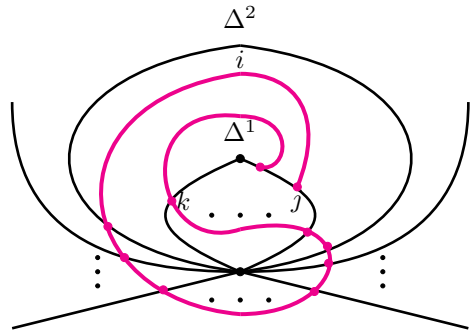

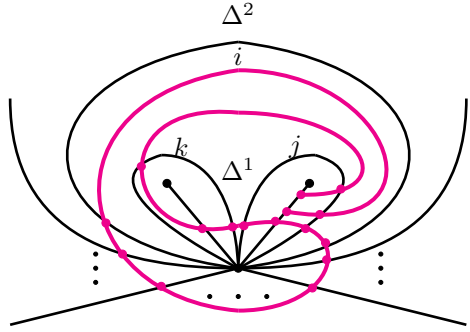
\begin{figure}[H]
\begin{center}

\begin{tikzpicture}[scale=0.75]
\filldraw [black] (0,0) circle (2pt)

(-1.25,1.5) circle (2pt)
(1.25,1.5) circle (2pt)

(0,-0.5) circle (1pt)
(-0.5,-0.5) circle (1pt)
(0.5,-0.5) circle (1pt)
(-2.5,-0.25) circle (1pt)
(-2.5,0) circle (1pt)
(-2.5,0.25) circle (1pt)
(2.5,-0.25) circle (1pt)
(2.5,0) circle (1pt)
(2.5,0.25) circle (1pt);
\draw (0,1.75)node{$\Delta^{1}$}; 
\draw (0,4.5)node{$\Delta^{2}$}; 
\draw (0,3.75)node{$i$}; 
\draw (1,2.15)node{$j$}; 
\draw (-1,2.15)node{$k$};

\draw[line width=1pt] (0,0) ..controls(-4,0.5) and (-4,3.5) .. (0, 4); 
\draw[line width=1pt] (0,0) ..controls(4,0.5) and (4,3.5) .. (0, 4);

\draw[line width=1pt] (0,0) -- (-1.25, 1.5); 
\draw[line width=1pt] (0,0) -- (1.25, 1.5);

\draw[line width=1pt] (0,0) ..controls (0,0) and (-3,1.5) .. (-1.35, 2);
\draw[line width=1pt] (0,0) ..controls (0,0) and (0,2) .. (-1.35, 2);

\draw[line width=1pt] (0,0) ..controls (0,0) and (3,1.5) .. (1.35, 2);
\draw[line width=1pt] (0,0) ..controls (0,0) and (0,2) .. (1.35, 2);

\draw[line width=1pt] (0,0) ..controls (-2,0) and (-4,0.5) .. (-4, 3);
\draw[line width=1pt] (0,0) -- (-4, -1);

\draw[line width=1pt] (0,0) ..controls (2,0) and (4,0.5) .. (4, 3);
\draw[line width=1pt] (0,0) -- (4, -1); 
\filldraw [magenta] (0.85,1) circle (2pt)
(1.43,0.95) circle (2pt)
(-2.32,0.8) circle (2pt)
(-2.03,0.25) circle (2pt)
(-1.35,-0.32) circle (2pt)
(1.3,-0.35) circle (2pt)
(1.55,0.15) circle (2pt)
(1.53,0.45) circle (2pt)
(1.18,0.7) circle (2pt)
(0.65,0.83) circle (2pt)
(0.1,0.75) circle (2pt)
(-0.1,0.73) circle (2pt)
(-0.63,0.7) circle (2pt)
(-1.7,1.8) circle (2pt)
(1.8,1.4) circle (2pt)
(1.1,1.3) circle (2pt);
\draw[color=magenta][line width=1.5pt] (0, 3.5) .. controls (2.75,3.5) and (3.75,0.5) .. (0.85,1);
\draw[color=magenta][line width=1.5pt] (0,-0.75) ..controls(-3,-0.5) and (-3.5,3) .. (0, 3.5);
\draw[color=magenta][line width=1.5pt] (0,-0.75) ..controls (2,-0.75) and (2.2,1.2) .. (0, 0.75);
\draw[color=magenta][line width=1.5pt] (0,0.75) ..controls(-2,0.25) and (-2.5,3) .. (0, 2.75); 
\draw[color=magenta][line width=1.5pt] (0, 2.75) ..controls (2.25,2.75) and (3,1.25) .. (1.1,1.3); 
\end{tikzpicture}

\end{center}
\caption{The arc $j$ is the not-folded side of a self-folded triangle $\Delta'$. \label{b121b}}

\end{figure}

Once drawn the curve $\gamma_{1}$, we notice that the cross point  of $\gamma_{1}$ with the arcs of $\tau^{\circ}$ denote   it by $r_{l}$ belong to the relative interior of the arc $j$. If the arc $j$ is not the non folded side of any folded side $\Delta'$, then we draw a based loop in the cross point $r_{l}$. This loop encloses only the puncture $p$, so that the cross point $r'_{1}$ of $\gamma_{2}$ with the arcs of $\tau^{\circ}$ belong to the relative interior of the arc $j$. Otherwise, if the arc $j$ is the non folded side of a self-folded triangle $\Delta'$, then we draw a based loop in the cross point $r_{l-2}$. This loop encloses only the puncture $p$ so that the segment $[r'_{0},r'_{2}]_{\gamma_{2}}$ is not contractible to the puncture $p$, with the homotopy that avoids $M$ and each of whose intermediate maps are segments with endpoints in the arcs of $\tau^{\circ}$ to which $r'_{0}$ and $r'_{2}$ belong.

Let's observe that if the arc $k$ is not a side of any self-folded triangle then the segments $[r_{l},r_{l+1}]_{\gamma_{1}}$ and $[r'_{0},r'_{1}]_{\gamma_{2}}$ are homotopics with the homotopy that avoids $M$ and each of whose intermediate maps are segments with endpoints in the arcs of $\tau^{\circ}$ to which $r'_{0}$,$r'_{2}$, $r_{l}$ and $r_{l+1}$ belong. Otherwise, if $k$ is the non folded side of a self-folded triangle $\Delta'$, then the segments $[r_{l-2},r_{l+1}]_{\gamma_{1}}$ and $[r'_{0},r'_{3}]_{\gamma_{2}}$ are homotopics with the homotopy that avoids $M$ and each of whose intermediate maps are segments with endpoints in the arcs of $\tau^{\circ}$ to which $r'_{0}$,$r'_{3}$, $r_{l}$, $r'_{l-2}$ belong. We denote by $\gamma_{3}$ the resulting curve of the next two steps:

\begin{itemize}
\item{}Joins the curves $\gamma_{1}$ and $\gamma_{2}$.
\item{}Identify the segments $[r_{l-2},r_{l+1}]_{\gamma_{1}}$ and $[r'_{0},r'_{3}]_{\gamma_{2}}$ in case of $k$ being not a side of any self-folded triangle. If $k$ is the non folded side of a self-folded triangle, then we identify $[r_{l-2},r_{l+1}]_{\gamma_{1}}$ and $[r'_{0},r'_{3}]_{\gamma_{2}}$. As is shown in Figure \ref{b121aa} and Figure \ref{b121bb}.

\end{itemize} 

\begin{figure}[H]
\begin{center}

\begin{tikzpicture}[scale=0.75]
\filldraw [black] (0,0) circle (2pt)
(0,2) circle (2pt)
(0,1) circle (1pt)
(-0.5,1) circle (1pt)
(0.5,1) circle (1pt)
(0,-0.5) circle (1pt)
(-0.5,-0.5) circle (1pt)
(0.5,-0.5) circle (1pt)
(-2.5,-0.25) circle (1pt)
(-2.5,0) circle (1pt)
(-2.5,0.25) circle (1pt)
(2.5,-0.25) circle (1pt)
(2.5,0) circle (1pt)
(2.5,0.25) circle (1pt);

\draw (0,2.5)node{$\Delta^{1}$}; 
\draw (0,4.5)node{$\Delta^{2}$}; 
\draw (0,3.75)node{$i$}; 
\draw (1,1.25)node{$j$}; 
\draw (-1,1.25)node{$k$}; 
\draw[line width=1pt] (0,0) ..controls(-4,0.5) and (-4,3.5) .. (0, 4); 
\draw[line width=1pt] (0,0) ..controls(4,0.5) and (4,3.5) .. (0, 4);

\draw[line width=1pt] (0,0) ..controls(-1.75,0.75) and (-1.75,1.25) .. (0, 2); 
\draw[line width=1pt] (0,0) ..controls(1.75,0.75) and (1.75,1.25) .. (0, 2);

\draw[line width=1pt] (0,0) ..controls (-2,0) and (-4,0.5) .. (-4, 3);
\draw[line width=1pt] (0,0) -- (-4, -1);

\draw[line width=1pt] (0,0) ..controls (2,0) and (4,0.5) .. (4, 3);
\draw[line width=1pt] (0,0) -- (4, -1);

\filldraw [magenta] (1,1.5) circle (2pt)
(-2.32,0.8) circle (2pt)
(-2.03,0.25) circle (2pt)
(-1.35,-0.32) circle (2pt)
(1.3,-0.35) circle (2pt)
(1.55,0.15) circle (2pt)
(1.53,0.45) circle (2pt)
(1.18,0.7) circle (2pt)
(-1.2,1.25) circle (2pt)
(0.35,1.85) circle (2pt);

\draw[color=magenta][line width=1.5pt] (0, 3.5) .. controls (1,3.5) and (1.5,2.75) .. (1,1.5);
\draw[color=magenta][line width=1.5pt] (0,-0.75) ..controls(-3,-0.5) and (-3.5,3) .. (0, 3.5);
\draw[color=magenta][line width=1.5pt] (0,-0.75) ..controls (2,-0.75) and (2.2,1.2) .. (0, 0.75);
\draw[color=magenta][line width=1.5pt] (0,0.75) ..controls(-1.5,0.25) and (-2,3) .. (0, 2.75); 
\draw[color=magenta][line width=1.5pt] (0, 2.75) ..controls (1,2.75) and (0.85,1.75) .. (0.35,1.85);

\filldraw [magenta] (-1.2,1.25) circle (2pt)
(1.87,0.55) circle (2pt)
(1.8,0.2) circle (2pt)
(1.51,-0.4) circle (2pt)
(-1.65,-0.4) circle (2pt)
(-1.61,0.15) circle (2pt)
(-1.5,0.35) circle (2pt)
(-1.2,0.7) circle (2pt);
\draw[color=magenta][line width=1.5pt] (0.35,1.85) ..controls (-2,0.5) and (-2.75,-1) .. (0.35,-1); 
\draw[color=magenta][line width=1.5pt] (0.35,-1) ..controls (2.5,-0.75) and (2.5,2.9) .. (0,3.1);
\draw[color=magenta][line width=1.5pt] (0,3.1) ..controls (-1.5,3.1) and (-2,1.5) .. (-1.2,1.25);

\end{tikzpicture}

\end{center}
\caption{The arc $j$ is not a side of any sefl-folded triangle. \label{b121aa}}

\end{figure}

\begin{figure}[H]
\begin{center}

\begin{tikzpicture}[scale=0.75]
\filldraw [black] (0,0) circle (2pt)

(-1.25,1.5) circle (2pt)
(1.25,1.5) circle (2pt)

(0,-0.5) circle (1pt)
(-0.5,-0.5) circle (1pt)
(0.5,-0.5) circle (1pt)
(-2.5,-0.25) circle (1pt)
(-2.5,0) circle (1pt)
(-2.5,0.25) circle (1pt)
(2.5,-0.25) circle (1pt)
(2.5,0) circle (1pt)
(2.5,0.25) circle (1pt);
\draw (0,1.75)node{$\Delta^{1}$}; 
\draw (0,4.5)node{$\Delta^{2}$}; 
\draw (0,3.75)node{$i$}; 
\draw (1,2.15)node{$j$}; 
\draw (-1,2.15)node{$k$};

\draw[line width=1pt] (0,0) ..controls(-4,0.5) and (-4,3.5) .. (0, 4); 
\draw[line width=1pt] (0,0) ..controls(4,0.5) and (4,3.5) .. (0, 4);

\draw[line width=1pt] (0,0) -- (-1.25, 1.5); 
\draw[line width=1pt] (0,0) -- (1.25, 1.5);

\draw[line width=1pt] (0,0) ..controls (0,0) and (-3,1.5) .. (-1.35, 2);
\draw[line width=1pt] (0,0) ..controls (0,0) and (0,2) .. (-1.35, 2);

\draw[line width=1pt] (0,0) ..controls (0,0) and (3,1.5) .. (1.35, 2);
\draw[line width=1pt] (0,0) ..controls (0,0) and (0,2) .. (1.35, 2);

\draw[line width=1pt] (0,0) ..controls (-2,0) and (-4,0.5) .. (-4, 3);
\draw[line width=1pt] (0,0) -- (-4, -1);

\draw[line width=1pt] (0,0) ..controls (2,0) and (4,0.5) .. (4, 3);
\draw[line width=1pt] (0,0) -- (4, -1); 
\filldraw [magenta] (0.85,1) circle (2pt)
(1.43,0.95) circle (2pt)
(-2.32,0.8) circle (2pt)
(-2.03,0.25) circle (2pt)
(-1.35,-0.32) circle (2pt)
(1.3,-0.35) circle (2pt)
(1.55,0.15) circle (2pt)
(1.53,0.45) circle (2pt)
(1.18,0.7) circle (2pt)
(0.65,0.83) circle (2pt)
(0.1,0.75) circle (2pt)
(-0.1,0.73) circle (2pt)
(-0.63,0.7) circle (2pt)
(-1.7,1.8) circle (2pt)
(1.8,1.4) circle (2pt)
(1.1,1.3) circle (2pt);
\draw[color=magenta][line width=1.5pt] (0, 3.5) .. controls (2.75,3.5) and (3.75,0.5) .. (0.85,1);
\draw[color=magenta][line width=1.5pt] (0,-0.75) ..controls(-3,-0.5) and (-3.5,3) .. (0, 3.5);
\draw[color=magenta][line width=1.5pt] (0,-0.75) ..controls (2,-0.75) and (2.2,1.2) .. (0, 0.75);
\draw[color=magenta][line width=1.5pt] (0,0.75) ..controls(-2,0.25) and (-2.5,3) .. (0, 2.75); 
\draw[color=magenta][line width=1.5pt] (0, 2.75) ..controls (2.25,2.75) and (3,1.25) .. (1.1,1.3);

\filldraw [magenta] (0.25,1.22) circle (2pt) 
(-0.25,1.13) circle (2pt)
(-0.77,0.93) circle (2pt)
(-1.17,0.7) circle (2pt)
(-1.45,0.4) circle (2pt)
(-1.6,0.18) circle (2pt)
(-1.55,-0.4) circle (2pt)
(2.08,-0.5) circle (2pt)
(3.02,0.75) circle (2pt)
(2.95,1.62) circle (2pt);

\draw[color=magenta][line width=1.5pt] (1.1,1.3) ..controls (-2,1.3) and (-2.75,-1) .. (0.35,-1); 
\draw[color=magenta][line width=1.5pt] (0.35,-1) ..controls (4,-0.75) and (4,3) .. (0,3.1);
\draw[color=magenta][line width=1.5pt] (0,3.1) ..controls(-0.15,3.1) and (-1.9,3) .. (-1.7,1.8) ;

\end{tikzpicture}

\end{center}
\caption{The arc $j$ is the non folded side of a self-folded triangle $\Delta'$. \label{b121bb}}
\end{figure}
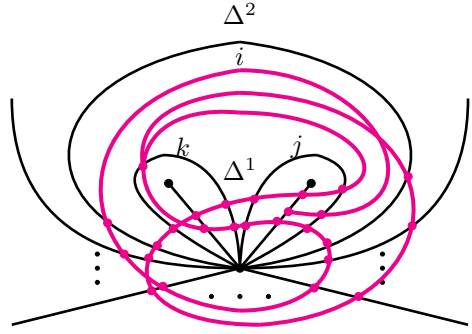
Now, if the arc $j$ is not the non folded side of a self-folded triangle, then we draw a curve  $\gamma_{4}$ from $r_{2}$ to $r'_{l'-1}$ so that the relative interior of $\gamma_{4}$ belongs to the interior of $\Delta_{2}$ and it is disjointed of $\gamma_{3}$. Otherwise, if the arc $j$ is the non folded side of a self-folded triangle $\Delta'$, then we draw a curve  $\gamma_{4}$ from $r_{3}$ to $r'_{l'-2}$ so that the relative interior of $\gamma_{4}$ belongs to the interior of $\Delta_{2}$ and is disjointed of $\gamma_{3}$.

The curve is defined by $i_{0}:=\gamma_{3}\cup \gamma_{4}$ and the curve $i':=i_{0}$, that is, the curve $i'$ is equal to the curve $i_{0}$. 

\end{myEnumerate}
\end{myEnumerate}

\item[$b_{2}$)] $\delta_{\tau}(p) = 0$.

\begin{myEnumerate}
\item[$b_{21}$)] The arc $i$ has a label "plain"at both ends.
We define $i_{0}:=i$ and it is assigned the label "plain" at both ends, as is shown in Figure \ref{i' b21)}.

\begin{figure}[H]
\centering 
$$ \begin{tikzpicture}[scale=0.8]

\filldraw [black] (-2,0) circle (2pt)
(-2,-2) circle (2pt)
(-2,1) circle (1pt)
(-2.25,1) circle (1pt)
(-1.75,1) circle (1pt);
\draw (-1.7,-0.25)node{$p$};

\draw[line width=1pt] (-2,-2) -- (-2, 0);
\draw[line width=1pt] (-2,-2) ..controls(-2,-2) and (-3.5,0.5) .. (-2, 0.5);
\draw[line width=1pt] (-2,-2) ..controls(-2,-2) and (-0.5,0.5) .. (-2, 0.5);

\filldraw [magenta] (-1.45,0.25) circle (2pt)
(-2.55,0.25) circle (2pt);

\draw[color=magenta][line width=1.5pt] (-2,0) ..controls(-3.5,0.5) and (-3.5,2) .. (-2, 2)
node[pos=0.5,left] {$i'$};
\draw[color=magenta][line width=1.5pt] (-2,0) ..controls(-0.5,0.5) and (-0.5,2) .. (-2, 2);

\end{tikzpicture}$$
\caption{Case $b_{21}$) of the Definition \ref{curva i'}.}\label{i' b21)}.
\end{figure}
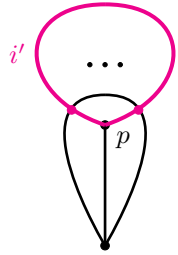 

\item[$b_{22}$)] The arc $i$ has a label $"\Bowtie"$ at both ends. We denote by $j$ the unique arc of $\tau^{\circ}$ which is incident at $p$ and we denote by $m$ the unique loop which encloses the puncture $p$.

According to this configuration we insert a point $p'$ in the arc $j$ such that $[p,p']_{j}\cap \tau^{\circ}=\emptyset$.

Otherwise, if we tour the arc $i$ from $p$ clockwise we denote by $x$ and $y$ the first and latest intersection points of $i$ with $m$ respectively. We define $i_{0}:= \gamma' \cup[x,y]_{i}\cup \gamma''$, where $\gamma'$ is the segment of $i$ which connect  $p$ with $x$ such that $y\notin \gamma'$ and $\gamma''$ is the segment which connects  $y$ with $p'$ such that  $x\notin \gamma''$. Let's observe that $[y,p']$ is homotopic to $[y,p]_{i}$ with the homotopy that avoids $M$ and each of whose intermediate maps are segments with endpoints in the arcs $j$ and $m$ respectively. It is assigned the label $\Bowtie$ at both ends of $i_{0}$, as is shown in Figure \ref{i' b22)}.

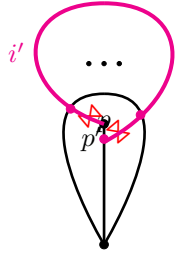
\begin{figure}[H]
\centering 
$$ \begin{tikzpicture}[scale=0.8]

\filldraw [black] (-2,0) circle (2pt)
(-2,-2) circle (2pt)
(-2,1) circle (1pt)
(-2.25,1) circle (1pt)
(-1.75,1) circle (1pt);
\draw (-2,0)node{$p$};
\draw (-2.2,-0.25)node{$p'$};

\draw[line width=1pt] (-2,-2) -- (-2, 0);
\draw[line width=1pt] (-2,-2) ..controls(-2,-2) and (-3.5,0.5) .. (-2, 0.5);
\draw[line width=1pt] (-2,-2) ..controls(-2,-2) and (-0.5,0.5) .. (-2, 0.5);

\filldraw [magenta] (-1.4,0.15) circle (2pt)
(-2,-0.25) circle (2pt)
(-2.55,0.25) circle (2pt);
\draw [color=red] (-1.8,-0.15)node {$\rotatebox{-50}{\textbf{\Bowtie}}$}; 
\draw [color=red] (-2.2,0.1)node {$\rotatebox{40}{\textbf{\Bowtie}}$}; 

\draw[color=magenta][line width=1.5pt] (-2,0) ..controls(-3.5,0.5) and (-3.5,2) .. (-2, 2)
node[pos=0.5,left] {$i'$};
\draw[color=magenta][line width=1.5pt] (-2,-0.3) ..controls(-0.5,0.5) and (-0.5,2) .. (-2, 2);

\end{tikzpicture}$$
\caption{Case $b_{22}$) of the Definition \ref{curva i'}.}\label{i' b22)}.
\end{figure}

\end{myEnumerate}
\end{myEnumerate}

\end{myEnumerate}

\end{definition}

The next remark (Remark \ref{ext i'}) is about the Definition \ref{curva i'} and it is useful to extend the construction $i'$ to any tagged triangulation.

\begin{remark}\label{ext i'}
Let $\tau$ be a tagged triangulation of a surface with marked points  $(\Sigma,M)$. Let's observe that from a tagged arc $i$, which does not belong to $\tau$, whose ends do not have negative signature, it is possible to construct the curve $i'$ followings the construction of the Definition \ref{curva i'}. 
\end{remark}

Now, let $\tau$ be a tagged triangulation of a surface with marked points  $(\Sigma,M)$. It follows from the Remark \ref{ext i'} that to construct the curve $i'$ in $(\Sigma,M)$ from a tagged arc $i$, which does not belong to $\tau$. It is enough to consider the cases when one end of $i$ (or both) has (have) negative signature. Then let's consider the three subcases mentioned, for each subcase, it is considered that if the label of the arc $i$ at the ends is "plain" or $\Bowtie$. We denote by $p$ and $q$ the ends of the tagged arc $i$, and we denote by $\tau'_{p}$ the tagged triangulation, which is different to $\tau$ only at the signature of the puncture $p$. That is, for every puncture which is different to $p$, the signature in $\tau$ and $\tau'_{p}$ is equal, otherwise, the signature of $p$ in $\tau$ is negative and the signature of $p$ in $\tau'_{p}$ is positive.

\begin{itemize}
\item[a)] The signature of $p$ is negative and the signature of $q$ is positive.
\begin{itemize}
\item[$a_{1}$)]The label of the arc $i$ at $p$ is "plain".

Let's denote by $i_{1}$ to the tagged arc which is obtained from $i$   changing the label "plain" at the end $p$  by the label $\Bowtie$. For the tagged arc $i_{1}$ it is possible to construct the curve $i'_{1}$ in $(\Sigma,M)$ with the triangulation $\tau'_{p}$. Under these conditions let's define the tagged curve $i'$ as $i'_{1}$. 

\item[$a_{2}$)]The label of the arc $i$ at $p$ is $\Bowtie$.

Let's denote by $i_{1}$ the tagged arc which is obtained  changing the label $\Bowtie$ at the end $p$  by the label "plain". For the tagged arc $i_{1}$ it is possible to construct the curve $i'_{1}$ in $(\Sigma,M)$ with the triangulation $\tau'_{p}$. Under these conditions let's define the tagged curve $i'$ as $i'_{1}$.

\end{itemize}

\item[b)] The signature of $q$ is negative and the signature of $p$ is positive.

\begin{itemize}
\item[$b_{1}$)]The label of the arc $i$ at $q$ is "plain".

Let's denote by $i_{1}$ the tagged arc which is obtained by changing  the label "plain" at the end $q$  by the label $\Bowtie$. For the tagged arc $i_{1}$ it is possible to construct the curve $i'_{1}$ in $(\Sigma,M)$ with the triangulation $\tau'_{q}$. Under these conditions let's define the tagged curve $i'$ as $i'_{1}$. 

\item[$b_{2}$)]The label of the arc $i$ at $q$ is $\Bowtie$.

Let's denote by $i_{1}$ the tagged arc which is obtained by changing the label $\Bowtie$ at the end $q$  by the label "plain". For the tagged arc $i_{1}$ it is possible to construct the curve $i'_{1}$ in $(\Sigma,,M)$ with the triangulation $\tau'_{q}$. Under these conditions let's define the tagged curve $i'$ as $i'_{1}$.

\end{itemize}

\item[c)] The signature of $p$ and $q$ is negative.

We denote by $\tau'_{pq}$ the triangulation that results from $\tau$ by changing the signature at the puncture $p$ and $q$  from negative to positive.

\begin{itemize}
\item[$c_{1}$)]The label of the arc $i$ at $p$ and $q$ is "plain".

Let's denote by $i_{1}$ the tagged arc which is obtained from $i$ by changing the label "plain" at the ends $p$ and $q$  by the label $\Bowtie$. For the tagged arc $i_{1}$ it is possible to construct the curve $i'_{1}$ in $(\Sigma,M)$ with the triangulation $\tau'_{pq}$. Under these conditions let's define the tagged curve $i'$ as $i'_{1}$. 

\item[$c_{2}$)]The label of the arc $i$ at $p$ is "plain" and the label at $q$ is $\Bowtie$.

Let's denote by $i_{1}$ the tagged arc which is obtained from $i$ by changing the labels "plain" and  $\Bowtie$  by the labels $\Bowtie$ and "plain" respectively. For the tagged arc $i_{1}$ it is possible to construct the curve $i'_{1}$ in $(\Sigma,M)$ with the triangulation $\tau'_{pq}$. Under these conditions let's define the tagged curve $i'$ as $i'_{1}$.

\item[$c_{3}$)]The label of the arc $i$ at $p$ is $\Bowtie$ and the label at $q$ is "plain".

Let's denote by $i_{1}$ the tagged arc which is obtained from $i$ by changing the labesl $\Bowtie$  and  "plain"  by the labels "plain" and $\Bowtie$ respectively. For the tagged arc $i_{1}$ it is possible to construct the curve $i'_{1}$ in $(\Sigma,M)$ with the triangulation $\tau'_{pq}$. Under this condition let's define the tagged curve $i'$ as $i'_{1}$.

\item[$c_{4}$)]The label of the arc $i$ at $p$ and $q$ is $\Bowtie$.

Let's denote by $i_{1}$ the tagged arc which is obtained from $i$ by changing the label $\Bowtie$ at the ends $p$ and $q$  by the label "plain". For the tagged arc $i_{1}$ it is possible to construct the curve $i'_{1}$ in $(\Sigma,M)$ with the triangulation $\tau'_{pq}$. Under these conditions let's define the tagged curve $i'$ as $i'_{1}$. 

\end{itemize}

\end{itemize}

\textbf{Examples of the derivative curve of an arc, Definition \ref{curva i'}.}

\begin{itemize}

\item[1)] Example 1.
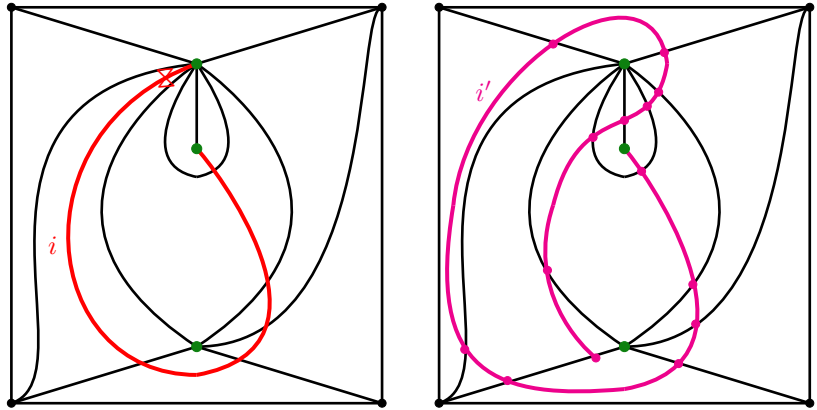
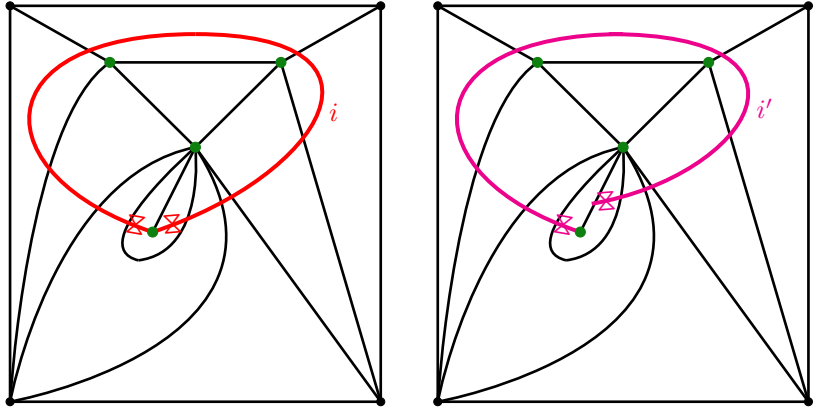
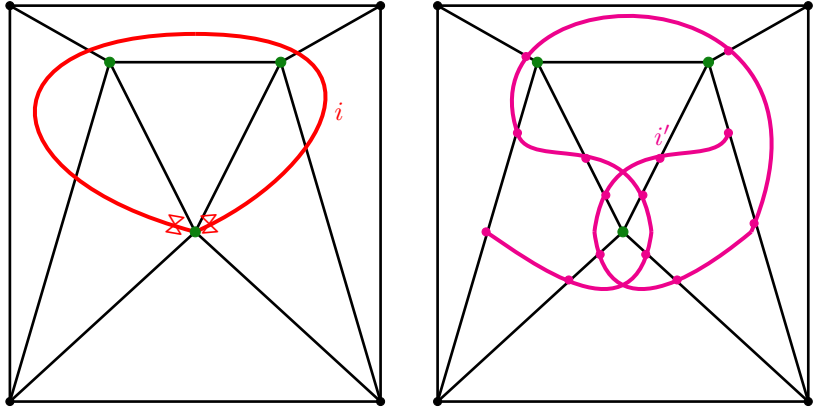
\begin{figure}[H]
\centering 
$$ 
$$
\end{figure} 

\end{itemize}

Using the Definition \ref{curva i'}, we are going to define the crossing points of the curve $i'$ with each arc of the triangulation $\tau^{\circ}$.

\begin{definition}(crossing point)\label{cruce}Let $\tau$ be a tagged triangulation of a surface with marked points $(\Sigma, M)$, $j$ an arc in $\tau^{\circ}$ and $i$ an arc which does not belong to $\tau^{\circ}$. Let's denote by $p$, $q$ the ends of the arc $j$, we say that $x\in i'\cap j$ is a crossing point of $i'$ with $j$ if it meets one of the three conditions described below. 

\begin{itemize}
\item[a)]The arc $j$ is not a side of any self-folded triangle.

$x\neq p$ and $x\neq q$.

\item[b)]The arc $j$ is the folded  side of a self-folded triangle $\Delta'$. Without loss of generallity we can suppose that $\delta_{\tau}(p)=0\neq \delta_{\tau}(q)$.

If the puncture $p$ is an end of $i'$ and the label of $i'$ at $p$ is $\Bowtie$ then $x\neq q$. Otherwise $q\neq x\neq p$.

\item[c)]The arc $j$ is the non folded side of a self-folded triangulo $\Delta'$, under this conditions the puncure $p$ is equal to $q$. Let's denote by $p'$ the unique puncture which it is  enclosed by the arc $j$ and we will denoted by $j'$ the unique arc in $\tau^{\circ}$ which is incident at $p'$.

If $t$ is an end of the curve $i'$ and it is labeled $\Bowtie$, then the intersection of the relative interior of $[x,t]_{i'}$ with $\tau^{\circ}$ is not empty. Otherwise $x$ belongs to the relative interior of the arc $j$.
\end{itemize}
\end{definition}

\begin{remark}\label{obcruce}
This remark is about the Definition \ref{cruce}.
\begin{myEnumerate}
\item[1)] Let $x$ and $x'$ be two crossing points as in $c)$ of the Definition \ref{cruce}. If the relative interior of the segment $[x,x']_{i'}$ cuts only the arc $j$ in one point, then we identify the crossing point $x$ with $x'$. We can denote by $x$ both crossing points and it will think  it as the same crossing point. 
\item[2)] If one end of the curve $i'$ is a puncture $p$ with signature zero and the label of $i'$ at the end $p$ is $\Bowtie$. Then the puncture $p$ will be considered as a point in the relative interior of the unique arc which is adjacent  to the puncture $p$.
\end{myEnumerate}

\end{remark}

With the help of the Definition \ref{cruce} and the Remark \ref{obcruce}, we are going to define the string representation denoted by $m(\tau,i)$.

Let $j$ be an arc in $\tau^{\circ}$ and $q_{j_{1}},...q_{j_{{\mathds A}(i',j)}}$ an enumeration of the different crossing points of the curve $i'$ with the arc $j$. Let's define the vector space associated to the arc $j$ in $m(\tau,i)$ as:

$$(m(\tau,i))_{j}:=K^{{\mathds A}(i',j)}=\displaystyle{\bigoplus_{r=1}^{{\mathds A}(i',j)}}Kq_{j,r},$$

where ${\mathds A}(i',j)$ denotes the number of different crossing points of the curve $i'$ with the arc $j$ and $K$ is a field.

Now we are going to define the linear tranformation associated to an arrow in $\widehat{Q}(\tau^{\circ})$. Given an arrow $\alpha: j\longrightarrow k$ in $\widehat{Q}(\tau^{\circ})$ there is only one puncture $p$ such that the arrow $\alpha$ is contractible to $p$ with the homotopy that avoids $M$ and each of whose intermediate maps are segments with endpoints in the arcs $j$ and $k$ respectively. Let $q_{j_{1}},...q_{j_{{\mathds A}(i',j)}}$ and $q_{k_{1}},...q_{k_{{\mathds A}(i',k)}}$ enumerating of the different crossing points of $i'$ with $j$ and $k$ respectively. We will denote by $K_{j,t}$ the copy of the field $K$ associated to the crossing point $q_{j,t}$. For $1 \leq s \leq {\mathds A}(i',j)$ and $1 \leq r \leq {\mathds A}(i',k)$ let's define $[(m(\tau,i'))_{\alpha}]_{r,s}: K_{j,s}\longrightarrow K_{k,r}$ as follows. We will consider all the possible cases for the arcs $j$ and $k$. The relative interior of the segment $[q_{j,s},q_{k,r}]_{i'}$ is disjointed of $\tau^{\circ}$ except when $j$ or $k$ (or both) is (are) the folded side(s) of a self-folded triangle $\Delta'$; in this case the segment $[q_{j,s},q_{k,r}]_{i'}$ intersects the non folded side of $\Delta'$. In addition we are going to distinguish the subcases when the segment $[q_{j,s},q_{k,r}]_{i'}$ is contractible to $p$ or if it is not.

\begin{itemize}

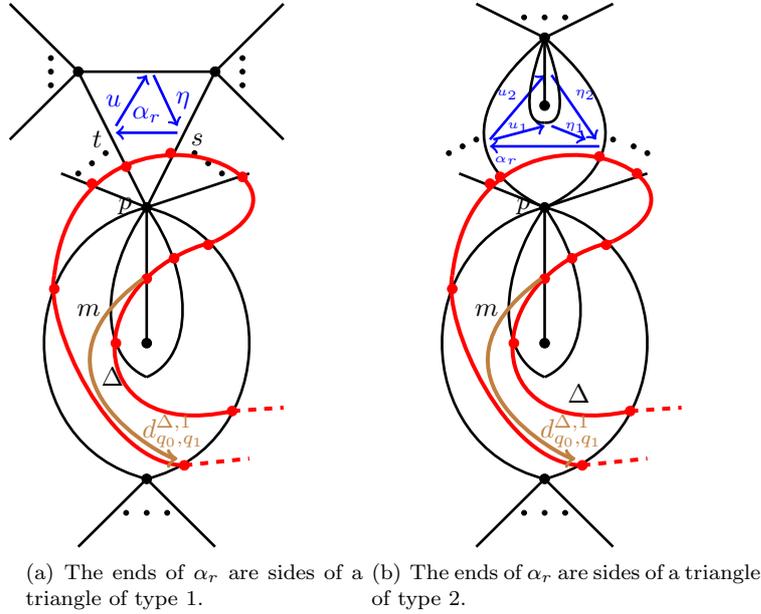
\item[1)]The arcs $j$ and $k$ are not  sides of any self-folded triangle.
\begin{itemize}
\item[a)] The segment $[q_{j,s},q_{k,r}]_{i'}$ is contractible to the puncture $p$.

$[(m(\tau,i'))_{\alpha}]_{r,s}=1$.

\item[b)] The segment $[q_{j,s},q_{k,r}]_{i'}$ is not contractible to the puncture $p$.

$[(m(\tau,i'))_{\alpha}]_{r,s}=0$.

\end{itemize}

\item[2)]The arc $j$ is not a side of any self-folded triangle and the arc $k$ is the folded side of a self-folded triangle.

\begin{itemize}
\item[a)] The segment $[q_{j,s},q_{k,r}]_{i'}$ is contractible to the puncture $p$.

$[(m(\tau,i'))_{\alpha}]_{r,s}=1$.

\item[b)] The segment $[q_{j,s},q_{k,r}]_{i'}$ is not contractible to the puncture $p$.

$[(m(\tau,i'))_{\alpha}]_{r,s}=0$.

\end{itemize}

\item[3)]The arc $j$ is not a side of any self-folded triangle and the arc $k$ is the non folded side of a self-folded triangle $\Delta'$, let's denote by $k'$ the folded side of $\Delta'$. If there are two crossing points denoted by $q_{k,r}$, then we are going to consider the copy of $q_{k,r}$ such that the segment $[q_{j,s},q_{k,r}]_{i'}$ does not intersect the folded side of $\Delta'$.

Under these conditions the segment $[q_{j,s},q_{k,r}]_{i'}$ is always contractible to the puncture $p$.

Then $[(m(\tau,i'))_{\alpha}]_{r,s}=-1$ if one of the two conditions below are met

\begin{itemize}
\item{}The curve $i'$ has one end at the opposite puncture to $p$ in the arc $k'$ (denoted by $x$) and the segment $[q_{k,r},x]_{i'}$ is disjointed to the triangulation $\tau^{\circ}$.

\item{}There is one crossing point $x$ of $i'$ with $k'$ such that the relative interior of the segment $[q_{k,r},x]_{i'}$ is disjointed of the triangulation $\tau^{\circ}$ and the segment $[q_{j,s},q_{k,r}]_{i'} \cup [q_{k,r},x]_{i'} $ is contractible to the puncture $p$.

\end{itemize}

On the onther hand $[(m(\tau,i'))_{\alpha}]_{r,s}=1$ if the following is met:

There is one crossing point $x$ of $i'$ with the arc $k$ such that the relative interior of the segment $[q_{k,r},x]_{i'}$ is disjointed of the triangulation $\tau^{\circ}$ and the segment $[q_{j,s},q_{k,r}]_{i'} \cup [q_{k,r},x]_{i'} $ is not contractible to $p$. 

and $[(m(\tau,i'))_{\alpha}]_{r,s}=0$ in any other case.

\item[4)]The arc $j$ is the folded side of a self-folded triangle and the arc $k$ is not a side of any self-folded triangle.

\begin{itemize}
\item[a)] The segment $[q_{j,s},q_{k,r}]_{i'}$ is contractible to the puncture $p$.

$[(m(\tau,i'))_{\alpha}]_{r,s}=1$.

\item[b)] The segment $[q_{j,s},q_{k,r}]_{i'}$ is not contractible to the puncture $p$.

$[(m(\tau,i'))_{\alpha}]_{r,s}=1$.
\end{itemize}

\item[5)]The arc $j$ is the non folded side of a self-folded triangle $\Delta'$ and the arc $k$ is not a side of any self-folded triangle, let's denote by $j'$ the folded side of the $\Delta'$. If there are two crossing points denoted by $q_{j,s}$ then we are going to choose the copy of $q_{j,s}$ such that the segment $[q_{j,s},q_{k,r}]_{i'}$ does not intersect the arc $j'$.

Under these conditions the segment $[q_{j,s},q_{k,r}]_{i'}$ is always contractible to the puncture $p$.

$[(m(\tau,i'))_{\alpha}]_{r,s}=1$ if one of the two conditions below are met

\begin{itemize}
\item{}The curve $i'$ has one end (denoted by $x$) at the opposite puncture to $p$ in the arc $j'$  and the segment $[x,q_{j,s}]_{i'}$ is disjointed to the triangulation $\tau^{\circ}$.

\item{}There is one crossing point $x$ of $i'$ with $j'$ such that the relative interior of the segment $[x,q_{j,s}]_{i'}$ is disjointed of the triangulation $\tau^{\circ}$ and the segment $[x,q_{j,s}]_{i'} \cup [q_{j,s},q_{k,r}]_{i'} $ is contractible to the puncture $p$.

\end{itemize} 

$[(m(\tau,i'))_{\alpha}]_{r,s}=0$ in any other case.

\item[6)] The arcs $j$ and $k$ are folded sides of different self-folded triangles $\Delta'$ and $\Delta''$, we denote by $p'$ and $p''$ the end of the arcs $j$ and $k$ with signature equals to zero respectively.

\begin{itemize}
\item[a)] The segment $[q_{j,s},q_{k,r}]_{i'}$ is contractible to the puncture $p$.

$[(m(\tau,i'))_{\alpha}]_{r,s}=1$.

\item[b)] The segment $[q_{j,s},q_{k,r}]_{i'}$ is not contractible to the puncture $p$.

$[(m(\tau,i'))_{\alpha}]_{r,s}=1$ if there is only one crossing point $q_{k,r'}$ of the curve $i'$ with the arc $k$ such that: 

\begin{itemize}
\item{} The segment $[q_{j,s},q_{k,r'}]_{i'}$ does not contain to the crossing point $q_{k,r}$ and it is contractible to the puncture $p$,
\item{}The curve $\gamma:=[q_{k,r'},q_{j,s}]_{i}\cup [q_{j,s},q_{k,r}]_{i'}$ surrounds the puncture $p'$ ,
\item{}$\gamma\cup [q_{k,r},q_{k,r'}]_{k}$ divides $\Sigma$ in two regions and one of those two regions is homeomorphic to a disk that contains one puncture namely $p'$. 
\end{itemize}

$[(m(\tau,i'))_{\alpha}]_{r,s}=0$ in any other case. 
\end{itemize}

\item[7)]The arc $j$ is the folded side of a self-folded triangle $\Delta'$ and the arc $k$ is the non folded side of a self-folded triangle $\Delta''$. We denote by $j'$ the non folded side of $\Delta'$, $k'$ to the folded side of $\Delta''$, $p'$ to the end of the arc $j$ with signature zero and $p''$ to the end of $k'$ with signature equal to zero.

If there are two crossing points denoted by $q_{k,r}$, then we are going to consider the copy of $q_{k,r}$ with the property that the segment $[q_{j,s},q_{k,r}]_{i'}$ does not intersect $k'$.

\begin{itemize}
\item[a)] The segment $[q_{j,s},q_{k,r}]_{i'}$ is contractible to the puncture $p$.

$[(m(\tau,i'))_{\alpha}]_{r,s}=1$ If there is one crossing point $q_{j,s'}$ of $i'$ with the arc $j$ such that $\gamma:= [q_{j,s},q_{k,r}]_{i'}\cup[q_{k,r},q_{j,s'}]_{i'}$ surrounds $p''$ and $\gamma \cup [q_{j,s'},q_{j,s}]_{j}$ divides $\Sigma$ in two regions. One of those two regions is homeomorphic to a disk that contains one puncture namely $p''$. 

Otherwise $[(m(\tau,i'))_{\alpha}]_{r,s}=-1$.
\item[b)] The segment $[q_{j,s},q_{k,r}]_{i'}$ is not contractible to puncture $p$.

$[(m(\tau,i'))_{\alpha}]_{r,s}=-1$ if one of the two conditions described below are met.

\begin{itemize}
\item{} The puncture $p''$ is an end of $i'$, the label of $i'$ at the end $p''$ is "plain" and there is one crossing point $q_{k',r'}$ of $i'$ with $k'$ such that:

\begin{itemize} 
\item{}The relative interior of the segment $[p'',q_{k,r}]_{i'}$ does not intersect any arc of $\tau^{\circ}$,
\item{} The union $[p'',q_{k,r}]_{i'}\cup [q_{k,r},q_{j,s}]_{i'}\cup [q_{j,s},q_{k',r'}]_{i'}\cup [q_{k',r'},p'']_{k'}$ divides $\Sigma$ in two regions and one of those two regions is homeomorphic to a monogon with puncture $p'$,

\item{}The marked point on the boundary of the monogon is puncture $p''$.
\end{itemize}
\item{} There are two crossing points of $i'$ with $k'$ denoted by $q_{k',r'}$ and $q_{k',r''}$  such that: 
\begin{itemize}
\item{}The relative interior of the segment $[q_{k',r'},q_{k,r}]$ is disjointed to the triangulation $\tau^{\circ}$,
\item{}The union $\gamma:=[q_{k',r'},q_{k,r}]_{i'}\cup[q_{k,r},q_{j,s}]_{i'}\cup[q_{j,s},q_{k',r''}]_{i'}$ surrounds the puncture $p'$,
\item{}$\gamma \cup [q_{k',r''},q_{k',r'}]_{k'}$ divides $\Sigma$ into two regions and one of those two regions is homeomorphic to a disk with the puncture $p'$,

\end{itemize} 

\end{itemize}

$[(m(\tau,i'))_{\alpha}]_{r,s}=1$ if one of the two conditions described below are met.

\begin{itemize}
\item{}There are two crossing points of $i'$ with $k'$ denoted by $q_{k',r'}$ and $q_{k',r''}$  such that: 
\begin{itemize}
\item{}The relative interior of the segment $[q_{k',r'},q_{k,r}]$ is disjointed to the triangulation $\tau^{\circ}$,
\item{}The union $\gamma:=[q_{k',r''},q_{j,s}]_{i'}\cup[q_{j,s},q_{k,r}]_{i'}\cup[q_{k,r},q_{k',r'}]_{i'}$ surrounds the puncture $p'$,
\item{}$\gamma \cup [q_{k',r'},q_{k',r''}]_{k'}$ divides $\Sigma$ in two regions and one of those two regions is homeomorphic to a disk with the punctures $p'$ and $p''$,

\end{itemize}

\item{} There is one crossing point of $i'$ with $j$ denoted by $q_{j,s'}$  such that $\gamma:=[q_{j,s},q_{k,r}]_{i'}\cup[q_{k,r},q_{j,s'}]_{i'}$ surrounds $p''$ and $\gamma \cup [q_{j,s'},q_{j,s}]_{i'}$ divides $\Sigma$ in two regions. One of those two regions is homeomorphic to a disk with the punctures $p'$ and $p''$,

$[(m(\tau,i'))_{\alpha}]_{r,s}=0$ in any other case. 

\end{itemize}

\end{itemize}

\item[8)]The arc $j$ is the non folded side of a self-folded triangle $\Delta'$ and the arc $k$ is the non folded side of a self-folded triangle $\Delta''$. We denote by $j'$ the folded side of $\Delta'$.

If there are two crossing points denoted by $q_{j,s}$ then we are going to consider the copy of $q_{j,s}$ with the property that the segment $[q_{j,s},q_{k,r}]_{i'}$ does not intersect $j'$.

\begin{itemize}
\item[a)] The segment $[q_{j,s},q_{k,r}]_{i'}$ is contractible to the puncture $p$.

Let $x$ be the crossing point of $i'$ with $j'$ with the property that $[x,q_{j,s}]_{i'}$ is disjointed to $\tau^{\circ}$.

$[(m(\tau,i'))_{\alpha}]_{r,s}=1$ if $[x,q_{j,s}]_{i'}\cup [q_{j,s},q_{k,r}]_{i'}$ is contractible to $p$ and $[(m(\tau,i'))_{\alpha}]_{r,s}=0$ in otherwise.

\item[b)] The segment $[q_{j,s},q_{k,r}]_{i'}$ is not contractible to the puncture $p$.

$[(m(\tau,i'))_{\alpha}]_{r,s}=0$.
\end{itemize}

\item[9)]The arcs $j$ and $k$ are non folded sides of self-folded triangles $\Delta'$ and $\Delta''$ respectively. We denote by $j'$ and $k'$ the folded sides of $\Delta'$ and $\Delta''$ respectively, in addition $p'$ and $p''$ will denote the ends of the arcs $j'$ and $k'$ with signature equal to zero respectively.

If there are two crossing points denoted by $q_{j,s}$ or $q_{k,r}$ (or both) in the arcs $j$ or $k$ (or both)  then we consider the copy(ies) of $q_{j,s}$ or $q_{k,r}$ (or both) such that the segment $[q_{j,s},q_{k,r}]_{i'} $ does not intersect the arc $j'$ neither $k'$.

Under these conditions the segment $[q_{j,s},q_{k,r}]_{i'}$ is always contractible to the puncture $p$.

$[(m(\tau,i'))_{\alpha}]_{r,s}=1$ if one of the two conditions described below are met. 
\begin{itemize}
\item{} There are crossing points of $i'$ with $j'$  denoted by $q_{j',s'}$ and $q_{j',s''}$ such that:
		\begin{itemize}
\item{}The relative interior of the segment $[q_{j',s'},q_{j,s}]_{i'}$ is disjointed to $\tau^{\circ}$,
\item{}$\gamma:= [q_{j',s'},q_{j,s}]_{i'}\cup[q_{j,s},q_{k,r}]_{i'}\cup[q_{k,r},q_{j',s''}]_{i'}$ surrounds the puncture $p''$,
\item{}$\gamma \cup [q_{j',s''},q_{j',s'}]_{j'}$ divides $\Sigma$ in two regions and one of those two regions is homeomorphic to a disk with the puncture $p''$.
		\end{itemize}

\item{}The puncture $p'$ is an end of the curve $i'$, the curve $i'$ is labeled "plain" at the end $p'$ and there is only one crossing point of $i'$ with $j'$ denoted by  $q_{j',s'}$  such that:

\begin{itemize}
\item{} The relative interior of the segment $[p',q_{j,s}]_{i'}$ is disjointed to $\tau^{\circ}$, 
\item{} The union $[p',q_{j,s}]_{i'}\cup[q_{j,s},q_{k,r}]_{i'}\cup[q_{k,r},q_{j',s'}]_{i'}\cup[q_{j',s'},p']_{j'}$ divides $\Sigma$ into two regions and one of those two regions is homeomorphic to a monogon with the punctures $p''$.
\item{} The marked point at the boundary is the puncture $p'$.
\end{itemize}
\end{itemize}

$[(m(\tau,i'))_{\alpha}]_{r,s}=-1$ if one of the four conditions described below are met.

\begin{itemize}

\item{} There are two crossing points $q_{k',r'}$, $q_{k',r''}$ in $k'$ and one crossing point $q_{j',s'}$ in $j'$ such that: 

		\begin{itemize}
\item{}The interior relative of the segments $[q_{j',s'},q_{j,s}]_{i'}$, $ [q_{k,r},q_{k',r'}]_{i'}$ are disjointed to $\tau^{\circ}$,
\item{} The curve $\gamma:=[q_{k',r'},q_{k,r}]_{i'}\cup[q_{k,r},q_{j,s}]_{i'}\cup[q_{j,s},q_{k',r''}]_{i'}$ surrounds the puncture $p'$,
\item{} $\gamma \cup [q_{k',r''},q_{k',r'}]_{k'}$ divides $\Sigma$ into two regions and one of those two regions is homeomorphic to a disk with the punctures $p'$.
		\end{itemize} 

\item{} There are two crossing points $q_{j',s'}$, $q_{j',s''}$ in $j'$ such that:

		\begin{itemize}
\item{}The interior relative of the segment $[q_{j',s'},q_{j,s}]_{i'}$ is disjointed to $\tau^{\circ}$,
\item{} The curve $\gamma:=[q_{j',s''},q_{k,r}]_{i'}\cup[q_{k,r},q_{j,s}]_{i'}\cup[q_{j,s},q_{j',s'}]_{i'}$ surrounds the puncture $p''$,
\item{} $\gamma \cup [q_{j',s'},q_{j',s''}]_{j'}$ divides $\Sigma$ into two regions and one of those two regions is homeomorphic to a disk with the puncture $p''$.
		\end{itemize}

\item{} The five conditions described below are met:

		\begin{itemize}
\item{}The puncture $p'$ is an end of $i'$,
\item{} The curve $i'$ has a label "plain" at the end $p'$,
\item{} There is one crossing point of $i'$ with $j'$ denoted by $q_{j',s'}$  such that the relative interior of the segment $[q_{j',s'},q_{j,s}]_{i'}$ is disjointed to $\tau^{\circ}$,
\item{} The union $[p',q_{k,r}]_{i'}\cup[q_{k,r},q_{j,s}]_{i'}\cup[q_{j,s},q_{j',s'}]_{i'}\cup[q_{j',s'},p']_{j'}$ divides $\Sigma$ into two regions and one of those two regions is homeomorphic to a monogon with the puncture $p''$,
\item{}The marked point on the boundary of the monogon is $p'$.
		\end{itemize} 

\item{} The five conditions described below are met:

\begin{itemize}
\item{}The puncture $p''$ is an end of the curve $i'$, 
\item{} The curve $i'$ has a label "plain" at the end $p''$,
\item{}There is one crossing point of $i'$ with $k'$ denoted by $q_{k',r'}$  such that the relative interior of the segment $[q_{k',r'},q_{k,r}]_{i'}$ is disjointed to $\tau^{\circ}$,
\item{} The union $[p'',q_{j,s}]_{i'}\cup[q_{j,s},q_{k,r}]_{i'}\cup[q_{k,r},q_{k',r'}]_{i'}\cup[q_{k',r'},p']_{k'}$ divides $\Sigma$ in two regions and one of those two regions is homeomorphic to a monogon with the puncture $p'$,
\item{} The marked point on the boundary of the monogon is $p''$.
\end{itemize}

\end{itemize}

$[(m(\tau,i'))_{\alpha}]_{r,s}=0$ in any other case. 

\end{itemize}

\textbf{EXAMPLES OF THE STRING REPRESENTATION $m(\tau,i)$}.
\begin{itemize}

\item[1)] Example 1.
\begin{figure}[H]
\centering 
$$ 
 }
\end{center}
\caption{Quiver.}

\end{figure}

\begin{figure}[H]
\begin{center}
\hspace*{0.6em}{%
}

\end{center}
\caption{String representation $m(\tau,i)$.}

\end{figure}

\item[2)] Example 2.
\begin{figure}[H]
\centering 
$$ %
 }

\end{center}
\caption{Quiver.}

\end{figure}

\begin{figure}[H]
\begin{center}

\hspace*{0.6em}{%
}
\end{center}
\caption{String representation $m(\tau,i)$.}

\end{figure}

\item[3)] Example 3.

\begin{figure}[H]
\centering 
$$ %
 }
\end{center}
\caption{Quiver.}

\end{figure}

\begin{figure}[H]
\begin{center}
\hspace*{0.6em}{%
}
\end{center}
\caption{String representation $m(\tau,i)$.}

\end{figure}

\item[4)] Example 4.

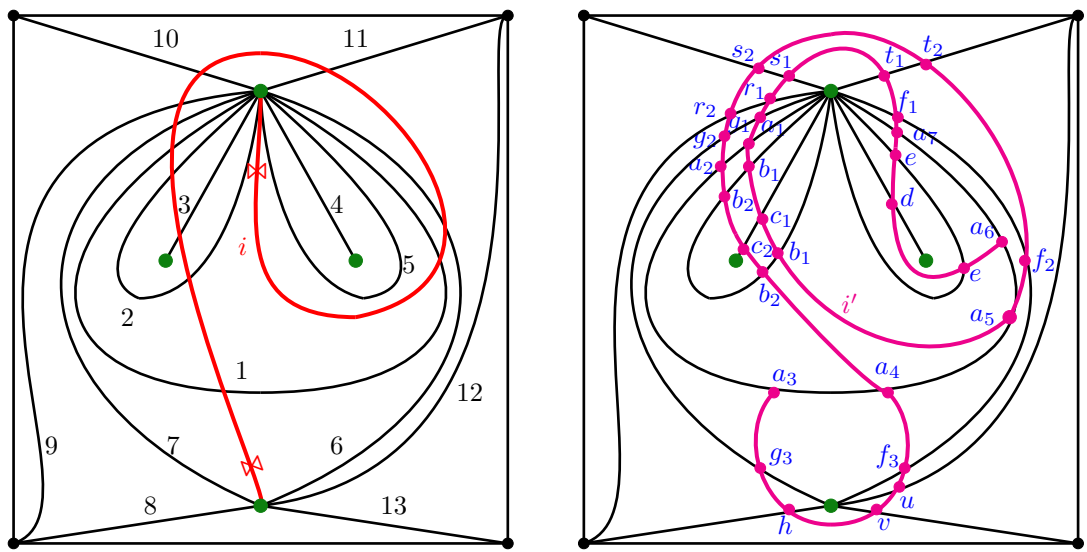
\begin{figure}[H]
\centering 
$$ %
}
\end{center}
\caption{String representation $m(\tau,i)$.}

\end{figure}

\item[5)] Example 5.
\begin{figure}[H]
\centering 
$$ %
}

\end{center}
\caption{String representation $m(\tau,i)$.}

\end{figure}

\item[6)] Example 6.
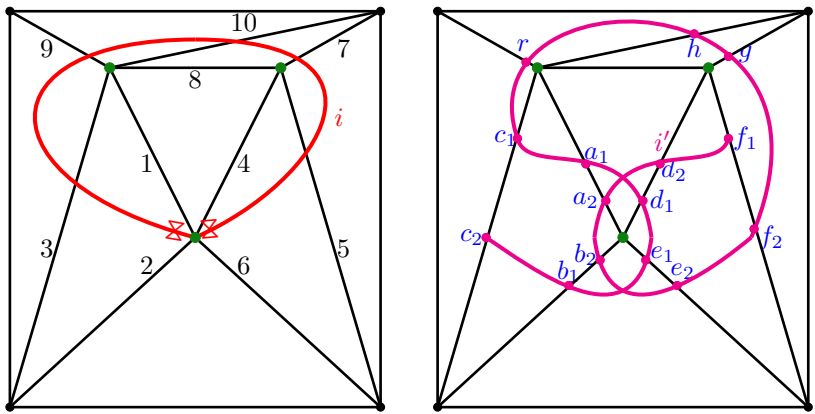
\begin{figure}[H]
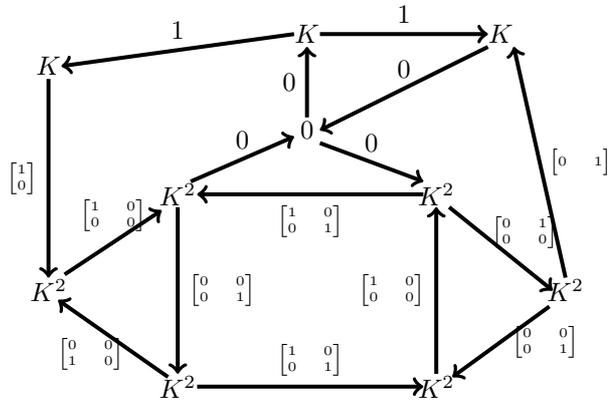

\centering 
$$ %
}

\end{center}
\caption{String representation $m(\tau,i)$.}

\end{figure}

\end{itemize}

An easy calculus verifies that the string representations do not satisfy the Jacobian relations, so in the next section, we are going to modify $m(\tau,i)$ for the Jacobian relations will meet.

\section{Arc representation $M(\tau,i)$.}

\setpapersize{A4} 
\setmargins{2.5cm} 
{1.5cm} 
{0.5cm} 
{23.42cm} 
{10pt} 
{1cm} 
{0pt} 
{2cm} 

\begin{definition}(Set ${\mathcal B}^{\Delta,1}_{i',j}$)\label{conjunto B})
Let $\tau$ be a tagged triangulation of a surface with marked points  $(\Sigma, M)$. Let's consider a self-folded triangle $\Delta$ in $\tau^{\circ}$, $j$ in $\tau^{\circ}$ and $i\notin \tau$ a tagged arc which minimizes the intersection points with the arcs of $\tau$. Using the curve $i'=i'(\tau,i)$ (Definition \ref{curva i'}) we are going to define the set ${\mathcal B}^{\Delta,1}_{i',j}$. The definition depends on whether the arc $j$ is the folded side of a self-folded triangle or if it is not:

\begin{itemize}
\item[a)]The arc $j$ is not the folded side of any self-folded triangle.

Let's consider the following local configuration,

\begin{itemize}
\item[i)]The arc $j$ is a side of the triangle $\Delta$,
\item[ii)]Let $p$ be a puncture with signature different from zero such that the puncture $p$ is opposite vertex to $j$ in $\Delta$,
\item[iii)]Let $q_{0}$,$q_{1}\in i'\cap j$ be two intersection points such that the segment $[q_{0},q_{1}]_{i'}$, surround the puncture $p$, 
\item[iv)]$l$ is the number of   intersection points of the relative interior of $[q_{0},q_{1}]_{i'}$ with the arcs of $\tau^{\circ}$,
\item[v)]$r_{s}$ with $s\in \{1,2,l\}$, is the $s^{th}$ intersection point of some arc of $\tau^{\circ}$ with the relative interior of $[q_{0},q_{1}]_{i'}$.
\end{itemize}

Every time we have the configuration described above $( i)-v) )$, we define the family of points $(q_{0},q_{1},r_{1}.r_{2},r_{l},p)$ as an element of the set ${\mathcal B}^{\Delta,1}_{i',j}$.

\begin{figure}[H]
\centering 
\begin{tikzpicture}[scale=0.75]

\filldraw [black] (-5,0) circle (2pt)
(-5,-4) circle (2pt)
(-8.5,-2) circle (2pt)
(-9.5,-2) circle (1pt)
(-9.5,-2.5) circle (1pt)
(-9.5,-1.5) circle (1pt);
\draw (-6,-2)node{$\Delta$};
\draw[line width=1pt] (-5, 0) -- (-5, -4)
node[pos=0.5,left] {$j$};
\draw[line width=1pt] (-5, 0) -- (-8.5, -2);
\draw[line width=1pt] (-5, -4) -- (-8.5, -2);
\draw[line width=1pt] (-8.5, -2) -- (-8.5,0);
\draw[line width=1pt] (-8.5, -2) -- (-10.5,-0.5);
\draw[line width=1pt] (-8.5, -2) -- (-10.5,-3.5);
\draw[line width=1pt] (-8.5, -2) -- (-8.5,-4)
node[pos=0.2,left] {$p$};

\filldraw [red] (-5,-1) circle (2pt)
(-5,-3) circle (2pt)
(-6.35,-3.2) circle (2pt)
(-8.5,-3.3) circle (2pt)
(-9.8,-3) circle (2pt)
(-9.8,-1) circle (2pt)
(-8.5,-0.45) circle (2pt)
(-6.05,-0.6) circle (2pt);

\draw[color=red][line width=1.5pt] (-5,-1) .. controls(-8, 0.5) and (-10.5,-1) .. (-10.5,-2)
node[pos=0,right] {$q_{0}$}
node[pos=0.11,above] {$r_{1}$}
node[pos=0.4,above] {$r_{2}$};
\draw[color=red][line width=1.5pt] (-5,-3) .. controls(-8, -3.5) and (-10.5,-3.5) .. (-10.5,-2)
node[pos=0,right] {$q_{1}$}
node[pos=0.15,below] {$r_{l}$};
\filldraw [black] (0,0) circle (2pt)
(-1.5,-2) circle (2pt)
(1.5,-2) circle (2pt)
(0,0.5) circle (1pt)
(0.3,0.5) circle (1pt)
(-0.3,0.5) circle (1pt);
\draw (0,-3)node{$\Delta$};
\draw[line width=1pt] (0,0) -- (-1, 1)
node[pos=0.3,right] {$p$};
\draw[line width=1pt] (0,0) -- (1, 1);
\draw[line width=1pt] (0,0) -- (-1.5,-2);
\draw[line width=1pt] (0,0) -- (1.5,-2);
\draw[line width=1pt] (0, 0) .. controls(-4,-0.1) and (-4,-3.9) .. (0, -4);
\draw[line width=1pt] (0, 0) .. controls(4,-0.1) and (4,-3.9) .. (0, -4)
node[pos=1,below] {$j$};
\draw[line width=1pt] (0, 0) .. controls(-1.5,-0.5) and (-2,-2) .. (-1.75, -2.25);
\draw[line width=1pt] (0, 0) .. controls(0.05,-0.5) and (-0.75,-2.5) .. (-1.75, -2.25)
node[pos=0.8,right] {$m$}
node[pos=0.8,left] {$\Delta'$};
\draw[line width=1pt] (0, 0) .. controls(0,-1) and (0.7,-2.5) .. (1.75, -2.25)
node[pos=0.8,left] {$m'$}
node[pos=0.5,right] {$\Delta''$};
\draw[line width=1pt] (0, 0) .. controls(0,0) and (2.25,-1.5) .. (1.75, -2.25);


\filldraw [red] (-1,-3.9) circle (2pt)
(-1.25,-2.15) circle (2pt)
(-1.3,-1.7) circle (2pt)
(-1.23,-0.8) circle (2pt)
(-1.05,-0.15) circle (2pt)
(-0.65,0.6) circle (2pt)
(0.65,0.6) circle (2pt)
(1.05,-0.15) circle (2pt)
(1.23,-0.95) circle (2pt)
(1.3,-1.7) circle (2pt)
(1.25,-2.2) circle (2pt)
(1,-3.9) circle (2pt);

\draw[color=red][line width=1.5pt] (-1,-3.9) .. controls(-2.5,2.5) and (2.5,2.5) .. (1,-3.9)
node[pos=1,below] {$q_{0}$}
node[pos=0,below] {$q_{1}$}
node[pos=0.08,left] {$r_{l}$}
node[pos=0.87,right] {$r_{2}$}
node[pos=0.92,right] {$r_{1}$};
\end{tikzpicture}
\end{figure}

\begin{figure}[H]
\centering 
$$ \begin{tikzpicture}[scale=0.75]

\filldraw [black] (0,0) circle (2pt)
(0,-4) circle (2pt)
(0,-2) circle (2pt)
(0,0.5) circle (1pt)
(0.3,0.5) circle (1pt)
(-0.3,0.5) circle (1pt)
(0,-4.5) circle (1pt)
(0.3,-4.5) circle (1pt)
(-0.3,-4.5) circle (1pt);
\draw (0,-3)node{$\Delta$};
\draw (4,-3.25)node{$\Delta$};
\draw[line width=1pt] (0,0) -- (-1, 1)
node[pos=0.3,right] {$p$};
\draw[line width=1pt] (0,0) -- (1, 1);
\draw[line width=1pt] (0,-4) -- (1, -5);
\draw[line width=1pt] (0,-4) -- (-1, -5);
\draw[line width=1pt] (0,0) -- (0, -2);
\draw[line width=1pt] (0,0) .. controls(-0.7,-1) and (-0.7,-2.2) .. (0, -2.5)
node[pos=0.5,left] {$m$};
\draw[line width=1pt] (0,0) .. controls(0.7,-1) and (0.7,-2.2) .. (0, -2.5);
\draw[line width=1pt] (0,0) .. controls(-2, -0.5) and (-2,-3.5) .. (0, -4)
node[pos=0.5,right] {$j$};
\draw[line width=1pt] (0,0) .. controls(2, -0.5) and (2,-3.5) .. (0, -4);

\filldraw [black] 
(4,0) circle (2pt)
(4,-4) circle (2pt)
(4,-2) circle (2pt)
(4,0.5) circle (1pt)
(4.3,0.5) circle (1pt)
(3.7,0.5) circle (1pt)
(4,-4.5) circle (1pt)
(4.3,-4.5) circle (1pt)
(3.7,-4.5) circle (1pt);
\draw[line width=1pt] (4,0) -- (3, 1)
node[pos=0.25,right] {$p$};
\draw[line width=1pt] (4,0) -- (5, 1);
\draw[line width=1pt] (4,-4) -- (5, -5);
\draw[line width=1pt] (4,-4) -- (3, -5);
\draw[line width=1pt] (4,0) -- (4, -2);
\draw[line width=1pt] (4,0) .. controls(2, -0.5) and (2,-3.5) .. (4, -4);
\draw[line width=1pt] (4,0) .. controls(6, -0.5) and (6,-3.5) .. (4, -4)
node[pos=0.5,left] {$j$};
\draw[line width=1pt] (4,0) .. controls(3.3,-1) and (3.3,-2.2) .. (4, -2.5)
node[pos=0.5,left] {$m$};
\draw[line width=1pt] (4,0) .. controls(4.7,-1) and (4.7,-2.2) .. (4, -2.5);



\filldraw [red] (4.55,-3.8) circle (2pt)
(5.25,-3) circle (2pt)
(3.55,-2) circle (2pt)
(4.4,-0.75) circle (2pt)
(4,-1.05) circle (2pt)
(4.9,-0.55) circle (2pt)
(4.8,0.75) circle (2pt)
(3.5,0.55) circle (2pt)
(2.65,-1.2) circle (2pt);
\draw[color=red][line width=1.5pt] (5.25,-3) .. controls(3,-3.5) and (3,-1) .. (5,-0.5)
node[pos=0,right] {$q_{0}$}
node[pos=0.5,right] {$r_{1}$}
node[pos=0.75,right] {$r_{2}$};
\draw[color=red][line width=1.5pt] (5,-0.5) .. controls(7,0.5) and (3,2) .. (2.65,-1)
node[pos=1,left] {$r_{l}$};
\draw[color=red][line width=1.5pt] (4.55,-3.8) .. controls(3.45,-3.8) and (2.5,-1.7) .. (2.65,-1)
node[pos=0,right] {$q_{1}$}; 

\filldraw [red] (-0.55,-3.8) circle (2pt)
(-1.25,-3) circle (2pt)
(1.3,-1.05) circle (2pt)
(0.7,0.7) circle (2pt)
(-0.75,0.75) circle (2pt)
(-0.9,-0.5) circle (2pt)
(-0.4,-0.75) circle (2pt)
(0,-1.05) circle (2pt)
(0,-1.05) circle (2pt)
(0.35,-2.1) circle (2pt);

\draw[color=red][line width=1.5pt] (-0.55,-3.8) .. controls(1,-3.5) and (1.5,-0.75) .. (1.25,-0.75)
node[pos=0,left] {$q_{0}$}
node[pos=0.9,right] {$r_{1}$};
\draw[color=red][line width=1.5pt] (1.3,-0.75) .. controls(1.25,2) and (-2.25,0.75) .. (-1,-0.5)
node[pos=0.25,right] {$r_{2}$};
\draw[color=red][line width=1.5pt] (-1.25,-3) .. controls(1,-3) and (0.75,-1) .. (-1,-0.5)
node[pos=0.3,right] {$r_{l}$}
node[pos=0,left] {$q_{1}$};

\end{tikzpicture}$$
\caption{Case $a$) of the Definition \ref{conjunto B}.}
\end{figure}

\item[b)]The arc $j$ is the folded side of a self-folded triangle $\Delta'$, let's denoted by $m$ to the non folded side of $\Delta'$.

\begin{itemize}
\item[i)] The arc $m$ is aside of $\Delta$,

\item[ii)]Let $p$ be a puncture with signature different from zero such that the puncture $p$ is opposite vertex to $m$ in $\Delta$,
\item[iii)] Let $q_{0}$,$q_{1}\in i'\cap j$ be two intersection points such that the segment $[q_{0},q_{1}]_{i'}$ surrounds the puncture $p$, 
\item[iv)]$l$ is the number  of intersection points of the relative interior of $[q_{0},q_{1}]_{i'}$ with the arcs of $\tau^{\circ}$,
\item[v)]$r_{s}$ with $s\in \{1,2,l\}$, is the $s^{th}$ intersection point of some arc of $\tau^{\circ}$ with the relative interior of $[q_{0},q_{1}]_{i'}ry$.

\item[vi)]Let $[q_{1},*]_{i}$ be the segment of $i'$ that follows $[q_{0},q_{1}]_{i'}$ and $q_{0}\notin [q_{1},*]_{i}$. With this notations $q_{2}$ is the first intersection point of $[q_{1},*]_{i'}$ with some arc in $\tau^{\circ} \setminus {m}$.
\end{itemize}

Every time we have the configuration described above $( i)-vi) )$, we define the family of points $(q_{0},q_{1},q_{2},r_{1}.r_{2},r_{3},r_{l},p)$ as an element of the set ${\mathcal B}^{\Delta,1}_{i',j}$.
\end{itemize}

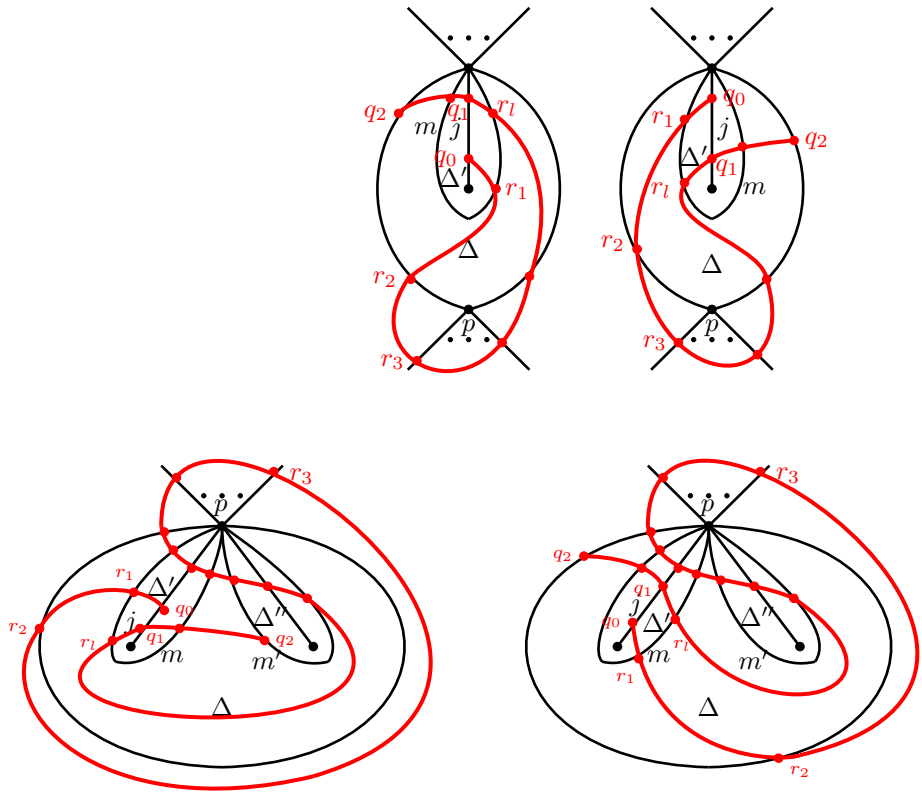
\begin{figure}[H]
\centering 
$$ \begin{tikzpicture}[scale=0.8]

\filldraw [black] (0,0) circle (2pt)
(0,-4) circle (2pt)
(0,-2) circle (2pt)
(0,0.5) circle (1pt)
(0.3,0.5) circle (1pt)
(-0.3,0.5) circle (1pt)
(0,-4.5) circle (1pt)
(0.3,-4.5) circle (1pt)
(-0.3,-4.5) circle (1pt);
\draw (0,-3)node{$\Delta$};
\draw (4,-3.25)node{$\Delta$};
\draw (3.7,-1.45)node{$\Delta'$};
\draw (-0.25,-1.8)node{$\Delta'$};
\draw (4.7,-2)node{$m$};
\draw (4.2,-1)node{$j$};
\draw (-0.7,-1)node{$m$};
\draw (-0.2,-1)node{$j$};
\draw (0,-4.3)node{$p$};
\draw[line width=1pt] (0,0) -- (-1, 1);
\draw[line width=1pt] (0,0) -- (1, 1);
\draw[line width=1pt] (0,-4) -- (1, -5);
\draw[line width=1pt] (0,-4) -- (-1, -5);
\draw[line width=1pt] (0,0) -- (0, -2);
\draw[line width=1pt] (0,0) .. controls(-0.7,-1) and (-0.7,-2.2) .. (0, -2.5);
\draw[line width=1pt] (0,0) .. controls(0.7,-1) and (0.7,-2.2) .. (0, -2.5);
\draw[line width=1pt] (0,0) .. controls(-2, -0.5) and (-2,-3.5) .. (0, -4);
\draw[line width=1pt] (0,0) .. controls(2, -0.5) and (2,-3.5) .. (0, -4);

\filldraw [black] 
(4,0) circle (2pt)
(4,-4) circle (2pt)
(4,-2) circle (2pt)
(4,0.5) circle (1pt)
(4.3,0.5) circle (1pt)
(3.7,0.5) circle (1pt)
(4,-4.5) circle (1pt)
(4.3,-4.5) circle (1pt)
(3.7,-4.5) circle (1pt);
\draw (4,-4.3)node{$p$};
\draw[line width=1pt] (4,0) -- (3, 1);

\draw[line width=1pt] (4,0) -- (5, 1);
\draw[line width=1pt] (4,-4) -- (5, -5);
\draw[line width=1pt] (4,-4) -- (3, -5);
\draw[line width=1pt] (4,0) -- (4, -2);
\draw[line width=1pt] (4,0) .. controls(2, -0.5) and (2,-3.5) .. (4, -4);
\draw[line width=1pt] (4,0) .. controls(6, -0.5) and (6,-3.5) .. (4, -4);
\draw[line width=1pt] (4,0) .. controls(3.3,-1) and (3.3,-2.2) .. (4, -2.5);
\draw[line width=1pt] (4,0) .. controls(4.7,-1) and (4.7,-2.2) .. (4, -2.5);



\filldraw [red] (4,-0.5) circle (2pt)
(3.55,-0.85) circle (2pt)
(2.77,-3) circle (2pt)
(3.45,-4.55) circle (2pt)
(4.75,-4.75) circle (2pt)
(4.9,-3.5) circle (2pt)
(3.55,-1.9) circle (2pt)
(4,-1.5) circle (2pt) 
(5.35,-1.2) circle (2pt)
(4.5,-1.3) circle (2pt);
\draw[color=red][line width=1.5pt] (4,-0.5) .. controls(3.25,-1) and (2.75,-2) .. (2.75,-2.9)
node[pos=0,right] {$q_{0}$}
node[pos=0.2,left] {$r_{1}$}
node[pos=1,left] {$r_{2}$};
\draw[color=red][line width=1.5pt] (2.75,-2.9) .. controls(3,-5.5) and (5.5,-5.5) .. (4.9,-3.5)
node[pos=0.3,left] {$r_{3}$};
\draw[color=red][line width=1.5pt] (4.9,-3.5) .. controls(4.8,-2.8) and (2.5,-2.5) .. (4,-1.5)
node[pos=0.8,left] {$r_{l}$}
node[pos=1.05,below] {$q_{1}$};
\draw[color=red][line width=1.5pt] (4,-1.5) ..controls(4.25,-1.3) and (5.25,-1.2) .. (5.35,-1.2)
node[pos=1,right] {$q_{2}$};

\filldraw [red] (-1.15,-0.75) circle (2pt)
(-0.3,-0.5) circle (2pt)
(0,-0.5) circle (2pt)
(0.4,-0.75) circle (2pt)
(1,-3.45) circle (2pt)
(0.55,-4.55) circle (2pt)
(-0.85,-4.85) circle (2pt)
(-0.95,-3.5) circle (2pt)
(0.45,-2) circle (2pt)
(0,-1.5) circle (2pt) ;

\draw[color=red][line width=1.5pt] (-1.15,-0.75) ..controls(-0.95,-0.5) and (-0.2,-0.4) .. (0,-0.5)
node[pos=0,left] {$q_{2}$};

\draw[color=red][line width=1.5pt] (0,-0.5) .. controls(1,-1) and (1.5,-2) .. (1,-3.5)
node[pos=-0.05,below] {$q_{1}$}
node[pos=0.1,right] {$r_{l}$};
\draw[color=red][line width=1.5pt] (1,-3.5) .. controls(0.5,-6) and (-2,-5) .. (-1,-3.5)
node[pos=0.6,left] {$r_{3}$}
node[pos=1,left] {$r_{2}$};
\draw[color=red][line width=1.5pt] (-1,-3.5) .. controls(-0.5,-3) and (1.25,-2.5) .. (0,-1.5)
node[pos=0.8,right] {$r_{1}$}
node[pos=1,left] {$q_{0}$};

\end{tikzpicture}$$
\end{figure}

\begin{figure}[H]
\centering 
\vspace{-8em}\hspace*{-6em}{\begin{tikzpicture}[scale=0.8]

\filldraw [black] (0,0) circle (2pt)
(-1.5,-2) circle (2pt)
(1.5,-2) circle (2pt)
(0,0.5) circle (1pt)
(0.3,0.5) circle (1pt)
(-0.3,0.5) circle (1pt);
\draw (0,-3)node{$\Delta$};
\draw (-1,-1)node{$\Delta'$};
\draw (-1.5,-1.6)node{$j$};
\draw (0.8,-1.5)node{$\Delta''$};
\draw[line width=1pt] (0,0) -- (-1, 1)
node[pos=0.3,right] {$p$};
\draw[line width=1pt] (0,0) -- (1, 1);
\draw[line width=1pt] (0,0) -- (-1.5,-2);
\draw[line width=1pt] (0,0) -- (1.5,-2);
\draw[line width=1pt] (0, 0) .. controls(-4,-0.1) and (-4,-3.9) .. (0, -4);
\draw[line width=1pt] (0, 0) .. controls(4,-0.1) and (4,-3.9) .. (0, -4);
\draw[line width=1pt] (0, 0) .. controls(-1.5,-0.5) and (-2,-2) .. (-1.75, -2.25);
\draw[line width=1pt] (0, 0) .. controls(0.05,-0.5) and (-0.75,-2.5) .. (-1.75, -2.25)
node[pos=0.8,right] {$m$};
\draw[line width=1pt] (0, 0) .. controls(0,-1) and (0.7,-2.5) .. (1.75, -2.25)
node[pos=0.8,left] {$m'$};
\draw[line width=1pt] (0, 0) .. controls(0,0) and (2.25,-1.5) .. (1.75, -2.25);


\filldraw [black] (8,0) circle (2pt)
(6.5,-2) circle (2pt)
(9.5,-2) circle (2pt)
(8,0.5) circle (1pt)
(8.3,0.5) circle (1pt)
(7.7,0.5) circle (1pt);
\draw (8,-3)node{$\Delta$};
\draw (7.15,-1.6)node{$\Delta'$};
\draw (6.8,-1.3)node{$j$};
\draw (8.8,-1.5)node{$\Delta''$};
\draw[line width=1pt] (8,0) -- (7, 1)
node[pos=0.3,right] {$p$};
\draw[line width=1pt] (8,0) -- (9, 1);
\draw[line width=1pt] (8,0) -- (6.5,-2);
\draw[line width=1pt] (8,0) -- (9.5,-2);
\draw[line width=1pt] (8, 0) .. controls(4,-0.1) and (4,-3.9) .. (8, -4);
\draw[line width=1pt] (8, 0) .. controls(12,-0.1) and (12,-3.9) .. (8, -4);
\draw[line width=1pt] (8, 0) .. controls(6.5,-0.5) and (6,-2) .. (6.25, -2.25);
\draw[line width=1pt] (8, 0) .. controls(8.05,-0.5) and (7.25,-2.5) .. (6.25, -2.25)
node[pos=0.8,right] {$m$};
\draw[line width=1pt] (8, 0) .. controls(8,-1) and (8.7,-2.5) .. (9.75, -2.25)
node[pos=0.8,left] {$m'$};
\draw[line width=1pt] (8, 0) .. controls(8,0) and (10.25,-1.5) .. (9.75, -2.25);


\filldraw [red] (-0.95,-1.4) circle (2pt)
(-1.45,-1.1) circle (2pt)
(-3,-1.7) circle (2pt)
(0.85,0.9) circle (2pt)
(-0.75,0.8) circle (2pt)
(-0.95,-0.1) circle (2pt)
(-0.8,-0.4) circle (2pt)
(-0.5,-0.7) circle (2pt)
(-0.2,-0.8) circle (2pt)
(0.2,-0.9) circle (2pt)
(0.75,-1) circle (2pt)
(1.4,-1.2) circle (2pt)
(-1.8,-1.9) circle (2pt)
(-1.35,-1.7) circle (2pt)
(-0.7,-1.7) circle (2pt)
(0.7,-1.9) circle (2pt);
\draw[color=red][line width=1.5pt] (-0.95,-1.4) ..controls(-1.1,-1) and (-2.5,-0.8) .. (-3,-1.7)
node[pos=-0,right] {{\scriptsize $q_{0}$}}
node[pos=0.4,above] {{\scriptsize $r_{1}$}}
node[pos=1,left] {{\scriptsize $r_{2}$}}; 
\draw[color=red][line width=1.5pt] (-3,-1.7) .. controls(-4,-3) and (-2,-5) .. (1.5,-4.15);
\draw[color=red][line width=1.5pt] (1.5,-4.15) .. controls(7,-2.5) and (-1,3.5) .. (-1,0)
node[pos=0.66,right] {$r_{3}$};
\draw[color=red][line width=1.5pt] (-1,0) .. controls(-0.8,-1) and (0.5,-0.8) .. (1.4,-1.2);
\draw[color=red][line width=1.5pt] (1.4,-1.2) .. controls(5,-4) and (-5.5,-3.5) .. (-1.35,-1.7)
node[pos=0.95,left] {{\scriptsize $r_{l}$}}
node[pos=1.02,below] {{\scriptsize $q_{1}$}};

\draw[color=red][line width=1.5pt] (-1.35,-1.7) .. controls(-1.25,-1.6) and (0.6,-1.8) .. (0.7,-1.9)
node[pos=1,right] {{\scriptsize $q_{2}$}};

\filldraw [red] (7.25,-1) circle (2pt)
(9.15,-3.85) circle (2pt)
(8.85,0.9) circle (2pt)
(7.25,0.8) circle (2pt)
(7.05,-0.1) circle (2pt)
(7.2,-0.4) circle (2pt)
(7.5,-0.7) circle (2pt)
(7.8,-0.8) circle (2pt)
(8.2,-0.9) circle (2pt)
(8.75,-1) circle (2pt)
(9.4,-1.2) circle (2pt)
(7.45,-1.55) circle (2pt)
(6.75,-1.6) circle (2pt)
(6.85,-2.2) circle (2pt)
(6.9,-0.7) circle (2pt)
(5.95,-0.5) circle (2pt);
\draw[color=red][line width=1.5pt] (6.75,-1.6) .. controls(7,-4) and (9.3,-3.9) .. (9.5,-3.8)
node[pos=-0,left] {{\scriptsize $q_{0}$}}
node[pos=0.15,left] {{\scriptsize $r_{1}$}}
node[pos=1,below] {{\scriptsize $r_{2}$}};
\draw[color=red][line width=1.5pt] (9.5,-3.8) .. controls(15,-2.5) and (7,3.5) .. (7,0) 
node[pos=0.66,right] {$r_{3}$};
\draw[color=red][line width=1.5pt] (7,0) .. controls(7.2,-1) and (8.5,-0.8) .. (9.4,-1.2);
\draw[color=red][line width=1.5pt] (9.4,-1.2) .. controls(12,-3) and (8,-3.7) .. (7.25,-1) 
node[pos=0.9,below] {{\scriptsize $r_{l}$}}
node[pos=0.99,left] {{\scriptsize $q_{1}$}};
\draw[color=red][line width=1.5pt] (7.25,-1) ..controls(7,-0.5) and (5.95,-0.5) .. (5.95,-0.5)
node[pos=1,left] {{\scriptsize $q_{2}$}};

\end{tikzpicture}}

\end{figure}

\begin{figure}[H]
\centering 
\vspace{-8em}\hspace*{-4em}{\begin{tikzpicture}[scale=0.8]

\filldraw [black] (0,0) circle (2pt)
(-1.5,-2) circle (2pt)
(1.5,-2) circle (2pt)
(0,0.5) circle (1pt)
(0.3,0.5) circle (1pt)
(-0.3,0.5) circle (1pt);
\draw (0,-3)node{$\Delta$};
\draw (-1,-1)node{$\Delta'$};
\draw (-1.5,-1.6)node{$j$};
\draw (0.8,-1.5)node{$\Delta''$};
\draw[line width=1pt] (0,0) -- (-1, 1)
node[pos=0.3,right] {$p$};
\draw[line width=1pt] (0,0) -- (1, 1);
\draw[line width=1pt] (0,0) -- (-1.5,-2);
\draw[line width=1pt] (0,0) -- (1.5,-2);
\draw[line width=1pt] (0, 0) .. controls(-4,-0.1) and (-4,-3.9) .. (0, -4);
\draw[line width=1pt] (0, 0) .. controls(4,-0.1) and (4,-3.9) .. (0, -4);
\draw[line width=1pt] (0, 0) .. controls(-1.5,-0.5) and (-2,-2) .. (-1.75, -2.25);
\draw[line width=1pt] (0, 0) .. controls(0.05,-0.5) and (-0.75,-2.5) .. (-1.75, -2.25)
node[pos=0.8,right] {$m$};
\draw[line width=1pt] (0, 0) .. controls(0,-1) and (0.7,-2.5) .. (1.75, -2.25)
node[pos=0.8,left] {$m'$};
\draw[line width=1pt] (0, 0) .. controls(0,0) and (2.25,-1.5) .. (1.75, -2.25);


\filldraw [black] (8,0) circle (2pt)
(6.5,-2) circle (2pt)
(9.5,-2) circle (2pt)
(8,0.5) circle (1pt)
(8.3,0.5) circle (1pt)
(7.7,0.5) circle (1pt);
\draw (8,-3)node{$\Delta$};
\draw (7.15,-1.6)node{$\Delta'$};
\draw (6.8,-1.3)node{$j$};
\draw (8.8,-1.5)node{$\Delta''$};
\draw[line width=1pt] (8,0) -- (7, 1)
node[pos=0.3,right] {$p$};
\draw[line width=1pt] (8,0) -- (9, 1);
\draw[line width=1pt] (8,0) -- (6.5,-2);
\draw[line width=1pt] (8,0) -- (9.5,-2);
\draw[line width=1pt] (8, 0) .. controls(4,-0.1) and (4,-3.9) .. (8, -4);
\draw[line width=1pt] (8, 0) .. controls(12,-0.1) and (12,-3.9) .. (8, -4);
\draw[line width=1pt] (8, 0) .. controls(6.5,-0.5) and (6,-2) .. (6.25, -2.25);
\draw[line width=1pt] (8, 0) .. controls(8.05,-0.5) and (7.25,-2.5) .. (6.25, -2.25)
node[pos=0.8,right] {$m$};
\draw[line width=1pt] (8, 0) .. controls(8,-1) and (8.7,-2.5) .. (9.75, -2.25)
node[pos=0.8,left] {$m'$};
\draw[line width=1pt] (8, 0) .. controls(8,0) and (10.25,-1.5) .. (9.75, -2.25);

\filldraw [red] (1.3,-1.8) circle (2pt)
(0.5,-1.75) circle (2pt)
(-1.6,-1.35) circle (2pt)
(-1,-1.25) circle (2pt)
(-0.3,-1.05) circle (2pt)
(0.15,-0.85) circle (2pt)
(0.45,-0.6) circle (2pt)
(0.6,-0.4) circle (2pt)
(0.7,-0.05) circle (2pt)
(0.6,0.6) circle (2pt)
(-0.8,0.8) circle (2pt)
(1.5,-1.25) circle (2pt)
(2.75,-1.2) circle (2pt)
(0.75,-1) circle (2pt);
(0.8,-1) 
\draw[color=red][line width=1.5pt] (1.3,-1.8) .. controls(1.2,-0.5) and (3.3,-1).. (3.3,-2)
node[pos=0,left] {$q_{1}$}
node[pos=0.3,below] {$r_{l}$}; 
\draw[color=red][line width=1.5pt] (-3.5,-2) ..controls(-3.6,-5) and (3.05,-5) .. (3.3,-2);
\draw[color=red][line width=1.5pt] (-3.5,-2) ..controls(-3.1,1) and (0.2,1.2) .. (0.6,0.6);

\draw[color=red][line width=1.5pt] (-1.75,-1.35) ..controls(-1.15,-1.35) and (1.3,-1) .. (0.6,0.6)
node[pos=0.3,below] {$r_{3}$}; 
\draw[color=red][line width=1.5pt] (-1.65,-1.35) ..controls (-3.5,-2) and (-1.2,-2.9) .. (-1,-2.75)
node[pos=0,above] {$r_{2}$};
\draw[color=red][line width=1.5pt] (-1,-2.75).. controls (0.45,-2.7) and (0.5,-1.5) .. (0.8,-1)
node[pos=1,right] {$q_{0}$}
node[pos=0.6,left] {$r_{1}$};
(1.3,-1.8)

\filldraw [red] (7.25,-1) circle (2pt)
(6.55,-1) circle (2pt)
(9.15,-3.85) circle (2pt)
(8.85,0.9) circle (2pt)
(7.25,0.8) circle (2pt)
(7.05,-0.1) circle (2pt)
(7.2,-0.4) circle (2pt)
(7.5,-0.7) circle (2pt)
(7.8,-0.8) circle (2pt)
(8.2,-0.9) circle (2pt)
(8.75,-1) circle (2pt)
(9.4,-1.2) circle (2pt)
(7.05,-2) circle (2pt)
(6.75,-1.6) circle (2pt);
\draw[color=red][line width=1.5pt] (7.25,-1) .. controls(3,-1) and (7.3,-3.9) .. (9.5,-3.8)
node[pos=-0,right] {{\scriptsize $q_{0}$}}
node[pos=0.06,above] {{\scriptsize $r_{1}$}}
node[pos=1,below] {{\scriptsize $r_{2}$}};
\draw[color=red][line width=1.5pt] (9.5,-3.8) .. controls(15,-2.5) and (7,3.5) .. (7,0) 
node[pos=0.66,right] {$r_{3}$}
node[pos=0.9,left] {{\scriptsize $r_{l-8}$}}
node[pos=1,left] {{\scriptsize$r_{l-7}$}};
\draw[color=red][line width=1.5pt] (7,0) .. controls(7.2,-1) and (8.5,-0.8) .. (9.4,-1.2)
node[pos=1.01,right] {{\scriptsize $r_{l-1}$}};
\draw[color=red][line width=1.5pt] (9.4,-1.2) .. controls(12,-3) and (8,-3.7) .. (6.75,-1.6) 
node[pos=0.9,below] {{\scriptsize $r_{l}$}}
node[pos=0.99,left] {{\scriptsize $q_{1}$}};

\end{tikzpicture}} 

\end{figure} 
\begin{figure}[H]
\centering 
\vspace{-8em}\hspace*{-4em}{\begin{tikzpicture}[scale=0.8]

\filldraw [black] (0,0) circle (2pt)
(-1.5,-2) circle (2pt)
(1.5,-2) circle (2pt)
(0,0.5) circle (1pt)
(0.3,0.5) circle (1pt)
(-0.3,0.5) circle (1pt);
\draw (0,-3)node{$\Delta$};
\draw (-1,-1)node{$\Delta'$};
\draw (-1.5,-1.6)node{$j$};
\draw (0.8,-1.5)node{$\Delta''$};
\draw[line width=1pt] (0,0) -- (-1, 1)
node[pos=0.3,right] {$p$};
\draw[line width=1pt] (0,0) -- (1, 1);
\draw[line width=1pt] (0,0) -- (-1.5,-2);
\draw[line width=1pt] (0,0) -- (1.5,-2);
\draw[line width=1pt] (0, 0) .. controls(-4,-0.1) and (-4,-3.9) .. (0, -4);
\draw[line width=1pt] (0, 0) .. controls(4,-0.1) and (4,-3.9) .. (0, -4);
\draw[line width=1pt] (0, 0) .. controls(-1.5,-0.5) and (-2,-2) .. (-1.75, -2.25);
\draw[line width=1pt] (0, 0) .. controls(0.05,-0.5) and (-0.75,-2.5) .. (-1.75, -2.25)
node[pos=0.8,right] {$m$};
\draw[line width=1pt] (0, 0) .. controls(0,-1) and (0.7,-2.5) .. (1.75, -2.25)
node[pos=0.8,left] {$m'$};
\draw[line width=1pt] (0, 0) .. controls(0,0) and (2.25,-1.5) .. (1.75, -2.25);


\filldraw [black] (8,0) circle (2pt)
(6.5,-2) circle (2pt)
(9.5,-2) circle (2pt)
(8,0.5) circle (1pt)
(8.3,0.5) circle (1pt)
(7.7,0.5) circle (1pt);
\draw (8,-3)node{$\Delta$};
\draw (7.15,-1.6)node{$\Delta'$};
\draw (6.8,-1.3)node{$j$};
\draw (8.8,-1.5)node{$\Delta''$};
\draw[line width=1pt] (8,0) -- (7, 1)
node[pos=0.3,right] {$p$};
\draw[line width=1pt] (8,0) -- (9, 1);
\draw[line width=1pt] (8,0) -- (6.5,-2);
\draw[line width=1pt] (8,0) -- (9.5,-2);
\draw[line width=1pt] (8, 0) .. controls(4,-0.1) and (4,-3.9) .. (8, -4);
\draw[line width=1pt] (8, 0) .. controls(12,-0.1) and (12,-3.9) .. (8, -4);
\draw[line width=1pt] (8, 0) .. controls(6.5,-0.5) and (6,-2) .. (6.25, -2.25);
\draw[line width=1pt] (8, 0) .. controls(8.05,-0.5) and (7.25,-2.5) .. (6.25, -2.25)
node[pos=0.8,right] {$m$};
\draw[line width=1pt] (8, 0) .. controls(8,-1) and (8.7,-2.5) .. (9.75, -2.25)
node[pos=0.8,left] {$m'$};
\draw[line width=1pt] (8, 0) .. controls(8,0) and (10.25,-1.5) .. (9.75, -2.25);

\filldraw [red] (1.3,-1.8) circle (2pt)
(0.7,-1.85) circle (2pt)
(-1.6,-1.35) circle (2pt)
(-1,-1.25) circle (2pt)
(-0.3,-1.05) circle (2pt)
(0.15,-0.85) circle (2pt)
(0.45,-0.6) circle (2pt)
(0.6,-0.4) circle (2pt)
(0.7,-0.05) circle (2pt)
(0.6,0.6) circle (2pt)
(-0.8,0.8) circle (2pt)
(1.1,-0.9) circle (2pt)
(2.75,-1.2) circle (2pt)
(0.75,-1) circle (2pt);
\draw[color=red][line width=1.5pt] (0.8,-1) .. controls(1.2,-0.5) and (3.3,-1).. (3.3,-2)
node[pos=0,left] {$q_{1}$}
node[pos=0.3,below] {$r_{l}$}; 
\draw[color=red][line width=1.5pt] (-3.5,-2) ..controls(-3.6,-5) and (3.05,-5) .. (3.3,-2);
\draw[color=red][line width=1.5pt] (-3.5,-2) ..controls(-3.1,1) and (0.2,1.2) .. (0.6,0.6);

\draw[color=red][line width=1.5pt] (-1.75,-1.35) ..controls(-1.15,-1.35) and (1.3,-1) .. (0.6,0.6)
node[pos=0.3,below] {$r_{3}$}; 
\draw[color=red][line width=1.5pt] (-1.65,-1.35) ..controls (-3.5,-2) and (-1.2,-2.9) .. (-1,-2.75)
node[pos=0,above] {$r_{2}$};
\draw[color=red][line width=1.5pt] (-1,-2.75).. controls (0.45,-2.7) and (0.5,-1.5) .. (1.3,-1.8)
node[pos=1,right] {$q_{0}$}
node[pos=0.6,left] {$r_{1}$};

\filldraw [red] (7.25,-1) circle (2pt)
(6.3,-1.55) circle (2pt)
(9.15,-3.85) circle (2pt)
(8.85,0.9) circle (2pt)
(7.25,0.8) circle (2pt)
(7.05,-0.1) circle (2pt)
(7.2,-0.4) circle (2pt)
(7.5,-0.7) circle (2pt)
(7.8,-0.8) circle (2pt)
(8.2,-0.9) circle (2pt)
(8.75,-1) circle (2pt)
(9.4,-1.2) circle (2pt)
(7.45,-1.5) circle (2pt)
(6.75,-1.6) circle (2pt);
\draw[color=red][line width=1.5pt] (6.75,-1.6) .. controls(3,-1) and (7.3,-3.9) .. (9.5,-3.8)
node[pos=-0,right] {{\scriptsize $q_{0}$}}
node[pos=0.06,above] {{\scriptsize $r_{1}$}}
node[pos=1,below] {{\scriptsize $r_{2}$}};
\draw[color=red][line width=1.5pt] (9.5,-3.8) .. controls(15,-2.5) and (7,3.5) .. (7,0) 
node[pos=0.66,right] {$r_{3}$}
node[pos=0.9,left] {{\scriptsize $r_{l-8}$}}
node[pos=1,left] {{\scriptsize$r_{l-7}$}};
\draw[color=red][line width=1.5pt] (7,0) .. controls(7.2,-1) and (8.5,-0.8) .. (9.4,-1.2)
node[pos=1.01,right] {{\scriptsize $r_{l-1}$}};
\draw[color=red][line width=1.5pt] (9.4,-1.2) .. controls(12,-3) and (8,-3.7) .. (7.25,-1)
node[pos=0.9,below] {{\scriptsize $r_{l}$}}
node[pos=0.99,left] {{\scriptsize $q_{1}$}};

\end{tikzpicture}} 
\caption{Case $b)$ of the Definition \ref{conjunto B}.}

\end{figure}

\end{definition}

The following Definition will be useful to define the detour matrix and the auxiliary matrix to modify the string representation $m(\tau,i)$ to obtain the required representation.

\begin{definition}\label{1-desviacion}(1-detour)For each element $(q_{0},q_{1},r_{1},r_{2},r_{l},p)$ or $(q_{0},q_{1},q_{2},r_{1},r_{2},r_{3},r_{l},p)$ in ${\mathcal B}^{\Delta,1}_{i',j}$ we are going to construct the 1-\textit{detour} $d_{(q_{0},q_{1})}^{\Delta,1}$ which is an oriented simple curve with ends $s(d_{(q_{0},q_{1})}^{\Delta,1})$ and $t(d_{(q_{0},q_{1})}^{\Delta,1})$ start and finish respectively. Only in case $b$) we will draw an auxiliary curve (1-auxiliar) $e_{(q_{0},q_{1})}^{\Delta,1}$ with ends $s(e_{(q_{0},q_{1})}^{\Delta,1})$ and $t(e_{(q_{0},q_{1})}^{\Delta,1})$. We are to consider the cases when the arc $j$ is the folded side of a self-folded triangle or if it is not.

\begin{itemize}
\item[a)] The arc $j$ is an arc in $\tau^{\circ}$ which is not the folded side of any self-folded triangle $\Delta'$. If $r_{1}$ does not belong to any side of any self-folded triangle then we draw the 1-\textit{detour} $d_{(q_{0},q_{1})}^{\Delta,1}$ with ends $s(d_{(q_{0},q_{1})}^{\Delta,1})=r_{1}$ and $t(d_{(q_{0},q_{1})}^{\Delta,1})=q_{1}$. Otherwise, we draw the 1-\textit{detour} $d_{(q_{0},q_{1})}^{\Delta,1}$ with ends $s(d_{(q_{0},q_{1})}^{\Delta,1})=r_{2}$ and $t(d_{(q_{0},q_{1})}^{\Delta,1})=q_{1}$.

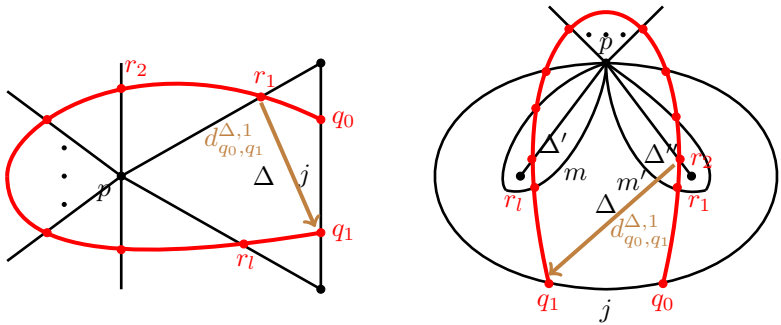
\begin{figure}[H]
\centering 
\begin{tikzpicture}[scale=0.75]

\filldraw [black] (-5,0) circle (2pt)
(-5,-4) circle (2pt)
(-8.5,-2) circle (2pt)
(-9.5,-2) circle (1pt)
(-9.5,-2.5) circle (1pt)
(-9.5,-1.5) circle (1pt);
\draw (-6,-2)node{$\Delta$};
\draw[line width=1pt] (-5, 0) -- (-5, -4)
node[pos=0.5,left] {$j$};
\draw[line width=1pt] (-5, 0) -- (-8.5, -2);
\draw[line width=1pt] (-5, -4) -- (-8.5, -2);
\draw[line width=1pt] (-8.5, -2) -- (-8.5,0);
\draw[line width=1pt] (-8.5, -2) -- (-10.5,-0.5);
\draw[line width=1pt] (-8.5, -2) -- (-10.5,-3.5);
\draw[line width=1pt] (-8.5, -2) -- (-8.5,-4)
node[pos=0.15,left] {$p$};

\filldraw [red] (-5,-1) circle (2pt)
(-5,-3) circle (2pt)
(-6.35,-3.2) circle (2pt)
(-8.5,-3.3) circle (2pt)
(-9.8,-3) circle (2pt)
(-9.8,-1) circle (2pt)
(-8.5,-0.45) circle (2pt)
(-6.05,-0.6) circle (2pt);

\draw[color=red][line width=1.5pt] (-5,-1) .. controls(-8, 0.5) and (-10.5,-1) .. (-10.5,-2)
node[pos=0,right] {$q_{0}$}
node[pos=0.11,above] {$r_{1}$}
node[pos=0.4,above] {$r_{2}$};
\draw[color=red][line width=1.5pt] (-5,-3) .. controls(-8, -3.5) and (-10.5,-3.5) .. (-10.5,-2)
node[pos=0,right] {$q_{1}$}
node[pos=0.15,below] {$r_{l}$};
\draw[->][color=brown][line width=1.5pt] (-6.05,-0.7) -- (-5.1,-2.9)
node[pos=0.3,left] {$d_{q_{0},q_{1}}^{\Delta , 1}$};

\filldraw [black] (0,0) circle (2pt)
(-1.5,-2) circle (2pt)
(1.5,-2) circle (2pt)
(0,0.5) circle (1pt)
(0.3,0.5) circle (1pt)
(-0.3,0.5) circle (1pt);
\draw (0,-2.5)node{$\Delta$};
\draw[line width=1pt] (0,0) -- (-1, 1)
node[pos=0.3,right] {$p$};
\draw[line width=1pt] (0,0) -- (1, 1);
\draw[line width=1pt] (0,0) -- (-1.5,-2);
\draw[line width=1pt] (0,0) -- (1.5,-2);
\draw[line width=1pt] (0, 0) .. controls(-4,-0.1) and (-4,-3.9) .. (0, -4);
\draw[line width=1pt] (0, 0) .. controls(4,-0.1) and (4,-3.9) .. (0, -4)
node[pos=1,below] {$j$};
\draw[line width=1pt] (0, 0) .. controls(-1.5,-0.5) and (-2,-2) .. (-1.75, -2.25);
\draw[line width=1pt] (0, 0) .. controls(0.05,-0.5) and (-0.75,-2.5) .. (-1.75, -2.25)
node[pos=0.7,right] {$m$}
node[pos=0.5,left] {$\Delta'$};
\draw[line width=1pt] (0, 0) .. controls(0,-1) and (0.7,-2.5) .. (1.75, -2.25)
node[pos=0.7,left] {$m'$}
node[pos=0.5,right] {$\Delta''$};
\draw[line width=1pt] (0, 0) .. controls(0,0) and (2.25,-1.5) .. (1.75, -2.25);


\filldraw [red] (-1,-3.9) circle (2pt)
(-1.25,-2.2) circle (2pt)
(-1.3,-1.7) circle (2pt)
(-1.23,-0.8) circle (2pt)
(-1.05,-0.15) circle (2pt)
(-0.65,0.6) circle (2pt)
(0.65,0.6) circle (2pt)
(1.05,-0.15) circle (2pt)
(1.23,-0.95) circle (2pt)
(1.3,-1.7) circle (2pt)
(1.25,-2.2) circle (2pt)
(1,-3.9) circle (2pt);

\draw[color=red][line width=1.5pt] (-1,-3.9) .. controls(-2.5,2.5) and (2.5,2.5) .. (1,-3.9)
node[pos=1,below] {$q_{0}$}
node[pos=0,below] {$q_{1}$}
node[pos=0.08,left] {$r_{l}$}
node[pos=0.87,right] {$r_{2}$}
node[pos=0.92,right] {$r_{1}$};

\draw[->][color=brown][line width=1.5pt] (1.2,-1.8) -- (-1,-3.75) 
node[pos=0.6,right] {$d_{q_{0},q_{1}}^{\Delta , 1}$};

\end{tikzpicture}
\end{figure}

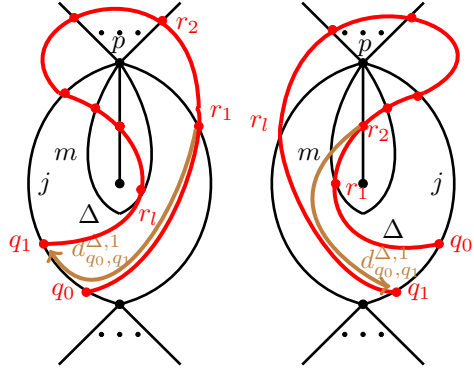
\begin{figure}[H]
\centering 
$$ \begin{tikzpicture}[scale=0.8]

\filldraw [black] (0,0) circle (2pt)
(0,-4) circle (2pt)
(0,-2) circle (2pt)
(0,0.5) circle (1pt)
(0.3,0.5) circle (1pt)
(-0.3,0.5) circle (1pt)
(0,-4.5) circle (1pt)
(0.3,-4.5) circle (1pt)
(-0.3,-4.5) circle (1pt);
\draw (-0.5,-2.5)node{$\Delta$};
\draw (4.5,-2.75)node{$\Delta$};
\draw[line width=1pt] (0,0) -- (-1, 1)
node[pos=0.3,right] {$p$};
\draw[line width=1pt] (0,0) -- (1, 1);
\draw[line width=1pt] (0,-4) -- (1, -5);
\draw[line width=1pt] (0,-4) -- (-1, -5);
\draw[line width=1pt] (0,0) -- (0, -2);
\draw[line width=1pt] (0,0) .. controls(-0.7,-1) and (-0.7,-2.2) .. (0, -2.5)
node[pos=0.5,left] {$m$};
\draw[line width=1pt] (0,0) .. controls(0.7,-1) and (0.7,-2.2) .. (0, -2.5);
\draw[line width=1pt] (0,0) .. controls(-2, -0.5) and (-2,-3.5) .. (0, -4)
node[pos=0.5,right] {$j$};
\draw[line width=1pt] (0,0) .. controls(2, -0.5) and (2,-3.5) .. (0, -4);

\filldraw [black] 
(4,0) circle (2pt)
(4,-4) circle (2pt)
(4,-2) circle (2pt)
(4,0.5) circle (1pt)
(4.3,0.5) circle (1pt)
(3.7,0.5) circle (1pt)
(4,-4.5) circle (1pt)
(4.3,-4.5) circle (1pt)
(3.7,-4.5) circle (1pt);
\draw[line width=1pt] (4,0) -- (3, 1)
node[pos=0.25,right] {$p$};
\draw[line width=1pt] (4,0) -- (5, 1);
\draw[line width=1pt] (4,-4) -- (5, -5);
\draw[line width=1pt] (4,-4) -- (3, -5);
\draw[line width=1pt] (4,0) -- (4, -2);
\draw[line width=1pt] (4,0) .. controls(2, -0.5) and (2,-3.5) .. (4, -4);
\draw[line width=1pt] (4,0) .. controls(6, -0.5) and (6,-3.5) .. (4, -4)
node[pos=0.5,left] {$j$};
\draw[line width=1pt] (4,0) .. controls(3.3,-1) and (3.3,-2.2) .. (4, -2.5)
node[pos=0.5,left] {$m$};
\draw[line width=1pt] (4,0) .. controls(4.7,-1) and (4.7,-2.2) .. (4, -2.5);



\filldraw [red] (4.55,-3.8) circle (2pt)
(5.25,-3) circle (2pt)
(3.55,-2) circle (2pt)
(4.4,-0.75) circle (2pt)
(4,-1.05) circle (2pt)
(4.9,-0.55) circle (2pt)
(4.8,0.75) circle (2pt)
(3.5,0.55) circle (2pt);
(2.65,-1.2) circle (2pt);
\draw[color=red][line width=1.5pt] (5.25,-3) .. controls(3,-3.5) and (3,-1) .. (5,-0.5)
node[pos=0,right] {$q_{0}$}
node[pos=0.5,right] {$r_{1}$}
node[pos=0.75,right] {$r_{2}$};
\draw[color=red][line width=1.5pt] (5,-0.5) .. controls(7,0.5) and (3,2) .. (2.65,-1)
node[pos=1,left] {$r_{l}$};
\draw[color=red][line width=1.5pt] (4.55,-3.8) .. controls(3.45,-3.8) and (2.5,-1.7) .. (2.65,-1)
node[pos=0,right] {$q_{1}$}; 
\draw[->][color=brown][line width=1.5pt] (3.95,-1.05) ..controls(2,-2. 5) and (4.3,-3.7) .. (4.45,-3.7)
node[pos=0.65,right] {$d_{q_{0},q_{1}}^{\Delta , 1}$}; 

\filldraw [red] (-0.55,-3.8) circle (2pt)
(-1.25,-3) circle (2pt)
(1.3,-1.05) circle (2pt)
(0.7,0.7) circle (2pt)
(-0.75,0.75) circle (2pt)
(-0.9,-0.5) circle (2pt)
(-0.4,-0.75) circle (2pt)
(0,-1.05) circle (2pt)
(0,-1.05) circle (2pt)
(0.35,-2.1) circle (2pt);
\draw[->][color=brown][line width=1.5pt] (1.25,-1.1) ..controls(1,-2.5) and (-0.25,-4.5) .. (-1.15,-3.1)
node[pos=0.7,above] {$d_{q_{0},q_{1}}^{\Delta , 1}$}; 

\draw[color=red][line width=1.5pt] (-0.55,-3.8) .. controls(1,-3.5) and (1.5,-0.75) .. (1.25,-0.75)
node[pos=0,left] {$q_{0}$}
node[pos=0.9,right] {$r_{1}$};
\draw[color=red][line width=1.5pt] (1.3,-0.75) .. controls(1.25,2) and (-2.25,0.75) .. (-1,-0.5)
node[pos=0.25,right] {$r_{2}$};
\draw[color=red][line width=1.5pt] (-1.25,-3) .. controls(1,-3) and (0.75,-1) .. (-1,-0.5)
node[pos=0.3,right] {$r_{l}$}
node[pos=0,left] {$q_{1}$};

\end{tikzpicture}$$
\caption{Case $a$) of the Definition \ref{1-desviacion}.}
\end{figure} 
\item[b)] The arc $j$ is the folded side of a self-folded triangle $\Delta'$ in $\tau^{\circ}$. We denote by $m'$ the non folded side of $\Delta'$.

Let's consider the cases when $r_{2}$ belongs to the non folded side of some self-folded triangle $\Delta''$ different to $\Delta'$ or $r_{2}$ does not belong to any self-folded triangle:

\begin{itemize}
\item[$b_{1}$)]The crossing point $r_{2}$ belongs to the non folded side of a self-folded triangle $\Delta''$ different than $\Delta'$. We denote by $m''$ the non folded side of $\Delta''$.

Under these conditions we draw the 1-\textit{detour} $d_{(q_{0},q_{1})}^{\Delta,1}$ with ends $s(d_{(q_{0},q_{1})}^{\Delta,1})=r_{3}$ and $t(d_{(q_{0},q_{1})}^{\Delta,1})=r_{l}$. In addition, if the segment $[q_{0},r_{2}]_{i'}$ is contractible to $p$ with the homotopy that preserves $q_{0}$ in $j$ and $r_{2}$ in $m''$, and the segment
$[q_{1},r_{l-1}]_{i'}$ is not contractible to $p$, then we draw the auxiliary curve $e_{(q_{0},q_{1})}^{\Delta,1}$ with ends $s(e_{(q_{0},q_{1})}^{\Delta,1})=r_{3}$ and $t(d_{(q_{0},q_{1})}^{\Delta,1})=q_{1}$.

\begin{figure}[H]
\centering 
\vspace{0em}\hspace*{-6em}{\begin{tikzpicture}[scale=0.8]

\filldraw [black] (0,0) circle (2pt)
(-1.5,-2) circle (2pt)
(1.5,-2) circle (2pt)
(0,0.5) circle (1pt)
(0.3,0.5) circle (1pt)
(-0.3,0.5) circle (1pt);
\draw (0,-2.5)node{$\Delta$};
\draw (-1,-0.8)node{$\Delta''$};
\draw (1.4,-1.4)node{$j$};
\draw (0.8,-1.5)node{$\Delta'$};
\draw[line width=1pt] (0,0) -- (-1, 1)
node[pos=0.3,right] {$p$};
\draw[line width=1pt] (0,0) -- (1, 1);
\draw[line width=1pt] (0,0) -- (-1.5,-2);
\draw[line width=1pt] (0,0) -- (1.5,-2);
\draw[line width=1pt] (0, 0) .. controls(-4,-0.1) and (-4,-3.9) .. (0, -4);
\draw[line width=1pt] (0, 0) .. controls(4,-0.1) and (4,-3.9) .. (0, -4);
\draw[line width=1pt] (0, 0) .. controls(-1.5,-0.5) and (-2,-2) .. (-1.75, -2.25);
\draw[line width=1pt] (0, 0) .. controls(0.05,-0.5) and (-0.75,-2.5) .. (-1.75, -2.25)
node[pos=0.8,right] {$m''$};
\draw[line width=1pt] (0, 0) .. controls(0,-1) and (0.7,-2.5) .. (1.75, -2.25)
node[pos=0.8,left] {$m'$};
\draw[line width=1pt] (0, 0) .. controls(0,0) and (2.25,-1.5) .. (1.75, -2.25);


\filldraw [black] (8,0) circle (2pt)
(6.5,-2) circle (2pt)
(9.5,-2) circle (2pt)
(8,0.5) circle (1pt)
(8.3,0.5) circle (1pt)
(7.7,0.5) circle (1pt);
\draw (8,-3.5)node{$\Delta$};
\draw (7.15,-1.6)node{$\Delta''$};
\draw (9.4,-1.4)node{$j$};
\draw (8.6,-1.2)node{$\Delta'$};
\draw[line width=1pt] (8,0) -- (7, 1)
node[pos=0.3,right] {$p$};
\draw[line width=1pt] (8,0) -- (9, 1);
\draw[line width=1pt] (8,0) -- (6.5,-2);
\draw[line width=1pt] (8,0) -- (9.5,-2);
\draw[line width=1pt] (8, 0) .. controls(4,-0.1) and (4,-3.9) .. (8, -4);
\draw[line width=1pt] (8, 0) .. controls(12,-0.1) and (12,-3.9) .. (8, -4);
\draw[line width=1pt] (8, 0) .. controls(6.5,-0.5) and (6,-2) .. (6.25, -2.25);
\draw[line width=1pt] (8, 0) .. controls(8.05,-0.5) and (7.25,-2.5) .. (6.25, -2.25)
node[pos=0.8,right] {$m$};
\draw[line width=1pt] (8, 0) .. controls(8,-1) and (8.7,-2.5) .. (9.75, -2.25);
\draw[line width=1pt] (8, 0) .. controls(8,0) and (10.25,-1.5) .. (9.75, -2.25)
node[pos=0.8,right] {$m'$};;

\filldraw [red] (1.3,-1.8) circle (2pt)
(1.75,-1.75) circle (2pt)
(-1.6,-1.35) circle (2pt)
(-1,-1.25) circle (2pt)
(-0.3,-1.05) circle (2pt)
(0.15,-0.85) circle (2pt)
(0.45,-0.6) circle (2pt)
(0.6,-0.4) circle (2pt)
(0.7,-0.05) circle (2pt)
(0.6,0.6) circle (2pt)
(-0.8,0.8) circle (2pt)
(1.1,-0.85) circle (2pt)
(2.75,-1.2) circle (2pt)
(0.75,-1) circle (2pt);
\draw[color=red][line width=1.5pt] (0.8,-1) .. controls(1.2,-0.5) and (3.3,-1).. (3.3,-2)
node[pos=0,left] {$q_{1}$}
node[pos=0.2,above] {$r_{l}$}; 
\draw[color=red][line width=1.5pt] (-3.5,-2) ..controls(-3.6,-5) and (3.05,-5) .. (3.3,-2);
\draw[color=red][line width=1.5pt] (-3.5,-2) ..controls(-3.1,1) and (0.2,1.2) .. (0.6,0.6);
\draw[color=red][line width=1.5pt] (-1.75,-1.35) ..controls(-1.15,-1.35) and (1.3,-1) .. (0.6,0.6)
node[pos=0.3,below] {$r_{3}$}; 
\draw[color=red][line width=1.5pt] (-1.65,-1.35) ..controls (-3.5,-2) and (-0.95,-3.5) .. (1.05,-2.75)
node[pos=0,above] {$r_{2}$};
\draw[color=red][line width=1.5pt] (1.05,-2.75).. controls (1.45,-2.7) and (3,-1.5) .. (1.3,-1.8)
node[pos=1,left] {$q_{0}$}
node[pos=0.87,above] {$r_{1}$};
\draw[color=brown][line width=1.5pt] (-1.05,-1.25) ..controls (-5,-1) and (-1.5,-4.5) .. (1.05,-3.25);
\draw[-<][color=brown][line width=1.5pt] (1.05,-3.25) ..controls(3.75,-2) and (1.1,-0.85) .. (1.4,-1.05)
node[pos=0.12,below] {$d_{q_{0},q_{1}}^{\Delta , 1}$};

\filldraw [red] (9.1,-1.5) circle (2pt)
(8.65,-1.85) circle (2pt)
(6.35,-1.35) circle (2pt)
(7,-1.3) circle (2pt)
(7.7,-1.1) circle (2pt)
(8.1,-0.9) circle (2pt)
(8.45,-0.6) circle (2pt)
(8.6,-0.45) circle (2pt)
(8.7,-0.05) circle (2pt)
(8.6,0.6) circle (2pt)
(7.2,0.8) circle (2pt)
(5,-2) circle (2pt)
(9.25,-2.2) circle (2pt)
(9.4,-1.9) circle (2pt)
(9.7,-1.6) circle (2pt)
(10.4,-0.8) circle (2pt);

\draw[color=red][line width=1.5pt] (5,-2) .. controls (5.2,-4.25) and (9,-3.5) .. (9.4,-1.9)
node[pos=1,right] {$q_{1}$}
node[pos=0.9,below] {$r_{l}$};
\draw[color=red][line width=1.5pt] (5,-2) ..controls(4.9,1) and (8.2,1.2) .. (8.6,0.6);
\draw[color=red][line width=1.5pt] (6.35,-1.35) ..controls(6.85,-1.35) and (9.3,-1) .. (8.6,0.6);
\draw[color=red][line width=1.5pt] (6.35,-1.35) ..controls (4.5,-2) and (7.05,-3.5) .. (9.05,-1.5)
node[pos=1,above] {$q_{0}$}
node[pos=0.87,above] {$r_{1}$}
node[pos=0.05,below] {$r_{2}$};
\draw[color=red][line width=1.5pt] (9.4,-1.9) -- (10.4,-0.8);

\draw[->][color=brown][line width=1.5pt] (7,-1.3) ..controls(3.1,-0.5) and (6,-5.25) .. (9.15,-2.2)
node[pos=0.8,below] {$d_{q_{0},q_{1}}^{\Delta , 1}$};

\draw[->][color=blue][line width=1.5pt] (7,-1.3) ..controls(3.8,-0.9) and (6,-4.8) .. (9.4,-1.9)
node[pos=0.8,above] {$e_{q_{0},q_{1}}^{\Delta , 1}$};

\end{tikzpicture}} 

\end{figure} 
\begin{figure}[H]
\centering 
\vspace{-4em}\hspace*{-2em}{\begin{tikzpicture}[scale=0.8]

\filldraw [black] (0,0) circle (2pt)
(-1.5,-2) circle (2pt)
(1.5,-2) circle (2pt)
(0,0.5) circle (1pt)
(0.3,0.5) circle (1pt)
(-0.3,0.5) circle (1pt);
\draw (0,-3)node{$\Delta$};
\draw (-1,-1)node{$\Delta''$};
\draw (-1.5,-1.6)node{$j$};
\draw (0.8,-1.5)node{$\Delta'$};
\draw[line width=1pt] (0,0) -- (-1, 1)
node[pos=0.3,right] {$p$};
\draw[line width=1pt] (0,0) -- (1, 1);
\draw[line width=1pt] (0,0) -- (-1.5,-2);
\draw[line width=1pt] (0,0) -- (1.5,-2);
\draw[line width=1pt] (0, 0) .. controls(-4,-0.1) and (-4,-3.9) .. (0, -4);
\draw[line width=1pt] (0, 0) .. controls(4,-0.1) and (4,-3.9) .. (0, -4);
\draw[line width=1pt] (0, 0) .. controls(-1.5,-0.5) and (-2,-2) .. (-1.75, -2.25);
\draw[line width=1pt] (0, 0) .. controls(0.05,-0.5) and (-0.75,-2.5) .. (-1.75, -2.25)
node[pos=0.8,right] {$m''$};
\draw[line width=1pt] (0, 0) .. controls(0,-1) and (0.7,-2.5) .. (1.75, -2.25)
node[pos=0.8,left] {$m'$};
\draw[line width=1pt] (0, 0) .. controls(0,0) and (2.25,-1.5) .. (1.75, -2.25);

\filldraw [black] (8,0) circle (2pt)
(6.5,-2) circle (2pt)
(9.5,-2) circle (2pt)
(8,0.5) circle (1pt)
(8.3,0.5) circle (1pt)
(7.7,0.5) circle (1pt);
\draw (8,-3)node{$\Delta$};
\draw (7,-1)node{$\Delta''$};
\draw (6.5,-1.6)node{$j$};
\draw (8.8,-1.5)node{$\Delta'$};
\draw[line width=1pt] (8,0) -- (7, 1)
node[pos=0.3,right] {$p$};
\draw[line width=1pt] (8,0) -- (9, 1);
\draw[line width=1pt] (8,0) -- (6.5,-2);
\draw[line width=1pt] (8,0) -- (9.5,-2);
\draw[line width=1pt] (8, 0) .. controls(4,-0.1) and (4,-3.9) .. (8, -4);
\draw[line width=1pt] (8, 0) .. controls(12,-0.1) and (12,-3.9) .. (8, -4);
\draw[line width=1pt] (8, 0) .. controls(6.5,-0.5) and (6,-2) .. (6.25, -2.25);
\draw[line width=1pt] (8, 0) .. controls(8.05,-0.5) and (7.25,-2.5) .. (6.25, -2.25)
node[pos=0.8,right] {$m''$};
\draw[line width=1pt] (8, 0) .. controls(8,-1) and (8.7,-2.5) .. (9.75, -2.25)
node[pos=0.8,left] {$m'$};
\draw[line width=1pt] (8, 0) .. controls(8,0) and (10.25,-1.5) .. (9.75, -2.25);

\filldraw [red] (1.3,-1.8) circle (2pt)
(0.5,-1.75) circle (2pt)
(-1.6,-1.35) circle (2pt)
(-1,-1.25) circle (2pt)
(-0.3,-1.05) circle (2pt)
(0.15,-0.85) circle (2pt)
(0.45,-0.6) circle (2pt)
(0.6,-0.4) circle (2pt)
(0.7,-0.05) circle (2pt)
(0.6,0.6) circle (2pt)
(-0.8,0.8) circle (2pt)
(1.5,-1.25) circle (2pt)
(2.75,-1.2) circle (2pt)
(0.75,-1) circle (2pt);

\draw[color=red][line width=1.5pt] (1.3,-1.8) .. controls(1.2,-0.5) and (3.3,-1).. (3.3,-2)
node[pos=0,left] {$q_{1}$}
node[pos=0.3,above] {$r_{l}$}; 
\draw[color=red][line width=1.5pt] (-3.5,-2) ..controls(-3.6,-5) and (3.05,-5) .. (3.3,-2);
\draw[color=red][line width=1.5pt] (-3.5,-2) ..controls(-3.1,1) and (0.2,1.2) .. (0.6,0.6);

\draw[color=red][line width=1.5pt] (-1.75,-1.35) ..controls(-1.15,-1.35) and (1.3,-1) .. (0.6,0.6)
node[pos=0.3,below] {$r_{3}$}; 
\draw[color=red][line width=1.5pt] (-1.65,-1.35) ..controls (-3.5,-2) and (-1.2,-2.9) .. (-1,-2.75)
node[pos=0,above] {$r_{2}$};
\draw[color=red][line width=1.5pt] (-1,-2.75).. controls (0.45,-2.7) and (0.5,-1.5) .. (0.8,-1)
node[pos=1,right] {$q_{0}$}
node[pos=0.6,left] {$r_{1}$};

\draw[color=brown][line width=1.5pt] (-1,-1.2) ..controls(-4.2,-0.9) and (-2,-4.8) .. (1.55,-2.88)
node[pos=0.8,above] {$d_{q_{0},q_{1}}^{\Delta , 1}$};

\draw[->][color=brown][line width=1.5pt] (1.5,-2.88) ..controls (1.8,-2.8) and (3,-1) .. (1.6,-1.25);
\filldraw [red] (9.3,-1.8) circle (2pt)
(8.7,-1.9) circle (2pt)
(6.4,-1.35) circle (2pt)
(7,-1.25) circle (2pt)
(7.7,-1.05) circle (2pt)
(8.15,-0.85) circle (2pt)
(8.45,-0.6) circle (2pt)
(8.6,-0.4) circle (2pt)
(8.7,-0.05) circle (2pt)
(8.6,0.6) circle (2pt)
(7.2,0.8) circle (2pt)
(9.1,-0.85) circle (2pt)
(10.75,-1.2) circle (2pt)
(8.75,-1) circle (2pt);

\draw[color=red][line width=1.5pt] (8.8,-1) .. controls(9.2,-0.5) and (11.3,-1).. (11.3,-2)
node[pos=0,left] {$q_{1}$}
node[pos=0.3,above] {$r_{l}$}; 
\draw[color=red][line width=1.5pt] (4.5,-2) ..controls(4.4,-5) and (11.05,-5) .. (11.3,-2);
\draw[color=red][line width=1.5pt] (4.5,-2) ..controls(4.9,1) and (8.2,1.2) .. (8.6,0.6);

\draw[color=red][line width=1.5pt] (6.25,-1.35) ..controls(6.85,-1.35) and (9.3,-1) .. (8.6,0.6)
node[pos=0.3,below] {$r_{3}$}; 
\draw[color=red][line width=1.5pt] (6.35,-1.35) ..controls (4.5,-2) and (6.8,-2.9) .. (7,-2.75)
node[pos=0,above] {$r_{2}$};
\draw[color=red][line width=1.5pt] (7,-2.75).. controls (8.45,-2.7) and (8.5,-1.5) .. (9.3,-1.8)
node[pos=1,right] {$q_{0}$}
node[pos=0.6,left] {$r_{1}$};

\draw[color=brown][line width=1.5pt] (7,-1.2) ..controls(3.8,-0.9) and (6,-4.8) .. (9.55,-2.88)
node[pos=0.8,above] {$d_{q_{0},q_{1}}^{\Delta , 1}$};

\draw[->][color=brown][line width=1.5pt] (9.5,-2.88) ..controls (9.8,-2.8) and (11,-1) .. (9.35,-1);

\end{tikzpicture}} 
\caption{Case $b_{1})$ of the Definition \ref{1-desviacion}.}

\end{figure} 

\item[$b_{2}$)]The crossing point $r_{2}$ does not belong to any side of any self-folded triangle.

We draw the 1-\textit{detour} $d_{(q_{0},q_{1})}^{\Delta,1}$ with ends $s(d_{(q_{0},q_{1})}^{\Delta,1})=r_{2}$ and $t(d_{(q_{0},q_{1})}^{\Delta,1})=q_{1}$ if the conditions of one of the two followings subcases are met.

\begin{myEnumerate}
\item[$b_{21}$)] The arc $m'$ belongs to a triangle of type 2.
\begin{myEnumerate}
\item[i)] The arc $m'$ is based on one punctured denoted by $q$,
\item[ii)] The segment $[q_{0},r_{2}]_{i'}$ is contractibe to $q$ with the homotopy that preserves the ends in the respective arc of $\tau^{\circ}$.
\end{myEnumerate}
\item[$b_{22}$)] The arc $m'$ belongs to a triangle of type 3.
\begin{myEnumerate}
\item[i)] The arc $m'$ is based on one punctured denoted by $p$,
\item[ii)] The segment $[q_{0},r_{2}]_{i'}$ is contractibe to $p$ with the homotopy that preserves the ends in the respective arc of $\tau^{\circ}$,
\item[iii)]The segment $[q_{1},r_{l-1}]_{i'}$ is not contractible to the puncture $p$ with the homotopy that preserves the ends in the respective arc of $\tau^{\circ}$.
\end{myEnumerate}
\end{myEnumerate}

Under the situation decribed in $b_{21}$) we draw the auxiliary curve $e_{(q_{0},q_{1})}^{\Delta,1}$ with ends $s(e_{(q_{0},q_{1})}^{\Delta,1})=r_{1}$ and $t(d_{(q_{0},q_{1})}^{\Delta,1})=q_{2}$. In addition, if the crossing point $q_{0}$ is an end of $i'$, then we draw another auiliary curve with ends $s(e_{(q_{0},q_{1})}^{\Delta,1})=r_{2}$ and $t(d_{(q_{0},q_{1})}^{\Delta,1})=r_{l}$. Otherwise, in the situatin described in $b_{22}$) we draw an auxiliary curve $e_{(q_{0},q_{1})}^{\Delta,1}$ with ends $s(e_{(q_{0},q_{1})}^{\Delta,1})=r_{2}$ y $t(d_{(q_{0},q_{1})}^{\Delta,1})=r_{l}$.

Now, if $b_{21}$) and $b_{22}$) are not met, then we draw the 1-\textit{detour} $d_{(q_{0},q_{1})}^{\Delta,1}$ with ends $s(d_{(q_{0},q_{1})}^{\Delta,1})=r_{2}$ and $t(d_{(q_{0},q_{1})}^{\Delta,1})=r_{l}$.

\end{itemize} 

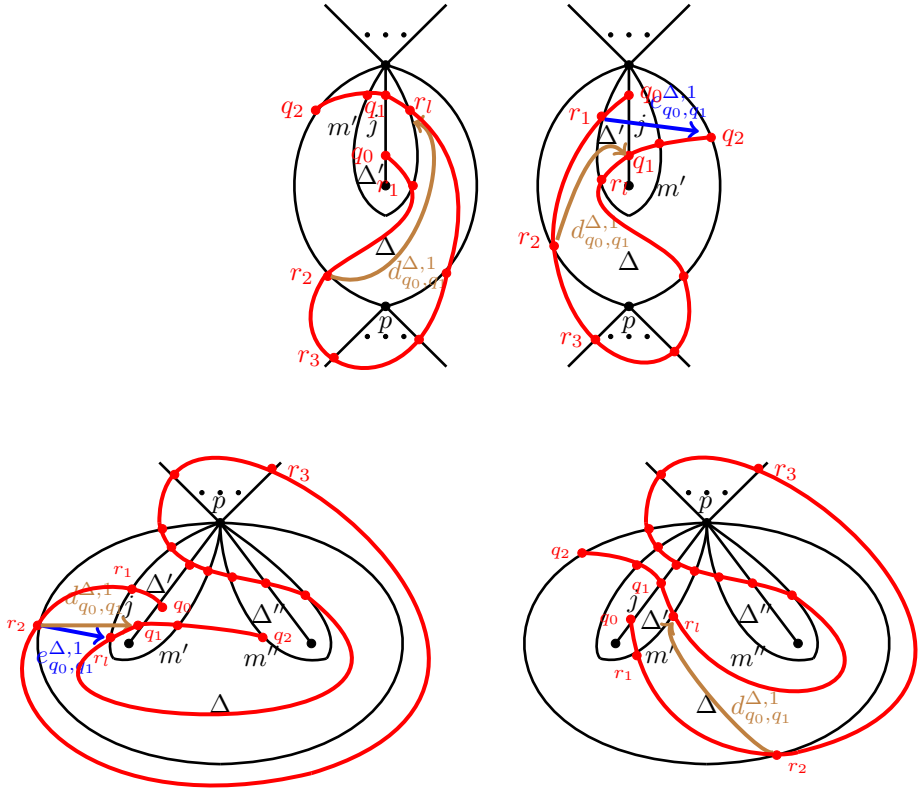
\begin{figure}[H]
\centering 
$$ \begin{tikzpicture}[scale=0.8]

\filldraw [black] (0,0) circle (2pt)
(0,-4) circle (2pt)
(0,-2) circle (2pt)
(0,0.5) circle (1pt)
(0.3,0.5) circle (1pt)
(-0.3,0.5) circle (1pt)
(0,-4.5) circle (1pt)
(0.3,-4.5) circle (1pt)
(-0.3,-4.5) circle (1pt);
\draw (0,-3)node{$\Delta$};
\draw (4,-3.25)node{$\Delta$};
\draw (3.7,-1.15)node{$\Delta'$};
\draw (-0.25,-1.8)node{$\Delta'$};
\draw (4.7,-2)node{$m'$};
\draw (4.2,-1)node{$j$};
\draw (-0.7,-1)node{$m'$};
\draw (-0.2,-1)node{$j$};
\draw (0,-4.3)node{$p$};
\draw[line width=1pt] (0,0) -- (-1, 1);
\draw[line width=1pt] (0,0) -- (1, 1);
\draw[line width=1pt] (0,-4) -- (1, -5);
\draw[line width=1pt] (0,-4) -- (-1, -5);
\draw[line width=1pt] (0,0) -- (0, -2);
\draw[line width=1pt] (0,0) .. controls(-0.7,-1) and (-0.7,-2.2) .. (0, -2.5);
\draw[line width=1pt] (0,0) .. controls(0.7,-1) and (0.7,-2.2) .. (0, -2.5);
\draw[line width=1pt] (0,0) .. controls(-2, -0.5) and (-2,-3.5) .. (0, -4);
\draw[line width=1pt] (0,0) .. controls(2, -0.5) and (2,-3.5) .. (0, -4);

\filldraw [black] 
(4,0) circle (2pt)
(4,-4) circle (2pt)
(4,-2) circle (2pt)
(4,0.5) circle (1pt)
(4.3,0.5) circle (1pt)
(3.7,0.5) circle (1pt)
(4,-4.5) circle (1pt)
(4.3,-4.5) circle (1pt)
(3.7,-4.5) circle (1pt);
\draw (4,-4.3)node{$p$}; 
\draw[line width=1pt] (4,0) -- (3, 1);
\draw[line width=1pt] (4,0) -- (5, 1);
\draw[line width=1pt] (4,-4) -- (5, -5);
\draw[line width=1pt] (4,-4) -- (3, -5);
\draw[line width=1pt] (4,0) -- (4, -2);
\draw[line width=1pt] (4,0) .. controls(2, -0.5) and (2,-3.5) .. (4, -4);
\draw[line width=1pt] (4,0) .. controls(6, -0.5) and (6,-3.5) .. (4, -4);
\draw[line width=1pt] (4,0) .. controls(3.3,-1) and (3.3,-2.2) .. (4, -2.5);
\draw[line width=1pt] (4,0) .. controls(4.7,-1) and (4.7,-2.2) .. (4, -2.5);



\filldraw [red] (4,-0.5) circle (2pt)
(3.55,-0.85) circle (2pt)
(2.77,-3) circle (2pt)
(3.45,-4.55) circle (2pt)
(4.75,-4.75) circle (2pt)
(4.9,-3.5) circle (2pt)
(3.55,-1.9) circle (2pt)
(4,-1.5) circle (2pt) 
(4.5,-1.3) circle (2pt) 
(5.35,-1.2) circle (2pt);
\draw[->][color=brown][line width=1.5pt] (2.82,-2.9) ..controls (3.1,-2) and (3.5,-1) .. (3.95,-1.5)
node[pos=0.05,right] {$d_{q_{0},q_{1}}^{\Delta , 1}$};
\draw[->][color=blue][line width=1.5pt] (3.6,-0.9) -- (5.15,-1.1)
node[pos=0.8,above] {$e_{q_{0},q_{1}}^{\Delta , 1}$};

\draw[color=red][line width=1.5pt] (4,-0.5) .. controls(3.25,-1) and (2.75,-2) .. (2.75,-2.9)
node[pos=0,right] {$q_{0}$}
node[pos=0.2,left] {$r_{1}$}
node[pos=1,left] {$r_{2}$};
\draw[color=red][line width=1.5pt] (2.75,-2.9) .. controls(3,-5.5) and (5.5,-5.5) .. (4.9,-3.5)
node[pos=0.3,left] {$r_{3}$};
\draw[color=red][line width=1.5pt] (4.9,-3.5) .. controls(4.8,-2.8) and (2.5,-2.5) .. (4,-1.5)
node[pos=0.8,right] {$r_{l}$}
node[pos=1.05,below] {$q_{1}$};
\draw[color=red][line width=1.5pt] (4,-1.5) ..controls(4.25,-1.3) and (5.25,-1.2) .. (5.35,-1.2)
node[pos=1,right] {$q_{2}$};

\filldraw [red] (-1.15,-0.75) circle (2pt) 
(-0.3,-0.5) circle (2pt)
(0,-0.5) circle (2pt)
(0.4,-0.75) circle (2pt)
(1,-3.45) circle (2pt)
(0.55,-4.55) circle (2pt)
(-0.85,-4.85) circle (2pt)
(-0.95,-3.5) circle (2pt)
(0.45,-2) circle (2pt)
(0,-1.5) circle (2pt);
\draw[->][color=brown][line width=1.5pt] (-0.9,-3.5) ..controls (0.5,-4) and (1.25,-1.1) .. (0.45,-0.95)
node[pos=0.2,right] {$d_{q_{0},q_{1}}^{\Delta , 1}$};

\draw[color=red][line width=1.5pt] (-1.15,-0.75) ..controls(-0.95,-0.5) and (-0.2,-0.4) .. (0,-0.5)
node[pos=0,left] {$q_{2}$};

\draw[color=red][line width=1.5pt] (0,-0.5) .. controls(1,-1) and (1.5,-2) .. (1,-3.5)
node[pos=-0.05,below] {$q_{1}$}
node[pos=0.1,right] {$r_{l}$};
\draw[color=red][line width=1.5pt] (1,-3.5) .. controls(0.5,-6) and (-2,-5) .. (-1,-3.5)
node[pos=0.6,left] {$r_{3}$}
node[pos=1,left] {$r_{2}$};
\draw[color=red][line width=1.5pt] (-1,-3.5) .. controls(-0.5,-3) and (1.25,-2.5) .. (0,-1.5)
node[pos=0.8,left] {$r_{1}$}
node[pos=1,left] {$q_{0}$};

\end{tikzpicture}$$
\end{figure}

\begin{figure}[H]
\centering 
\vspace{-8em}\hspace*{-2em}{\begin{tikzpicture}[scale=0.8]

\filldraw [black] (0,0) circle (2pt)
(-1.5,-2) circle (2pt)
(1.5,-2) circle (2pt)
(0,0.5) circle (1pt)
(0.3,0.5) circle (1pt)
(-0.3,0.5) circle (1pt);
\draw (0,-3)node{$\Delta$};
\draw (-1,-1)node{$\Delta'$};
\draw (-1.5,-1.4)node{$j$};
\draw (0.8,-1.5)node{$\Delta''$};
\draw[line width=1pt] (0,0) -- (-1, 1)
node[pos=0.3,right] {$p$};
\draw[line width=1pt] (0,0) -- (1, 1);
\draw[line width=1pt] (0,0) -- (-1.5,-2);
\draw[line width=1pt] (0,0) -- (1.5,-2);
\draw[line width=1pt] (0, 0) .. controls(-4,-0.1) and (-4,-3.9) .. (0, -4);
\draw[line width=1pt] (0, 0) .. controls(4,-0.1) and (4,-3.9) .. (0, -4);
\draw[line width=1pt] (0, 0) .. controls(-1.5,-0.5) and (-2,-2) .. (-1.75, -2.25);
\draw[line width=1pt] (0, 0) .. controls(0.05,-0.5) and (-0.75,-2.5) .. (-1.75, -2.25)
node[pos=0.8,right] {$m'$};
\draw[line width=1pt] (0, 0) .. controls(0,-1) and (0.7,-2.5) .. (1.75, -2.25)
node[pos=0.8,left] {$m''$};
\draw[line width=1pt] (0, 0) .. controls(0,0) and (2.25,-1.5) .. (1.75, -2.25);


\filldraw [black] (8,0) circle (2pt)
(6.5,-2) circle (2pt)
(9.5,-2) circle (2pt)
(8,0.5) circle (1pt)
(8.3,0.5) circle (1pt)
(7.7,0.5) circle (1pt);
\draw (8,-3)node{$\Delta$};
\draw (7.15,-1.6)node{$\Delta'$};
\draw (6.8,-1.3)node{$j$};
\draw (8.8,-1.5)node{$\Delta''$};
\draw[line width=1pt] (8,0) -- (7, 1)
node[pos=0.3,right] {$p$};
\draw[line width=1pt] (8,0) -- (9, 1);
\draw[line width=1pt] (8,0) -- (6.5,-2);
\draw[line width=1pt] (8,0) -- (9.5,-2);
\draw[line width=1pt] (8, 0) .. controls(4,-0.1) and (4,-3.9) .. (8, -4);
\draw[line width=1pt] (8, 0) .. controls(12,-0.1) and (12,-3.9) .. (8, -4);
\draw[line width=1pt] (8, 0) .. controls(6.5,-0.5) and (6,-2) .. (6.25, -2.25);
\draw[line width=1pt] (8, 0) .. controls(8.05,-0.5) and (7.25,-2.5) .. (6.25, -2.25)
node[pos=0.8,right] {$m'$};
\draw[line width=1pt] (8, 0) .. controls(8,-1) and (8.7,-2.5) .. (9.75, -2.25)
node[pos=0.8,left] {$m''$};
\draw[line width=1pt] (8, 0) .. controls(8,0) and (10.25,-1.5) .. (9.75, -2.25);


\filldraw [red] (-0.95,-1.4) circle (2pt)
(-1.45,-1.1) circle (2pt)
(-3,-1.7) circle (2pt)
(0.85,0.9) circle (2pt)
(-0.75,0.8) circle (2pt)
(-0.95,-0.1) circle (2pt)
(-0.8,-0.4) circle (2pt)
(-0.5,-0.7) circle (2pt)
(-0.2,-0.8) circle (2pt)
(0.2,-0.9) circle (2pt)
(0.75,-1) circle (2pt)
(1.4,-1.2) circle (2pt)
(-1.8,-1.9) circle (2pt)
(-1.35,-1.7) circle (2pt)
(-0.7,-1.7) circle (2pt)
(0.7,-1.9) circle (2pt);
\draw[->][color=blue][line width=1.5pt] (-2.95,-1.7) -- (-1.9,-1.9)
node[pos=0.4,below] {$e_{q_{0},q_{1}}^{\Delta , 1}$};
\draw[->][color=brown][line width=1.5pt] (-2.95,-1.7) -- (-1.45,-1.7)
node[pos=0.6,above] {$d_{q_{0},q_{1}}^{\Delta , 1}$};

\draw[color=red][line width=1.5pt] (-0.95,-1.4) ..controls(-1.1,-1) and (-2.5,-0.8) .. (-3,-1.7)
node[pos=-0,right] {{\scriptsize $q_{0}$}}
node[pos=0.4,above] {{\scriptsize $r_{1}$}}
node[pos=1,left] {{\scriptsize $r_{2}$}}; 
\draw[color=red][line width=1.5pt] (-3,-1.7) .. controls(-4,-3) and (-2,-5) .. (1.5,-4.15);
\draw[color=red][line width=1.5pt] (1.5,-4.15) .. controls(7,-2.5) and (-1,3.5) .. (-1,0)
node[pos=0.66,right] {$r_{3}$};
\draw[color=red][line width=1.5pt] (-1,0) .. controls(-0.8,-1) and (0.5,-0.8) .. (1.4,-1.2);
\draw[color=red][line width=1.5pt] (1.4,-1.2) .. controls(5,-4) and (-5.5,-3.5) .. (-1.35,-1.7)
node[pos=0.94,below] {{\scriptsize $r_{l}$}}
node[pos=1.02,below] {{\scriptsize $q_{1}$}};

\draw[color=red][line width=1.5pt] (-1.35,-1.7) .. controls(-1.25,-1.6) and (0.6,-1.8) .. (0.7,-1.9)
node[pos=1,right] {{\scriptsize $q_{2}$}};

\filldraw [red] (5.95,-0.5) circle (2pt)
(6.85,-0.7) circle (2pt)
(7.25,-1) circle (2pt)
(9.15,-3.85) circle (2pt)
(8.85,0.9) circle (2pt)
(7.25,0.8) circle (2pt)
(7.05,-0.1) circle (2pt)
(7.2,-0.4) circle (2pt)
(7.5,-0.7) circle (2pt)
(7.8,-0.8) circle (2pt)
(8.2,-0.9) circle (2pt)
(8.75,-1) circle (2pt)
(9.4,-1.2) circle (2pt)
(7.45,-1.55) circle (2pt)
(6.75,-1.6) circle (2pt)
(6.85,-2.2) circle (2pt);
\draw[->][color=brown][line width=1.5pt] (9.1,-3.8) ..controls (8.9,-3.95) and (7.15,-2) .. (7.45,-1.65)
node[pos=0.45,right] {$d_{q_{0},q_{1}}^{\Delta , 1}$};

\draw[color=red][line width=1.5pt] (6.75,-1.6) .. controls(7,-4) and (9.3,-3.9) .. (9.5,-3.8)
node[pos=-0,left] {{\scriptsize $q_{0}$}}
node[pos=0.15,left] {{\scriptsize $r_{1}$}}
node[pos=1,below] {{\scriptsize $r_{2}$}};
\draw[color=red][line width=1.5pt] (9.5,-3.8) .. controls(15,-2.5) and (7,3.5) .. (7,0) 
node[pos=0.66,right] {$r_{3}$};
\draw[color=red][line width=1.5pt] (7,0) .. controls(7.2,-1) and (8.5,-0.8) .. (9.4,-1.2);
\draw[color=red][line width=1.5pt] (9.4,-1.2) .. controls(12,-3) and (8,-3.7) .. (7.25,-1) 
node[pos=0.85,above] {{\scriptsize $r_{l}$}}
node[pos=0.99,left] {{\scriptsize $q_{1}$}};
\draw[color=red][line width=1.5pt] (7.25,-1) ..controls(7,-0.5) and (5.95,-0.5) .. (5.95,-0.5)
node[pos=1,left] {{\scriptsize $q_{2}$}};

\end{tikzpicture}}

\end{figure} 
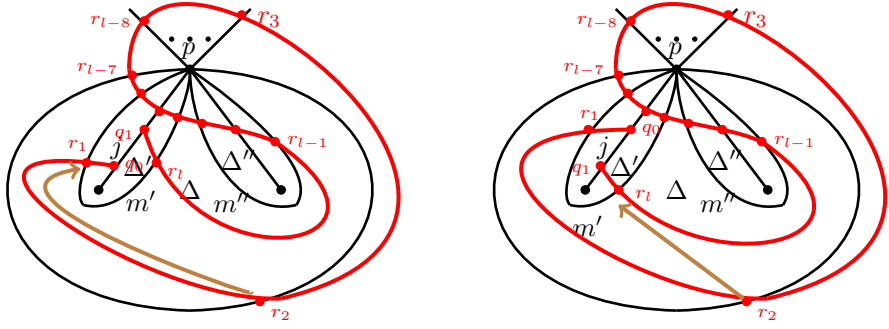
\begin{figure}[H]
\centering 
\vspace{-8em}\hspace*{-2em}{\begin{tikzpicture}[scale=0.8]


\filldraw [black] (0,0) circle (2pt)
(-1.5,-2) circle (2pt)
(1.5,-2) circle (2pt)
(0,0.5) circle (1pt)
(0.3,0.5) circle (1pt)
(-0.3,0.5) circle (1pt);
\draw (0,-2)node{$\Delta$};
\draw (-0.85,-1.6)node{$\Delta'$};
\draw (-1.2,-1.3)node{$j$};
\draw (0.8,-1.5)node{$\Delta''$};
\draw[line width=1pt] (0,0) -- (-1, 1)
node[pos=0.3,right] {$p$};
\draw[line width=1pt] (0,0) -- (1, 1);
\draw[line width=1pt] (0,0) -- (-1.5,-2);
\draw[line width=1pt] (0,0) -- (1.5,-2);
\draw[line width=1pt] (0, 0) .. controls(-4,-0.1) and (-4,-3.9) .. (0, -4);
\draw[line width=1pt] (0, 0) .. controls(4,-0.1) and (4,-3.9) .. (0, -4);
\draw[line width=1pt] (0, 0) .. controls(-1.5,-0.5) and (-2,-2) .. (-1.75, -2.25);
\draw[line width=1pt] (0, 0) .. controls(0.05,-0.5) and (-0.75,-2.5) .. (-1.75, -2.25)
node[pos=0.65,below] {$m'$};
\draw[line width=1pt] (0, 0) .. controls(0,-1) and (0.7,-2.5) .. (1.75, -2.25)
node[pos=0.8,left] {$m''$};
\draw[line width=1pt] (0, 0) .. controls(0,0) and (2.25,-1.5) .. (1.75, -2.25);


\filldraw [black] (8,0) circle (2pt)
(6.5,-2) circle (2pt)
(9.5,-2) circle (2pt)
(8,0.5) circle (1pt)
(8.3,0.5) circle (1pt)
(7.7,0.5) circle (1pt);
\draw (8,-2)node{$\Delta$};
\draw (7.15,-1.6)node{$\Delta'$};
\draw (6.8,-1.3)node{$j$};
\draw (8.8,-1.5)node{$\Delta''$};
\draw[line width=1pt] (8,0) -- (7, 1)
node[pos=0.3,right] {$p$};
\draw[line width=1pt] (8,0) -- (9, 1);
\draw[line width=1pt] (8,0) -- (6.5,-2);
\draw[line width=1pt] (8,0) -- (9.5,-2);
\draw[line width=1pt] (8, 0) .. controls(4,-0.1) and (4,-3.9) .. (8, -4);
\draw[line width=1pt] (8, 0) .. controls(12,-0.1) and (12,-3.9) .. (8, -4);
\draw[line width=1pt] (8, 0) .. controls(6.5,-0.5) and (6,-2) .. (6.25, -2.25);
\draw[line width=1pt] (8, 0) .. controls(8.05,-0.5) and (7.25,-2.5) .. (6.25, -2.25)
node[pos=0.9,below] {$m'$};
\draw[line width=1pt] (8, 0) .. controls(8,-1) and (8.7,-2.5) .. (9.75, -2.25)
node[pos=0.8,left] {$m''$};
\draw[line width=1pt] (8, 0) .. controls(8,0) and (10.25,-1.5) .. (9.75, -2.25);

\filldraw [red] (-0.75,-1) circle (2pt)
(-1.7,-1.55) circle (2pt)
(1.15,-3.85) circle (2pt)
(0.85,0.9) circle (2pt)
(-0.75,0.8) circle (2pt)
(-0.95,-0.1) circle (2pt)
(-0.8,-0.4) circle (2pt)
(-0.5,-0.7) circle (2pt)
(-0.2,-0.8) circle (2pt)
(0.2,-0.9) circle (2pt)
(0.75,-1) circle (2pt)
(1.4,-1.2) circle (2pt)
(-0.55,-1.55) circle (2pt)
(-1.25,-1.6) circle (2pt);
\draw[color=red][line width=1.5pt] (-1.25,-1.6) .. controls(-5,-1) and (-0.7,-3.9) .. (1.5,-3.8)
node[pos=-0,right] {{\scriptsize $q_{0}$}}
node[pos=0.06,above] {{\scriptsize $r_{1}$}}
node[pos=1,below] {{\scriptsize $r_{2}$}};
\draw[color=red][line width=1.5pt] (1.5,-3.8) .. controls(7,-2.5) and (-1,3.5) .. (-1,0) 
node[pos=0.66,right] {$r_{3}$}
node[pos=0.9,left] {{\scriptsize $r_{l-8}$}}
node[pos=1,left] {{\scriptsize$r_{l-7}$}};
\draw[color=red][line width=1.5pt] (-1,0) .. controls(-0.8,-1) and (0.5,-0.8) .. (1.4,-1.2)
node[pos=1.01,right] {{\scriptsize $r_{l-1}$}};
\draw[color=red][line width=1.5pt] (1.4,-1.2) .. controls(4,-3) and (0,-3.7) .. (-0.75,-1)
node[pos=0.85,above] {{\scriptsize $r_{l}$}}
node[pos=0.99,left] {{\scriptsize $q_{1}$}};

\draw[<-][color=brown][line width=1.5pt] (-1.8,-1.65) ..controls (-4,-2) and (1.1,-3.8) .. (1,-3.7);


\filldraw [red] (7.25,-1) circle (2pt)
(6.55,-1) circle (2pt)
(9.15,-3.85) circle (2pt)
(8.85,0.9) circle (2pt)
(7.25,0.8) circle (2pt)
(7.05,-0.1) circle (2pt)
(7.2,-0.4) circle (2pt)
(7.5,-0.7) circle (2pt)
(7.8,-0.8) circle (2pt)
(8.2,-0.9) circle (2pt)
(8.75,-1) circle (2pt)
(9.4,-1.2) circle (2pt)
(7.05,-2) circle (2pt)
(6.75,-1.6) circle (2pt);
\draw[color=red][line width=1.5pt] (7.25,-1) .. controls(3,-1) and (7.3,-3.9) .. (9.5,-3.8)
node[pos=-0,right] {{\scriptsize $q_{0}$}}
node[pos=0.06,above] {{\scriptsize $r_{1}$}}
node[pos=1,below] {{\scriptsize $r_{2}$}};
\draw[color=red][line width=1.5pt] (9.5,-3.8) .. controls(15,-2.5) and (7,3.5) .. (7,0) 
node[pos=0.66,right] {$r_{3}$}
node[pos=0.9,left] {{\scriptsize $r_{l-8}$}}
node[pos=1,left] {{\scriptsize$r_{l-7}$}};
\draw[color=red][line width=1.5pt] (7,0) .. controls(7.2,-1) and (8.5,-0.8) .. (9.4,-1.2)
node[pos=1.01,right] {{\scriptsize $r_{l-1}$}};
\draw[color=red][line width=1.5pt] (9.4,-1.2) .. controls(12,-3) and (8,-3.7) .. (6.75,-1.6) 
node[pos=0.85,above] {{\scriptsize $r_{l}$}}
node[pos=0.99,left] {{\scriptsize $q_{1}$}};

\draw[<-][color=brown][line width=1.5pt] (7.05,-2.25) -- (9.1,-3.8);


\end{tikzpicture}} 
\caption{Case $b_{2})$ of the Definition \ref{1-desviacion}.}

\end{figure} 
\end{itemize}

\end{definition}

Now, we will assign a label in the set $\{1,2\}$ to the 1-detours $d_{(q_{0},q_{1})}^{\Delta,1}$ this is going to be useful to define the detour matrix. 

we are going to assign the label "1" to the 1-detour $d_{(q_{0},q_{1})}^{\Delta,1}$ which satisfies the conditions of one of the two cases described below.

\begin{itemize}
\item[$a_{1}$)] The arc $j$ is not the folded side of any self-folded triangle and one of the two conditions described below are met.

	\begin{itemize}
		\item[i)]The arc $j$ is a side of a triangle of type 3.
		\item[ii)]$\Delta$ is a triangle of type 2 and the intersection point $r_{l}$ of $\gamma$ with the arcs of $\tau^{\circ}$ belong to the unique arc of $\Delta$ which is a side of a self-folded triangle.
	\end{itemize}  

\item[$a_{2}$)] The arc $j$ is the folded side of a self-folded triangle $\Delta'$ with non folded side $m$.

\begin{itemize}
\item[i)] $\Delta$ is a triangle of type 3, where $m'$ denote the unique side of $\Delta$ which is not a side of any self-folded triangle.

\item[ii)] $r_{2}$ belong to arc $m'$.
\end{itemize}
\end{itemize}

If the conditions of $a_{1}$) and $a_{2}$) are not met then we assign to $d_{(q_{0},q_{1})}^{\Delta,1}$ the label "2".

We call detours of type 1 to the detours with label "1" and the detours of type 2 will be the detours with labels 2.

Let's consider $n\geq 1$, once all the $n$-detours and $n$-auxiliary are drawn for ($\tau^{\circ}$,$i'$), we take an arc $j$ and fix one triangle $\Delta$. Now, we are goint to define sets $\mathcal{B}^{\Delta , n+1}_{i',j}$.

\begin{definition}(Set ${\mathcal B}^{\Delta,n+1}_{i',j}$)\label{conjunto B_{n+1}}
Let $\tau$ be a tagged triangulation of a surface with marked points $(\Sigma,M)$, let's consider a triangle $\Delta$ in $\tau^{\circ}$, $j$ and arc of $\tau^{\circ}$ and $i\notin \tau$ a tagged arc which minimizes the intersection points with the arc of $\tau$. We are going to define the set ${\mathcal B}^{\Delta,n+1}_{i',j}$ with the help of $i'$ and the $n$-detours constructed before. The construcction depends on if the arc $j$ is the folded side of a self-folded triangle, or if it is not:

\begin{itemize}
\item[a)] The arc $j$ is not the folded side of any self-folded triangle.

Under this condition the element $(q_{0},q_{1},s(d^{\Delta',n}),r_{2},r_{l},p)$ belongs to ${\mathcal B}^{\Delta,n+1}_{i',j}$ if $s(d^{\Delta',n})=r_{1}$ does not belong to any self-folded triangle and the conditions of $a$) Definition \ref{conjunto B} are met. The curve $\gamma$ is obtained as the union of the segments $[q_{0},r_{1}]_{i'}$, $d^{\Delta',n}$ y $[r_{2},q_{1}]_{i'}$ and is oriented from $q_{0}$ to $q_{1}$.

\begin{figure}[H]
\centering 
\begin{tikzpicture}[scale=0.75]

\filldraw [black] (-5,0) circle (2pt)
(-5,-4) circle (2pt)
(-8.5,-2) circle (2pt)
(-9.5,-2) circle (1pt)
(-9.5,-2.5) circle (1pt)
(-9.5,-1.5) circle (1pt);
\draw (-6,-2)node{$\Delta$};
\draw[line width=1pt] (-5, 0) -- (-5, -4)
node[pos=0.5,left] {$j$};
\draw[line width=1pt] (-5, 0) -- (-8.5, -2);
\draw[line width=1pt] (-5, -4) -- (-8.5, -2);
\draw[line width=1pt] (-8.5, -2) -- (-8.5,0);
\draw[line width=1pt] (-8.5, -2) -- (-10.5,-0.5);
\draw[line width=1pt] (-8.5, -2) -- (-10.5,-3.5);
\draw[line width=1pt] (-8.5, -2) -- (-8.5,-4)
node[pos=0.15,left] {$p$};

\filldraw [red] 
(-5,-3) circle (2pt)
(-6.35,-3.2) circle (2pt)
(-8.5,-3.3) circle (2pt)
(-9.8,-3) circle (2pt)
(-9.8,-1) circle (2pt)
(-8.5,-0.45) circle (2pt)
(-6.05,-0.6) circle (2pt)
(-5,-1) circle (2pt);
\draw[color=red][line width=1.5pt] (-5,-1) -- (-6.05,-0.6) 
node[pos=0,right] {$q^{0}$}
node[pos=1,above] {$r_{1}$};

\draw[->][color=brown][line width=1.5pt] (-6.15,-0.6) ..controls(-6.65,-0.5) and (-7.9,-0.35) .. (-8.4,-0.45)
node[pos=0.5,above] {$d^{\Delta',n}$};

\draw[color=red][line width=1.5pt] (-8.5,-0.45) .. controls (-8.5,-0.45) and (-10.5,-1) .. (-10.5,-2)
node[pos=0.2,above] {$r_{2}$};
\draw[color=red][line width=1.5pt] (-5,-3) .. controls(-8, -3.5) and (-10.5,-3.5) .. (-10.5,-2)
node[pos=0,right] {$q_{1}$}
node[pos=0.15,below] {$r_{l}$};

\filldraw [black] (0,0) circle (2pt)
(0,-4) circle (2pt)
(0,-2) circle (2pt)
(0,0.5) circle (1pt)
(0.3,0.5) circle (1pt)
(-0.3,0.5) circle (1pt)
(0,-4.5) circle (1pt)
(0.3,-4.5) circle (1pt)
(-0.3,-4.5) circle (1pt);
\draw (0,-3)node{$\Delta$};

\draw[line width=1pt] (0,0) -- (-1, 1)
node[pos=0.3,right] {$p$};
\draw[line width=1pt] (0,0) -- (1, 1);
\draw[line width=1pt] (0,-4) -- (1, -5);
\draw[line width=1pt] (0,-4) -- (-1, -5);
\draw[line width=1pt] (0,0) -- (0, -2);
\draw[line width=1pt] (0,0) .. controls(-0.7,-1) and (-0.7,-2.2) .. (0, -2.5)
node[pos=0.5,left] {$m$};
\draw[line width=1pt] (0,0) .. controls(0.7,-1) and (0.7,-2.2) .. (0, -2.5);
\draw[line width=1pt] (0,0) .. controls(-2, -0.5) and (-2,-3.5) .. (0, -4)
node[pos=0.5,right] {$j$};
\draw[line width=1pt] (0,0) .. controls(2, -0.5) and (2,-3.5) .. (0, -4);


\filldraw [red] (-0.55,-3.8) circle (2pt)
(-1.25,-3) circle (2pt)
(1.3,-1.05) circle (2pt)
(0.7,0.7) circle (2pt)
(-0.75,0.75) circle (2pt)
(-0.95,-0.5) circle (2pt)
(-0.4,-0.75) circle (2pt)
(0,-1.05) circle (2pt)
(0,-1.05) circle (2pt)
(0.35,-2.1) circle (2pt);

\draw[color=red][line width=1.5pt] (-0.55,-3.8) .. controls(1,-3.5) and (1.5,-0.75) .. (1.3,-1.05) 
node[pos=0,left] {$q_{0}$}
node[pos=0.9,right] {$r_{1}$};
\draw[->][color=brown][line width=1.5pt] (1.3,-0.95) ..controls(1.5,-0.8) and (0.9,0.5) .. (0.8,0.65)
node[pos=0.5,right] {$d^{\Delta',n}$};

\draw[color=red][line width=1.5pt] (0.7,0.7) .. controls(0.25,1.25) and (-2.25,0.75) .. (-1,-0.5)
node[pos=0.1,above] {$r_{2}$};
\draw[color=red][line width=1.5pt] (-1.25,-3) .. controls(1,-3) and (0.75,-1) .. (-1,-0.5)
node[pos=0.3,right] {$r_{l}$}
node[pos=0,left] {$q_{1}$};

\end{tikzpicture}
\caption{Case $a$) of the Definition \ref{conjunto B_{n+1}}.}
\end{figure} 

\item[b)]The arc $j$ is the folded side of a self-folded triangle $\Delta'$. We denote by $m$ the non folded side of $\Delta'$.

Under these conditions the element $(q_{0},q_{1},q_{2},r_{1},s(d^{\Delta',n}),r_{3},r_{l},p)$ belongs to ${\mathcal B}^{\Delta,n+1}_{i',j}$  if  the conditions of $b_{21}$) Definition \ref{conjunto B} are met. The curve $\gamma$ is obtained as the union of the segments $[q_{0},r_{2}]_{i'}$, $d^{\Delta',n}$ and $[r_{3},q_{1}]_{i'}$, and is oriented from $q_{0}$ to $q_{1}$.

\begin{figure}[H]
\centering 
$$ \begin{tikzpicture}[scale=0.75]

\filldraw [black] 
(4,0) circle (2pt)
(4,-4) circle (2pt)
(4,-2) circle (2pt)
(4,0.5) circle (1pt)
(4.3,0.5) circle (1pt)
(3.7,0.5) circle (1pt)
(4,-4.5) circle (1pt)
(4.3,-4.5) circle (1pt)
(3.7,-4.5) circle (1pt);
\draw (4,-4.3)node{$p$}; 
\draw[line width=1pt] (4,0) -- (3, 1);
\draw[line width=1pt] (4,0) -- (5, 1);
\draw[line width=1pt] (4,-4) -- (5, -5);
\draw[line width=1pt] (4,-4) -- (3, -5);
\draw[line width=1pt] (4,0) -- (4, -2);
\draw[line width=1pt] (4,0) .. controls(2, -0.5) and (2,-3.5) .. (4, -4);
\draw[line width=1pt] (4,0) .. controls(6, -0.5) and (6,-3.5) .. (4, -4);
\draw[line width=1pt] (4,0) .. controls(3.3,-1) and (3.3,-2.2) .. (4, -2.5);
\draw[line width=1pt] (4,0) .. controls(4.7,-1) and (4.7,-2.2) .. (4, -2.5)
node[pos=0.5,right] {$m$};



\filldraw [red] (4,-0.5) circle (2pt)
(3.55,-0.85) circle (2pt)
(2.77,-3) circle (2pt)
(3.45,-4.55) circle (2pt)
(4.75,-4.75) circle (2pt)
(4.9,-3.5) circle (2pt)
(3.55,-1.9) circle (2pt)
(4,-1.5) circle (2pt) 
(4.5,-1.3) circle (2pt) 
(5.35,-1.2) circle (2pt);

\draw[color=red][line width=1.5pt] (4,-0.5) .. controls(3.25,-1) and (2.75,-2) .. (2.77,-3)
node[pos=0,right] {$q_{0}$}
node[pos=0.2,left] {$r_{1}$}
node[pos=1,left] {$r_{2}$};
\draw[->][color=brown][line width=1.5pt] (2.77,-3.1) ..controls (2.67,-3.6) and (3,-4.05).. (3.4,-4.45)
node[pos=0.5,left] {$d^{\Delta',n}$};
\draw[color=red][line width=1.5pt] (3.45,-4.55) .. controls(3.55,-4.65) and (5.5,-5.5) .. (4.9,-3.5)
node[pos=0.1,below] {$r_{3}$};
\draw[color=red][line width=1.5pt] (4.9,-3.5) .. controls(4.8,-2.8) and (2.5,-2.5) .. (4,-1.5)
node[pos=0.8,right] {$r_{l}$}
node[pos=1.05,below] {$q_{1}$};
\draw[color=red][line width=1.5pt] (4,-1.5) ..controls(4.25,-1.3) and (5.25,-1.2) .. (5.35,-1.2)
node[pos=1,right] {$q_{2}$};

\end{tikzpicture}$$
\caption{Case $b)$ of the Definition \ref{conjunto B_{n+1}}.}

\end{figure} 
\end{itemize}

\end{definition}

Using the sets ${\mathcal B}^{\Delta,n+1}_{i',j}$ we will define the $(n+1)$-detours and the $(n+1)$-auxiliaries.

\begin{definition}\label{n+1-desviacion}($(n+1)$-detour) For each element $(q_{0},q_{1},s(d^{\Delta',n}),r_{2},r_{l},p)$ \'o $(q_{0},q_{1},q_{2},r_{1},s(d^{\Delta',n}),r_{3},r_{l},p)$ which belong to ${\mathcal B}^{\Delta,n+1}_{i',j}$, we will construct the $(n+1)$-\textit{detour} $d_{(q_{0},q_{1})}^{\Delta,n+1}$ which is a simple oriented curve with ends $s(d_{(q_{0},q_{1})}^{\Delta,n+1})$ and $t(d_{(q_{0},q_{1})}^{\Delta,n+1})$, which start and finish respectively. Only in the case $b$) we draw a curve ($(n+1)$-auxiliar) $e_{(q_{0},q_{1})}^{\Delta,n+1}$ with ends $s(e_{(q_{0},q_{1})}^{\Delta,n+1})$ and $t(e_{(q_{0},q_{1})}^{\Delta,n+1})$. We will consider the cases when arc $j$ is the folded side, or if it is not:

\begin{myEnumerate}

\item[a)] $j$ is not the folded side of any self-folded triangle.

For each element $(q_{0},q_{1},s(d^{\Delta',n}),r_{2},r_{l},p)$ in ${\mathcal B}^{\Delta,n+1}_{i',j}$ we draw the $(n+1)$-detour $d_{(q_{0},q_{1})}^{\Delta,n+1}$ contained in $\Delta$ with ends $s(d_{(q_{0},q_{1})}^{\Delta,n+1})=s(d_{(q_{0},q_{1})}^{\Delta',n})=r_{1}$ and $t(d_{(q_{0},q_{1})}^{\Delta,n+1})=q_{1}$.

\begin{figure}[H]
\centering 
\begin{tikzpicture}[scale=0.75]

\filldraw [black] (-5,0) circle (2pt)
(-5,-4) circle (2pt)
(-8.5,-2) circle (2pt)
(-9.5,-2) circle (1pt)
(-9.5,-2.5) circle (1pt)
(-9.5,-1.5) circle (1pt);
\draw (-6,-2)node{$\Delta$};
\draw[line width=1pt] (-5, 0) -- (-5, -4)
node[pos=0.5,left] {$j$};
\draw[line width=1pt] (-5, 0) -- (-8.5, -2);
\draw[line width=1pt] (-5, -4) -- (-8.5, -2);
\draw[line width=1pt] (-8.5, -2) -- (-8.5,0);
\draw[line width=1pt] (-8.5, -2) -- (-10.5,-0.5);
\draw[line width=1pt] (-8.5, -2) -- (-10.5,-3.5);
\draw[line width=1pt] (-8.5, -2) -- (-8.5,-4)
node[pos=0.15,left] {$p$};

\filldraw [red] 
(-5,-3) circle (2pt)
(-6.35,-3.2) circle (2pt)
(-8.5,-3.3) circle (2pt)
(-9.8,-3) circle (2pt)
(-9.8,-1) circle (2pt)
(-8.5,-0.45) circle (2pt)
(-6.05,-0.6) circle (2pt)
(-5,-1) circle (2pt);
\draw[->][color=brown][line width=1.5pt] (-6.05,-0.7) -- (-5.1,-2.9)
node[pos=0.3,left] {$d_{q_{0},q_{1}}^{\Delta , n+1}$};
\draw[color=red][line width=1.5pt] (-5,-1) -- (-6.05,-0.6) 
node[pos=0,right] {$q_{0}$}
node[pos=1,above] {$r_{1}$};

\draw[->][color=brown][line width=1.5pt] (-6.15,-0.6) ..controls(-6.65,-0.5) and (-7.9,-0.35) .. (-8.4,-0.45)
node[pos=0.5,above] {$d^{\Delta',n}$};

\draw[color=red][line width=1.5pt] (-8.5,-0.45) .. controls (-8.5,-0.45) and (-10.5,-1) .. (-10.5,-2)
node[pos=0.2,above] {$r_{2}$};
\draw[color=red][line width=1.5pt] (-5,-3) .. controls(-8, -3.5) and (-10.5,-3.5) .. (-10.5,-2)
node[pos=0,right] {$q_{1}$}
node[pos=0.15,below] {$r_{l}$};

\filldraw [black] (0,0) circle (2pt)
(0,-4) circle (2pt)
(0,-2) circle (2pt)
(0,0.5) circle (1pt)
(0.3,0.5) circle (1pt)
(-0.3,0.5) circle (1pt)
(0,-4.5) circle (1pt)
(0.3,-4.5) circle (1pt)
(-0.3,-4.5) circle (1pt);
\draw (-0.5,-2.5)node{$\Delta$};

\draw[line width=1pt] (0,0) -- (-1, 1)
node[pos=0.3,right] {$p$};
\draw[line width=1pt] (0,0) -- (1, 1);
\draw[line width=1pt] (0,-4) -- (1, -5);
\draw[line width=1pt] (0,-4) -- (-1, -5);
\draw[line width=1pt] (0,0) -- (0, -2);
\draw[line width=1pt] (0,0) .. controls(-0.7,-1) and (-0.7,-2.2) .. (0, -2.5)
node[pos=0.5,left] {$m$};
\draw[line width=1pt] (0,0) .. controls(0.7,-1) and (0.7,-2.2) .. (0, -2.5);
\draw[line width=1pt] (0,0) .. controls(-2, -0.5) and (-2,-3.5) .. (0, -4)
node[pos=0.5,right] {$j$};
\draw[line width=1pt] (0,0) .. controls(2, -0.5) and (2,-3.5) .. (0, -4);


\filldraw [red] (-0.55,-3.8) circle (2pt)
(-1.25,-3) circle (2pt)
(1.3,-1.05) circle (2pt)
(0.7,0.7) circle (2pt)
(-0.75,0.75) circle (2pt)
(-0.95,-0.5) circle (2pt)
(-0.4,-0.75) circle (2pt)
(0,-1.05) circle (2pt)
(0,-1.05) circle (2pt)
(0.35,-2.1) circle (2pt);
\draw[->][color=brown][line width=1.5pt] (1.25,-1.1) ..controls(1,-2.5) and (-0.25,-4.5) .. (-1.15,-3.1)
node[pos=0.7,above] {$d_{q_{0},q_{1}}^{\Delta , n+1}$}; 

\draw[color=red][line width=1.5pt] (-0.55,-3.8) .. controls(1,-3.5) and (1.5,-0.75) .. (1.3,-1.05) 
node[pos=0,left] {$q_{0}$}
node[pos=0.9,right] {$r_{1}$};
\draw[->][color=brown][line width=1.5pt] (1.3,-0.95) ..controls(1.5,-0.8) and (0.9,0.5) .. (0.8,0.65)
node[pos=0.5,right] {$d^{\Delta',n}$};

\draw[color=red][line width=1.5pt] (0.7,0.7) .. controls(0.25,1.25) and (-2.25,0.75) .. (-1,-0.5)
node[pos=0.1,above] {$r_{2}$};
\draw[color=red][line width=1.5pt] (-1.25,-3) .. controls(1,-3) and (0.75,-1) .. (-1,-0.5)
node[pos=0.3,right] {$r_{l}$}
node[pos=0,left] {$q_{1}$};

\end{tikzpicture}
\caption{Case $a$) of the Definition \ref{n+1-desviacion}.}
\end{figure} 
\item[b)] The arc $j$ is the folded side of a self-folded triangle $\Delta'$. We will denote by $m$ the non folded side of $\Delta'$.
For each element $(q_{0},q_{1},q_{2},r_{1},s(d^{\Delta',n}),r_{3},r_{l},p)$ in ${\mathcal B}^{\Delta,n+1}_{i',j}$ we will draw the $(n+1)$-detour $d_{(q_{0},q_{1})}^{\Delta,n+1}$ with ends $s(d_{(q_{0},q_{1})}^{\Delta,n+1})=s(d_{(q_{0},q_{1})}^{\Delta',n})$ and $t(d_{(q_{0},q_{1})}^{\Delta,n+1})=q_{1}$. In this case we also draw the $(n+1)$ auxiliary $e_{(q_{0},q_{1})}^{\Delta,n+1}$ with ends $s(e_{(q_{0},q_{1})}^{\Delta,n+1})=r_{1}$ and $t(e_{(q_{0},q_{1})}^{\Delta,n+1})=q_{2}$.

\begin{figure}[H]
\centering 
$$ \begin{tikzpicture}[scale=0.85]

\filldraw [black] 
(4,0) circle (2pt)
(4,-4) circle (2pt)
(4,-2) circle (2pt)
(4,0.5) circle (1pt)
(4.3,0.5) circle (1pt)
(3.7,0.5) circle (1pt)
(4,-4.5) circle (1pt)
(4.3,-4.5) circle (1pt)
(3.7,-4.5) circle (1pt);
\draw (4,-4.3)node{$p$}; 
\draw[line width=1pt] (4,0) -- (3, 1);
\draw[line width=1pt] (4,0) -- (5, 1);
\draw[line width=1pt] (4,-4) -- (5, -5);
\draw[line width=1pt] (4,-4) -- (3, -5);
\draw[line width=1pt] (4,0) -- (4, -2);
\draw[line width=1pt] (4,0) .. controls(2, -0.5) and (2,-3.5) .. (4, -4);
\draw[line width=1pt] (4,0) .. controls(6, -0.5) and (6,-3.5) .. (4, -4);
\draw[line width=1pt] (4,0) .. controls(3.3,-1) and (3.3,-2.2) .. (4, -2.5);
\draw[line width=1pt] (4,0) .. controls(4.7,-1) and (4.7,-2.2) .. (4, -2.5);



\filldraw [red] (4,-0.5) circle (2pt)
(3.55,-0.85) circle (2pt)
(2.77,-3) circle (2pt)
(3.45,-4.55) circle (2pt)
(4.75,-4.75) circle (2pt)
(4.9,-3.5) circle (2pt)
(3.55,-1.9) circle (2pt)
(4,-1.5) circle (2pt) 
(4.5,-1.3) circle (2pt) 
(5.35,-1.2) circle (2pt);
\draw[->][color=brown][line width=1.5pt] (2.82,-2.9) ..controls (3.1,-2) and (3.5,-1) .. (3.95,-1.5)
node[pos=0.05,right] {$d_{q_{0},q_{1}}^{\Delta ,n+ 1}$};
\draw[->][color=blue][line width=1.5pt] (3.6,-0.9) -- (5.15,-1.1)
node[pos=0.95,above] {$e_{q_{0},q_{1}}^{\Delta ,n+ 1}$};

\draw[color=red][line width=1.5pt] (4,-0.5) .. controls(3.25,-1) and (2.75,-2) .. (2.77,-3)
node[pos=0,right] {$q_{0}$}
node[pos=0.2,left] {$r_{1}$}
node[pos=1,left] {$r_{2}$};
\draw[->][color=brown][line width=1.5pt] (2.77,-3.1) ..controls (2.67,-3.6) and (3,-4.05).. (3.4,-4.45)
node[pos=0.5,left] {$d^{\Delta',n}$};
\draw[color=red][line width=1.5pt] (3.45,-4.55) .. controls(3.55,-4.65) and (5.5,-5.5) .. (4.9,-3.5)
node[pos=0.1,below] {$r_{3}$};
\draw[color=red][line width=1.5pt] (4.9,-3.5) .. controls(4.8,-2.8) and (2.5,-2.5) .. (4,-1.5)
node[pos=0.8,right] {$r_{l}$}
node[pos=1.05,below] {$q_{1}$};
\draw[color=red][line width=1.5pt] (4,-1.5) ..controls(4.25,-1.3) and (5.25,-1.2) .. (5.35,-1.2)
node[pos=1,right] {$q_{2}$};

\end{tikzpicture}$$
\caption{Case $b)$ of the Definition \ref{n+1-desviacion}.}

\end{figure}

\end{myEnumerate}

\end{definition}

Let $d_{(q_{0},q_{1})}^{\Delta,m}$ be a $m$-detour with $1 \leq m \leq n$, we assign the label "1" to $d_{(q_{0},q_{1})}^{\Delta,m}$ if the label of $d_{(q_{0},q_{1})}^{\Delta,1}$ is $"1"$. Otherwise we assign the label "2" to $d_{(q_{0},q_{1})}^{\Delta,m}$.

Now, we are going to analyze the arcs $j$ of $\tau^{\circ}$ to determine the number of detour matrixes and auxiliary matrixes which it will assign to the arc $j$.

\begin{myEnumerate}
\item[a)] The arc $j$ is a side of a self-folded triangle $\Delta'$, let's denote by $m$ the non folded side of $\Delta'$.

If the arc $m$ is a side of a triangle $\Delta$ of type 3 with label $l_{3}$ in $\Delta$, we assign to arcs $j$ and $m$ two detour matrixes and two auxiliary matrixes. Otherwise, we assign only one detour matrix and one auxiliary matrix.

\item[b)] The arc $j$ is not a side of any self-folded triangle.
Under this condition, the arc $j$ is a side of exactly two triangles which are not self-folded triangles. We will denote by $\Delta_{1}$ and $\Delta_{2}$ the not self-folded triangles which share the arc $j$.

\begin{myEnumerate}
\item[$b_{1}$)] If the conditions of one of the two subcases described below are met, we assign two detour matrixes and two auxiliary matrixes to arc $j$.

\begin{myEnumerate}
\item[$b_{11}$)] $\Delta_{1}$ and $\Delta_{2}$ are triangles of type 1.
\item[$b_{12}$)] $\Delta_{1}$ is of type 1, $\Delta_{2}$ is of type 2 and $j$ has a label $l_{3}$ in $\Delta_{2}$.

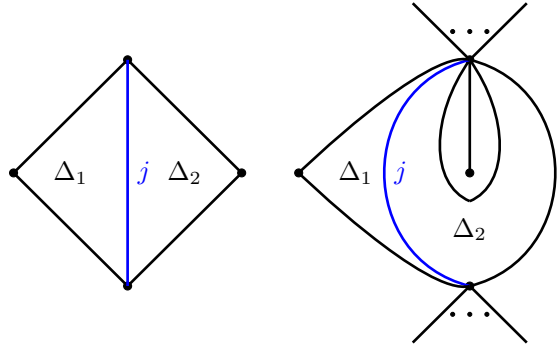
\begin{figure}[H]
\centering 
\begin{tikzpicture}[scale=0.75]

\filldraw [black] (-6,0) circle (2pt)
(-8, -2) circle (2pt)
(-4,-2) circle (2pt)
(-6,-4) circle (2pt);
\draw (-5,-2)node{$\Delta_{2}$};

\draw (-7,-2)node{$\Delta_{1}$};

\draw[line width=1pt] (-6,0) -- (-8, -2);
\draw[line width=1pt] (-4,-2) -- (-6,0) ;

\draw[line width=1pt] (-6,-4) -- (-8, -2);
\draw[line width=1pt] (-4,-2) -- (-6,-4) ;

\draw[color=blue][line width=1pt] (-6,0) -- (-6, -4)
node[pos=0.5,right] {$j$};


\filldraw [black] (0,0) circle (2pt)
(-3,-2) circle (2pt)
(0,-4) circle (2pt)
(0,-2) circle (2pt)
(0,0.5) circle (1pt)
(0.3,0.5) circle (1pt)
(-0.3,0.5) circle (1pt)
(0,-4.5) circle (1pt)
(0.3,-4.5) circle (1pt)
(-0.3,-4.5) circle (1pt);
\draw (0,-3)node{$\Delta_{2}$};
\draw (-2,-2)node{$\Delta_{1}$};

\draw[line width=1pt] (0,0) -- (-1, 1);
\draw[line width=1pt] (0,0) -- (1, 1);
\draw[line width=1pt] (0,-4) -- (1, -5);
\draw[line width=1pt] (0,-4) -- (-1, -5);
\draw[line width=1pt] (0,0) -- (0, -2);
\draw[line width=1pt] (0,0) .. controls(-0.7,-1) and (-0.7,-2.2) .. (0, -2.5);
\draw[line width=1pt] (0,0) .. controls(0.7,-1) and (0.7,-2.2) .. (0, -2.5);
\draw[color=blue][line width=1pt] (0,0) .. controls(-2, -0.5) and (-2,-3.5) .. (0, -4)
node[pos=0.5,right] {$j$};
\draw[line width=1pt] (0,0) .. controls(2, -0.5) and (2,-3.5) .. (0, -4);

\draw[line width=1pt] (0,0) ..controls (-0.75,0.2 ) and (-3,-2) .. (-3, -2);

\draw[line width=1pt] (0,-4) ..controls (-0.75,-4.2 ) and (-3,-2) .. (-3, -2);

\end{tikzpicture}
\end{figure} 
\end{myEnumerate}

\item[$b_{2}$)] If the conditions of one of the four subcases described below are met, we assign three detour matrixes and three auxiliaries matrixes to the arc $j$.

\begin{myEnumerate}
\item[$b_{21}$)] $\Delta_{1}$ is of type 3 y $\Delta_{2}$ is of type 1.
\item[$b_{22}$)] $\Delta_{1}$ is of type 2, $\Delta_{2}$ is of type and $j$ has a label $l_{2}$ in $\Delta_{1}$. 
\item[$b_{23}$)] $\Delta_{1}$ and $\Delta_{2}$ are of type 2, $j$ has a label $l_{2}$ in $\Delta_{1}$ and $l_{3}$ in $\Delta_{2}$.
\item[$b_{24}$)] $\Delta_{1}$ is of type 3 and $\Delta_{2}$ is of type 2, $j$ has a label $l_{1}$ in $\Delta_{1}$ and $l_{3}$ in $\Delta_{2}$.
\end{myEnumerate}
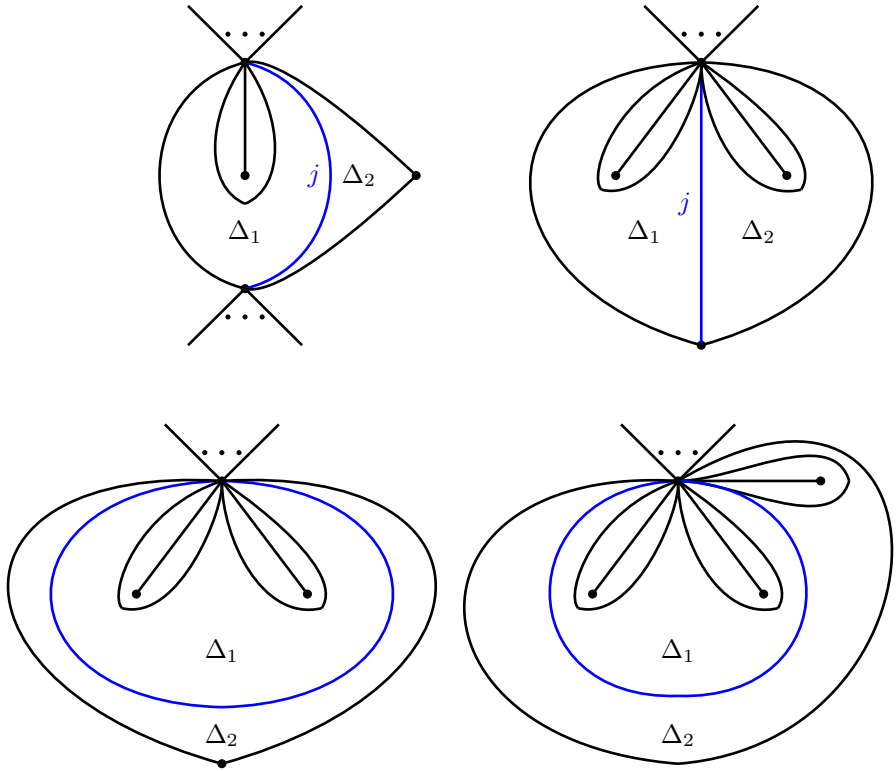
\begin{figure}[H]
\centering 
\begin{tikzpicture}[scale=0.75]

\filldraw [black] (8,0) circle (2pt)
(6.5,-2) circle (2pt)
(9.5,-2) circle (2pt)
(8, -5) circle (2pt)
(8,0.5) circle (1pt)
(8.3,0.5) circle (1pt)
(7.7,0.5) circle (1pt);
\draw (7,-3)node{$\Delta_{1}$};

\draw (9,-3)node{$\Delta_{2}$};
\draw[line width=1pt] (8,0) -- (7, 1);

\draw[color=blue][line width=1pt] (8,0) -- (8, -5)
node[pos=0.5,left] {$j$};
\draw[line width=1pt] (8,0) -- (9, 1);
\draw[line width=1pt] (8,0) -- (6.5,-2);
\draw[line width=1pt] (8,0) -- (9.5,-2);
\draw[line width=1pt] (8, 0) .. controls(4,-0.1) and (4,-3.9) .. (8, -5);
\draw[line width=1pt] (8, 0) .. controls(12,-0.1) and (12,-3.9) .. (8, -5);
\draw[line width=1pt] (8, 0) .. controls(6.5,-0.5) and (6,-2) .. (6.25, -2.25);
\draw[line width=1pt] (8, 0) .. controls(8.05,-0.5) and (7.25,-2.5) .. (6.25, -2.25);
\draw[line width=1pt] (8, 0) .. controls(8,-1) and (8.7,-2.5) .. (9.75, -2.25);
\draw[line width=1pt] (8, 0) .. controls(8,0) and (10.25,-1.5) .. (9.75, -2.25);


\filldraw [black] (0,0) circle (2pt)
(3,-2) circle (2pt)
(0,-4) circle (2pt)
(0,-2) circle (2pt)
(0,0.5) circle (1pt)
(0.3,0.5) circle (1pt)
(-0.3,0.5) circle (1pt)
(0,-4.5) circle (1pt)
(0.3,-4.5) circle (1pt)
(-0.3,-4.5) circle (1pt);
\draw (0,-3)node{$\Delta_{1}$};
\draw (2,-2)node{$\Delta_{2}$};

\draw[line width=1pt] (0,0) -- (-1, 1);
\draw[line width=1pt] (0,0) -- (1, 1);
\draw[line width=1pt] (0,-4) -- (1, -5);
\draw[line width=1pt] (0,-4) -- (-1, -5);
\draw[line width=1pt] (0,0) -- (0, -2);
\draw[line width=1pt] (0,0) .. controls(-0.7,-1) and (-0.7,-2.2) .. (0, -2.5);
\draw[line width=1pt] (0,0) .. controls(0.7,-1) and (0.7,-2.2) .. (0, -2.5);
\draw[line width=1pt] (0,0) .. controls(-2, -0.5) and (-2,-3.5) .. (0, -4);
\draw[color=blue][line width=1pt] (0,0) .. controls(2, -0.5) and (2,-3.5) .. (0, -4)
node[pos=0.5,left] {$j$};

\draw[line width=1pt] (0,0) ..controls (0.75,0.2 ) and (3,-2) .. (3, -2);

\draw[line width=1pt] (0,-4) ..controls (0.75,-4.2 ) and (3,-2) .. (3, -2);

\end{tikzpicture}
\end{figure} 
\begin{figure}[H]
\centering 
\vspace{-4em}\hspace*{-6em}{\begin{tikzpicture}[scale=0.75]

\filldraw [black] (0,0) circle (2pt)
(-1.5,-2) circle (2pt)
(1.5,-2) circle (2pt)
(0, -5) circle (2pt)
(0,0.5) circle (1pt)
(0.3,0.5) circle (1pt)
(-0.3,0.5) circle (1pt);
\draw (0,-3)node{$\Delta_{1}$};
\draw (0,-4.5)node{$\Delta_{2}$};
\draw[line width=1pt] (0,0) -- (-1, 1);
\draw[line width=1pt] (0,0) -- (1, 1);
\draw[line width=1pt] (0,0) -- (-1.5,-2);
\draw[line width=1pt] (0,0) -- (1.5,-2);
\draw[color=blue][line width=1pt] (0, 0) .. controls(-4,-0.1) and (-4,-3.9) .. (0, -4);
\draw[color=blue][line width=1pt] (0, 0) .. controls(4,-0.1) and (4,-3.9) .. (0, -4);
\draw[line width=1pt] (0, 0) .. controls(-1.5,-0.5) and (-2,-2) .. (-1.75, -2.25);
\draw[line width=1pt] (0, 0) .. controls(0.05,-0.5) and (-0.75,-2.5) .. (-1.75, -2.25);
\draw[line width=1pt] (0, 0) .. controls(0,-1) and (0.7,-2.5) .. (1.75, -2.25);
\draw[line width=1pt] (0, 0) .. controls(0,0) and (2.25,-1.5) .. (1.75, -2.25);

\draw[line width=1pt] (0, 0) .. controls(-5,0.3) and (-5,-3.6) .. (0, -5);
\draw[line width=1pt] (0, 0) .. controls(5,0.3) and (5,-3.6) .. (0, -5);


\filldraw [black] (8,0) circle (2pt)
(6.5,-2) circle (2pt)
(9.5,-2) circle (2pt)
(10.5,0) circle (2pt)
(8,0.5) circle (1pt)
(8.3,0.5) circle (1pt)
(7.7,0.5) circle (1pt);
\draw (8,-3)node{$\Delta_{1}$};
\draw (8,-4.5)node{$\Delta_{2}$};
\draw[line width=1pt] (8,0) -- (7, 1); 
\draw[line width=1pt] (8,0) -- (9, 1);

\draw[line width=1pt] (8,0) -- (6.5,-2);
\draw[line width=1pt] (8,0) -- (9.5,-2);

\draw[line width=1pt] (8,0) -- (10.5,0);
\draw[line width=1pt] (8,0) ..controls (9,0) and (10.7,-1) .. (11,0);
\draw[line width=1pt] (8,0) ..controls (9,0) and (10.7,1) .. (11,0);

\draw[color=blue][line width=1pt] (8, 0) .. controls(5,-0.1) and (5,-3.9) .. (8, -3.8); 
\draw[color=blue][line width=1pt] (8, 0) .. controls(11,-0.1) and (11,-3.9) .. (8, -3.8);

\draw[line width=1pt] (8, 0) .. controls(6.5,-0.5) and (6,-2) .. (6.25, -2.25);
\draw[line width=1pt] (8, 0) .. controls(8.05,-0.5) and (7.25,-2.5) .. (6.25, -2.25);
\draw[line width=1pt] (8, 0) .. controls(8,-1) and (8.7,-2.5) .. (9.75, -2.25);
\draw[line width=1pt] (8, 0) .. controls(8,0) and (10.25,-1.5) .. (9.75, -2.25);

\draw[line width=1pt] (8, 0) .. controls(3,0.3) and (3,-4.6) .. (8, -5);
\draw[line width=1pt] (8, 0) .. controls(13,3) and (13,-4.6) .. (8, -5);

\end{tikzpicture}} 

\end{figure}

\item[$b_{3}$)] If the conditions of one of the two subcases described below are met, we assign four detour matrixes and four auxiliary matrixes to arc $j$.
\begin{myEnumerate}
\item[$b_{31}$)] $\Delta_{1}$ and $\Delta_{2}$ are triangles of type 2, and the arc $j$ has a label $l_{2}$ at both triangles.
\item[$b_{32}$)] $\Delta_{1}$ is  triangle of type 3, $\Delta_{2}$ is  tringle of type 2, $j$ has a label $l_{1}$ in $\Delta_{1}$ and a label $l_{2}$ in $\Delta_{2}$.
\end{myEnumerate}

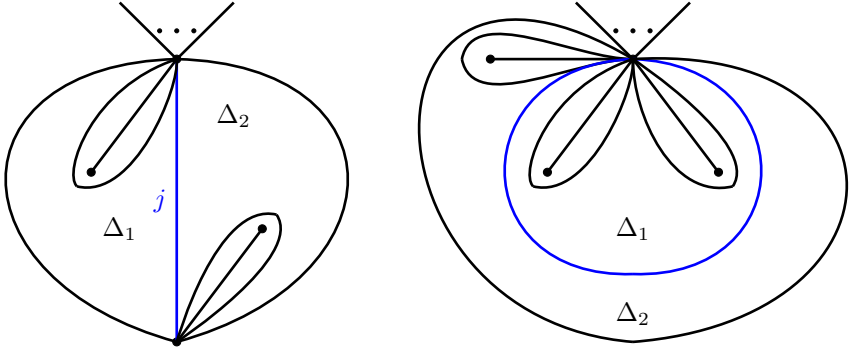
\begin{figure}[H]
\centering 
\vspace{-4em}\hspace*{-6em} \begin{tikzpicture}[scale=0.75]

\filldraw [black] (0,0) circle (2pt)
(-1.5,-2) circle (2pt)
(1.5,-3) circle (2pt)
(0, -5) circle (2pt)
(0,0.5) circle (1pt)
(0.3,0.5) circle (1pt)
(-0.3,0.5) circle (1pt);
\draw (-1,-3)node{$\Delta_{1}$};
\draw (1,-1)node{$\Delta_{2}$};
\draw[line width=1pt] (0,0) -- (-1, 1);
\draw[line width=1pt] (0,0) -- (1, 1); 

\draw[color=blue][line width=1pt] (0,0) -- (0, -5)
node[pos=0.5,left] {$j$};

\draw[line width=1pt] (0,0) -- (-1.5,-2);
\draw[line width=1pt] (0,-5) -- (1.5,-3);

\draw[line width=1pt] (0, 0) .. controls(-4,-0.1) and (-4,-3.9) .. (0, -5); 
\draw[line width=1pt] (0, 0) .. controls(4,-0.1) and (4,-3.9) .. (0, -5);

\draw[line width=1pt] (0, 0) .. controls(-1.5,-0.5) and (-2,-2) .. (-1.75, -2.25);
\draw[line width=1pt] (0, 0) .. controls(0.05,-0.5) and (-0.75,-2.5) .. (-1.75, -2.25);
\draw[line width=1pt] (0, -5) .. controls(0,-5) and (0.7,-2.5) .. (1.75, -2.75); 
\draw[line width=1pt] (0, -5) .. controls(0,-5) and (2.25,-3.5) .. (1.75, -2.75);


\filldraw [black] (8,0) circle (2pt)
(6.5,-2) circle (2pt)
(9.5,-2) circle (2pt)
(5.5,0) circle (2pt)
(8,0.5) circle (1pt)
(8.3,0.5) circle (1pt)
(7.7,0.5) circle (1pt);
\draw (8,-3)node{$\Delta_{1}$};
\draw (8,-4.5)node{$\Delta_{2}$};
\draw[line width=1pt] (8,0) -- (7, 1); 
\draw[line width=1pt] (8,0) -- (9, 1);

\draw[line width=1pt] (8,0) -- (6.5,-2);
\draw[line width=1pt] (8,0) -- (9.5,-2);

\draw[line width=1pt] (8,0) -- (5.5,0);
\draw[line width=1pt] (8,0) ..controls (7,0) and (5.2,-1) .. (5,0);
\draw[line width=1pt] (8,0) ..controls (7,0) and (5.2,1) .. (5,0);

\draw[color=blue][line width=1pt] (8, 0) .. controls(5,-0.1) and (5,-3.9) .. (8, -3.8); 
\draw[color=blue][line width=1pt] (8, 0) .. controls(11,-0.1) and (11,-3.9) .. (8, -3.8);

\draw[line width=1pt] (8, 0) .. controls(6.5,-0.5) and (6,-2) .. (6.25, -2.25);
\draw[line width=1pt] (8, 0) .. controls(8.05,-0.5) and (7.25,-2.5) .. (6.25, -2.25);
\draw[line width=1pt] (8, 0) .. controls(8,-1) and (8.7,-2.5) .. (9.75, -2.25);
\draw[line width=1pt] (8, 0) .. controls(8,0) and (10.25,-1.5) .. (9.75, -2.25);

\draw[line width=1pt] (8, 0) .. controls(3,3) and (3,-4.6) .. (8, -5); 
\draw[line width=1pt] (8, 0) .. controls(13,0.3) and (13,-4.6) .. (8, -5);

\end{tikzpicture}
\end{figure} 
\end{myEnumerate}

\end{myEnumerate}

Then we are going to define the detour matrixes.

\begin{definition}(Detour matrix)\label{matriz} Let $\tau$ be a tagged triangulation of a surface with marked points $(\Sigma,M)$. Using the detour curves of $(\tau^{\circ},i')$ and considering the cases when the arc $j$ is a side of a self-folded triangle or if it is not. We define the detour matrix for each arc $j\in \tau^{\circ}$, the arrows and columns are indexed by the different crossing points of $i'$ with $j$.

Let's define the $q$-th column (associated to the $q$-th crossing point) of the detour matrix through the followings rules:
\end{definition}

\begin{myEnumerate}
\item[a)] The arc $j$ is a side of a self-folded triangle.
Under this condition we consider if the arc $j$ is a folded side of a self-folded triangle or if it is not:

\begin{myEnumerate}
\item[$a_{1}$)] The arc $j$ is the folded side of a self-folded triangle $\Delta'$. Let's denote by $\Delta$ the unique not self-folded triangle in $\tau^{\circ}$ which share one arc with $\Delta'$.

We are going to consider if we assign one or two detour matrixes to the arc $j$.

\begin{myEnumerate}
\item[$a_{11}$)] We assign to the arc $j$ only one detour matrix denoted by $D^{\overline{\Delta}}_{i',j}$.

\begin{myEnumerate}
\item[i)] The $q$-th coordinate is 1.
\item[ii)] If $\tilde{q}$ is a crossing point which is a final point of a detour 
$d^{\overline{\Delta},n}_{(q,\tilde{q})}$ and there is one element $(q,\tilde{q},r_{1},r_{2},r_{l},p)$ or $(q,\tilde{q},q_{1},r_{1},r_{2},r_{3},r_{l},p)$ in ${\mathcal B}^{\Delta,n}_{i',j}$, then the $\tilde{q}$-th coordinate is $\delta_{\tau}(p)$. In addition, if $d^{\Delta,n}_{(q,\tilde{q})}$ is of type 1 then the coordinate associated to the crossing point $r_{l-5}$ of $\gamma$ with $j$ is also $\delta_{\tau}(p)$. 
\item[iii)] All other coordinates of $q$ are zero.
\end{myEnumerate}

\item[$a_{12}$)] We assign two detour matrixes to the arc $j$, let's denote by $D^{\Delta}_{i',j}$ and $D^{\overline{\Delta}}_{i',j}$ the detour matrixes.

The $q$-th coordinate is 1 for both matrixes and zero  all other coordinates are zero.

\end{myEnumerate}
\item[$a_{2}$)] The arc $j$ is the non folded side of a self-folded triangle $\Delta'$. Let's denote by  $j'$ the folded side of $\Delta'$.

\begin{myEnumerate}
\item[$a_{21}$)]We assign only one matrix to the arc $j$ denoted by $D^{\Delta}_{i',j}$. 
	\begin{myEnumerate}
		\item[i)] The $q$-th coordinate is 1.
		\item[ii)] If $\tilde{q}$ is a crossing point which is a final point of a detour $d^{\Delta,n}_{(q,\tilde{q})}$ and there is one element $(q_{0},q_{1},q_{2},q,r_{2},r_{3}, \tilde{q},p)$ in ${\mathcal B}^{\Delta,n}_{i',j'}$, then:

 The $\tilde{q}$-th coordinate is $\delta_{\tau}(p)$ if one of the two conditions described below are met; 

			\begin{itemize}
				\item{} $\Delta$ is a triangle of type 2.
				\item{} $\Delta$ is a triangle of type 3 and $\gamma \cup [q_{1},q_{0}]_{j}$ is a closed curve which divides $\Sigma$ into two regions. One region is homeomorphic to a disk with only one puncture $p$, where $\gamma$ is the curve which corresponds to the element $(q_{0},q_{1},q_{2},q,r_{2},r_{3},\tilde{q},p)$. In addition, if $d^{\Delta}$ is of type 1 then the coordinate which corresponds to the $(r-4)$ -th crossing point of $\gamma$ with the arcs of $\tau^{\circ}$ is also $\delta_{\tau}(p)$.

			\end{itemize}

The $\tilde{q}$-th coordinate is -$\delta_{\tau}(p)$  if  the curve $\gamma \cup [q_{1},q_{0}]_{j}$ divides $\Sigma$ into two regions. One region is homeomorphic to a disk with only two punctures being one puncture $p$, where $\gamma$ is the curve which corresponds to the element $(q_{0},q_{1},q_{2},q,r_{2},r_{3},\tilde{q},p)$. In addition, if $d^{\Delta}$ is of type 1, then the coordinate which corresponds to the $(r-4)$-th crossing point of $\gamma$ with the arcs of $\tau^{\circ}$ is also -$\delta_{\tau}(p)$.

\item[iii)] All other coordinates are zero.
\end{myEnumerate}

\item[$a_{22}$)] We assign two matrixes to arc $j$, let's denote by $D^{\Delta}_{i',j}$ and $D^{\overline{\Delta}}_{i',j}$ the detour matrixes.

\begin{myEnumerate}
\item[i)] The $q$-th coordiate for both matrixes is 1,
\item[ii)] If $\tilde{q}$ is a crossing point which is a final point of a detour $d^{\Delta,n}_{(q,\tilde{q})}$ and there is an element $(q_{0},q_{1},q_{2},q,r_{2},r_{3},\tilde{q},p)$ in ${\mathcal B}^{\Delta,n}_{i',j'}$ then:

\begin{itemize}
\item{}The $\tilde{q}$-th coordinate of matrix $D^{\overline{\Delta}}_{i',j}$ is $\delta_{\tau}(p)$ if $\gamma \cup [q_{1},q]_{j'}$  is a closed curve which divides $\Sigma$ in two regions. One region is homeomorphic to a disk with only one puncture being the puncture $p$.

\item{} The $\tilde{q}$-th coordinate of the matrix $D^{\overline{\Delta}}_{i',j}$ is -$\delta_{\tau}(p)$ if $\gamma \cup [q_{1},q]_{j'}$ is a closed curve which divides $\Sigma$ in two regions. One region is homeomorphic to a disk with only two punctures being one puncture $p$.

\end{itemize}

\item[iii)] All other coordinates of $q$ are zero for both matrixes.

\end{myEnumerate}

\end{myEnumerate}

\end{myEnumerate}

\item[b)] The arc $j$ is not a side of any self-folded triangle.

If the arc $j$ is not a side of any self-folded triangle, then $j$ is a side of exactly two triangles which are not self-folded. Let's denote by $\Delta_{1}$ and $\Delta_{2}$ those triangles. We are going to consider the cases when the arc $j$ is assigned 2,3 or 4 detour matrixes.

\begin{myEnumerate}
\item[$b_{1}$)] We assign two detour matrixes to arc $j$.
We associate two matrixes $D^{\Delta}_{i',j}$ one by each triangle $\Delta_{k}$ with $k=1,2$.

\begin{myEnumerate}
\item[i)] The $q$-th coordinate is 1.
\item[ii)]If $\tilde{q}$ is a crossing point, which is a final point of $d^{\Delta,n}_{(q,\tilde{q})}$, and there is an element $(q,\tilde{q},r_{1},r_{2},r_{l},p)$ in ${\mathcal B}^{\Delta,n}_{i',j}$, then 
the $\tilde{q}$-th coordinate is $\delta_{\tau}(p)$. In addition, if $d^{\Delta,n}_{(q,\tilde{q})}$ is of type 1, then the coordinate associated to the crossing point $r_{l-3}$ of $\gamma$ with $j$ is also $\delta_{\tau}(p)$. 
\item[iii)] All other coordinates of $q$ are zero for both matrixes.
\end{myEnumerate}

\item[$b_{2}$)] We assign three detour matrixes to arc $j$.
We associate three matrixes, two of them denoted by $D^{\Delta}_{i',j}$, one by each triangle $\Delta_{k}$ with $k=1,2$. In addition, we assign the third matrix denoted by $D^{\overline{\Delta}}_{i',j}$ where $\Delta$ is a triangle of type 3 or $\Delta$ is a triangle of type 2.

\begin{myEnumerate}
\item[i)] The $q$-th coordinate of the three matrixes is 1.
\item[ii)] If $\tilde{q}$ is a crossing point, which is final point of a detour $d^{\Delta,n}_{(q,\tilde{q})}$ and there is an element $(q,\tilde{q},r_{1},r_{2},r_{l},p)$ in ${\mathcal B}^{\Delta,n}_{i',j}$, then the $\tilde{q}$ -th coordinate of $D^{\Delta}_{i',j}$ is $\delta_{\tau}(p)$ if the relative interior of the curve $d^{\Delta}_{(q,\tilde{q})} $ that does not intersect any arc of $\tau^{\circ}$. Otherwise the $\tilde{q}$-th coordinate of $D^{\overline{\Delta}}_{i´,j}$ is $\delta_{\tau}(p)$. In addition, if $d^{\Delta}_{(q,\tilde{q})}$ is of type 1, then the coordinate corresponding to the crossing point $r_{l-3}$ of the curve $\gamma$ with the arc $j$ is also $\delta_{\tau}(p)$. 
\item[iii)] All other coordinates of $q$ are zero for the three matrixes.
\end{myEnumerate}

\item[$b_{3}$)] We assign four detour matrixes to the arc $j$.
We associate four matrixes two of them denoted by $D^{\Delta}_{i',j}$ one for each triangle $\Delta_{k}$ with $k=1,2$ and two additional denoted by $D^{\overline{\Delta}}_{i',j}$.

\begin{myEnumerate}
\item[i)] The $q$-th coordinate of the four matrixes is 1.
\item[ii)] If $\tilde{q}$ is the final point of a detour $d^{\Delta,n}_{(q,\tilde{q})}$ and there is an element $(q,\tilde{q},r_{1},r_{2},r_{l},p)$ in ${\mathcal B}^{\Delta,n}_{i',j}$, then the $\tilde{q}$-th coordinate of $D^{\Delta}_{i',j}$ is $\delta_{\tau}(p)$ if the relative interior of the curve $d^{\Delta}_{(q,\tilde{q})}t$ does not intersect any arc of $\tau^{\circ}$. Otherwise the $\tilde{q}$-th coordinate of $D^{\overline{\Delta}}$ is $\delta_{\tau}(p)$. In addition, if $d^{\Delta}_{(q,\tilde{q})}$ is of type 1, then the coordinate corresponding to the crossing point $r_{l-3}$ of the curve $\gamma$ with the arc $j$ is also $\delta_{\tau}(p)$.

\item[iii)] All other coordinates of $q$ are zero for the four matrixes.
\end{myEnumerate}

\end{myEnumerate}

\end{myEnumerate}

Now, we are going to define the auxiliary matrixes.

\begin{definition}(Auxiliary matrix )\label{matriz auxiliar} Let $\tau$ be a tagged triangulation of a surface with marked points $(\Sigma,M)$. Using the auxiliary curves $(\tau^{\circ},i')$ and considering the cases when arc $j$ is a side of a self-folded triangle or it is not. We define the auxiliary matrix for each arc $j\in \tau^{\circ}$. The arrows and columns are indexed by the different crossing points of $i'$ with $j$.

Let's define the $q$-th column (associated to the $q$-th crossing point) of the auxiliary matrix through the followings rules:
\end{definition}

\begin{myEnumerate}
\item[a)] The arc $j$ is a side of a self-folded triangle $\Delta'$. Let's denote by $\Delta$ the unique triangle in $\tau^{\circ}$ which shares one side with $\Delta$.

Since $j$ is a side of a self-folded triangle we are going to consider the cases when $j$ is the folded side or $j$ is the non folded side.

\begin{myEnumerate}
\item[$a_{1}$)] The arc $j$ is the folded side of $\Delta'$.
We are going to consider the cases when arc $j$ has been assigned one or two matrixes.

\begin{myEnumerate}
\item[$a_{11}$)]We assign only one matrix to the arc $j$ denoted by $E^{\overline{\Delta}}_{i',j}$. 

The $q$-th coordinate is 1 and all others are zero. 

\item[$a_{12}$)] We assign two matrixes to the arc $j$. Let's denote by $E^{\Delta}_{i',j}$ and $E^{\overline{\Delta}}_{i',j}$ both auxiliary matrixes.
\begin{myEnumerate}
\item[i)] The $q$-th coordinate of both matrixes is 1.
\item[ii)] If $\tilde{q}\in i'\cap j$ is a crossing point which is a final point of an auxiliary curve $e^{\Delta,n}_{(q,\tilde{q})}$ and there is an element $(q,\tilde{q},q_{2},r_{1},r_{2},r_{3},r_{l},p)$ in ${\mathcal B}^{\Delta,n}_{i',j}$, then the $\tilde{q}$-th coordinate of $E^{\Delta}_{i',j}$ is $\delta_{\tau}(p)$.

\item[iii)]All other coordinates of $q$ are zero.
\end{myEnumerate}

\end{myEnumerate}

\item[$a_{2}$)] The arc $j$ is the non folded side of $\Delta'$, let's denote by $j'$ the folded side of $\Delta'$.

We are going to consider the cases when we assign one or two matrixes to arc $j$.
\begin{myEnumerate}
\item[$a_{21}$)] We assign only one matrix to $j$ denoted by $E^{\Delta}_{i',j}$.
\begin{myEnumerate}
\item[i)] The $q$-th coordinate is 1.
\item[ii)] If $\tilde{q}\in i'\cap j$ is a crossing point which is a final point of a auxiliary curve $e^{\Delta,n}_{(q,\tilde{q})}$ and there is an element $(q_{0},q_{1},q_{2},q,r_{2},r_{3},\tilde{q},p)$ in ${\mathcal B}^{\Delta,n}_{i',j'}$, then the $\tilde{q}$-th coordinate is $\delta_{\tau}(p)$.
\item[iii)] All other coordinates of $q$ are zero.
\end{myEnumerate}

\item[$a_{22}$)] We assign two matrixes to arc $j$. Let's denote by $E^{\Delta}_{i',j}$ and $E^{\overline{\Delta}}_{i',j}$ the auxiliary matrixes assigned to arc $j$.

The $q$-th coordinate for both matrixes is 1 and the others are 0.
\end{myEnumerate}

\end{myEnumerate}
\item[b)] The arc $j$ is not a side of any self-folded triangle. 

We are going to consider the cases when we assign 2,3 o 4 matrixes to the arc $j$.
\begin{myEnumerate}
\item[$b_{1}$)] We assign 2 matrixes to $j$.

Assign one matrix for each non self-folded triangle which has arc $j$ as one of its sides. Let's denote by $E^{\Delta}_{i',j}$ the auxiliary matrixes.

The $q$-th coordinate for both matrixes is 1 and the others are 0.

\item[$b_{2}$)] We assign 3 matrixes to arc $j$.

Assign one matrix for each non-delf-folded triangle which share arc $j$ as one of its sides and let's denote it by $E^{\Delta}_{i',j}$. In addition, we define one additional matrix denoted by $E^{\overline{\Delta}}_{i',j}$, in this case $\Delta$ is a triangle of type 2 such that arc $j$ has label $l_{2}$ in $\Delta$.

We are going to define the $q$-th column in the matrix $E^{\overline{\Delta}}_{i',j}$.

The $q$-th coordinate is 1 and the others are 0.

Now let's define the matrixes $E^{\Delta}_{i',j}$.

We will denote by $\Delta'$ the unique self-folded triangle which shares one side with $\Delta$, being $m$ and $m'$ the sides of $\Delta'$, where $m'$ is the folded side of $\Delta'$. Let's consider the cases when $q$ is a crossing point, which is a final point of an auxiliary curve, or it is not:

\begin{myEnumerate}
\item[i)] $q$ is not a final point of any auxiliary curve.
\begin{itemize}
\item[$\bullet$] The $q$-th coordinate is 1.
\item[$\bullet$] If $\tilde{q}$ is the final point of an auxiliary curve $e^{\Delta,n}_{i',j}$, there are crossing points denoted by $x$,$x'$ of $i'$ with $m'$, the relative interior of the segment $[q,x]_{i'}$ intersects the arc $m$ in only one point and there is $(x,x',\tilde{q},r_{1},r_{2},r_{3},r_{l},p)\in {\mathcal B}^{\Delta,n}_{i',m'}$, then the $\tilde{q}$-th coordinate is $\delta_{\tau}(p)$.

\end{itemize}

\item[ii)] $q$ is a final point of an auxiliary curve $e^{\Delta,n}_{i',j}$.
\begin{itemize}
\item[$\bullet$] The $q$-th coordinate is 0.
\item[$\bullet$] If $x$,$x'$ are crossing points of $i'$ with $m'$, $q'$ is a crossing point of $i'$ with $j$, the relative interior of the segment $[q',x]_{i'}$ intersects the arc $m$ in only one point, $q''$ is the final point of the auxiliary curve $e^{\Delta,n-1}_{i',j}$ and there is one element $(x,x',q,r_{1},r_{2},r_{3},r_{l},p)\in {\mathcal B}^{\Delta,n}_{i',m'}$, then the $q''$-th coordinate is 1 when $n\textgreater 1$ and the $q'$-th coordinate is 1 when $n=1$.

\end{itemize} 
\end{myEnumerate}

\item[$b_{3}$)] We assign four matrixes to the arc $j$.

Let's assign one matrix for each non self-folded triangle which has $j$ as a side, let's denote by $E^{\Delta}_{i',j}$ the auxiliary matrixes. In addition we will define two additional matrixes which we denote by $E^{\overline{\Delta}}_{i',j}$.

We are going to define the $q$-th column of $E^{\overline{\Delta}}_{i',j}$.

The $q$-th coordinate is 1 and all the others are zero.

Now let's define the matrixes $E^{\Delta}_{i',j}$.

We will denote by $\Delta'$ the unique self-folded triangle which share one side with $\Delta$, being $m$ and $m'$ the sides of $\Delta'$, where $m'$ is the folded side of $\Delta'$. Let's consider the cases when $q$ is a crossing point, which is a final point of an auxiliary curve or if it is not:

\begin{myEnumerate}
\item[i)] $q$ is not a final point of of any auxiliary curve.
\begin{itemize}
\item[$\bullet$] The $q$-th coordinate is 1.
\item[$\bullet$] If $\tilde{q}$ is the final point of an auxiliary curve $e^{\Delta,n}_{i',j}$, there are crossing points denoted by $x$,$x'$ of $i'$ with $m'$, the relative interior of the segment $[q,x]_{i'}$ intersects the arc $m$ in only one point and there is $(x,x',\tilde{q},r_{1},r_{2},r_{3},r_{l},p)\in {\mathcal B}^{\Delta,n}_{i',m'}$, then the $\tilde{q}$-th coordinate is $\delta_{\tau}(p)$.

\end{itemize}

\item[ii)] $q$ is a final point of an auxiliary curve $e^{\Delta,n}_{i',j}$.
\begin{itemize}
\item[$\bullet$] The $q$-th coordinate is 0.
\item[$\bullet$] If $x$,$x'$ are crossing points of $i'$ with $m'$, $q'$ is a crossing point of $i'$ with $j$, the relative interior of the segment $[q',x]_{i'}$ intersects the arc $m$ in only one point, $q''$ is the final point of the auxiliary curve $e^{\Delta,n-1}_{i',j}$ and there is one element $(x,x',q,r_{1},r_{2},r_{3},r_{l},p)\in {\mathcal B}^{\Delta,n}_{i',m'}$, then the $q''$-th coordinate is 1 when $n\textgreater 1$ and the $q'$-th coordinate is 1 when $n=1$.

\end{itemize} 
\end{myEnumerate}

\end{myEnumerate}
\end{myEnumerate}

Once we have defined the detour matrixes, auxiliary matrixes and the string representations we are ready to define the arc representation $M(\tau,i)$.

\begin{definition}\label{rep de arco}
Let $\{p,q\}$ be the set of ends of the arc $j$, first let's define the vertex of the arc representation $M(\tau,i)$. Let $j$ be an arc of $\tau^{\circ}$, we consider the cases when the puntures $p$ or $q$ (or both) are marked points or if they are not.

\begin{myEnumerate}
\item[a)] $p$ and $q$ are marked points.

$(M(\tau,i))_{j}=(m(\tau,i))_{j}$.

\item[b)] Only one end of $i'$ is a marked point. Without loss of generality we supposed that the puncture $q$ is not a marked point and belongs to an arc $j_{1}$ of $\tau^{\circ}$.

Under these conditions we consider the following two subcases:

\begin{myEnumerate}
\item[$b_{1}$)] The arc $j_{1}$ is the folded side of self-folded triangle $\Delta'$. Let's denote by $m$ the non folded side of $\Delta'$ and by $q'$ the crossing point of $i'$ with $m$ such that the relative interior of the segment $[q',q]_{i'}$ does not intersect any arc of $\tau^{\circ}$.

\begin{myEnumerate}
\item[$b_{11}$)] If $j=j_{1}$, then $(M(\tau,i))_{j}=(m(\tau,i))_{j}/K_{j,q}$.
\item[$b_{12}$)] $j=m$.
$(M(\tau,i))_{j}=(m(\tau,i))_{j}/K_{j,q'}$.
\item[$b_{13}$)] If $m \neq j \neq j_{1}$ then $(M(\tau,i))_{j}=(m(\tau,i))_{j}$.
\end{myEnumerate}

\item[$b_{2}$)] The arc $j_{1}$ is not a side of any self-folded triangle.

\begin{myEnumerate}
\item[$b_{21}$)] If $j=j_{1}$ then:

$(M(\tau,i))_{j}=(m(\tau,i))_{j}/K_{j,q}$.

\item[$b_{22}$)] If $ j \neq j_{1}$, then $(M(\tau,i))_{j}=(m(\tau,i))_{j}$.
\end{myEnumerate}

\end{myEnumerate}
\item[c)]Any end of the curve $i'$ is a marked point. Let'¿s suppose that the puncture $p$ belongs to the arc $j_{1}$ and the puncture $q$ belongs to $j_{2}$, where $j_{1}$ and $j_{2}$ are arcs of $\tau^{\circ}$.
We are going to consider the subcases described below:

\begin{myEnumerate}
\item[$c_{1}$)] The arc $j_{k}$ is the folded side of $\Delta_{k}$ for $k=1,2$. Let's denote by $m_{k}$ the non folded side of $\Delta_{k}$ and by $q_{k}$ the crossing point of $i'$ with $m_{k}$ such that the relative interior of the segments $[q_{1},p]_{i'}$, $[q_{2},q]_{i'}$ do not intersect any arc of $\tau^{\circ}$.

\begin{myEnumerate}
\item[$c_{11}$)] If $j=j_{1}$, then $(M(\tau,i))_{j}=(m(\tau,i))_{j}/K_{j,p}$.
\item[$c_{12}$)] If $j=j_{2}$, then $(M(\tau,i))_{j}=(m(\tau,i))_{j}/K_{j,q}$.
\item[$c_{13}$)] $j=m_{k}$.
$(M(\tau,i))_{j}=(m(\tau,i))_{j}/K_{j,q_{k}}$.
\item[$c_{14}$)] If $m_{k} \neq j \neq j_{k}$, then $(M(\tau,i))_{j}=(m(\tau,i))_{j}$.
\end{myEnumerate}

\item[$c_{2}$)]The arc $j_{1}$ is the folded side of a self-folded triangle $\Delta'$ and $j_{2}$ is not a side of any self-folded triangle. We denote by $m$ the non folded side of $\Delta'$ and by $p'$ the crossing point of $i'$ with $m$ such that the relative interior of the segment $[p',p]_{i'}$ does not intersect any arc of $\tau^{\circ}$.

\begin{myEnumerate}
\item[$c_{21}$)] If $j=j_{1}$, then $(M(\tau,i))_{j}=(m(\tau,i))_{j}/K_{j,p}$.
\item[$c_{22}$)] $j=m$.
$(M(\tau,i))_{j}=(m(\tau,i))_{j}/K_{j,p'}$.
\item[$c_{23}$)] $j=j_{2}$.

$(M(\tau,i))_{j}=(m(\tau,i))_{j}/K_{j,q}$.

\item[$c_{24}$)] If $m \neq j \neq j_{r}$ for $r=1,2$, then $(M(\tau,i))_{j}=(m(\tau,i))_{j}$.
\end{myEnumerate}

The case when $j_{2}$ is the folded side of a self-folded triangle $\Delta'$ and $j_{1}$ is not a side of any self-folded triangle is symmetric.

\item[$c_{3}$)] For $k=1,2$, arc $j_{k}$ is not the folded side of any self-folded triangle.
\begin{myEnumerate}
\item[$c_{31}$)] $j=j_{1}$.

$(M(\tau,i))_{j}=(m(\tau,i))_{j}/K_{j,p}$.
\item[$c_{32}$)] $j=j_{2}$.

$(M(\tau,i))_{j}=(m(\tau,i))_{j}/K_{j,q}$.
\item[$c_{33}$)] If $ j \neq j_{r}$, then $(M(\tau,i))_{j}=(m(\tau,i))_{j}$.
\end{myEnumerate}

\end{myEnumerate}

\end{myEnumerate}

Now, we are going to define the linear transformations of each arrow $\alpha:j \rightarrow k$ in $\widehat{Q}(\tau^{\circ})$, let's denote by $\Delta^{\alpha}$ the unique self-folded triangle which contains a segment of $\alpha$ in the interior. For each arrow $\alpha:j \rightarrow k$ we split the matrix $m(\tau,i)_{\alpha}$ as the sum of two matrixes denoted by $m(\tau,i)_{\alpha}^{+}$ and $m(\tau,i)_{\alpha}^{-}$, where $m(\tau,i)_{\alpha}^{+}$ is the matrix with positive entries of $m(\tau,i)_{\alpha}$ and $m(\tau,i)_{\alpha}^{-}$ is the matrix with negative entries of $m(\tau,i)_{\alpha}$. We will consider the cases when $p$ or $q$ (or both) are marked points or if they are not:

\begin{myEnumerate}
\item[a)] The points $p$ and $q$ are marked points.

Remembering how we draw $(\widehat{Q}(\tau^{\circ}),\widehat{S}(\tau^{\circ}))$ in $(\Sigma,M)$ (See the section 3), we will consider the following cases: 

\begin{myEnumerate}
\item[$a_{1}$)] $\alpha$ intersect exactly one arc of $\tau^{\circ}$. 
If both matrixes $m(\tau,i)_{\alpha}^{+}$ and $m(\tau,i)_{\alpha}^{-}$ are not the matrix zero, then
$$(M(\tau,i))_{\alpha}= E^{\overline{\Delta}^{\alpha}}_{i',k}(D^{\overline{\Delta}^{\alpha}}_{i',k}(m(\tau,i))^{+}_{\alpha} + (m(\tau,i))^{-}_{\alpha}).$$ 

Otherwise $(M(\tau,i))_{\alpha}=E^{\Delta^{\alpha}}_{i',k}D^{\Delta^{\alpha}}_{i',k}m(\tau,i)_{\alpha}$. 

\item[$a_{2}$)] $\alpha$ does not intersect any arc of $\tau^{\circ}$ or $\alpha$ intersects exactly two arcs of $\tau^{\circ}$.
If both matrixes $m(\tau,i)_{\alpha}^{+}$ and $m(\tau,i)_{\alpha}^{-}$ are not the matrix zero, then 
$$(M(\tau,i))_{\alpha}=E^{\Delta^{\alpha}}_{i',k}(D^{\Delta^{\alpha}}_{i',k}(m(\tau,i))^{+}_{\alpha}+(m(\tau,i))^{-}_{\alpha}).$$ 

Otherwise $(M(\tau,i))_{\alpha}=E^{\Delta^{\alpha}}_{i',k}D^{\Delta^{\alpha}}_{i',k}m(\tau,i)_{\alpha}$. 
\end{myEnumerate}

\item[b)] Only one end of the curve $i'$ is a marked point. Without loss of generality, we suppose that $q$ is not a marked point and belongs to an arc $j_{1}$ of $\tau^{\circ}$.
We will consider the cases when $j_{1}$ is the folded side of a self-folded triangle or it is not.

\begin{myEnumerate}
\item[$b_{1}$)] The arc $j_{1}$ is the folded side of a self-folded triangle$\Delta'$. Let's denote by $m$ the non folded side of $\Delta'$.

If any end of $\alpha$ is $j_{1}$ or $m$, then $M(\tau,i)_{\alpha}$ is defined as in the case $a)$. Otherwise we consider the followings subcases: 

\begin{myEnumerate}
\item[$b_{11}$)] $j_{1}=j$ or $j=m$.
Let's consider the followings subcases:

\begin{myEnumerate}
\item[$b_{111}$)] $\alpha$ intersects exactly one arc of $\tau^{\circ}$.

If both matrixes $m(\tau,i)_{\alpha}^{+}$ and $m(\tau,i)_{\alpha}^{-}$ are not the matrix zero, then
$$(M(\tau,i))_{\alpha}= (E^{\overline{\Delta}^{\alpha}}_{i',k}(D^{\overline{\Delta}^{\alpha}}_{i',k}(m(\tau,i))^{+}_{\alpha} + (m(\tau,i))^{-}_{\alpha}))\iota.$$ 

Otherwise $(M(\tau,i))_{\alpha}= E^{\overline{\Delta}^{\alpha}}_{i',k}D^{\overline{\Delta}^{\alpha}}_{i',k}m(\tau,i)_{\alpha}\iota $.

Where $\iota:M(\tau,i)_{j}\longrightarrow m(\tau,i)_{j}$ is the canonical inclusion.

\item[$b_{112}$)] $\alpha$ does not intersect any arc of $\tau^{\circ}$ or $\alpha$ intersects exactly two arcs of $\tau^{\circ}$.

If both matrixes $m(\tau,i)_{\alpha}^{+}$ and $m(\tau,i)_{\alpha}^{-}$ are not the matrix zero, then $(M(\tau,i))_{\alpha}=(E^{\Delta^{\alpha}}_{i',k}( D^{\Delta^{\alpha}}_{i',k}m(\tau,i)^{+}_{\alpha} + m(\tau,i)^{-}_{\alpha} ))\iota$. Otherwise $(M(\tau,i))_{\alpha}=E^{\Delta^{\alpha}}_{i',k}D^{\Delta^{\alpha}}_{i',k}m(\tau,i)_{\alpha} \iota$.

\end{myEnumerate} 
Where $\iota:M(\tau,i)_{j}\longrightarrow m(\tau,i)_{j}$ is the canonical inclusion. 
\item[$b_{12}$)] $j_{1}=k$ or $m=k$.
\begin{myEnumerate}
\item[$b_{121}$)] $\alpha$ intersects exactly one arc of $\tau^{\circ}$.

If both matrixes $m(\tau,i)_{\alpha}^{+}$ and $m(\tau,i)_{\alpha}^{-}$ are not the matrix zero, then
$$(M(\tau,i))_{\alpha}= \pi(E^{\overline{\Delta}^{\alpha}}_{i',k}(D^{\overline{\Delta}^{\alpha}}_{i',k}(m(\tau,i))^{+}_{\alpha} + (m(\tau,i))^{-}_{\alpha})).$$ 

Otherwise, $(M(\tau,i))_{\alpha}= \pi E^{\overline{\Delta}^{\alpha}}_{i',k}D^{\overline{\Delta}^{\alpha}}_{i',k}m(\tau,i)_{\alpha} $ .

Where $\pi: m(\tau,i)_{j_{1}}\longrightarrow M(\tau,i)_{j_{1}}$ is the canonical projection.

\item[$b_{122}$)] $\alpha$ does not intersect any arc of $\tau^{\circ}$ or $\alpha$ intersects exactly two arcs of $\tau^{\circ}$.
If both matrixes $m(\tau,i)_{\alpha}^{+}$ and $m(\tau,i)_{\alpha}^{-}$ are not the matrix zero, then $(M(\tau,i))_{\alpha}=\pi (E^{\Delta^{\alpha}}_{i',k}(D^{\Delta^{\alpha}}_{i',k}(m(\tau,i))^{+}_{\alpha} + (m(\tau,i))^{-}_{\alpha}))$. Otherwise $(M(\tau,i))_{\alpha}=\pi E^{\Delta^{\alpha}}_{i',k}D^{\Delta^{\alpha}}_{i',k}m(\tau,i)_{\alpha} $.
\end{myEnumerate}
\end{myEnumerate} 

Where $\pi: m(\tau,i)_{j}\longrightarrow M(\tau,i)_{j}$ is the canonical projection. 
\item[$b_{2}$)] $j_{1}$ is not the folded side of any self-folded triangle.
If any end of $\alpha$ is $j_{1}$, then $M(\tau,i)_{\alpha}$ is defined as in $a)$. Otherwise let's consider the next two cases:

\begin{myEnumerate}
\item[$b_{21}$)] $j_{1}=j$.
Under this situation we will consider the following two cases:

\begin{myEnumerate}
\item[$b_{211}$)] $\alpha$ intersects exactly one arc of $\tau^{\circ}$.

If both matrixes $m(\tau,i)_{\alpha}^{+}$ and $m(\tau,i)_{\alpha}^{-}$ are not the matrix zero, then $(M(\tau,i))_{\alpha}=(E^{\overline{\Delta}^{\alpha}}_{i',k}( D^{\overline{\Delta}^{\alpha}}_{i',k}m(\tau,i)^{+}_{\alpha} + m(\tau,i)^{-}_{\alpha} ))\iota$. Otherwise $(M(\tau,i))_{\alpha}=E^{\overline{\Delta}^{\alpha}}_{i',k}D^{\overline{\Delta}^{\alpha}}_{i',k}m(\tau,i)_{\alpha} \iota$.
\item[$b_{212}$)] $\alpha$ does not intersect any arc of $\tau^{\circ}$ or $\alpha$ intersects exactly two arcs $\tau^{\circ}$.
If both matrixes $m(\tau,i)_{\alpha}^{+}$ and $m(\tau,i)_{\alpha}^{-}$ are not the matrix zero, then $(M(\tau,i))_{\alpha}=(E^{\Delta^{\alpha}}_{i',k}( D^{\Delta^{\alpha}}_{i',k}m(\tau,i)^{+}_{\alpha} + m(\tau,i)^{-}_{\alpha} ))\iota$. Otherwise $(M(\tau,i))_{\alpha}=E^{\Delta^{\alpha}}_{i',k}D^{\Delta^{\alpha}}_{i',k}m(\tau,i)_{\alpha} \iota$.
\end{myEnumerate}
Where $\iota:M(\tau,i)_{j}\longrightarrow m(\tau,i)_{j}$ is the canonical inclusion.

\item[$b_{22}$)] $j_{1}=k$.
\begin{myEnumerate}
\item[$b_{221}$)] $\alpha$ intersects exactly one arc of $\tau^{\circ}$. 
If both matrixes $m(\tau,i)_{\alpha}^{+}$ and $m(\tau,i)_{\alpha}^{-}$ are not the matrix zero, then $(M(\tau,i))_{\alpha}=\pi (E^{\overline{\Delta}^{\alpha}}_{i',k}( D^{\overline{\Delta}^{\alpha}}_{i',k}m(\tau,i)^{+}_{\alpha} + m(\tau,i)^{-}_{\alpha} ))$. Otherwise $(M(\tau,i))_{\alpha}=\pi E^{\overline{\Delta}^{\alpha}}_{i',k}D^{\overline{\Delta}^{\alpha}}_{i',k}m(\tau,i)_{\alpha} $.
\item[$b_{222}$)] $\alpha$ does not intersect any arc of $\tau^{\circ}$ or $\alpha$ intersects exactly two arcs of $\tau^{\circ}$.
If both matrixes $m(\tau,i)_{\alpha}^{+}$ y $m(\tau,i)_{\alpha}^{-}$ are not the matrix zero, then $(M(\tau,i))_{\alpha}=\pi (E^{\Delta^{\alpha}}_{i',k}( D^{\Delta^{\alpha}}_{i',k}m(\tau,i)^{+}_{\alpha} + m(\tau,i)^{-}_{\alpha} ))$. Otherwise $(M(\tau,i))_{\alpha}=\pi E^{\Delta^{\alpha}}_{i',k}D^{\Delta^{\alpha}}_{i',k}m(\tau,i)_{\alpha} $.
\end{myEnumerate}
\end{myEnumerate}
Where $\pi: m(\tau,i)_{j}\longrightarrow M(\tau,i)_{j}$ is the canonical projection.
\end{myEnumerate}

\item[c)] Any end of $i'$ is a marked point. We suppose that the puncture $p$ belongs to $j_{1}$ and $q$ belongs to $j_{2}$, where $j_{1}$ and $j_{2}$ belong to the triangulation $\tau^{\circ}$.

In this situation we apply the rules in $b)$ in the arcs $j_{1}$ and $j_{2}$ to define the linear transformation $(M(\tau,i))_{\alpha}$.

\end{myEnumerate}

\end{definition}

\textbf{EXAMPLES OF THE ARC REPRESENTATION $M(\tau,i)$} (Definition \ref{rep de arco}).
\begin{itemize}

\item[1)] Example 1.
\begin{figure}[H]
\centering 
$$ 
 }
\end{center}
\caption{Quiver.}

\end{figure}

\begin{figure}[H]
\begin{center}
\hspace*{-2em}{%
}

\end{center}
\caption{Arc representation $M(\tau,i)$.}

\end{figure}

\item[2)] Example 2.

\begin{figure}[H]
\centering 
$$ %
 }

\end{center}
\caption{Quiver.}

\end{figure}

\begin{figure}[H]
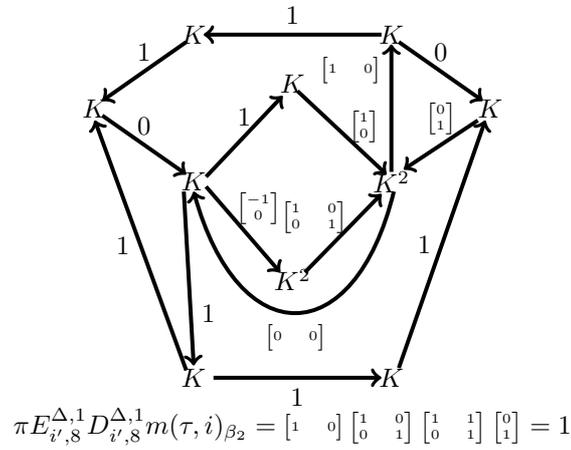

\begin{center}
\hspace*{0.6em}{%
}

\end{center}
\caption{Arc representation $M(\tau,i)$.}

\end{figure}

\end{itemize}

\include{sec5}

\end{document}